\newif\ifarxiv
\arxivtrue

\PassOptionsToPackage{pagebackref}{hyperref} 
\documentclass{acta_2022}

\ifarxiv
  \makeatletter
  \def\ps@titlepage{\ps@empty}
  \makeatother
\fi

\usepackage[T1]{fontenc}
\usepackage[utf8]{inputenc} 
\usepackage{newtxtext}\usepackage{newtxmath}
\usepackage{mathtools}
\usepackage{amsmath}
\usepackage{nicefrac}

\ifarxiv
\else
\fi
\usepackage[longnamesfirst,sort]{natbib} \usepackage{har2nat}
\usepackage{citet-and} 
\usepackage{cite-yearsort} 
\setcitestyle{comma,aysep={}}

\renewcommand*{\backref}[1]{}
\renewcommand*{\backrefalt}[4]{%
	\ifcase #1 %
	(No citations.)
	\or
	(Cited page~#2.)
	\else
	(Cited pages~#2.)
	\fi
}
\shortcites{samplepreprint} 
\ifarxiv
  
\else
  
\fi

\DeclareSymbolFont{CMlargesymbols}{OMX}{cmex}{m}{n}
\DeclareMathDelimiter{(}{\mathopen} {operators}{"28}{CMlargesymbols}{"00}
\DeclareMathDelimiter{)}{\mathclose}{operators}{"29}{CMlargesymbols}{"01}
\DeclareMathAlphabet\mathcal{OMS}{cmsy}{m}{n}
\SetMathAlphabet\mathcal{bold}{OMS}{cmsy}{b}{n}
\DeclareMathAlphabet{\mathsf}{T1}{lmss}{m}{n}
\pagerange{\pageref{firstpage}--\pageref{lastpage}}
\allowdisplaybreaks
\usepackage[UKenglish]{babel} 
\newcommand{\ignore}[1]{}
\usepackage[twoside,
	paperheight=247mm,
	paperwidth=174mm,
	marginratio={3:5,2:3},
	left=18mm,
	top=20mm]{geometry}
\usepackage[hang]{subfigure}
\numberwithin{figure}{section}
\numberwithin{table}{section}
\usepackage{subfigure}
\usepackage{graphics,graphicx}
\usepackage[cmyk]{xcolor}
\usepackage{enumitem} 

\usepackage{calc,tikz,textcomp}
\usetikzlibrary{calc}
\usetikzlibrary{tikzmark}
\usetikzlibrary{positioning}

\usepackage{subcaption}
\usepackage{enumitem}
\usepackage{soul}
\usepackage{algorithm, algpseudocode, algorithmicx}
\usepackage{comment}
\usepackage{xspace}
\usepackage{mleftright}
\usepackage{booktabs}
\usepackage[nameinlink,capitalise]{cleveref}
\crefname{equation}{}{}

\crefname{definition}{Definition}{Definitions}

\makeatletter
\@ifpackageloaded{xcolor}{}{\usepackage[dvipsnames]{xcolor}}
\makeatother

\usepackage[most]{tcolorbox}

\colorlet{algcolor}{blue}
\algrenewcommand\alglinenumber[1]{\sf\scriptsize\color{algcolor}{#1}}
\algrenewcommand\algorithmiccomment[1]{\hfill{\color{algcolor}\(\triangleright\) #1}}

\renewcommand{\R}{\mathbb{R}}
\newcommand{\C}{\mathbb{C}}

\newcommand{\onevec}{\mathbf{1}}

\DeclareMathOperator{\tr}{tr}
\DeclareMathOperator{\offdiag}{offdiag}
\DeclareMathOperator{\orth}{\mathbf{Orth}}
\DeclareMathOperator{\Orth}{\mathbf{Orth}}
\DeclareMathOperator{\vectorization}{vec}
\newcommand{\mat}[1]{\boldsymbol{#1}}
\newcommand{\tensor}[1]{\boldsymbol{\mathcal{#1}}}
\renewcommand{\vec}[1]{\boldsymbol{#1}}
\newcommand{\lowrank}[1]{\mleft\llbracket #1 \mright\rrbracket}
\newcommand{\norm}[1]{\mleft\| #1 \mright\|}
\newcommand{\ip}[2]{\mleft\langle #1, \, #2\mright\rangle}
\DeclareMathOperator{\range}{range}
\DeclareMathOperator{\rank}{rank}
\newcommand{\uinorm}[1]{{\left\vert\kern-0.25ex\left\vert\kern-0.25ex\left\vert #1 
		\right\vert\kern-0.25ex\right\vert\kern-0.25ex\right\vert}}

\DeclareMathOperator{\spn}{span}
\renewcommand{\top}{\tp}

\newcommand{\expmat}[1]{\begin{bmatrix} #1 \end{bmatrix}}
\newcommand{\twobytwo}[4]{\expmat{#1 & #2 \\ #3 & #4}}
\newcommand{\twobyone}[2]{\expmat{#1 \\ #2}}
\newcommand{\onebytwo}[2]{\expmat{#1 & #2}}
\newcommand{\flatonebytwo}[2]{[\begin{matrix} #1 & #2 \end{matrix}]}
\newcommand{\flatmatrix}[1]{[\begin{matrix} #1 \end{matrix}]}

\newcommand{\Id}{\mathbf{I}}
\newcommand{\evec}{\mathbf{e}}
\newcommand{\tp}{\mathsf{T}}

\DeclareMathOperator{\diag}{diag}

\DeclareMathOperator{\Var}{Var}

\DeclareMathOperator{\expect}{\mathbb{E}}
\DeclareMathOperator{\prob}{\mathbb{P}}
\DeclareMathOperator{\Unif}{\textnormal{\textsc{Unif}}}
\DeclareMathOperator{\Normal}{\textnormal{\textsc{Normal}}}

\newcommand{\order}{\mathcal{O}}
\newcommand{\orderish}{\widetilde{\mathcal{O}}}
\DeclareMathOperator{\poly}{poly}

\DeclareMathOperator*{\argmax}{argmax}
\DeclareMathOperator*{\argmin}{argmin}
\DeclareMathOperator*{\nrank}{nrank}
\newcommand{\set}[1]{\mathsf{#1}}
\newcommand{\e}{\mathrm{e}}

\renewcommand{\hat}[1]{\widehat{#1}}
\renewcommand{\tilde}[1]{\widetilde{#1}}
\newcommand{\actionbox}[1]{\begin{tcolorbox}[colback=white,colframe=black,width=\columnwidth,boxsep=5pt,arc=4pt]
    #1
\end{tcolorbox}}

\renewcommand{\d}{\mathrm{d}}

\tcbuselibrary{skins, breakable}

\newtcolorbox{infobox}{
  enhanced,
  breakable,
  colback=gray!5,
  colframe=black,
  fonttitle=\bfseries,
  title={},
  sharp corners=south, 
  boxrule=0.4pt,
  left=6pt,
  right=6pt,
  top=6pt,
  bottom=6pt,
  borderline west={2pt}{0pt}{black},
  before skip=10pt,
  after skip=10pt
}

\usepackage{listings}
\definecolor{bg}{rgb}{0.97,0.97,0.97}
\definecolor{mlgreen}{rgb}{0.20,0.65,0.28}
\lstdefinelanguage{Matlab}{
    keywords={break, case, catch, continue, else, elseif, end, for, function,
              global, if, otherwise, persistent, return, switch, try, while},
    keywordstyle=\color{blue},
    identifierstyle=\color{black},
    commentstyle=\color{mlgreen},
    stringstyle=\color{magenta},
    backgroundcolor=\color{bg},
    morecomment=[l][\color{mlgreen}]{\%},
    morestring=[b]"
}
\lstdefinelanguage{C}{
    keywords={auto, break, case, char, const, continue, default, do, double, else, enum, extern, float, for, goto, if, inline, int, long, register, restrict, return, short, signed, sizeof, static, struct, switch, typedef, union, unsigned, void, volatile, while},
    keywordstyle=\color{blue},
    identifierstyle=\color{black},
    commentstyle=\color{mlgreen},    
    stringstyle=\color{magenta},
    backgroundcolor=\color{bg},      
    morecomment=[l][\color{mlgreen}]{//},  
    morecomment=[s][\color{mlgreen}]{/*}{*/},  
    morestring=[b]"
}
\crefname{program}{Program}{Programs}
\crefname{conjecture}{Conjecture}{Conjectures}

\newfloat{program}{htbp}{lop}[section]
\floatname{program}{Program}

\ExplSyntaxOn
\NewDocumentCommand{\mytexttt}{m}
 {
  \texttt
   {
    \tl_set:Nn \l_tmpa_tl { #1 }
    \tl_replace_all:Nen \l_tmpa_tl { \char_generate:nn { `_ } { 8 } } { \_ }
    \tl_use:N \l_tmpa_tl
   }
 }
\cs_generate_variant:Nn \tl_replace_all:Nnn { Ne }
\ExplSyntaxOff

\usepackage{epigraph}
    \usepackage{etoolbox}
    \makeatletter
    \newlength\epitextskip
    \pretocmd{\@epitext}{\em}{}{}
    \apptocmd{\@epitext}{\em}{}{}
    \patchcmd{\epigraph}{\@epitext{#1}\\}{\@epitext{#1}\\[\epitextskip]}{}{}
    \makeatother

\newcommand{\hatbold}[1]{\skew{4}\widehat{\smash{\boldsymbol{#1}}\mathstrut}}
\newcommand{\tildebold}[1]{\skew{4}\widetilde{\smash{\boldsymbol{#1}}\mathstrut}}
\newcommand{\Ahat}{\smash{\hatbold{A}}}
\newcommand{\Bhat}{\smash{\hatbold{B}}}
\newcommand{\bhat}{\smash{\widehat{\vec{b}}}}
\newcommand{\yhat}{\smash{\widehat{\vec{y}}}}
\newcommand{\normal}[2]{\Normal(#1,#2)}
\newcommand{\Atilde}{\smash{\tildebold{A}}}
\newcommand{\atilde}{\smash{\widetilde{\vec{a}}}}

\newcommand{\xtilde}{\widetilde{\boldsymbol{x}}}
\newcommand{\xhat}{\widehat{\boldsymbol{x}}}
\newcommand{\Xtilde}{\smash{\tildebold{X}}}
\newcommand{\Btilde}{\smash{\tildebold{B}}}
\newcommand{\Utilde}{\smash{\tildebold{U}}}
\newcommand{\outprod}[1]{#1^{\vphantom{\tp}}#1^\tp}
\newcommand{\Qhat}{\smash{\hatbold{Q}}}
\newcommand{\Hhat}{\smash{\hatbold{H}}}
\newcommand{\Rhat}{\smash{\hatbold{R}}}
\newcommand{\simiid}{\stackrel{\rm iid}{\sim}}
\newcommand{\pinv}{\dagger}
\newcommand{\plusplus}{\textsf{++}\xspace}
\newcommand{\RPCholesky}{RPCholesky\xspace}
\newcommand{\XTrace}{XTrace\xspace}
\newcommand{\XDiag}{XDiag\xspace}
\newcommand{\XNysDiag}{XNysDiag\xspace}
\newcommand{\XNysTrace}{XNysTrace\xspace}
\newcommand{\err}{\mathrm{err}}

\def\Nystrom{Nystr{\"o}m\xspace}

\newcommand\subsetsim{\mathrel{%
  \ooalign{\raise0.2ex\hbox{$\subset$}\cr\hidewidth\raise-0.8ex\hbox{\scalebox{0.9}{$\sim$}}\hidewidth\cr}}}

\newcommand*{\Scale}[2][4]{\scalebox{#1}{\ensuremath{#2}}}%

       \newcommand{\rowsize}{.7cm}
       \newcommand{\rrowsize}{.5cm}
       \newcommand{\colsize}{1.8cm}

\newcommand{\iset}{\set{I}}
\newcommand{\isettwo}{\set{J}}
\newcommand{\isetthree}{\set{K}}
\DeclareMathOperator{\Median}{Median}
\DeclareMathOperator{\Lip}{Lip}
\DeclareMathOperator{\vol}{vol}
\newcommand{\psd}{psd\xspace}
\newcommand{\pd}{pd\xspace}
\usetikzlibrary{decorations.pathreplacing,calligraphy}
\newcommand{\surveytopic}[1]{\textcolor{gray}{[#1]}}

\makeatletter
\renewenvironment{proof}[1][Proof]{%
  \rm\trivlist
  \item[\hskip\labelsep{\it #1.\hspace*{0.5em}}]%
}{%
  \hspace*{\fill}{$\actabox$}\endtrivlist
}
\makeatother

\usepackage{thmtools}
\declaretheorem[name=Theorem,numberwithin=section]{theorem}

\declaretheorem[name=Proposition,sibling=theorem]{proposition}
\declaretheorem[name=Corollary,sibling=theorem]{corollary}
\declaretheorem[name=Definition,sibling=theorem]{definition}
\declaretheorem[name=Remark,sibling=theorem]{remark}
\declaretheorem[name=Warning,sibling=theorem]{warning}

\title[Randomized Linear Algebra in Scientific Computing]{Randomized algorithms in numerical linear algebra for scientific computing}

\author[Epperly, Kressner, Nakatsukasa]{
Ethan N. W. Epperly${}^1$, Daniel Kressner${}^2$, and Yuji Nakatsukasa${}^3$\\\\
${}^1$ Department of Mathematics, University of California Berkeley,\\\vspace{0.3em} Berkeley, California 94720, USA E-mail: \href{mailto:eepperly@berkeley.edu}{eepperly@berkeley.edu}\\
${}^2$ Institute of Mathematics, École Polytechnique Fédérale de Lausanne (EPFL),\\\vspace{0.3em} CH-1015 Lausanne, Switzerland E-mail: \href{mailto:daniel.kressner@epfl.ch}{daniel.kressner@epfl.ch}\\
${}^3$ Mathematical Institute, University of Oxford,\\ Oxford OX2 6GG, United Kingdom E-mail: \href{mailto:yuji.nakatsukasa@maths.ox.ac.uk}{yuji.nakatsukasa@maths.ox.ac.uk}
}

\begin{document}

\label{firstpage}
\maketitle
\vspace{5em}

\begin{abstract}
Randomized numerical linear algebra (RNLA) has fundamentally 
changed the landscape
of scientific computing over the last twenty to thirty years.
While early breakthroughs focused primarily on low-rank approximation and least-squares problems, randomized techniques are now used to solve a broad range of linear algebra problems, supporting applications spanning data science, optimization, and computational engineering.
This paper surveys some of the main ideas of RNLA, focusing on methods that are effective in practice and on simple theoretical bounds that can be used to understand them. 
We highlight the core benefits of randomization: improved scalability, algorithmic simplicity, and---contrary to what the name may suggest---remarkable reliability and numerical stability. 
Topics include sketching, low-rank approximation, least squares, orthogonalization and approximate orthogonalization, randomized Krylov subspace methods, preconditioning, trace and spectral density estimation, row and column selection algorithms, posterior error estimation, randomized regularization, and tensor-structured randomness.
\end{abstract}

\tableofcontents

\section{Introduction}

Now is a time for tremendous change in computational science, and computation more generally.
The size and scale of problems that we wish to solve are larger than ever, and machine learning and artificial intelligence are fundamentally changing the computational landscape.
While much remains uncertain, what is not in doubt is that linear algebra plays a fundamental role in the next era of computing, and numerical linear algebra (NLA)---the mathematical study of algorithms for linear algebra problems---remains as relevant as ever.

Traditionally, NLA has focused on several core problems: Matrix multiplication $\mat{A}\mat{B}$, solution to systems of equations $\mat{A}\vec{x} = \vec{b}$, least-squares fitting $\min_{\vec{x}} \norm{\mat{A}\vec{x} - \vec{b}}$, eigenvalue problems $\mat{A}\vec{x} = \lambda \vec{x}$, and singular value problems.
Alongside these classical tasks, low-rank approximation, trace and spectral
density estimation, and computations involving matrix functions have become
increasingly prominent.
In all cases, the goal is to solve these problems as rapidly as possible while achieving a satisfactory level of accuracy.

Standard NLA algorithms typically fall into one of two categories: direct methods and iterative methods. 
Direct methods are designed to be exceptionally robust and reliable and function largely as black-box solvers.
These led to some of the highest-quality software packages, most prominently LAPACK~\cite{anderson1999lapack}, which stands among the greatest contributions by the NLA community. 
However, without exploiting additional structure such as sparsity, direct methods typically run in cubic time in the matrix dimension $n$, making them infeasible for large-scale problems (say $n>10^5$) on standard hardware. 
Iterative methods~\cite{saad2003iterative,saad2011numerical}, by contrast, produce a sequence of approximate solutions to a linear algebra problem.
Each iteration is typically cheap, and the goal is to design an iteration that converges to a desired level of accuracy in as few steps as possible.
When they work, iterative methods can be very fast.
But these methods are also somewhat fragile, with convergence behavior that depends sensitively on the properties of the matrix.
Often, a problem must be transformed or ``preconditioned'' to obtain favorable results, and the search for good preconditioners remains an active field of research.

NLA has been undergoing a renaissance over the last two or three decades, establishing \emph{randomized} methods as a third category, complementing both direct and iterative solvers. 
Its roots are much older; early examples include the von Neumann--Ulam Monte Carlo method for matrix inversion~\citeyear{forsythe1950matrix} and Wilkinson's recommendation~\citeyear{wilkinson1965algebraic} to  use random starting vectors in the Lanczos method.
Notwithstanding examples such as these, NLA researchers have traditionally tried to use randomness in as few places as possible.
\emph{Randomized numerical linear algebra (RNLA)} inverts this design philosophy, actively incorporating randomness in algorithm design and using probabilistic arguments to study the performance of algorithms.
As this survey highlights, randomization can be a powerful tool, potentially offering speed, flexibility, and (probabilistic) guarantees that deterministic methods cannot provide.
While randomized methods are often understood to be inherently less accurate than classical approaches, this need not be the case: 
Some randomized algorithms yield solutions of modest accuracy at lightning speed \cite{sarlos2006improved,halko2011finding}, and others can match the level of accuracy of traditional methods while still being faster \cite{epperly2026fasta}.
RNLA has been an active research area for several decades, and the community has been productive: A number of important algorithms now have a substantial
body of supporting theory, and new methods continue to emerge every year.

\subsection{What does randomization do for you?}

On its face, injecting randomness into linear algebra may seem silly.
After all, the inputs to these problems are matrices and vectors, arrays of numbers with no inherent probabilistic content.
How could randomness possibly help?
What can it do for us?
To begin to get a sense of an answer, we can examine different ways that randomness has been used productively in linear algebra computations.

\begin{enumerate}
    \item \textbf{Fast, approximate solutions.}
    Randomized algorithms can be used to give fast approximate solutions to linear algebra problems such as highly overdetermined least-squares problems (\cref{sec:sas}), low-rank matrix approximation (\cref{sec:lowrank}), and computing the distribution of eigenvalues (\cref{sec:sde}).
    In effect, such algorithms allow us to trade the accuracy and consistent behavior of traditional deterministic algorithms for speed.
    One approach to producing such fast approximate solutions is \emph{sketching}, where one reduces the size of a matrix $\mat{A} \in \R^{m\times n}$ by multiplying it against a matrix $\mat{\Omega} \in \R^{n\times s}$ with $s \ll n$, or $\mat{\Psi} \in \R^{m\times s}$ with $s\ll m$. Choosing the test matrices $\mat{\Omega}$ or $\mat{\Psi}$ suitably at random
allows us, with overwhelming probability, to ensure nonnegligible overlap with
the ``important'' directions in the row or column space of $\mat{A}$,
\emph{without prior knowledge of these directions.}
The resulting \emph{sketch} $\mat{A}\mat{\Omega}$ or $\mat{\Psi}^\tp\mat{A}$ can then be used to reveal low-dimensional structure in $\mat{A}$.
Using sketching, we can solve linear algebra problems with runtime scaling beyond any known deterministic method.
    For instance, we can compute near-optimal rank-$k$ approximations to $\mat{A}$ in roughly $\order(mn + (m+n)k^2)$ operations (\cref{sec:lowrank-generalized-nystrom}) and we can solve overdetermined linear least-squares problems to high accuracy in roughly $\order(mn+n^3)$ operations (\cref{sec:sketch-and-precondition}); compare with $\order(mnk)$ and $\order(mn^2)$ for known deterministic methods \cite{ghashami2016frequent,golub1965numerical}.
    
    \item \textbf{Turning every problem instance into an \emph{average-case} instance.}
    Many computational tasks in NLA are easy when applied to a ``generic'' (i.e., random) matrix.
    For instance, for many nice random matrix models, we can safely perform Gaussian elimination without pivoting, any set of $r$ columns spans a near-optimal rank-$r$ approximation, and all of the eigenvalues are well-separated, guaranteeing the success of eigenvalue algorithms such as spectral bisection and small block-size Lanczos methods.
    Unfortunately, as the saying goes, ``a random matrix is not a typical matrix'', and many matrices in applications look nothing like a random matrix.
    However, by injecting randomness in the algorithm, we can perturb, rotate, or transform a matrix to ``look more random'', thus enabling us to use approaches that would be unsafe otherwise.
    Sketching algorithms, discussed above, already provide one example, as the sketch $\mat{A}\mat{\Omega}$ or $\mat{\Psi}^\tp \mat{A}$ can be seen as extracting a subset of columns or rows of $\mat{A}$, after a random change of basis.
    But there are many more applications.
    In \cref{sec:regularization}, we shall see several examples of how random regularization can be used to solve, e.g., difficult eigenvalue problems. 
    
    \item \textbf{Matrix-free linear algebra and Monte Carlo estimates.}
    In many applications, we want to understand matrices $\mat{A}$ that can be accessed only via matrix--vector multiplications $\vec{z} \mapsto \mat{A}\vec{z}$.
    Examples include Jacobian and Hessian matrices in machine learning, where $\mat{A}$ is applied via automatic differentiation, and the matrix inverse $\mat{A} = \mat{M}^{-1}$ or other matrix functions $\mat{A} = f(\mat{M})$, where we can use a fast solver or iterative algorithms to apply $\mat{A}$.
    Using randomized algorithms, we can develop \emph{unbiased random estimates} of important quantities associated with $\mat{A}$, including its trace (\cref{sec:trace}), its diagonal (\cref{sec:diagonal}), and even the distribution of its eigenvalues (\cref{sec:sde}).
    These estimates are often only accurate to a few digits, but that is often enough.
    Moreover, by using variance reduction techniques (\cref{sec:variance-reduction,sec:xtrace}), the accuracy of these algorithms can often be improved by orders of magnitude.
    These types of matrix-free algorithms are, in some sense, impossible to realize efficiently in a deterministic way without prior information about $\mat{A}$:
    A deterministic choice of probe vectors can always be defeated
by a matrix whose important directions have little or no overlap with those
probes. 
    
    \item \textbf{Solving problems at extreme scale.}
    Many problems in modern scientific computing and machine learning are of such enormous scale that we cannot hope to read, store, and process all of the problem data using traditional methods \cite{meanti2020kernel,greene2022approximating}.
    Randomization offers a powerful way of solving such problems.
    Using randomized matrix approximation techniques, we can compute accurate approximations to large matrices without ever needing to read or generate all of the matrix entries (\cref{sec:colum-nystrom-linear}), and we can approximately solve least-squares regression problems $\min_{\vec{x}} \norm{\mat{A}\vec{x} - \vec{b}}$ where we never inspect all entries of the right-hand side $\vec{b}$~(\cref{sec:activelsq}).
    Using tensor-structured random embeddings, we can rapidly solve analysis tasks involving large datasets and scientific problems involving exponentially high-dimensional systems from quantum physics (\cref{sec:tensors}).
    And, taken to the extreme, RNLA ideas can be used to give principled solutions to linear algebraic problems in \emph{infinite dimensions} \cite{persson2025randomized,boulle2022generalization,boulle2023elliptic,adcock2025optimal,shustin2022semi,zvonek2025conthutch}.
\end{enumerate}

\subsection{Are randomized algorithms reliable?}

All of the benefits of randomization may seem great, but a basic concern may linger:
If an algorithm uses randomness, won't it fail some of the time?
In principle, yes.
With the exception of a few special randomized algorithms that succeed with 100\% probability, most randomized algorithms do have some probability of producing unsatisfactory results.
But this failure probability can be made negligibly small.
If one flips a coin 1000 times, one would be shocked if it came up heads every time or even 900 times.
In much the same way, the error of well-designed randomized algorithms tends to concentrate around its typical value up to a small standard deviation.
The error may occasionally be two or three standard deviations above the typical value during ordinary use, but the probability of being fifty standard deviations above is typically negligible.
For many randomized algorithms, the failure probability can be made arbitrarily small (say, $10^{-10}$ or even $10^{-100}$) by paying only a small price in efficiency.
Moreover, we need not simply trust the output of a randomized algorithm blindly, as we can use a posteriori error estimation techniques (\cref{sec:error-estimation}) to monitor the quality of the computed output at runtime.

Nonetheless, there is a real tradeoff.
Most of the randomized algorithms we discuss in this survey give approximate solutions, and the quality of the approximation will vary somewhat depending on the random seed.
By using such randomized algorithms, we trade some of the predictability and reproducibility of deterministic methods for speed. 
In many applications, this is a bargain we are happy to make.

Related to this point, let us highlight that our focus on randomized algorithms is pragmatic, not dogmatic.
While this survey is focused on randomized methods, we do not take the view that randomized methods are always or even often superior to deterministic methods. 
If a deterministic method achieves the same accuracy and runtime scaling as a randomized one, then it is natural to prefer the deterministic method.
Nonetheless, there are good reasons to believe that there are some problems for which the best deterministic algorithm is fundamentally worse than the best randomized algorithm.
Consider, for instance, the problem of computing the top eigenvalue of a symmetric positive semi-definite matrix.
Using randomness, we can simply generate a Gaussian random initial vector and obtain a satisfactory estimate from a few steps of the power method~\cite{kuczyski1992estimating,martinsson2020randomized}.
For any fixed, deterministic choice of the initial vector, there are matrices for which
the dominant eigenvector is orthogonal to this choice, causing the power method to not converge to the top eigenvalue.
To estimate the top eigenvalue deterministically, we must either give up on the power method or we must run some additional calculation to find a good initialization.
Despite decades of searching, no deterministic algorithm for top eigenvalue estimation is known that comes close to matching the runtime of the power method with random initialization.

\subsection{What is in this survey, and what is not}

This survey presents some of the most significant success stories in RNLA, in as accessible a manner as possible for a broad readership. 
We assume little background besides basic numerical linear algebra and probability theory.
For references on the former, we recommend~\cite{trefethen1997numerical,demmel1997applied,golub2013matrix}.
For the latter, we recommend \cite{grimmett2020probability,blitzstein2019introduction,tropp2023cms}.
Major topics in the survey include sketching (\cref{sec:dimension-reduction,sec:embeddings}), low-rank approximation (\cref{sec:lowrank,sec:cur,sec:column-nystrom}), (approximate) orthogonalization and speeding up Krylov subspace methods (\cref{sec:orthognalization,sec:krylov}), randomized preconditioning (\cref{sec:precond}), trace estimation and spectral density estimation (\cref{sec:trace,sec:sde}), selecting representative sets of rows and columns (\cref{sec:selection,sec:cur,sec:column-nystrom}), posterior error estimation (\cref{sec:error-estimation}), regularization by randomization (\cref{sec:regularization}), and tensor problems (\cref{sec:tensors}).

Our emphasis is on algorithms, concepts, and analyses relevant to scientific computing. While our choice of topics is influenced by applications and computational experience, we provide only limited discussion of implementations and specific scientific use cases. We have attempted to present the simplest theory we are aware of for each method, with a focus on simple bounds that are directly useful for understanding the algorithms. The theory of RNLA is deep, and many important results are mentioned only briefly or not at all.

Linear algebra is present in almost every field of numerical computation, and randomized algorithms are proliferating throughout the computing sciences.
As such, there is a very large class of topics that fit under the broad umbrella of ``randomized numerical linear algebra''.
Even with this survey's length, we only have room to cover a subset of topics.

We omit or only briefly discuss some topics already covered thoroughly in the fairly recent survey of \citet{martinsson2020randomized}, including single-vector power and Krylov iteration for extremal eigenvalues, single-pass Schatten-norm estimation, (full) rank-revealing matrix factorizations, hierarchical low-rank approximation, and graph Laplacian linear solvers.
Another major omission is the field of stochastic iterative methods, such as
stochastic gradient, randomized Kaczmarz, and randomized coordinate descent methods.
A thorough discussion of these methods could easily fill another survey of this length; we recommend \cite{strohmer2009randomized,needell2016stochastic,ma2015convergence,gower2015randomized,ferreira2025survey} as jumping off points for the interested reader.
Also, we do not systematically treat streaming models, in which matrix data
arrive sequentially and storage is restricted.
While some of the methods for low-rank approximation covered in this survey have been motivated by streaming (\cref{sec:lowrank-nystrom,sec:lowrank-generalized-nystrom}), we have chosen to not emphasize this topic; the interested reader can consult, e.g., \cite{clarkson2009numerical,ghashami2016frequent,tropp2019streaming}.
RNLA approaches have also proven effective for solving linear algebra problems at massive scale, where the dimension is either extremely large (say, $>10^{50}$) or infinite.
While some aspects of our work touch on this area (e.g., our discussion of active least-squares regression in \cref{sec:selection} and tensors in \cref{sec:tensors}), we largely do not focus on this setting.
See \citet{weare2026randomly,braverman2022sublinear,jin2023moments,greene2022approximating} for extreme-scale RNLA.

\subsection{Other surveys and expositions}

The field of RNLA is already gifted with several survey articles and other expository references.
Here we provide a partial list.
For some surveys, we highlight distinctive features in brackets.

\paragraph{Randomized numerical linear algebra: A perspective on the field with an eye to software \cite{murray2023randomized}.}
This monograph is written to provide an overview of RNLA, with a focus on core algorithms to be incorporated in the RandBLAS and RandLAPACK libraries.
The survey is particularly strong in its discussion of sketching, leverage score sampling, tensor-structured sketching, and all matters of implementation.

\citet{murray2023randomized} provide an impressively broad overview of RNLA.
Our survey provides a complementary perspective, placing additional emphasis on sketching for Krylov subspace methods (\cref{sec:krylov}), spectral density estimation (\cref{sec:sde}), row and column selection (\cref{sec:selection,sec:cur,sec:column-nystrom}), and randomized regularization (\cref{sec:regularization}).
A revision of this book is forthcoming, so this summary of differences may not reflect the final version.

\paragraph{Randomized numerical linear algebra: Foundations and algorithms \cite{martinsson2020randomized}.} This survey remains a foundational reference for modern RNLA.
It already covers a substantial portion of modern RNLA and, as mentioned above,  we have deliberately chosen to avoid extended discussion of several topics treated there. Our survey complements it by emphasizing subsequent developments,
especially orthogonalization (\cref{sec:orthognalization}), sketching for Krylov methods (\cref{sec:krylov}), variance-reduced trace estimation (\cref{sec:trace}), randomized regularization (\cref{sec:regularization}), column \Nystrom approximation (\cref{sec:column-nystrom}), and tensor-structured randomness (\cref{sec:tensors}).
As the importance of index selection methods, including CUR, has grown in the past years, we
provide a much broader discussion on these topics (\cref{sec:selection,sec:cur,sec:column-nystrom}).

\paragraph{Surveys on randomized low-rank approximation.}
Low-rank approximation is one of the defining successes of randomization in linear algebra. 
Several surveys are dedicated to this topic, including \citet{halko2011finding} \surveytopic{randomized SVD}, \citet{martinsson2018randomized} \surveytopic{column selection, adaptive stopping}, \citet{tropp2023randomized} \surveytopic{randomized block Krylov iteration}, and \citet{pearce2026randomized} \surveytopic{CP/Tucker tensor decompositions}.

\paragraph{Surveys from a computer science perspective.}
There is an alternative perspective on RNLA, especially popular in (theoretical) computer science, that places greater emphasis on approximate matrix multiplication, sketching, and index selection by norm-squared and leverage-score sampling.
Surveys on this perspective include \citet{woodruff2014sketching} \surveytopic{sketching and complexity lower bounds}, \citet{kannan2017randomized} \surveytopic{norm-squared sampling}, \citet{drineas2018lectures} \surveytopic{leverage-score sampling and structural error bounds}, and \citet{dereziski2024recent} \surveytopic{sharp analysis of sketching methods}.

\paragraph{Surveys on individual topics.}
To analyze randomized algorithms, knowledge of high-dimensional probability and (non-asymptotic) random matrix theory is fundamental.
Good references include \citet{vanhandel2014probability,vershynin2026high} \surveytopic{high-dimensional probability theory} and \citet{tropp2015introduction,tropp2026applied} \surveytopic{matrix concentration inequalities}.
\citet{derezinski2021determinantal} provide an introduction to determinantal point processes and their role in linear algebra computations.
\citet{chen2024lanczos} provides an introduction to the Lanczos algorithm, which is a fundamental ingredient in many randomized matrix algorithms.
Epperly's thesis \citeyearpar{epperly2025make} provides an overview of many aspects of RNLA, with a focus on index selection, stochastic trace estimation, and sketching for least squares.
\citet{dedamas2026randomized} discuss randomized
orthogonalization and its use in Krylov subspace methods. 

\subsection{Numerical experiments}

Throughout this survey, we shall present numerical examples to illustrate and compare different methods.
The experiments are small-scale and far from comprehensive, and we have deliberately selected examples that highlight certain properties of the algorithms we are discussing.
Notwithstanding these caveats, we believe the behavior of the algorithms in these experiments is consistent with our experience using them.
Code for the experiments is provided in
\actionbox{\centering\url{github.com/eepperly/RNLA-for-Scientific-Computing.git}}
\noindent All runtime experiments were run on the first author's MacBook Pro, which has an M3 Pro chip and 18 GB of memory.

\subsection{Notation and terminology}

Our choice of notation and terminology is largely standard (or at least not non-standard), and we comment only on potential sources of confusion.
We mostly focus on the real field $\R$ for concreteness, though essentially everything we do generalizes in a reasonably natural way to the complex numbers $\C$.
The symbol ${}^\pinv$ denotes the pseudoinverse.
The identity matrix is $\Id$, and the vector of ones is $\onevec$.
We occasionally use subscripts to indicate the dimension.

We frequently index matrices using MATLAB notation,  and we use the standard letters $\iset$ and $\isettwo$ to denote row and column index sets.
Occasionally, this will result in unfortunate situations where we encounter expressions such as $\Id(:,\set{I})$ where we index the identity matrix $\Id$ using the index set $\set{I}$.

A matrix $\mat{Q} \in \R^{m\times n}$ is \emph{orthonormal} if it satisfies $\mat{Q}^\top\mat{Q} = \Id$, and $\Orth(\mat{A})$ denotes an orthonormal matrix $\mat{Q}$ whose columns form an orthonormal basis for $\range(\mat{A})$.
We write $\lambda_i(\cdot)$ and $\sigma_i(\cdot)$ for the $i$th decreasingly ordered eigenvalue or singular value, respectively.
A \emph{rank-$r$ matrix} is a matrix of rank at most $r$, and we let $\lowrank{\mat{A}}_r$ denote any simultaneous best rank-$r$ approximation with respect to any unitarily invariant norm.
Despite the possibility of nonuniqueness, we shall generally refer to \emph{the} best rank-$r$ approximation.
We use \emph{psd} for a positive semidefinite matrix and \emph{pd} for a positive definite matrix; symmetry is tacitly assumed. 

The unadorned norm symbol $\norm{\cdot}$ denotes the $\ell_2$ norm of a vector or the spectral norm of a matrix.
For an $m\times n$ matrix $\mat{A}$ with $m\geq n$, 
the condition number is $\kappa(\mat{A}):=\sigma_1({\mat A})/\sigma_n(\mat A)$; 
when $m=n$,  $\kappa(\mat{A})=\norm{\mat{A}}\norm{\smash{\mat{A}^{-1}}}$.
The Frobenius norm is $\norm{\cdot}_{\rm F}$, the nuclear norm is $\norm{\cdot}_*$, and the Schatten $p$-norm is $\norm{\cdot}_{(p)}$.
We will occasionally use the Frobenius inner product $\ip{ \mat{A}}{\mat{B}}_{\rm F} \coloneqq \tr(\mat{A}^\top\mat{B})$.

We work frequently with random scalars, vectors, and matrices and do \emph{not} use any special notation to distinguish random quantities from non-random ones.
We use the convention that nonlinear operations bind before expectations, so that $\expect x^2 = \expect [x^2]$ not $(\expect x)^2$.
The multivariate normal distribution with mean $\vec{\mu}$ and covariance matrix $\mat{\Sigma}$ is $\Normal(\vec{\mu},\mat{\Sigma})$, and a uniform distribution on a set $\set{S}$ is denoted by $\Unif (\set{S})$.
The distribution $\Unif \{\pm1\}$ is called the \emph{Rademacher} distribution. 
The phrase ``independent and identically distributed'' carries its standard abbreviation ``iid'', and we write $i_1,i_2,\ldots \simiid \textsc{Distribution}$ to denote iid samples from $\textsc{Distribution}$.

We make use of standard asymptotic notations $\order(\cdot)$, $o(\cdot)$, $\Omega(\cdot)$, etc., which follow their standard definitions in computer science.
The symbols $\lesssim$ and $\gtrsim$ denote inequality up to an absolute constant.

Throughout this survey, we shall work with embedding or test matrices that are used to reduce the dimension of a matrix.
For consistency, all embeddings are denoted by Greek letters ($\mat{\Omega},\mat{\Psi} \in \R^{n\times s}$) and are tall and thin, by convention.

\subsection{AI usage}

This survey was written by the authors, and it reflects our personal knowledge, experience, and viewpoints doing research in RNLA.
AI tools were used in several supportive ways: to identify and correct typos and grammatical mistakes, to make global changes to terminology and notation, to generate TikZ diagrams, and to identify additional references that we may have missed.
The authors made extensive use of AI coding assistants in writing code and running numerical experiments.
The authors take full responsibility for the manuscript and its contents.

\section{Sketching and fundamental algorithms} \label{sec:dimension-reduction}

In modern scientific computing, we frequently encounter matrices of enormous size.
Yet often, it seems that the amount of ``information'' in a matrix is much less than is suggested by its nominal dimensions.
Perhaps the matrix is much taller than it is wide or wider than it is tall, or perhaps the matrix is approximately low-rank, with most of the ``energy'' of the matrix spanned by a small number of directions.
Such cases arise frequently in scientific computing. Very tall matrices arise naturally in subspace methods, fitting data with smooth functions, and model reduction (\cref{sec:krylov,sec:modred}). Low-rank matrices appear in kernel machine learning algorithms (\cref{sec:column-nystrom-applications}), tensor computations (\cref{sec:tensors}), signal processing \cite{xia2012superfast}, direct solvers for integral equations \cite{martinsson2019fast}, and numerous other areas \cite{beckermann2019bounds,udell2019why}.

When presented with such a very tall or approximately low-rank matrix, it stands to reason that we should be able to avoid directly processing the entire matrix and instead work with a compressed summary of its information.
This intuition motivates the notion of \emph{sketching}.
There are three different ways in which we will use sketching. Given a matrix
$\mat{A}\in\R^{m\times n}$, we may form
\begin{equation} \label{eq:sketching}
     \underbrace{\mat{A}\mat{\Omega}}_{\text{right sketch}}
    \in \R^{m\times s},
    \qquad
    \underbrace{\mat{\Psi}^{T}\mat{A}}_{\text{left sketch}}
    \in \R^{\ell\times n},
    \qquad
    \underbrace{\mat{\Psi}^{T}\mat{A}\mat{\Omega}}_{\text{two-sided sketch}}
    \in \R^{\ell\times s},
\end{equation}
where $\mat{\Omega}\in\R^{n\times s}$ and
$\mat{\Psi}\in\R^{m\times \ell}$ are tall-thin random matrices,
with $s\ll n$ and $\ell\ll m$, called \emph{random embeddings}, \emph{sketching matrices}, or
\emph{random test matrices}. 
For example, the sketch $\mat{A}\mat{\Omega}$ provides a compressed summary of $\mat{A}$ that is much smaller in size while containing useful information about the column space of $\mat{A}$.
By processing such sketches, we are able to obtain approximate solutions to various linear algebra problems.

We begin our discussion of randomized dimension reduction in \cref{sec:random-embeddings}, where we discuss random embeddings and how they can be used to compress a matrix.
Then we discuss two applications.
In \cref{sec:rangefinder}, we see how dimension reduction can be used to find a low-dimensional subspace that approximately contains a family of vectors, and in \cref{sec:sas}, we use dimension reduction to obtain an approximate solution to a large least-squares regression problem.

\subsection{The sketching paradigm and Gaussian embeddings} \label{sec:random-embeddings} 

The core ingredient in randomized dimension reduction is a tall-thin random
embedding. As illustrated in~\cref{eq:sketching}, \emph{sketching}---that is,
the process of applying a random embedding---can be used to reduce the number
of columns or rows of a matrix. Applied from the left to a collection of high-dimensional vectors, a suitable random embedding significantly reduces their dimension while approximately preserving
their norms. This idea is closely related to the celebrated Johnson--Lindenstrauss
lemma~\citeyear{johnson1984extensions}, see~\Cref{sec:row-norm-estimation}, and has been studied and utilized extensively
in theoretical computer science, see, e.g., \cite{achlioptas2003database,dasgupta2003elementary}.
In this survey, we take a somewhat different perspective, emphasizing the preservation of geometry on entire linear subspaces.

A good random embedding combines or samples coordinates while preserving the most 
relevant geometric information. This notion will be made more precise in \cref{sec:metric}.
One type of random embedding meeting these criteria is the Gaussian embedding.

\begin{definition}[Gaussian embedding] \label{def:gaussian}
    A Gaussian embedding is a random matrix $\mat{\Omega} \in \R^{n\times s}$
    populated with iid $\Normal(0,1/s)$ entries. 
    The parameter $s$ is called the \emph{embedding dimension}.
\end{definition}

The Gaussian embedding is perhaps the most widely used embedding in randomized linear algebra, and its mathematical properties are very well-understood \cite[\S8]{martinsson2020randomized}. However, as we will discuss in \cref{sec:desiderata},
there are good reasons to use other embeddings for certain use cases.
The full power of random embeddings, with Gaussian embeddings serving as a
canonical model, will soon become clear.
For now, let us highlight a few useful properties of Gaussian embeddings.

A first important  property of Gaussian embeddings is that they preserve squared
Euclidean and Frobenius norms in expectation. Specifically, for an
appropriately sized Gaussian embedding $\mat{\Omega}$,
\[
    \expect \norm{\mat{A}\mat{\Omega}}_{\rm F}^2
    = \norm{\mat{A}}_{\rm F}^2,
    \qquad
    \expect \norm{\mat{\Omega}^\tp\mat{A}}_{\rm F}^2
    = \norm{\mat{A}}_{\rm F}^2,
    \qquad
    \expect \norm{\mat{\Omega}^\tp \vec{a}}^2
    = \norm{\vec{a}}^2.
\]
This property is called \emph{isotropy}; we formalize it and discuss its role
for general random embeddings in \cref{sec:isotropy}.

As another consequential property, a Gaussian embedding is \emph{rotationally invariant}.
A Gaussian embedding $\mat{\Omega}$ has the same random distribution as any rotation $\mat{Q}^\tp \cdot \mat{\Omega}$ by orthogonal $\mat{Q}$.
As a consequence of this property, a right sketch $\mat{A} \mat{\Omega}$ alone carries no 
information about the co-range or right singular vectors of $\mat{A}$.
Indeed, given a singular value decomposition (SVD) $\mat{A} = \mat{U}\mat{\Sigma}\mat{V}^\tp$, rotational invariance implies that
\begin{equation*}
    \mat{A}\mat{\Omega} = \mat{U}\mat{\Sigma} \cdot (\mat{V}^\tp \mat{\Omega}) \quad \text{has the same distribution as} \quad \mat{U}\mat{\Sigma} \cdot \mat{\Omega}.
\end{equation*}
To learn about both the range and co-range of a matrix using sketching,
standard approaches either make multiple passes over \(\mat{A}\) or collect an additional left sketch.

\begin{remark}[The role of randomness]
    Right from the top, it is worth asking: What is randomness doing for us here?
    Why not compress a matrix using some fixed, non-random test matrix $\mat{\Omega}$?
    Most of the time, this will work perfectly well, but it can fail badly in some cases.
    What if we are interested in learning information about the range of $\mat{A}$ but are unfortunate enough to choose a test matrix $\mat{\Omega}$ lying in its nullspace?
    Given fixed $\mat{\Omega}$, it is always possible to find such a pathological $\mat{A}$.
    But given any fixed nonzero $\mat{A}$, a Gaussian $\mat{\Omega}$ avoids this exact pathology with probability one.
    If nearly any test matrix will do, why not pick one at random?
    Randomness allows us to hedge: a deterministic test matrix works well for most matrices $\mat{A}$, but fails on some inputs, whereas a random embedding lets us design procedures that work reliably for \emph{any} matrix $\mat{A}$, except with vanishingly small probability of failure.
\end{remark}

\subsection{Example: Randomized rangefinder} \label{sec:rangefinder}

As an initial example of randomized dimension reduction in linear algebra, consider the following task:
\emph{Given a collection of vectors $\vec{a}_1,\ldots,\vec{a}_n \in \R^m$, find a low-dimensional subspace $\set{V} \subseteq \R^m$ that nearly contains them.} 
Computing such a subspace can be useful for data analysis, as it describes a small number of degrees of freedom explaining most of the variation in the data. For centered data, this is the fundamental idea of principal component analysis (PCA, \citet{jolliffe2002principal}).

To formalize this problem, begin by assembling the vectors $\vec{a}_i$ as columns of a matrix $\mat{A} \in \R^{m\times n}$, and let the columns of $\mat{Q} \in \R^{m\times s}$ constitute an orthonormal basis for an $s$-dimensional subspace $\set{V}$.
We wish to choose $\mat{Q}$ so that the original data $\vec{a}_i$ are close to their orthogonal projections $\mat{Q}\mat{Q}^\tp \vec{a}_i$ onto the subspace $\set{V}$.
That is, in matrix notation, we wish to find an orthonormal matrix $\mat{Q}$ for which $\mat{A} \approx \mat{Q}\mat{Q}^\tp \mat{A}$.

An optimal subspace $\set{V}$ is spanned by the $s$ dominant left singular vectors of $\mat{A}$.
But computing these singular vectors can be expensive, particularly using classical direct methods, which require $\order(mn \min \{m,n\})$ operations.
To go faster, we can employ randomized dimension reduction. The
\emph{randomized rangefinder} in~\cref{alg:rangefinder} constructs the sketch
$
    \mat{Y}=\mat{A}\mat{\Omega}\in\R^{m\times s}
$
using a random embedding $\mat{\Omega}$, with the aim of capturing the range
of $\mat{A}$ by the range of $\mat{Y}$. A thin QR factorization is used 
to compute an orthonormal basis $\mat{Q}$ of the latter.
\begin{algorithm}[t]
\caption{Randomized rangefinder} \label{alg:rangefinder}
\begin{algorithmic}[1]
\Require Matrix $\mat{A} \in \R^{m \times n}$; sketch size $s \ll m,n$
\Ensure Orthonormal basis $\mat{Q} \in \R^{m \times s}$ for an approximate range of $\mat{A}$
\State Generate a random embedding $\mat{\Omega} \in \R^{n\times s}$
\State Form the sketch $\mat{Y} \gets \mat{A}\mat{\Omega} \in \R^{m\times s}$
\State Orthonormalize by QR factorization $(\mat{Q},\sim) \gets \Call{QR}{\mat{Y},\texttt{'econ'}}$
\end{algorithmic}
\end{algorithm}The randomized rangefinder was isolated and systematized by \citet*{halko2011finding}.
Closely related ideas appeared earlier in~\cite{papadimitriou1998latent,sarlos2006improved,liberty2007randomized,woolfe2008fast}
and in subsequent work by~\citet{martinsson2011randomized}.

\paragraph{Uses.}
The randomized rangefinder can be used to compress a data set.
Instead of storing the full $m$-dimensional vectors $\vec{a}_1,\ldots,\vec{a}_n$, we store the $s$-dimensional embeddings $\atilde_1 \coloneqq \mat{Q}^\tp \vec{a}_1,\ldots,\atilde_n =  \mat{Q}^\tp \vec{a}_n$.
One performs data analysis directly on the reduced data points $\atilde_i$.
Or, as needed, the full-dimensional data points can be approximately regenerated using the relation $\vec{a}_i \approx \mat{Q}\atilde_i$.

The randomized rangefinder is an important subroutine in randomized algorithms for low-rank approximation, which we will discuss in \cref{sec:lowrank}.

\paragraph{Analysis.}
The performance of the randomized rangefinder is quantified by the following bound due to \citet[Thm.~10.5]{halko2011finding}.

\begin{theorem}[Randomized rangefinder: Gaussian] \label{thm:rangefinder}
    Let $\mat{A} \in \R^{m\times n}$ be a matrix, let $\mat{\Omega} \in \R^{n\times s}$ with $2 \le s \le \min\{m,n\}$
    be a Gaussian embedding, and let $\mat{Q}$ be the output of the randomized rangefinder (\Cref{alg:rangefinder}). Then
    \begin{equation*}
        \expect \norm{\mat{A} - \mat{Q}\mat{Q}^\tp \mat{A}}_{\rm F}^2 \le \min_{0 \le r \le s-2} \left( 1 + \frac{r}{s - r - 1} \right) \norm{\mat{A} - \outprod{\mat{U}_r} \mat{A}}_{\rm F}^2,
    \end{equation*}
     where $\mat{U}_r \in \R^{m \times r}$ contains the $r$ dominant left singular vectors of $\mat{A}$.
\end{theorem}

This result shows that, on average, the $s$-dimensional subspace identified by the randomized rangefinder is comparable in quality to the optimal $r$-dimensional subspace, provided that $s$ is somewhat bigger than $r$.
In particular, to satisfy the assumptions of the theorem, the rangefinder must be implemented with at least two columns of oversampling, that is, $s \ge r+2$.
With $s = r+2$, the root-mean-square Frobenius norm error is within a factor
$\sqrt{r+1}$ of the optimal rank-$r$ error. With $s=2r+1$, this factor
improves to $\sqrt{2}$.
A proof of this result will appear in \cref{sec:analysisrandsvd}.

\begin{remark}[High-probability bounds]
    There are two natural ways of quantifying the error of a randomized approximation:
    bounds on the \emph{expected error}, such as \cref{thm:rangefinder}, or error bounds that hold except for some failure probability.
    These two types of error bounds often require different techniques to prove, and consequently algorithms sometimes come equipped with one type of bound and not the other. 
    Fortunately, many well-designed randomized linear algebra algorithms have errors that concentrate, being within a small factor of their expected value except with vanishingly small probability.
    For this reason, expectation or median bounds are usually sufficient to justify an algorithm, while high-probability bounds serve to provide extra confidence.
    A posteriori error estimators (\cref{sec:error-estimation}) can always be used to certify the quality of the computed output at runtime.
    We will make a point of highlighting any situations for which an algorithm is not known to concentrate. A notorious example is provided by algorithms utilizing high-dimensional rank-one random vectors; see~\Cref{sec:rankone}.
    For those interested, high-probability error bounds for rangefinders are given in \cite[Thm.~10.7]{halko2011finding} and \cite[Thm.~8.1]{tropp2023randomized}.
\end{remark}

\paragraph{Parameter selection.}

\Cref{thm:rangefinder} tells us how we can set the embedding dimension $s$ to produce a subspace that nearly captures the same information as the best $r$-dimensional subspace.
In practice, however, the target dimension $r$ is often not known in advance.
Instead, we either choose a fixed dimension $s$ (``I want to compress my data into 100 dimensions'') or choose the dimension $s$ adaptively to meet a tolerance $\tau$ (``I want the smallest $s$ for which $\norm{\mat{A} - \mat{Q}\mat{Q}^\tp \mat{A}}_{\rm F} \le \tau$ or $\norm{\vec{a}_i - \mat{Q}\mat{Q}^\tp \vec{a}_i} \le \tau$ for all $i$'').
To support these workflows, we generally desire \emph{a posteriori} bounds on the error.
These are discussed in \cref{sec:error-estimation}, and adaptive determination of the rank is discussed in \cref{sec:rankest}.

\paragraph{Computational cost.}
For a dense unstructured matrix $\mat{A}$, the randomized rangefinder expends $\order(mns)$ operations with a Gaussian embedding, which can be faster than the $\order(mn \min \{m,n\})$ cost of computing the optimal subspace by a dense SVD.
The computational cost can be reduced to roughly $\order(mn + ms^2)$ operations by using a structured random embedding; see \cref{sec:embeddings}.

\subsection{Example: Sketch-and-solve}\label{sec:sas}

As another example of how randomized dimension reduction can help in solving a linear algebra problem, we consider the linear least-squares problem
\begin{equation} \label{eq:least-squares}
    \min_{\mat{X} \in \R^{n\times p}} \norm{\mat{A}\mat{X}-\mat{B}}_{\rm F}^2 \quad \text{for } \mat{A} \in\R^{m\times n} \text{ and }\mat{B} \in \R^{m\times p}.
\end{equation}
We consider the highly overdetermined setting where $m \gg n$.
The case $p=1$ is common, for which the right-hand side $\mat{B} = \vec{b} \in \R^m$ becomes a single vector and $\norm{\cdot}_{\rm F} = \norm{\cdot}$ is the $\ell_2$ norm. For the sake of generality and later purposes, we will consider the general case featuring a matrix right-hand side.

Since $\mat{A}$ has many more rows than columns, it is natural to surmise that many of its rows contribute little.
Randomized dimension reduction allows us to compress this least-squares problem to a smaller size, at the cost of producing an approximate solution to the original problem \cref{eq:least-squares}.
Sketch-and-solve~\cite{drineas2006sampling,sarlos2006improved}, see 
\cref{alg:sketch-and-solve}, carries out this idea by first forming the sketches $\Atilde \coloneqq \mat{\Omega}^\tp \mat{A} \in \R^{s\times n}$ and $\Btilde \coloneqq \mat{\Omega}^\tp \mat{B} \in \R^{s\times p}$ for a random embedding $\mat{\Omega} \in \R^{m\times s}$. It then obtains an approximate least-squares solution using these sketches:
    \begin{equation*}
        \Xtilde \coloneqq \argmin_{\mat{Z} \in \R^{n\times p}} \norm{\Atilde \mat{Z}-\Btilde}_{\rm F}^2 = \Atilde^\pinv \Btilde.
    \end{equation*}
\begin{algorithm}[t]
\caption{Sketch-and-solve least squares} \label{alg:sketch-and-solve}
\begin{algorithmic}[1]
\Require Matrices $\mat{A} \in \R^{m\times n}$, $\mat{B} \in \R^{m\times p}$; sketch size $s$ with ($n \le s \ll m$)
\Ensure Approximate least-squares solution $\Xtilde \in \R^{n\times p}$
\State Generate a random embedding $\mat{\Omega} \in \R^{m\times s}$
\State Form the sketches $\Atilde \gets \mat{\Omega}^\tp\mat{A} \in \R^{s\times n}$ and $\Btilde \coloneqq \mat{\Omega}^\tp\mat{B} \in \R^{s\times p}$
\State Factorize the sketch $(\mat{Q},\mat{R}) \gets \Call{QR}{\Atilde,\texttt{'econ'}}$
\State Form the sketch-and-solve solution $\Xtilde \gets \mat{R}^{-1}(\mat{Q}^\tp \Btilde)$
\end{algorithmic}
\end{algorithm}Pictorially, 
\begin{equation}\label{eq:saspicture}
\min_{\mat{X}} \left\|\    \begin{tikzpicture}[node distance=1mm,baseline=-1mm]
       \node[draw, rectangle, solid, fill=gray!10, minimum width=\rrowsize*1.5, minimum height=\colsize] (A) {$\Scale[1]{\mat{A}}$};
   \end{tikzpicture}\hspace{0.3mm}
   \begin{tikzpicture}[node distance=1mm,baseline=-1mm]
       \node[draw, rectangle, solid, fill=gray!10, minimum width = 0.3cm, minimum height=\rowsize] (x) {$\Scale[.8]{\mat{X}}$};
   \end{tikzpicture}
-  \begin{tikzpicture}[node distance=1mm,baseline=-1mm]
       \node[draw, rectangle, solid, fill=gray!10, minimum width=0.3cm, minimum height=\colsize] (b) {$\Scale[.8]{B}$};
   \end{tikzpicture}
\ \right\|^2_{\rm F}
\quad \Rightarrow \quad 
\min_{\Xtilde}\left\|\    \begin{tikzpicture}[node distance=1mm,baseline=-1mm]
       \node[draw, rectangle, solid, fill=gray!10, minimum width=\rrowsize, minimum height=\colsize/2] (SA) {$\Scale[1]{\widetilde{\mat{A}}}$};
   \end{tikzpicture}\hspace{0.3mm}
   \begin{tikzpicture}[node distance=1mm,baseline=-1mm]
       \node[draw, rectangle, solid, fill=gray!10, minimum width = 0.3cm, minimum height=\rowsize] (x) {$\Scale[.7]{\widetilde{\mat{X}}}$};
   \end{tikzpicture}
-  \begin{tikzpicture}[node distance=1mm,baseline=-1mm]
       \node[draw, rectangle, solid, fill=gray!10, minimum width=0.3cm, minimum height=\colsize/2] (Sb) {$\Scale[.6]{\widetilde{\mat{B}}}$};
   \end{tikzpicture}
\ \right\|^2_{\rm F}. 
\end{equation}
Assuming that the sketched matrix $\Atilde$ has full column rank, the solution 
can be computed as 
$\Xtilde = \mat{R}^{-1}(\mat{Q}^\tp \Btilde)$ using a thin QR factorization 
$\Atilde = \mat Q \mat R$.

\paragraph{Analysis.}
The expected squared residual for a Gaussian embedding is known exactly.

\begin{theorem}[Sketch-and-solve: Gaussian] \label{thm:sketch-solve-gaussian}
    Let $\mat{A} \in \R^{m\times n}$ have full column rank, let $\mat{B} \in \R^{m\times p}$, and let $\mat{\Omega} \in \R^{m \times s}$ be a Gaussian embedding.
    Then the sketch-and-solve solution $\Xtilde = (\mat{\Omega}^\tp \mat{A})^\pinv (\mat{\Omega}^\tp \mat{B})$ satisfies
    \begin{equation*}
        \expect \norm{\mat{A}\Xtilde-\mat{B}}_{\rm F}^2 = \left(1 + \frac{n}{s-n-1}\right) \min_{\mat{X} \in \R^{n\times p}} \norm{\mat{A}\mat{X}-\mat{B}}_{\rm F}^2
    \end{equation*}
    for $s > n+1$.
    If $s\in\{n,n+1\}$ and $\range(\mat{B}) \not\subseteq \range(\mat{A})$, the squared expected value is infinite.
\end{theorem}

\Cref{thm:sketch-solve-gaussian} appears as \cite[Lem.~2.1]{bartan2020distributed}; see also earlier work by \citet[Eq.~(3)]{ahfock2021statistical}.
This theorem certifies that sketch-and-solve produces an expected squared residual that is within a small constant factor of the minimal squared least-squares residual, provided that the embedding dimension $s$ is modestly larger than the width $n$ of $\mat{A}$.
In particular, 
\begin{equation*}
    \expect \norm{\mat{A}\Xtilde - \mat{B}}_{\rm F}^2 = 2 \min_{\mat{X} \in \R^{n\times p} } \norm{\mat{A}\mat{X} - \mat{B}}_{\rm F}^2 \quad \text{for } s =2n+1.
\end{equation*}
In particular, if the least-squares problem is nearly consistent 
then the sketch-and-solve solution can be expected to have a tiny residual.

\paragraph{Level of accuracy and parameter selection.}
To achieve an expected squared residual that is within a factor $1+\varepsilon$ of optimal,
\begin{equation} \label{eq:ls-1+eps}
    \expect \norm{\mat{A}\Xtilde-\mat{B}}_{\rm F}^2 \le \left(1 + \varepsilon \right) \min_{\mat{X} \in \R^{n\times p}} \norm{\mat{A}\mat{X}-\mat{B}}_{\rm F}^2,
\end{equation}
requires embedding dimension $s \ge n/\varepsilon + n + 1$.
The rate of convergence is \emph{algebraic} in $s$, making it punishing to obtain a small suboptimality parameter $\varepsilon$ using sketch-and-solve.
Thus, the sketch-and-solve algorithm is primarily used with an embedding dimension between $s=2n$ and $s=4n$, resulting in accuracy parameters $\varepsilon$ asymptotically in the range $[1/3,1]$. 
To solve least-squares problems to higher accuracy, the sketch-and-precondition method from~\cref{sec:precond} may be preferable.

Another way of quantifying the error of the sketch-and-solve solution $\Xtilde$ is its distance $\norm{\Xtilde - \mat{X}_\star}_{\rm F}$ to the optimal least-squares solution $\mat{X}_\star$.
To this end, we note that~\cref{eq:ls-1+eps} and Pythagoras imply
\begin{equation*}
    \expect \norm{\smash{\Xtilde - \mat{X}_\star}}_{\rm F}
    \le \big( \expect \norm{\smash{\Xtilde - \mat{X}_\star}}^2_{\rm F} \big)^{1/2} \le  \kappa(\mat{A}) \cdot \sqrt{\varepsilon} \cdot \frac{\norm{\mat{A}\mat{X}_\star - \mat{B}}_{\rm F}}{\norm{\mat{A}}},
\end{equation*}
where $\kappa(\mat{A}) = \sigma_{\rm max}(\mat{A}) / \sigma_{\rm min}(\mat{A})$ denotes the condition number; see \cite[Lem.~1]{drineas2011faster}.
Since $\varepsilon$ cannot be efficiently made small with sketch-and-solve, this bound guarantees small solution error when the optimal residual is
small. The distinction between residual and solution accuracy has been
studied in depth by~\citet{pilanci2016iterative}; see also~\cite{sridhar2020lower}.
We will revisit the analysis of sketch-and-solve in \cref{subsec:sketchsolve_injection}. 

\actionbox{\textbf{Maxim.} Sketch-and-solve should generally be used in applications where achieving a squared residual within a factor two of optimal is sufficient.
In particular, in contexts where the distance to optimal solution $\norm{\Xtilde - \mat{X}_\star}_{\rm F}$ is important, sketch-and-solve should only be used when the least-squares residual $\norm{\mat{A}\mat{X}_\star-\mat{B}}_{\rm F}$ is small and $\kappa(\mat{A})$ is modest.}

\paragraph{Computational cost.}
With the generally recommended parameter setting $2n \le s \le 4n$, the cost of the sketch-and-solve algorithm with a Gaussian embedding is $\order(mn^2)$ operations for a single right-hand side  ($p=1$), matching the asymptotic cost of standard direct solvers for least-squares problems such as the normal equations or (Householder-based) QR factorization.
Despite its cost being comparable with direct deterministic methods, sketch-and-solve with a Gaussian embedding can still be useful in distributed or out-of-core computing environments or as a subroutine in a randomized algorithm for low-rank approximation (\cref{sec:lowrank}).
With a structured random embedding (\cref{sec:embeddings}), the cost of sketch-and-solve can be brought to roughly $\order(mn + n^3)$ operations, much faster than standard direct solvers when $1 \ll n \ll m$.

\section{Random embeddings: Desiderata and constructions}
\label{sec:embeddings}

In the previous section, we saw two algorithms---the randomized rangefinder and sketch-and-solve---which use sketching to approximately solve a linear algebra problem.
In this section, we discuss the choice of \emph{random embedding matrix} $\mat{\Omega}$ to use in these and other algorithms.

This section is roughly divided into three topics.
First, in \cref{sec:metric,sec:desiderata}, we discuss what properties we want from a random embedding to be useful in algorithms.
Then, in \cref{sec:dense-sketching,sec:sparse-embeddings,sec:laso,sec:srfts,sec:tensor-embeddings-intro,sec:index-sketching-intro}, we discuss constructions for embeddings and evaluate whether they meet these criteria.
Finally, in \cref{sec:which-embedding}, we summarize and provide recommendations on which embedding to use.

\subsection{Desiderata for random embeddings 1: Metric properties}

\label{sec:metric}

For a random embedding to be useful, it should respect features of the data being reduced.
In particular, it is essential that embeddings preserve the lengths of ``important directions'' in the row or column space of a matrix being reduced.

\subsubsection{Isotropy} \label{sec:isotropy}
The most elementary metric property for a random embedding is \emph{isotropy}.
\begin{definition}[Isotropic embedding] \label{def:isotropic}
    A random embedding $\mat{\Omega} \in \R^{n\times s}$ is \emph{isotropic} if it preserves the squared Euclidean norm of any vector, in expectation:
    \begin{equation*}
        \expect \norm{\smash{\mat{\Omega}^\tp\vec{x}}}^2 = \norm{\vec{x}}^2 \quad \text{for any } \vec{x} \in \R^n.
    \end{equation*}
\end{definition}
Equivalently, $\expect[ \mat{\Omega} \mat{\Omega}^\tp ] = \Id$.
Essentially every random embedding in common use satisfies the isotropy property, perhaps after scaling by an appropriate constant.
In one sense, isotropy is a strong property: it requires norm preservation for \emph{every}
vector, on average.  At the same time, it is
a weak probabilistic guarantee: it provides no control on the fluctuations for
an individual vector, let alone on the maximum and minimum dilation over a
\emph{family} of vectors. For example, a Gaussian vector (Gaussian embedding with $s=1$) is isotropic, but $s=1$ is usually not enough to acquire enough information.

\subsubsection{How does a Gaussian random embedding act on a subspace?}

To understand how random embeddings behave in linear algebra computations, we need to understand how the embedding acts on \textit{all} vectors in a \textit{linear subspace}.
Let us begin investigating this topic by first studying the special case of Gaussian embeddings.
Let $\set{V}$ denote an $r$-dimensional linear subspace, let $\mat{Q} \in \R^{n\times r}$ be an orthonormal matrix whose columns form a basis for $\set{V}$, and let $\mat{\Omega} \in \R^{n\times s}$ be a Gaussian embedding.

\begin{figure}
    \centering
    \includegraphics[width=0.99\linewidth]{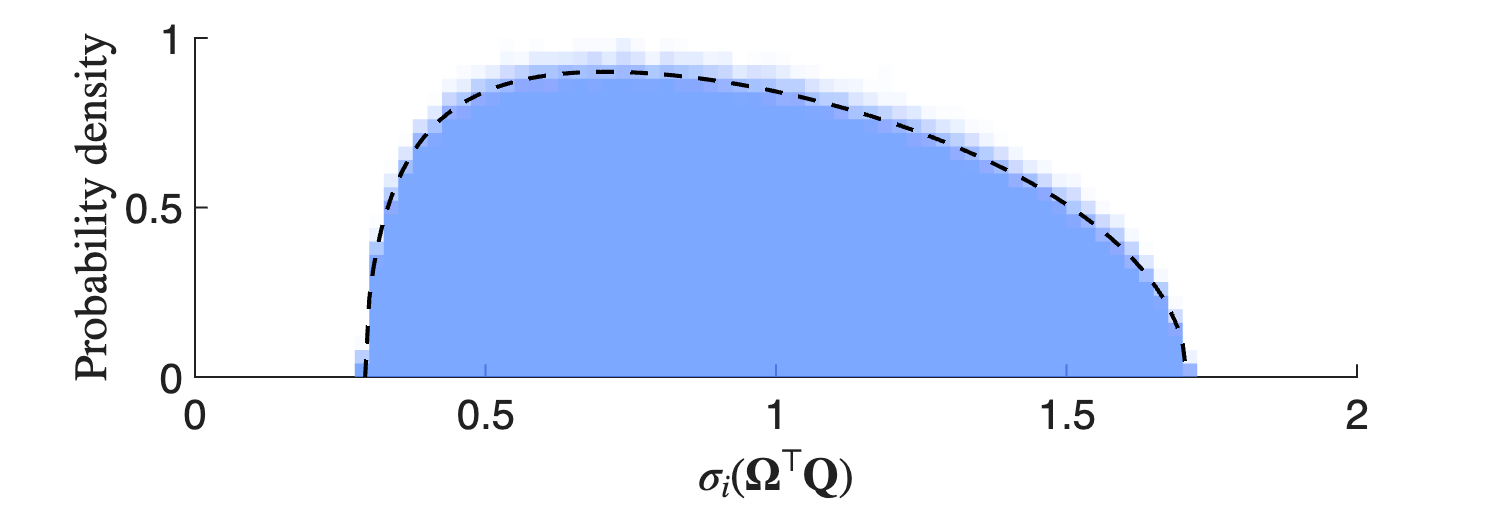}
    \caption{Histogram of singular values for product $\mat{\Omega}^\tp \mat{Q} \in \R^{s\times r}$ of a Gaussian embedding and an orthonormal matrix (100 superimposed trials).
    The subspace dimension is $r=1000$, and the embedding dimension is $s=2r=2000$.
    Across all trials, we see excellent agreement with the Marčenko–Pastur distribution, shown as a black dashed line.}
    \label{fig:mp}
\end{figure}

How much does the sketching map $\mat{\Omega}^\tp : \set{V} \to \R^s$ distort the sizes of vectors in the subspace $\set{V}$?
To answer this question, we can compute the singular values of the matrix $\mat{\Omega}^\tp \mat{Q}$.
The maximum stretch factor is the largest singular value, and the minimum stretch factor is the smallest singular value:
\begin{equation*}
    \max_{\vec{v} \in \set{V}\setminus\{\vec 0\}} \frac{\norm{\smash{\mat{\Omega}^\tp \vec{v}}}}{\norm{\vec{v}}} = \sigma_{\rm max}(\mat{\Omega}^\tp \mat{Q}) \quad \text{and} \quad \min_{\vec{v} \in \set{V}\setminus\{\vec 0\}} \frac{\norm{\smash{\mat{\Omega}^\tp \vec{v}}}}{\norm{\vec{v}}} = \sigma_{\rm min}(\mat{\Omega}^\tp \mat{Q})
\end{equation*}
The distribution of singular values is also informative.
In particular, the more tightly the singular values cluster around \(1\), the closer the embedding is to an isometry on the subspace.

For Gaussian embeddings, we can understand the distribution of singular values using random matrix theory.
When the dimensions $r$ and $s$ grow large with limiting ratio $r/s \to \varrho \in (0,1)$, the distribution of singular values approaches a continuous distribution known as the Marčenko–Pastur law \cite{marenko1967distribution}. 
It has a continuous density
\begin{equation*}
    f_{\rm mp}(\sigma) = \frac{\sqrt{(\sigma_+^2 - \sigma^2)(\sigma^2 - \sigma_-^2)}}{\pi \varrho \, \sigma} \quad \text{for } \sigma \in [\sigma_-,\sigma_+],
\end{equation*}
and it is supported on the interval between $\sigma_- = 1 - \sqrt{\varrho}$ and $\sigma_+ = 1+\sqrt{\varrho}$.
In addition, the largest and smallest singular values converge (almost surely) to the edges of this distribution as $r,s \to\infty$ \cite{bai1993limit}:
\begin{equation*} 
    \sigma_{\rm max}(\mat{\Omega}^\tp \mat{Q}) \to 1 + \sqrt{\varrho} \quad \text{and} \quad \sigma_{\rm min}(\mat{\Omega}^\tp \mat{Q}) \to 1 - \sqrt{\varrho}.
\end{equation*}
See \cref{fig:mp} for an illustration, and see \cref{thm:gaussian-ose} below.
For the practitioner, the asymptotic conclusion is that we should expect a Gaussian embedding of dimension $s$ to preserve the lengths of all vectors in an $r$-dimensional subspace up to a $1 \pm \sqrt{r/s}$ multiplicative factor.

\subsubsection{Universality}

What about other embeddings?
How do they act on a linear subspace?
Is the Marčenko–Pastur law specific to Gaussian embeddings,
or does similar behavior persist for other random embeddings?

Fortunately, a large class of random embeddings is observed to exhibit \emph{universal}
behavior, with important geometric properties matching those of Gaussian
embeddings~\cite{oymak2018universality,dobriban2019asymptotics,dereziski2023algorithmic}. For such embeddings,
the distribution of singular values of the matrix $\mat{\Omega}^\tp \mat{Q} \in \R^{s\times r}$ tends to follow the Marčenko–Pastur law, and the maximum and minimum stretch factors are approximately $1 \pm \sqrt{r/s}$.
Moreover, when we use these embeddings in algorithms such as the randomized rangefinder or sketch-and-solve, the accuracy of the resulting approximations often closely matches predictions from Gaussian theory (\cref{thm:sketch-solve-gaussian,thm:rangefinder}); see \cite{epperly2026sharp} for experiments.

For the practitioner, this universality phenomenon is very appealing.
When using an embedding with Gaussian-like behavior, we can 
achieve essentially the same accuracy as with a Gaussian embedding, potentially
at significantly reduced computational cost (\cref{sec:fast-embedding}).
Moreover, sharp theory for Gaussian embeddings (e.g., \cref{thm:rangefinder,thm:sketch-solve-gaussian}) can provide useful predictions for the behavior of these other embeddings in practice and guide the choice of parameters. We therefore generally recommend using random embeddings that exhibit Gaussian-like  universal behavior when implementing randomized algorithms.

We should emphasize that universality is not a law of nature; some types of random embeddings exhibit near-Gaussian behavior, while others do not.
And aside from a few settings where the universality principle has been justified rigorously \cite{oymak2018universality,dobriban2019asymptotics,dereziski2023algorithmic,chenakkod2024optimal}, Gaussian-like behavior for highly structured random embeddings remains largely an empirical observation, albeit one supported by a substantial body of evidence \cite{dong2023simpler,melnichenko2025choleskyqr,camao2025faster,epperly2025make,chen2026gpu,epperly2026sharp}.
Be careful!
Some structured embeddings, such as the SRFT with one round of randomization (\cref{sec:srfts}) and the CountSketch embedding (\cref{sec:sparse-embeddings}) can behave systematically differently (and worse) from Gaussian embeddings, whereas other embeddings can behave (slightly) better \cite{dobriban2019asymptotics,lacotte2020optimal,epperly2026sharp}.

\subsubsection{Subspace embeddings and injections.} \label{sec:embedinject}

Universality offers an appealing way to understand the behavior of a large class of random embeddings in algorithms.
Nonetheless, it is also worth developing analytical approaches to analyzing randomized algorithms that require weaker conditions on the embedding $\mat{\Omega}$.
This motivates the following definition.

\begin{definition}[Subspace embedding and injection] \label{def:embedinject}
    Let $\set{V} \subseteq \R^n$ be a linear subspace.
    A matrix $\mat{\Omega} \in \R^{n\times s}$ is said to be a \emph{subspace embedding} with \emph{distortion $\eta \ge 0$} for $\set{V}$ if
    \begin{equation} \label{eq:subpace-embedding}
        (1+\eta)^{-1} \norm{\vec{v}}^2 \le \norm{\smash{\mat{\Omega}^\tp\vec{v}}}^2 \le (1+\eta) \norm{\vec{v}}^2 \quad \text{for any } \vec{v} \in \set{V}.
    \end{equation}
    If only the lower bound holds, $\mat{\Omega}$ is called a \emph{subspace injection}.
\end{definition}
A subspace embedding simultaneously preserves the squared norms of all vectors in a subspace, up to scaling by a factor of $1+\eta$.
The subspace injection property is weaker, only requiring that the embedding not come close to annihilating any vector in the subspace.
The subspace embedding property is essentially due to \citet{sarlos2006improved}; the subspace injection property, as such, was introduced by \citet{camao2025faster}, with conceptually close precedents appearing earlier in, e.g., \cite{drineas2011faster}.

For most uses of sketching, the subspace injection property---together with isotropy---is sufficient to justify the performance of algorithms \cite[\S2]{camao2025faster}.
The stronger subspace embedding property becomes important when one needs two-sided, uniform control of the distortion over an entire subspace; notable examples include approximate orthogonalization (\cref{sec:orthognalization}) and preconditioning (\cref{sec:sketch-and-precondition}).

\begin{warning}[Competing definitions]
    Subspace embeddings are defined in several ways in the literature.
    Popular variations include not squaring the norms, using a factor of $1-\eta$ instead of $(1+\eta)^{-1}$ in the lower bound or allowing the $\eta$ parameters to differ for the upper and lower bounds.
    Fortunately, when $\eta \in [0,1/2]$, each of these conventions is equivalent up to rescaling $\eta$ by a small constant.
\end{warning}

For a Gaussian embedding, we can control the distortion $\eta$ by bounding the probability that a stray singular value escapes from the bounds of the Marčenko–Pastur region.
Here is a non-asymptotic result; see \cite[\S8.7.2]{martinsson2020randomized} and \cite[Thm.~II.13]{davidson2001local}.

\begin{theorem}[Gaussian embeddings are subspace embeddings] \label{thm:gaussian-ose}
    Fix a subspace $\set{V} \subseteq \R^n$ of dimension $r$, and let $\mat{\Omega} \in \R^{n\times s}$ be a Gaussian embedding.
    Then
    \begin{equation*}
        \max_{\vec{v} \in \set{V}\setminus\{\vec 0\}} \left| \frac{\norm{\mat{\Omega}^\tp \vec{v}}}{\norm{\vec{v}}} - 1\right| \le \sqrt{\frac{r}{s}} + \beta \quad \text{with probability at least } 1-2\e^{-\beta^2 s/2}.
    \end{equation*}
    In particular, selecting $s = \order(r/\eta^2)$ ensures that, with 99\% probability,  $\mat{\Omega}$ acts as a subspace embedding for $\set{V}$.
\end{theorem}

By universality, we expect that other dense random embeddings should 
also be subspace embeddings when $s = \order(r/\eta^2)$. This includes, for example, Rademacher embeddings (\Cref{sec:dense-sketching}). Subspace embedding guarantees for
other, more structured embeddings are more subtle and discussed below.

\subsubsection{Using the subspace embedding and injection properties}\label{subsec:sketchsolve_injection}

The subspace embedding and injection properties are powerful because they allow one to decouple the analysis of a randomized algorithm into two stages.
First, one proves the algorithm works when implemented with an embedding satisfying the subspace embedding or subspace injection property.
Then, one proves that the embedding used in the algorithm satisfies the necessary property.

Here, we develop a representative result for the sketch-and-solve algorithm; the randomized rangefinder will be treated in~\cref{sec:lowrank-error-analysis}.

\begin{theorem}[Sketch-and-solve: Subspace embedding or injection] \label{thm:sketch-and-solve-ose-osi}
    Let $\mat{\Omega} \in \R^{m\times s}$ be a matrix or random matrix, let $\mat{A} \in \R^{m\times n}$ be full-rank, let $\mat{B} \in \R^{m\times p}$, and $\Xtilde \coloneqq (\mat{\Omega}^\tp \mat{A})^\dagger (\mat{\Omega}^\tp\mat{B})$ be the sketch-and-solve solution.
    The following hold:
    \begin{enumerate}[label=(\alph*)]
        \item \label{item:sas-ose} \textbf{Subspace embedding.} Suppose $\mat{\Omega}$ is a subspace embedding for the range of the augmented matrix $\flatonebytwo{\mat{A}}{\mat{B}}$ with distortion $\eta$.
        Then 
        \begin{equation} \label{eq:ose-sketch-solve}
            \norm{\mat{A}\Xtilde-\mat{B}}_{\rm F}^2 \le (1+\eta)^2 \cdot \min_{\mat{X} \in \R^{n\times p}} \norm{\mat{A}\mat{X}-\mat{B}}_{\rm F}^2.
        \end{equation}
        
        \item \label{item:sas-osi-uniform} \textbf{Subspace injection, uniform bound.}
        Let $\mat{Q} = \orth(\mat{A})$, and suppose that $\mat{\Omega}$ is a subspace injection for $\range(\mat{A})$ with distortion $\eta$.
        Then
        \begin{equation} \label{eq:sketch-solve-universal2}
            \begin{split}
            \norm{\smash{\mat{A}\widetilde{\mat{X}}-\mat{B}}}_{\rm F}^2 
            &\le \frac{\norm{\mat{\Omega}}^2}{\sigma_{\rm min}^2(\mat{\Omega}^\top\mat{Q})}\cdot \min_{\mat{X} \in \R^{n\times p}}
    \norm{\mat{A}\mat{X}-\mat{B}}_{\rm F}^2 \\
            &\leq
    \norm{\mat{\Omega}}^2(1+\eta)
    \cdot \min_{\mat{X} \in \R^{n\times p}}
    \norm{\mat{A}\mat{X}-\mat{B}}_{\rm F}^2.
            \end{split}
        \end{equation}
        \item \label{item:sas-osi-probability} \textbf{Subspace injection, probabilistic bound.} Suppose that $\mat{\Omega}$ is an isotropic random embedding that, with probability $3/4$, satisfies the subspace injection property for $\range(\mat{A})$ with injectivity parameter $\eta$.
        Then 
        \begin{equation} \label{eq:osi-sketch-solve}
            \Median\!\big(\norm{\mat{A}\Xtilde-\mat{B}}_{\rm F}^2\big) \le 5(1+\eta) \cdot \min_{\mat{X} \in \R^{n\times p}} \norm{\mat{A}\mat{X}-\mat{B}}_{\rm F}^2.
        \end{equation}
    \end{enumerate}
\end{theorem}

Before we prove this result, let us make a few comments.
The first two results \ref{item:sas-ose} and \ref{item:sas-osi-uniform} are structural bounds that hold for any embedding, deterministic or randomized, whereas the third result is inherently probabilistic \ref{item:sas-osi-probability}.
Among the three results, only the subspace embedding result \ref{item:sas-ose} shows that $\Xtilde$ approaches optimality in the limit $\eta \downarrow 0$.
For the curious reader, we note that the rate of convergence is not sharp in the limit $\eta \to 0$, and the bound \cref{eq:ose-sketch-solve} can be improved so that the prefactor takes the form $1 + \order(\eta^2)$ \cite[Thm.~C.2]{epperly2025make}.
Both results \ref{item:sas-osi-uniform} and \ref{item:sas-osi-probability} use the subspace injection hypothesis, but the results have a fundamentally different character.
The bound \ref{item:sas-osi-uniform} holds even for a deterministic embedding $\mat{\Omega}$ and holds uniformly and simultaneously over all possible choices of $\mat{B}$.
However, for most embeddings, the norm $\norm{\mat{\Omega}}^2 = \order(m/s)$, so the prefactor in \cref{eq:sketch-solve-universal2} scales with the dimension.
The second bound \ref{item:sas-osi-probability} is probabilistic and controls the median error for a single right-hand side $\mat{B}$.
On the other hand, the dimension-dependent quantity $\norm{\mat{\Omega}}^2$ is replaced by a small, absolute constant.
All three bounds are informative in different ways, and each has its own uses.

\begin{proof}[Proof of \cref{thm:sketch-and-solve-ose-osi}]
    The proof of \cref{eq:ose-sketch-solve} is one line:
    For any $\mat{X} \in \R^{n\times p}$,
    \begin{align*}
        (1+\eta)^{-1}\norm{\mat{A}\Xtilde-\mat{B}}_{\rm F}^2 &\le \norm{\mat{\Omega}^\tp(\mat{A}\Xtilde-\mat{B})}_{\rm F}^2 \\ &\le \norm{\mat{\Omega}^\tp(\mat{A}\mat{X}-\mat{B})}_{\rm F}^2\le (1+\eta)\norm{\mat{A}\mat{X}-\mat{B}}_{\rm F}^2.
    \end{align*}
    The first and third bounds are the subspace embedding property, and the second is the optimality of $\Xtilde$ as a solution to the least-squares problem $\min_{\mat{X}} \norm{\mat{\Omega}^\tp(\mat{A}\mat{X}-\mat{B})}_{\rm F}^2$.

    Next, we prove \cref{eq:osi-sketch-solve}.
    By optimality of $\mat{X}_\star$, the residual $\mat{A}\mat{X}_\star - \mat{B}$ is orthogonal to the column space of $\mat{A}$.
    Therefore, by Pythagoras,
    \begin{equation} \label{eq:sketch-solve-universal1}
        \norm{\mat{A}\Xtilde - \mat{B}}_{\rm F}^2 = \norm{\mat{A}\mat{X}_\star - \mat{B}}_{\rm F}^2 + \norm{\mat{A}(\mat{X}_\star-\Xtilde)}_{\rm F}^2
    \end{equation}
    Thus, our task is to bound $\norm{\mat{A}(\mat{X}_\star-\Xtilde)}_{\rm F}^2$.
    To do so, observe that we may write
    \begin{equation} \label{eq:blubber}
        \mat{A}(\mat{X}_\star-\Xtilde) = \mat{A} (\mat{\Omega}^\tp \mat{A})^\pinv \mat{\Omega}^\tp (\mat{A}\mat{X}_\star - \mat{B}) = \mat{P} (\mat{A}\mat{X}_\star - \mat{B}),
    \end{equation}
    where
    \[
        \mat{P} \coloneqq \mat{A} (\mat{\Omega}^\tp \mat{A})^\pinv \mat{\Omega}^\tp = \mat{Q}(\mat{\Omega}^\tp \mat{Q})^\pinv\mat{\Omega}^\top.
    \]
    When $\mat{\Omega}$ is a subspace injection for $\range(\mat{A})$, then $\mat{\Omega}^\tp \mat{A}$ has full column-rank and $\mat{P}$
    is an oblique projector onto $\range(\mat{A})$.
    This decomposition implies the bound
    \begin{equation*}
        \norm{\mat{A}(\mat{X}_\star-\Xtilde)}^2_{\rm F} \le \norm{\smash{(\mat{\Omega}^\tp \mat{Q})^\pinv}}^2 \cdot \norm{\smash{\mat{\Omega}^\tp(\mat{A}\mat{X}_\star - \mat{B})}}_{\rm F}^2.
    \end{equation*}
    With probability 3/4, the subspace injection property holds and the first factor is bounded as $\norm{\smash{(\mat{\Omega}^\tp \mat{Q})^\pinv}}^2 \le 1+\eta$.
    To control the second factor, observe that isotropy implies that $\expect \norm{\smash{\mat{\Omega}^\tp(\mat{A}\mat{X}_\star - \mat{B})}}_{\rm F}^2 = \norm{\smash{\mat{A}\mat{X}_\star - \mat{B}}}_{\rm F}^2$.
    Therefore, by Markov's inequality, $\norm{\smash{\mat{\Omega}^\tp(\mat{A}\mat{X}_\star - \mat{B})}}_{\rm F}^2 \le 4\norm{\smash{\mat{A}\mat{X}_\star - \mat{B}}}_{\rm F}^2$ except with probability 1/4.
    Thus, the Pythagorean decomposition~\cref{eq:sketch-solve-universal1} gives
    \begin{align*}
        \norm{\mat{A}\Xtilde - \mat{B}}_{\rm F}^2 &= \norm{\mat{A}\mat{X}_\star - \mat{B}}_{\rm F}^2 + \norm{\mat{A}(\mat{X}_\star-\Xtilde)}_{\rm F}^2 \\ &\le \norm{\mat{A}\mat{X}_\star - \mat{B}}_{\rm F}^2 + \norm{\smash{\mat{A} (\mat{\Omega}^\tp \mat{A})^\pinv}}^2 \cdot \norm{\smash{\mat{\Omega}^\tp(\mat{A}\mat{X}_\star - \mat{B})}}_{\rm F}^2 \\
        &\le \norm{\mat{A}\mat{X}_\star - \mat{B}}_{\rm F}^2 + (1+\eta) \cdot 4 \norm{\mat{A}\mat{X}_\star - \mat{B}}_{\rm F}^2 \le 5(1+\eta) \norm{\mat{A}\mat{X}_\star - \mat{B}}_{\rm F}^2.
    \end{align*}
    This is \cref{eq:osi-sketch-solve}.

    Finally, we prove \cref{eq:sketch-solve-universal2}.
    This result is trivially satisfied when $m = n$, so we may assume that $m > n$ and, hence, $\mat{P} \ne \Id$.
    This ensures that 
    the identity $\norm{\Id-\mat{P}}=\norm{\mat{P}}$ holds~\cite{szyld2006many}.
    Using~\cref{eq:blubber} proves the result:
    \begin{align*}
        \norm{\smash{\mat{A}\widetilde{\mat{X}}-\mat{B}}}^2_{\rm F}
        &= \norm{\smash{(\Id-\mat{P})( \mat{A}\mat{X}_\star - \mat{B} )}}_{\rm F}^2 
        \leq \norm{\mat{P}}^2 \,
        \norm{\mat{A}\mat{X}_\star - \mat{B}}_{\rm F}^2 \\
        &\leq
        \norm{\smash{(\mat{\Omega}^\tp \mat{Q})^\pinv}}^2\,
        \norm{\mat{\Omega}}^2\,
        \norm{\mat{A}\mat{X}_\star - \mat{B}}_{\rm F}^2
        \le (1+\eta) \norm{\mat{\Omega}}^2
        \norm{\mat{A}\mat{X}_\star - \mat{B}}_{\rm F}^2.
    \end{align*}
    This is the final desired result \cref{eq:sketch-solve-universal2}.
\end{proof}

\subsubsection{Beyond the subspace embedding and injection properties} \label{sec:beyondsubspacemebddding}

So far, our governing principle has been that a good random embedding is one that acts on a subspace like a Gaussian embedding.
In the absence of a detailed analysis of Gaussian-like behavior for a particular embedding, we can instead analyze an embedding by using the subspace embedding and subspace injection properties, which merely require that the embedding respect the $\ell_2$ norms of all vectors within a subspace.
For those so inclined, there is a much greater array of mathematical properties to utilize that have been developed in the literature to understand sketching algorithms.
Let us briefly mention a few more properties and note connections.

\begin{itemize}
    \item A random embedding satisfies the \emph{Johnson--Lindenstrauss (JL) property} with \emph{size} $p$ if it controls the distortion~\cref{eq:subpace-embedding} over any finite set $\set{V}$ with $p$ elements \cite[see, e.g.,][]{dasgupta2003elementary}.
    Compare to the subspace embedding property, which requires \cref{eq:subpace-embedding} to hold for the infinitely many points comprising a linear subspace.
    One can deduce the oblivious subspace embedding property with dimension $n$ from the JL property with size $p = \e^{\order(n)}$ by covering the sphere with a net and applying a union bound \cite[Ch.~4]{vershynin2026high}.
    See \cref{sec:jl} for more discussion on the JL property.
    \item A random embedding is said to satisfy the \emph{approximate matrix multiplication} (AMM) property with distortion $\varepsilon$ and failure probability $\delta$ if
    \[
        \norm{\smash{\mat{A}^\tp\mat{\Omega}\mat{\Omega}^\tp\mat{B}
               -\mat{A}^\tp\mat{B}}}_{\rm F}
        \leq \varepsilon\norm{\mat{A}}_{\rm F}\norm{\mat{B}}_{\rm F} \quad \text{with probability }1-\delta
    \]
    holds for all matrices $\mat{A}$ and $\mat{B}$ with $n$ rows \cite{sarlos2006improved}.
    One can derive improved error bounds for the sketch-and-solve algorithm and randomized low-rank approximation algorithms by combining the AMM and subspace injection properties; see \cite[Thm.~4.2]{woodruff2014sketching} and \cite[Lem.~1]{drineas2011faster}.
    Also, by taking $\mat{A} = \mat{B} = \mat{Q}$, the AMM property implies that $\norm{\mat{Q}^\tp\mat{\Omega}\mat{\Omega}^\tp\mat{Q}-\Id} \leq\varepsilon r$, which implies the subspace embedding property when $\varepsilon \le \frac{\eta}{r(1+\eta)}$. 
    \item The quantities $\expect \big| \norm{\smash{\mat{\Omega}^\top\vec{v}}}^2 - 1\big|^p$ are called the \emph{JL moments} of a random embedding.
    They are the polynomial moments of the distortion of the embedding for a single unit-norm vector $\vec{v}$.
    Uniform control on the JL moments can be used to prove the JL, AMM, and subspace embedding properties; see \cite{kane2012sparser,ahle2020oblivious}.
    To establish the subspace injection property, control on all the JL moments is more than is needed, and control on the low-order moments is often enough; see~\cref{sec:tensors}.
\end{itemize}

Other properties for random embeddings include \emph{projection cost preservation} \cite{musco2020projection}, \emph{$(\varepsilon,\delta)$-unbiasedness} \cite{derezinski2021sparse}, and \emph{regularized spectral approximation} \cite{chenakkod2026accelerating}.
Lastly, let us mention that we have been solely concerned with the norm-preservation properties of random embeddings with respect to the $\ell_2$ norm.
For other $\ell_p$ norms, see \cite{woodruff2013subspace} and \cite[Ch.~3]{woodruff2014sketching}.

\subsection{Desiderata for random embeddings 2: Computational considerations} \label{sec:desiderata}

In addition to preserving the geometry of data being reduced, a random embedding must be efficient to construct, store, and apply to a matrix.
In this section, we review such computational considerations for random embeddings.

\subsubsection{Fast apply operation} \label{sec:fast-embedding}

The most fundamental operation for a random embedding is the apply operation $(\mat{A},\mat{\Omega}) \mapsto \mat{A}\mat{\Omega}$.
If both the 
embedding $\mat{\Omega}\in\R^{n\times s}$ and the matrix $\mat{A} \in \R^{m\times n}$ are stored as dense arrays, the apply step expends $\order(mns)$ arithmetic operations using conventional matrix multiplication.
We can go faster by using a \emph{fast embedding}.
Loosely speaking, a random embedding $\mat{\Omega}\in\R^{n\times s}$ is said to be \emph{fast} if the product $\mat{A}\mat{\Omega}$ with a dense array $\mat{A}\in\R^{m\times n}$ can be computed in $\order(mn \log (n))$ operations or fewer.\footnote{As we have already discussed, randomized algorithms may require applying embeddings both to the left and right (i.e., $\mat{\Omega}^\tp\mat{A}$ and $\mat{A}\mat{\Omega}$).
Algebraically, applying to the left is merely applying to the right on the transpose $\mat{\Omega}^\tp \mat{A} = (\mat{A}^\tp \mat{\Omega})^\tp$. Fast left and right apply operations are therefore formally equivalent.
Be warned, the cost of left and right sketching can be quite different in practice, and optimizing the storage format for arrays can lead to meaningful speedups in sketching algorithms.}

Examples of fast random embeddings will be presented in \cref{sec:sparse-embeddings,sec:srfts}.
Using a fast embedding is essential to obtain the best-possible computational scaling for many randomized algorithms.
For instance, sketch-and-solve for a single right-hand side runs in $\order(mn\log(s) + s n^2)$ operations when using a one-round SRFT embedding (\Cref{sec:srfts}) but requires $\order(mns)$ operations when using a Gaussian embedding.

\subsubsection{Other computational considerations} \label{sec:othercomp}

The apply operation is not the only source of expense when computing with random embeddings. 
We must also expend work to generate the random numbers comprising the embedding, and we must have space available to store it.\footnote{Even dense random embeddings need not necessarily be stored explicitly: they can be generated on demand, reducing memory and communication costs; see~\cite[§§2.1.1, 2.3]{murray2023randomized}.}
These additional sources of cost can be a real concern; see, e.g., \cite[p.~A1207]{minster2024parallel}. 
Fortunately, the best-known fast embeddings are also efficient to generate and store.

\subsubsection{Obliviousness}

In many cases, we want to use a random embedding that preserves distances in a subspace $\set{V}$ that is not known to us in advance.
For instance, we may be interested in an embedding that respects the $r$-dimensional dominant eigenspace of a symmetric matrix $\mat{A}$.
In this setting, it takes significant computational effort to identify the subspace, and it is precisely this effort that we seek to \emph{avoid} by using a randomized method!

A random embedding is \emph{oblivious} if its distribution is chosen independently of the problem data. More informally, an embedding is oblivious if it may be constructed without ever looking at the problem data. Gaussian random embeddings are oblivious, as are several more constructions discussed below.
It is, of course, too much to demand that an oblivious random embedding is guaranteed to be a subspace embedding or injection for \emph{every} $r$-dimensional subspace, but for any \emph{fixed}
$r$-dimensional subspace we expect this property to hold with high probability.

\begin{definition}[Oblivious subspace embedding and injection]
    Let $\mat{\Omega} \in \R^{n\times s}$ be a random matrix, and fix a \emph{dimension} $r$, \emph{distortion} $\eta$, and \emph{failure probability} $\delta$.
    \begin{itemize}
        \item The matrix $\mat{\Omega}$ is an \textbf{oblivious subspace embedding} if, for any $r$-dimensional subspace $\set{V}$, the matrix $\mat{\Omega}$ satisfies the subspace embedding property for $\set{V}$ with distortion $\eta$ with probability at least $1-\delta$.
        \item The matrix $\mat{\Omega}$ is an \textbf{oblivious subspace injection} if it is \emph{both} isotropic and, for any $r$-dimensional subspace $\set{V}$, the matrix $\mat{\Omega}$ satisfies the subspace injection property for $\set{V}$ with parameter $\eta$ with probability at least $1-\delta$.
    \end{itemize}
\end{definition}
\Cref{thm:gaussian-ose} shows that a Gaussian embedding is a subspace embedding provided the embedding dimension satisfies $s \gtrsim (r+\log(1/\delta))/\eta^2$.
We note that the isotropy property is essential in the definition of an oblivious subspace injection, as it fixes a scale.

\subsection{Dense, unstructured embeddings: Gaussian, Rademacher, et al.} \label{sec:dense-sketching}

As we have already seen, Gaussian embeddings (\cref{def:gaussian}) are very effective, possessing crisp theoretical analysis (\cref{thm:rangefinder,thm:sketch-solve-gaussian,thm:gaussian-ose}).
In fact, they enjoy (nearly) optimal error guarantees in many sketching-based algorithms; see \cite{witten2015randomized,sridhar2020lower,epperly2026sharp} for examples.

Other types of dense, unstructured embeddings can be useful in some cases.
\emph{Rademacher embeddings} are populated with iid $\pm 1/\sqrt{s}$ random values.
They behave similarly to Gaussian embeddings in practice, while potentially reducing the cost for random number generation. More generally, Gaussian and Rademacher embeddings belong to the broader class of \emph{sub-Gaussian embeddings}, consisting of random matrices that have
independent, suitably normalized entries with uniformly bounded sub-Gaussian norms.
Any such embedding is an oblivious subspace
embedding provided $s = \order((r+\log(1/\delta))/\eta^2)$;
see~\cite[Cor.~5.4]{dirksen2016dimensionality}.

As another model, we can use a scaled random orthonormal matrix 
\begin{equation*}
    \mat{\Omega} = \sqrt{\frac{n}{s}} \cdot \mat{Q} \quad \text{for } \mat{Q} \sim \Unif \{ \mat{U} \in \R^{n\times s} : \mat{U}^\tp \mat{U} = \Id \}
\end{equation*}
as an embedding, where $\mat{Q}$ is uniformly distributed on the Stiefel manifold.
For  sketch-and-solve, random orthonormal embeddings can achieve somewhat better accuracy than Gaussian embeddings and enjoy certain optimality properties~\cite{dobriban2019asymptotics,epperly2026sharp}.
However, in the typical setting where $s \ll n$, the differences between random partial isometries and Gaussian embeddings are slight and generating them requires extra computational cost~\cite{stewart1980efficient}.
We generally advise against random orthonormal embeddings for practical use.

As discussed in \Cref{sec:othercomp}, dense unstructured embeddings can be expensive to generate, store, and apply; the structured constructions below address these costs.

\subsection{Sparse embeddings} \label{sec:sparse-embeddings}

To perform sketching more rapidly, we can use \emph{sparse} random embedding matrices.
One useful exemplar from this class is the SparseStack embedding.

\begin{definition}[SparseStack]
    A \emph{CountSketch} $\mat{\Phi} \in \R^{n\times b}$ is a sparse random matrix with independent rows, where a single nonzero entry with a uniform $\pm 1$ value is placed at a uniformly random position in each row.
    For a given \emph{block size} $b \in \mathbb N$ and \emph{sparsity level} $\zeta \in \mathbb N$, a \emph{SparseStack} is a random matrix
    \begin{equation*}
        \mat{\Omega}  = \frac{1}{\sqrt{\zeta}} \begin{bmatrix}
            \mat{\Phi}_1 & \cdots & \mat{\Phi}_\zeta
        \end{bmatrix} \in \R^{n\times b\zeta} \quad \text{for iid CountSketches } \mat{\Phi}_1,\ldots,\mat{\Phi}_\zeta.
    \end{equation*}
\end{definition}

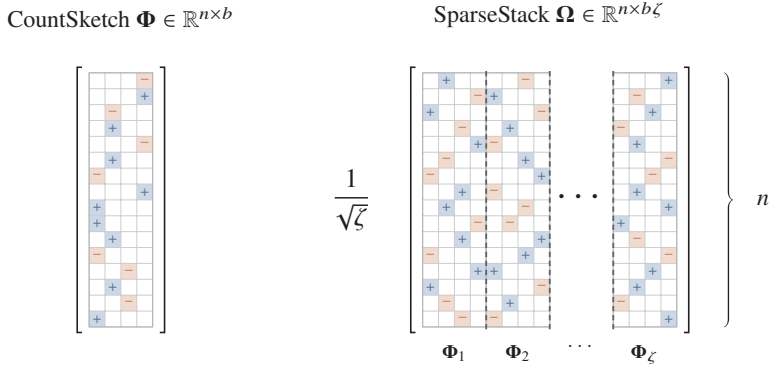
\begin{figure}[t]
  \centering
  \def\nrows{16}   
  \def\bcols{4}    
  \def\zt{5}       
  \def\nslots{4}   

  \definecolor{poscol}{RGB}{31,110,180}
  \definecolor{negcol}{RGB}{200,70,50}

  \newcommand{\sketch}[3]{\begin{scope}[shift={(#1,0)}]
      \foreach \r/\c/\s in {#3} {
        \pgfmathsetmacro{\y}{\nrows-\r}
        \ifnum\s=1
          \fill[poscol!22] (\c-1,\y) rectangle (\c,\y+1);
          \node[font=\tiny,poscol!80!black] at (\c-0.5,\y+0.5) {$+$};
        \else
          \fill[negcol!22] (\c-1,\y) rectangle (\c,\y+1);
          \node[font=\tiny,negcol!80!black] at (\c-0.5,\y+0.5) {$-$};
        \fi
      }
      \draw[gray!45,line width=0.15pt,step=1] (0,0) grid (#2,\nrows);
      \draw[gray!70,line width=0.4pt] (0,0) rectangle (#2,\nrows);
    \end{scope}
  }

  \newcommand{\brackets}[2]{\begin{scope}[shift={(#1,0)}]
      \draw[line width=0.6pt] (-0.45,-0.35) -- (-0.75,-0.35) -- (-0.75,\nrows+0.35) -- (-0.45,\nrows+0.35);
      \draw[line width=0.6pt] (#2+0.45,-0.35) -- (#2+0.75,-0.35) -- (#2+0.75,\nrows+0.35) -- (#2+0.45,\nrows+0.35);
    \end{scope}
  }

  \begin{tikzpicture}[
      x=0.21cm, y=0.21cm,
      brc/.style={decorate,decoration={calligraphic brace,mirror,amplitude=3pt},
                  line width=0.5pt,gray!80!black},
    ]

    \pgfmathsetmacro{\wide}{\nslots*\bcols}   
    \pgfmathsetmacro{\off}{\bcols+17}         

    \sketch{0}{\bcols}{1/4/-1, 2/4/1, 3/2/-1, 4/2/1, 5/4/-1, 6/2/1, 7/1/-1, 8/4/1,
      9/1/1, 10/1/1, 11/2/1, 12/1/-1, 13/3/-1, 14/2/1, 15/3/-1, 16/1/1}
    \brackets{0}{\bcols}
    \node[anchor=south,font=\footnotesize] at (\bcols/2,\nrows+2.2) {CountSketch $\mat{\Phi}\in\R^{n\times b}$};

    \sketch{\off}{\bcols}{1/2/1, 2/4/-1, 3/1/1, 4/3/-1, 5/4/1, 6/2/-1, 7/1/-1, 8/3/1,
      9/2/1, 10/4/-1, 11/3/1, 12/1/-1, 13/4/1, 14/1/1, 15/2/-1, 16/3/-1}
    \sketch{\off+\bcols}{\bcols}{1/3/-1, 2/1/1, 3/4/-1, 4/2/1, 5/1/-1, 6/3/1, 7/4/1, 8/1/-1,
      9/3/-1, 10/2/-1, 11/4/1, 12/3/1, 13/1/1, 14/4/-1, 15/2/1, 16/1/-1}
    \sketch{\off+3*\bcols}{\bcols}{1/4/1, 2/2/-1, 3/3/1, 4/1/-1, 5/2/1, 6/4/-1, 7/3/-1, 8/2/1,
      9/4/1, 10/1/1, 11/2/-1, 12/4/-1, 13/3/1, 14/2/1, 15/1/-1, 16/4/1}

    \brackets{\off}{\wide}

    \node[font=\Large] at (\off+2*\bcols+\bcols/2,\nrows/2) {$\cdots$};

    \foreach \k in {1,...,\the\numexpr\nslots-1\relax} {
      \pgfmathsetmacro{\xd}{\off+\k*\bcols}
      \draw[black!75,line width=0.7pt,dash pattern=on 2pt off 1.4pt] (\xd,-0.15) -- (\xd,\nrows+0.15);
    }

    \node[anchor=east,font=\small] at (\off-2.6,\nrows/2) {$\dfrac{1}{\sqrt{\zeta}}$};
    \node[anchor=south,font=\footnotesize] at (\off+\wide/2,\nrows+2.2)
          {SparseStack $\mat{\Omega}\in\R^{n\times b\zeta}$};

    \foreach \k/\lab in {0/{$\mat{\Phi}_1$}, 1/{$\mat{\Phi}_2$}, 2/{$\cdots$}, 3/{$\mat{\Phi}_\zeta$}} {
      \pgfmathsetmacro{\xc}{\off+\k*\bcols+\bcols/2}
      \node[anchor=north,font=\scriptsize] at (\xc,-0.5) {\lab};
    }

    \draw[brc] (\off+\wide+3.0,0) -- (\off+\wide+3.0,\nrows);
    \node[anchor=west,font=\footnotesize] at (\off+\wide+4.4,\nrows/2) {$n$};

  \end{tikzpicture}

  \caption{A CountSketch $\mat{\Phi}\in\R^{n\times b}$ (left) and a SparseStack
    $\mat{\Omega}\in\R^{n\times b\zeta}$ (right), shown for $n = 16$, $b = 4$,
    and $\zeta = 3$.  Each row of each CountSketch contains exactly one nonzero
    entry, equal to $\pm 1$ with equal probability, in a uniformly random column.
    The dashed lines mark the boundaries between the $\zeta$ independent
    CountSketches comprising $\mat{\Omega}$.}
  \label{fig:sparsestack}
\end{figure}

In other words, a SparseStack chunks every row into $\zeta$ equally sized pieces, in which we place a single nonzero entry in a uniformly random position.
See \cref{fig:sparsestack} for an illustration. By construction, a SparseStack $\mat{\Omega}$ is isotropic and
requires $\order(n\zeta)$ storage.
For a dense matrix $\mat{A} \in \R^{m\times n}$,
forming $\mat{A}\mat{\Omega}$ requires $\order(\zeta mn)$
operations.
For sparse $\mat{A}$, this reduces to
$\order(\zeta \cdot \mathrm{nnz}(\mat{A}))$ operations,
where $\mathrm{nnz}$ denotes the number of nonzero entries of a matrix.

\paragraph{Parameter settings and metric properties, empirically.}

With appropriate parameters, SparseStacks tend to achieve near-Gaussian levels of accuracy in sketching algorithms, while being significantly faster.
Empirically, we reliably obtain excellent performance for SparseStack embeddings when the sparsity level is very small, just $\zeta = 4$.
Especially for applications requiring only modest conditioning or a constant-factor approximation,
the embedding dimension $s$ can be set as for a Gaussian embedding, say $2r \le s \le 4r$.
When not stated otherwise, we use parameter values $\zeta = 4$ and $s = 2r$ in this work.

Some researchers, particularly in algorithms theory, utilize the plain CountSketch matrix as a test matrix in algorithms.
While there are good reasons to use the CountSketch in algorithms theory \cite{clarkson2017low}, one should exercise caution when deploying these test matrices in practice.
CountSketch embeddings can fail catastrophically when applied to a \emph{coherent}\footnote{A coherent subspace is one that contains vectors whose mass is concentrated on a small number of coordinates \cite[Def.~1]{drineas2012fast}.} subspace $\set{V}$ unless the embedding dimension is quadratically large, $s = \Omega(r^2)$.
By stacking together $\zeta$ copies of a CountSketch into a SparseStack, we obtain dramatic improvements in accuracy.
The intermediate parameter setting $\zeta = 2$ is more reliable, but still can be problematic.
For high-quality behavior on general subspaces, numerical experience suggests that one should use at least sparsity level $\zeta \ge 3$.
See \cite[Fig.~2]{camao2025faster} for a numerical demonstration of these effects.

\paragraph{Parameter settings and metric properties, theoretically.}

Theoretical analysis of SparseStacks and related sparse random embeddings has occupied the randomized linear algebra theory community for over a decade.
Here is one result due to~\citet{chenakkod2026optimal}.

\begin{theorem}[SparseStack: Subspace embedding] \label{thm:sparsestack}
    Let $\mat{\Omega} \in \R^{n\times s}$ be a SparseStack with sparsity level $\zeta$, and fix a subspace dimension $r$ and distortion parameter $\eta \in (0,1)$.
    Up to subpolylogarithmic factors in $r$, the conditions
$s \gtrsim r/\eta^2$ and 
$\zeta \gtrsim \log(r) / \eta$
    suffice for $\mat{\Omega}$ to be an oblivious subspace embedding of dimension $r$, distortion $\eta$, and failure probability $0.01$.
\end{theorem}
Here, ``up to subpolylogarithmic factors'' means additional multiplicative factors of the form $\log^{o(1)}(r)$ and the cited result assumes $\eta \ge r^{-\order(1)}$. Note that the result of \Cref{thm:sparsestack} continues to hold when the column positions and signs of the nonzero entries are only logarithmic-wise independent.
See \citet{nelson2013osnap,bourgain2015toward,cohen2016nearly,chenakkod2024optimal,chenakkod2025optimal,tikhomirov2026level,tropp2026comparisonb,tropp2026comparisona} for other analyses of sparse random embeddings.
\citet{nelson2014lower} showed that a sparsity level of $\zeta = \Omega(\log(r) / \log \log r)$ is necessary to achieve the subspace embedding property with embedding dimension $s = \order(r)$ and constant distortion $\eta$, so \cref{thm:sparsestack} is nearly sharp.
Very recent results \cite{heidary2026sparsestack,mai2026nelson} claim that SparseStacks are oblivious subspace embeddings with $s = \order(r/\eta^2)$ and $\zeta = \order(\log(r)/\eta)$, thereby removing the subpolylogarithmic factors in~\cref{thm:sparsestack}.

Subspace injection results from \citet{tropp2026comparisonb,tropp2026comparisona} show that $s = \order(r)$ and $\zeta = \order(\log r)$ suffice to achieve an oblivious subspace injection with constant distortion $\eta$; see also \cite{camao2025faster}.
\citet{huang2026well} recently showed that constant-sparsity 
SparseStack embeddings with $s=\order(r)$ cannot provide a worst-case subspace-injection guarantee with a dimension-independent injectivity parameter as $r \to \infty$.
Experiments reported in~\cite{camao2025faster} indicate that 
the value $\zeta = 4$ appears sufficient to ensure the subspace injection property for subspaces of dimension up to $r = 10^7$.

\paragraph{History and variants.}
The first work advocating sparse random matrices for fast dimension reduction may be the work of \citet{achlioptas2003database}.
The influential CountSketch (the $\zeta = 1$ case of SparseStacks) was introduced by \citet{charikar2004finding} and popularized for linear algebra computations by \citet{clarkson2009numerical,clarkson2017low}, \citet{nelson2013osnap}, and \citet{meng2013low}.
Building on earlier work on sparse Johnson--Lindenstrauss transforms, including \citet{dasgupta2010sparse}, \citet{kane2012sparser} introduced the block construction that we call SparseStack. They also proposed a similar SparseUniform embedding construction that places $\zeta$ nonzero randomly signed coordinates in uniformly random positions in each row.
Both the SparseStack and SparseUniform matrices appear to behave similarly in practice, though the former are easier to generate; see \cite{murray2023randomized,chen2026gpu} for a discussion of efficient generation of SparseUniform matrices.

\subsection{Subsampled randomized fast trigonometric transforms} \label{sec:srfts}

Another class of random embeddings is based on fast trigonometric transforms.

\begin{definition}[Subsampled randomized fast trig transform] \label{def:srtt}
    Let $\mat{F} \in \R^{n\times n}$ be an orthogonal matrix for which the matrix-vector products $\mat{F}\vec{z}$ and $\mat{F}^\tp\vec{z}$ can be evaluated in $\order(n\log n)$ operations, and let $\set{J} \subseteq \{1,\ldots,n\}$ be a uniformly random set of $s$ indices.
    A \emph{\textbf{S}ubsampled \textbf{R}andomized \textbf{F}ast trig \textbf{T}ransform} (SRFT) with $t$ rounds of randomization is defined as
    \begin{equation*}
        \mat{\Omega}^\tp = \sqrt{\frac{n}{s}} \Id(\set{J},:)\mat{F}\mat{D}_t\mat{F} \mat{D}_{t-1}  \cdots \mat{F} \mat{D}_1,
    \end{equation*}
    where $\mat{D}_i$ are independent diagonal matrices, each with independent Rademacher random variables on the diagonal.
\end{definition}

SRFTs can be implemented with any fast trig transform.
Popular choices include the discrete cosine transform, the Hadamard--Walsh transform (if $n$ is a power of two), or the discrete Fourier transform (which requires a complex extension of~\Cref{def:srtt} and complex arithmetic).
In the case where $\mat{F}$ is a Hadamard--Walsh transform, the embedding $\mat{\Omega}$ is often called a subsampled randomized \emph{Hadamard} transform (SRHT).
These constructions are closely related to the fast Johnson–Lindenstrauss transform of~\citet{ailon2006approximate}.

\paragraph{Parameter settings and metric properties, empirically.}

Most commonly, SRFTs are implemented with a single round of randomization ($t=1$). Two rounds ($t=2$) are necessary to reliably obtain a subspace embedding (or even a subspace injection) of dimension $s = \order(r)$ for certain \emph{coherent} subspaces, spanned by vectors with most of their mass concentrated in a few coordinates. In a different setting, $t=3$ rounds of randomization have been found
beneficial~\cite{andoni2015practical}.
For use cases in scientific computing, however, the available evidence suggests that $t=2$ rounds of randomization are sufficient.
With $t=2$ rounds of randomization, SRFTs are empirically observed to exhibit near-Gaussian behavior.
Therefore, to obtain an embedding for an $r$-dimensional subspace, one can generally use $2r \le s \le 4r$.

\paragraph{Parameter settings and metric properties, theoretically.}

For $t = 1$, subspace embedding properties of SRFTs are theoretically well understood. 

\begin{theorem}[SRFT: Subspace embedding] \label{thm:srht}
Let $\mat{\Omega} \in \R^{n \times s}$ be a one-round SRFT whose transform $\mat{F}$
has bounded coherence,
 $n \max_{i,j} |\mat{F}(i,j)|^2 = \order(1)$, and fix
    $r \le s$. If
    \[
        s \gtrsim \frac{(r + \log n)\log r}{\eta^2},
    \]
    then $\mat{\Omega}$ is an oblivious subspace embedding with distortion $\eta$ and fixed failure probability.
\end{theorem}

\Cref{thm:srht} is due to~\citet{tropp2011improved}; see~\cite{boutsidis2013improved} for related work.
Its coherence assumption is satisfied by 
all common SRFTs, including SRHT. \citet{tropp2011improved} also showed
that without multiple rounds of randomization, an SRHT \textit{requires} an embedding dimension of $s = \Omega(r \log r)$ to achieve even the oblivious subspace injection property with nonvanishing probability.
Very recent results~\cite{yang2026subspace,heidary2026rerandomized} claim that the SRHT with $t=2$ rounds is an oblivious subspace embedding with $s=\order(r/\eta^2)$ for fixed failure probability, confirming the empirical discussion above.

\paragraph{Resource usage.}
An SRFT is described by the $nt$ diagonal entries of the matrices $\mat{D}_i$ plus the $s$ randomly selected coordinates $\isettwo$.
As such, SRFTs are very economical to generate and store.
To apply an SRFT to a vector, we interleave multiplication by diagonal matrices $\mat{D}_i$ and applications of the fast trig transform $\mat{F}$, followed by subsampling to the index set $\isettwo$ at the last step.
The total runtime is $\order(t \cdot nm\log n)$ to compute the action of $\mat{\Omega}^\tp$ on a matrix $\mat{A} \in \R^{n\times m}$.

If only $t=1$ round of randomization is used, the product $\mat{\Omega}^\tp \mat{A}$ can be computed in a slightly faster $\order(nm\log s)$ operations by making use of a fast \emph{subsampled} trig transform algorithm \cite{sorensen1993efficient,ailon2009fast} to compute the action of $\Id(\isettwo,:)\mat{F}$.
This variant, faster in principle, is rarely implemented, and it is believed 
to lead only to minor speedups in practice.

As a downside, SRFTs do not leverage structure in the matrix $\mat{A}$ such as sparsity.
As such, SRFTs are not the appropriate dimension reduction map for sketching sparse matrices.

\paragraph{History and variants.}
Random embeddings similar to SRFTs were first used by \citet{ailon2006approximate}, and refined in subsequent works by \citet{woolfe2008fast} and \citet{drineas2011faster}.
Some researchers advocate using variants of SRFTs where signed permutations are used in place of random diagonal matrices \cite{tropp2019streaming,martinsson2020randomized}.
One-round SRFTs can also be made more reliable by replacing the final coordinate subsampling by a sparse random embedding; see \cite{cartis2021hashing,chenakkod2025optimal,camao2025faster}.
\citet{balabanov2023block} proposed a block SRHT that admits subspace embedding guarantees comparable to the standard SRHT while being better suited to parallel, distributed computation.

\subsection{A few more options: LASO and multisketching} \label{sec:laso}

We briefly comment on a few additional proposals for fast random embeddings:
\begin{itemize}
    \item \emph{Long-axis-sparse operators} (LASOs, \citet[\S2.4]{murray2023randomized})
    is a catch-all term for sparse random embeddings that have a small number of nonzeros per \emph{column} (long axis); see, e.g., \cite{derezinski2021sparse,dereziski2021newton}.
    LASOs are appealing because they are very fast.
    However, to be used safely as an embedding for an arbitrary subspace $\set{V}$, the sparsity pattern of the embedding must be adapted to the leverage scores of $\set{V}$ (\cref{def:leverage}), which are typically somewhat expensive to compute or estimate.
    Therefore, we do not recommend LASOs as a sketching technique for general use.
    \item \emph{Multisketching \cite{clarkson2017low,kapralov2016how,higgins2025analysis}.}
    The idea of multisketching is to build a random embedding as a composition
    $\mat{\Omega} = \mat{\Omega}_1\mat{\Omega}_2$ of two constituent random embeddings $\mat{\Omega}_1 \in \R^{n\times s_1}$ and $\mat{\Omega}_2 \in \R^{s_1 \times s_2}$.
    The idea is to use a very fast embedding $\mat{\Omega}_1$ with suboptimal embedding dimension to perform an initial phase of compression (e.g., a CountSketch or traditional SRFT), then use a slower, more reliable embedding to perform a final level of compression.
    Typically, we find that standalone sparse embeddings and SRFTs are cheaper than multisketching, but there is some evidence that highly optimized CountSketch--Gaussian multisketching may be faster for sketching very tall-thin matrices~\cite{higgins2025high,higgins2025analysis}.
\end{itemize}

\subsection{Tensor-structured embeddings} \label{sec:tensor-embeddings-intro}

One of the most exciting ongoing areas of development in randomized linear algebra is the development of specialized methods for \emph{tensors}, i.e., multilinear arrays~\cite{ballard2025tensor}.
For tensor problems, it is often desirable to use specialized embeddings which respect or exploit tensor-product structure.
These embeddings will be discussed in \cref{sec:tensors}.

\subsection{Index selection} \label{sec:index-sketching-intro}

So far, we have explored sketching techniques that seek to efficiently mix together the rows or columns of a matrix.
On the opposite extreme, we can also try to sketch a vector $\vec{a} \in \R^m$ or matrix $\mat{A} \in \R^{m\times n}$ by subsampling a collection of entries $\vec{a}(\iset)$ or rows $\mat{A}(\iset,:)$, characterized by some index collection $\iset \subseteq \{1,\ldots,m\}$ of small cardinality $s \coloneqq |\iset| \ll m$.
This is an appealing idea, but it must be applied with care.
Often, the rows of a matrix $\mat{A}$ are of non-uniform importance, and it is essential to choose the indices $\iset$ in a manner biased towards selecting the important rows.
For this reason, reliable dimension reduction techniques based on index selection
are typically \emph{non-oblivious}, requiring a sampling distribution adapted
to the data being processed.
Index selection techniques will be discussed in \cref{sec:selection}.

\subsection{Which embedding should I use?} \label{sec:which-embedding}

The question of which embedding to use continues to be actively discussed by researchers in randomized linear algebra, and there is no consensus.
Nonetheless, we believe the available numerical evidence is sufficient to justify some recommendations.
In this subsection, we present a representative numerical experiment (\cref{sec:which-embedding-setup,sec:which-embedding-results}), then discuss conclusions (\cref{sec:which-embedding-discussion}).

\subsubsection{Numerical evaluation: Setup} \label{sec:which-embedding-setup}

\Cref{fig:sketching-comparison} shows a comparison of nine different sketching techniques: Gaussian, SRFT (with both 1 and 2 rounds), CountSketch, SparseStack (with $\zeta = 4$), LASO, multisketching, uniform sampling, and leverage score sampling.
We test on two matrices, both of size $m = 10^6$ by $n = 10^2$.
The first matrix (``Arnoldi'') is the tall-thin matrix $\mat{A}$ obtained from the truncated Arnoldi process during the execution of the sGMRES method (\citet{nakatsukasa2024fast}, \cref{sec:sgmres}) on an advection--diffusion linear system; see \cref{fig:sgmres}.
The second matrix (``coherent'') is the adversarial matrix $\mat{A} =\flatonebytwo{\Id}{\mat{0}}^\tp$, which is a hard instance for certain sketching techniques.
Both test matrices $\mat{A}$ are stored as dense arrays.
For each test matrix $\mat{A}$ and random embedding $\mat{\Omega}$, we report two metrics: (1) the total time to generate the embedding $\mat{\Omega}$ and to apply it to compute $\mat{\Omega}^\tp \mat{A}$ and (2) the minimum singular value of $\mat{\Omega}^\tp \mat{Q}$ for $\mat{Q} = \orth(\mat{A})$, which measures the quality of $\mat{\Omega}$ as a subspace injection for $\set{V} = \range(\mat{A})$.

All sketching techniques are implemented using a plain MATLAB implementation.
The trigonometric transform for the SRFT is the discrete cosine transform.
For multisketching, we use a CountSketch with embedding dimension $s_1 = 12n^2$ composed with a Gaussian embedding.
Uniform and leverage-score sampling use $s$ iid samples with replacement, with the standard rescaling described in \cref{sec:leverage}.
For a LASO embedding, we test the LESS embedding \cite{derezinski2021sparse}, a non-oblivious sketching matrix that requires estimates of the leverage scores to construct. 
We use $\lceil \log m \rceil$ nonzeros \emph{per column} for the LESS embedding.
To implement leverage score sampling and LESS sketching, we estimate the leverage scores using the procedure of \citet{drineas2012fast}, as described in \cref{sec:leverage-estimation}.
Notably, this procedure uses a fast embedding to estimate the leverage scores, so using leverage score-based sketching techniques for a general matrix requires sketching twice: once to estimate the leverage scores and again to use them.

\subsubsection{Numerical evaluation: Results} \label{sec:which-embedding-results}

\begin{figure}
    \centering
    Arnoldi
    \includegraphics[width=0.99\linewidth]{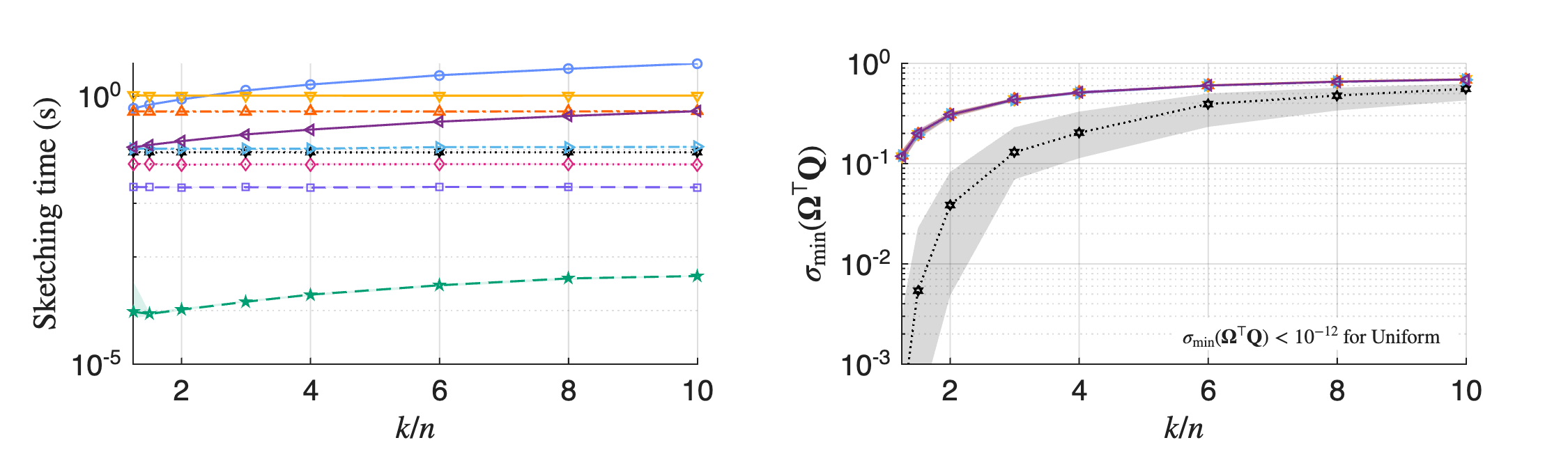}

    Coherent
    \includegraphics[width=0.99\linewidth]{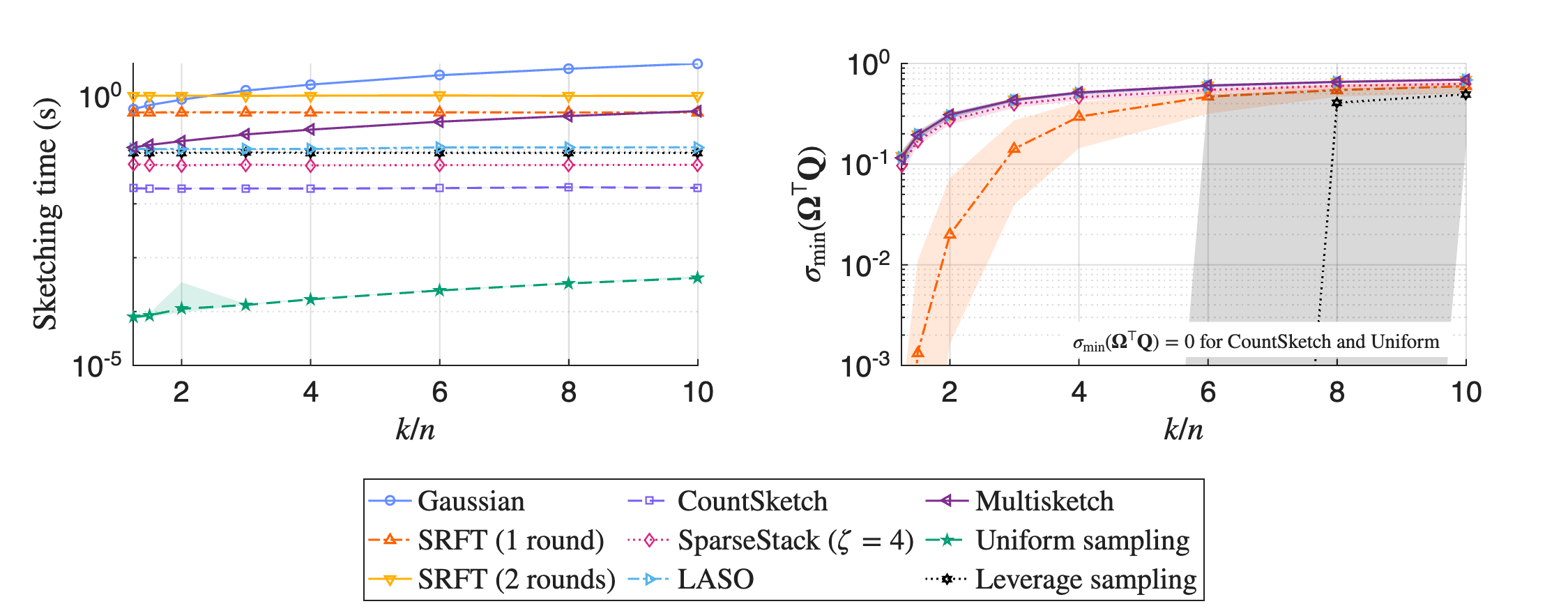}
    
    \caption{Comparison of runtime and minimum stretch $\sigma_{\rm min}(\mat{\Omega}^\tp\mat{Q})$ for nine different sketching techniques on two different test matrices; see text for details.
    Lines show medians, and shaded regions show 90\% and 10\% quantiles.
    The minimum-stretch line for uniform sampling is too small to be visible on both examples, as is the line for the CountSketch embedding on the coherent example.}
    \label{fig:sketching-comparison}
\end{figure}

The results are shown in \cref{fig:sketching-comparison}.
We see that several embeddings produce very similar results for the minimum stretch factor $\sigma_{\rm min}(\mat{\Omega}^\tp\mat{Q})$, consistent with near-universal behavior. 
Exceptions include uniform and leverage score sampling, SRFT with 1 round, and CountSketch, which have a substantially smaller minimum stretch factor than the other embeddings on at least one of the examples.
Among the embeddings with acceptable metric properties, the SparseStack is the fastest embedding in these experiments, followed by LASO, multisketching, SRFT, and Gaussian embeddings.
In particular, SparseStack is $10\times$ faster than SRFT with 1 round of randomization, $17\times$ faster than the Gaussian embedding, and $20\times$ faster than SRFT with 2 rounds of randomization at $s=2n$.
We make a few more specific comments.

\paragraph{Safe parameter choices.}
For the SparseStack and SRFT embeddings, the choice of parameter values is essential to obtain good results on hard inputs.
Indeed, a SparseStack with sparsity level $\zeta = 1$ (i.e., a CountSketch) annihilates a nonzero vector in the subspace $\set{V} = \range(\mat{A})$ on the coherent example with $s\le 10n$.
SRFT with 1 round of randomization does not fail quite so catastrophically, but the minimum stretch $\sigma_{\rm min}(\mat{\Omega}^\tp\mat{Q})$ is still $15\times$ smaller than for other types of embeddings at $s=2n$ on the coherent example.

\paragraph{Subsampling methods.}
The subsampling methods---uniform sampling and leverage score sampling---do not seem very appealing in these experiments.
The uniform sampling method is very cheap, but it fails to be a subspace injection on both examples for all $s\le 10n$.
The leverage score sampling method is better, but it still is less effective than the SparseStack embedding---being both more expensive and having a smaller minimum stretch $\sigma_{\rm min}(\mat{\Omega}^\tp\mat{Q})$.
This behavior is consistent with its worst-case $\Theta(n\log n)$
sample complexity; see \cref{thm:leverage}.
On the basis of these behaviors, we do not recommend subsampling methods as a \emph{general-purpose} sketching technique \emph{for arbitrary subspaces}.
However, there are settings where subsampling techniques are amazingly powerful; see \cref{sec:selection}.

\subsubsection{Conclusions} \label{sec:which-embedding-discussion}

On the basis of the available evidence \cite[see also][]{dong2023simpler,camao2025faster,chen2026gpu,murray2023randomized}, we recommend the SparseStack embedding for general-purpose use, provided that an efficient implementation is available.
There is also evidence that high-performance implementations of SRFTs and multisketching can be effective \cite{balabanov2023block,higgins2025high}.

\section{Low-rank approximation} \label{sec:lowrank}

Across many applications, we encounter matrices with rapidly decaying singular values ~\cite{beckermann2019bounds,udell2019why,webb2026low}. Such matrices are well-approximated by low-rank matrices, and it is often desirable to compute such an  approximation efficiently.
Randomized algorithms are often the method of choice for this task. This section will discuss such methods.

In modern RNLA, there is a common strategy for computing a low-rank approximation. First, to acquire information about the matrix $\mat{A}$, we apply 
it to one or more random test matrices. 
Then, we use this information to compute a low-rank approximation to $\mat{A}$.
To instantiate this paradigm, we must make two choices: \emph{(1) which test matrix or test matrices to use}
and \emph{(2) how to postprocess the resulting information to obtain a low-rank
approximation, that is, which approximation format to use.}
Several useful algorithms arise by combining different answers to questions (1) and (2); \cref{tab:lowrank-algorithms} summarizes the main examples
discussed in this and subsequent sections.

\begin{table}
\centering
\caption{Common algorithms for low-rank approximation, organized by choice of test matrices and approximation formats}
\label{tab:lowrank-algorithms}
\begin{tabular}{lccc}
\toprule
\textbf{Test matrices} & \textbf{Projection} & \textbf{\Nystrom} & \textbf{Gen.\ \Nystrom} \\
\midrule
Embedding & Randomized SVD & (Single-pass) \Nystrom & Gen.\ \Nystrom \\[0.05cm]
Powered & RSI  & NysSI &  \\[0.05cm]
Block Krylov & RBKI & NysBKI &  \\ \midrule
Index selection & Interpolative & Column \Nystrom & CUR \textbf{OR} cross \\
& decomposition & (partial Cholesky) & approximation \\
& (partial QR) & & (partial LU) \\
\bottomrule
\end{tabular}
\end{table}

We begin in \cref{sec:lowrank-why} by discussing the motivation for computing a low-rank approximation.
Then, in \cref{sec:lowrank-formats}, we address question (2) by discussing
three approximation formats:
projection, \Nystrom, and generalized \Nystrom approximations. Question (1) is addressed in 
\Cref{sec:lowrank-test-matrix}, discussing different choices of test
matrices and how these choices affect the performance of the resulting algorithms.
In \cref{sec:lowrank-error-analysis}, we develop representative error bounds,
and \cref{sec:lowrank-gram-correspondence} explains the important Gram correspondence
between projection and \Nystrom approximations.
\Cref{sec:lowrank-functions} discusses extensions to matrix functions, and \cref{sec:rankest} discusses estimation of the rank.
We conclude in \cref{sec:lowrank-discussion} with a comparison of the
different methods.

\subsection{Low-rank approximation and why to compute one} \label{sec:lowrank-why}

Let $\mat{A} \in \R^{m\times n}$ be a matrix with SVD $\mat{A} = \mat{U}\mat{\Sigma}\mat{V}^\tp$, and let $1\le r<\min\{m,n\}$.
By the Eckart--Young--Mirsky theorem, the \emph{truncated SVD} gives rise to a matrix
\[
\lowrank{\mat{A}}_r \coloneqq \mat{U}(:,1:r) \mat{\Sigma}(1:r,1:r) \mat{V}(:,1:r)^\tp
\]
that is a best rank-$r$ approximation to $\mat{A}$ in any unitarily invariant norm~\cite[\S 7.4.9]{horn2012matrix}. For the spectral and Frobenius norms, the error of this approximation is
\begin{equation*}
    \norm{\mat{A} - \lowrank{\mat{A}}_r} = \sigma_{r+1}(\mat{A}) \quad \text{and} \quad \norm{\mat{A} - \lowrank{\mat{A}}_r}_{\rm F}^2 = \sum_{j=r+1}^{\min\{m,n\}} \sigma_j^2(\mat{A}).
\end{equation*}
In particular, $\mat{A}$ has an accurate rank-$r$ approximation if and only if its trailing singular values $\sigma_{r+1}(\mat{A}), \sigma_{r+2}(\mat{A}), \ldots$ are small.

Constructing the optimal low-rank approximation by computing the full SVD and then truncating it is expensive, requiring $\order(mn\cdot \min(m,n))$ operations, and we are often willing to settle for a suboptimal---but cheap to compute---rank-$r$ approximation $\Ahat$.
To economize storage, we typically want to represent $\Ahat$ in factored form.
The simplest is to use an outer product representation $\Ahat = \mat{F}\mat{G}^\tp$ for some $\mat{F} \in \R^{m\times r}$ and $\mat{G} \in \R^{n\times r}$.
At modest additional cost, we can instead use a compact SVD representation $\Ahat = \mat{U}\mat{\Sigma}\mat{V}^\tp$ for orthonormal matrices $\mat{U} \in \R^{m\times r}$, $\mat{V} \in \R^{n\times r}$, and a diagonal psd matrix 
$\mat{\Sigma} \in \R^{r\times r}$.
When $m = n$ and the matrix $\mat{A}$ is psd, it is sensible to choose a psd low-rank approximation $\Ahat$, for which one  can use a symmetric outer product representation $\Ahat = \mat{F}\mat{F}^\tp$ for $\mat{F} \in \R^{m\times r}$ or a compact eigendecomposition representation $\Ahat = \mat{U}\mat{\Lambda}\mat{U}^\tp$ for orthonormal $\mat{U} \in \R^{m\times r}$ and diagonal psd $\mat{\Lambda}$.

Why would we want to compute a low-rank approximation to a matrix?
There are several reasons:

\begin{enumerate}[label=(\alph*)]
    \item \textbf{Compression.} To store the rank-$r$ approximation $\Ahat = \mat{F}\mat{G}^\tp$, we only need to store $r(m+n)$ values.
    This is much less than the $mn$ values needed to store the full matrix $\mat{A}$ when $r\ll \min\{m,n\}$. \label{item:lra-compression}
    \item \textbf{Speeding up matrix computations.}
Several common operations with a dense and unstructured matrix $\mat{A}$ become
cheaper when they are instead performed with a low-rank approximation $\Ahat$
represented in factored form.
For instance, the matrix-vector product $\mat{A}\vec{x}$ requires $\order(mn)$ operations, while the matrix-vector product $\Ahat\vec{x} = \mat{F}(\mat{G}^\tp \vec{x})$ requires only $\order(r(m+n))$ operations.
    Thus, when the approximation $\Ahat \approx \mat{A}$ is sufficiently accurate, using $\Ahat$ as a proxy for $\mat{A}$ can yield large computational savings. \label{item:lra-speeding}
    \item \textbf{Approximate SVD (or eigendecomposition).} Given the compact SVD representation of an approximation $\Ahat$, we can read off its singular values and singular vectors. If $\Ahat$ is a good approximation to $\mat{A}$, then the \emph{largest} singular values
    of $\Ahat$ will be good approximations to the matching singular values of
    $\mat{A}$. Under suitable assumptions on singular value gaps, the same holds for the corresponding dominant singular vectors and, more robustly, for  the singular subspace spanned by them~\cite{wedin1972perturbation}.
    As such, computing a low-rank approximation can provide a fast way to obtain an approximation to the partial SVD of a matrix.
    Similarly, if $\mat{A}$ and $\Ahat$ are \psd (or just symmetric), then one can use a low-rank approximation to obtain approximations to the dominant eigenvalues and, under suitable assumptions on eigenvalue gaps, to the eigenvectors or invariant subspaces of $\mat{A}$~\cite{davis1970rotation}. \label{item:lra-singular}
    \item \textbf{Dimension reduction and denoising in downstream tasks.} A central idea in PCA and model reduction is to identify a
low-dimensional subspace that approximately contains the data vectors or
solution snapshots~\cite[e.g.,][]{jolliffe2002principal,benner2015survey,sirovich1987turbulence}.
The row or column space of a low-rank approximation provides such a subspace for dimension reduction, which often can be leveraged for substantial speedups in subsequent computations. Moreover, focusing on the dominant singular vector components can have denoising effects~\cite{gavish2014optimal}. \label{item:lra-dimension}
\end{enumerate}
The level of accuracy demanded in low-rank approximation depends on the use case.
For tasks \ref{item:lra-compression} and \ref{item:lra-speeding}, we generally care about the error $\mat{A} - \Ahat$ measured in some norm.
For \ref{item:lra-singular} and \ref{item:lra-dimension}, we may care about the accuracy of individual
singular values and the corresponding dominant singular subspaces.
As a general rule, attaining a good approximation to a dominant singular subspace requires
a highly accurate low-rank approximation; the error needs to be small relative to the relevant singular value gap.

\subsection{Low-rank approximation formats} \label{sec:lowrank-formats}

Having motivated the problem of computing a low-rank approximation, we now discuss different formats for constructing such an approximation from a right sketch $\mat{A}\mat{\Omega}$ and, optionally, a left sketch $\mat{\Psi}^\tp \mat{A}$.
The choice of test matrices $\mat{\Omega}, \mat{\Psi}$ will be discussed below in \cref{sec:lowrank-test-matrix}.

\subsubsection{Projection approximation} \label{sec:lowrank-projection}

When we probe a matrix $\mat{A}$ with a test matrix $\mat{\Omega}$, the range of the resulting matrix $\mat{Y} = \mat{A}\mat{\Omega}$ is a subspace of the range of $\mat{A}$.
If we pick a ``good'' test matrix $\mat{\Omega}$, then $\range(\mat{Y})$ will be closely aligned with the dominant left singular vectors of $\mat{A}$, and thus the range of $\mat{Y}$ will form a low-dimensional subspace that captures most of the action of $\mat{A}$.
To form a low-rank approximation $\Ahat$ to $\mat{A}$, we simply project $\mat{A}$ onto the range of $\mat{Y}$:
\begin{equation*}
    \Ahat \coloneqq \mat{Q}\mat{Q}^\tp\mat{A}
    =     
    \mat{Q}\mat{B}^\tp \quad \text{where } \mat{Q} = \orth(\mat{Y}) \text{ and } \mat{B}^\tp = \mat{Q}^\tp \mat{A}.
\end{equation*}
We call this an \emph{(orthogonal) projection approximation} to $\mat{A}$, since it is obtained by orthogonally projecting $\mat{A}$ onto the range of $\mat{Y}$.
The name \emph{QB approximation} is also often used.
For simplicity, we shall often drop the modifier ``orthogonal'' as this is the most natural way of performing the projection.
\emph{Oblique} projections lead to the
(generalized) \Nystrom approximations discussed in \Cref{sec:lowrank-nystrom,sec:lowrank-generalized-nystrom}.
See \cref{alg:rsvd} for a stable implementation.

The orthogonal projection approximation is natural as it gives the optimal approximation to $\mat{A}$ among all approximations whose range is contained in the range of $\mat{Y}$, measured in 
any unitarily invariant norm $\uinorm{\cdot}$.

\begin{algorithm}[t]
\caption{Projection approximation (randomized SVD)} \label{alg:rsvd}
\begin{algorithmic}[1]
\Require Matrix $\mat{A} \in \R^{m\times n}$ and test matrix $\mat{\Omega} \in \R^{n\times s}$
\Ensure Approximation $\Ahat = \mat{Q}\mat{B}^\tp$ or $\Ahat = \mat{U}\mat{\Sigma}\mat{V}^\tp$
\State Form the sketch $\mat{Y} \coloneqq \mat{A}\mat{\Omega}$
\State Orthonormalize by QR factorization $(\mat{Q},\sim) \gets \Call{QR}{\mat{Y},\texttt{'econ'}}$
\State Form $\mat{B} \gets \mat{A}^\tp\mat{Q}$
\State \emph{Optional:} Compute an economy-size SVD $\mat{B}^\tp = \Utilde \mat{\Sigma}\mat{V}^\tp$ and set $\mat{U} \gets \mat{Q}\Utilde$
\end{algorithmic}
\end{algorithm}

\begin{proposition}[Optimality of projection approximation] \label{prop:projection}
    Let $\mat{A} \in \R^{m\times n}$ and $\mat{\Omega} \in \R^{n\times s}$ be matrices, and let $\Ahat = (\mat{A}\mat{\Omega})(\mat{A}\mat{\Omega})^\pinv \mat{A}$ be the associated projection approximation.
    Then, for any unitarily invariant norm $\uinorm{\cdot}$,
        \begin{equation} \label{eq:projection-approximation-optimality}
        \uinorm{\Ahat - \mat{A}} = \min_{\mat{X} : \range(\mat{X}) \subseteq \range(\mat{A}\mat{\Omega})} \uinorm{\mat{A} - \mat{X}}.
    \end{equation}
    However, the projection approximation $\Ahat$ need not be the unique minimizer.
\end{proposition}

\begin{proof}
    Introduce $\mat{Q} = \orth(\mat{A}\mat{\Omega})$, and complete $\mat{Q}$ to a square orthogonal matrix $\mat{U} = [\mat{Q},\mat{Q}_\perp] \in \R^{m\times m}$.
    Then for any matrix $\mat{X} = \mat{Q}\mat{C}^\top$, it holds that
    \begin{equation*}
        \uinorm{\mat{A} - \mat{X}} = \uinorm{\smash{\mat{U}^\top(\mat{A} - \mat{Q}\mat{C}^\top)}} = \uinorm{ \twobyone{\mat{Q}^\top\mat{A} - \mat{C}^\top}{\mat{Q}^\top_\perp\mat{A}} } \ge \uinorm{\mat{Q}_\perp^\top\mat{A}} = \uinorm{\mat{A} - \Ahat}.
    \end{equation*}
    The inequality follows from the fact that deleting rows of a matrix decreases its singular values and thus any unitarily invariant norm.
    The example $\mat{A} = \diag(2,1)$, $\mat{\Omega} = \evec_1$, and $\mat{X} = \diag(1,0)$ demonstrates non-uniqueness for the spectral norm.
\end{proof}

Given the matrix $\mat{Y} = \mat{A}\mat{\Omega}$, computing the projection approximation $\Ahat$ requires two steps: (1) computing an orthonormal basis $\mat{Q}$ for the range of $\mat{Y}$, and (2) computing the matrix $\mat{B} = \mat{A}^\tp \mat{Q}$.
The former step takes $\order(ms^2)$ operations, while the latter step requires applying the transpose  of $\mat{A}$ (i.e., the adjoint operator) to $\mat{Q}$, which takes $\order(mn s)$ operations when $\mat{A}$ is dense and unstructured.
The optional conversion to compact SVD form requires an additional $\order(ns^2)$ operations.
Observe that forming such a projection approximation always requires at least two sequential passes over the matrix $\mat{A}$: one pass to compute $\mat{Y} = \mat{A}\mat{\Omega}$, and another pass to compute $\mat{B} = \mat{A}^\tp \mat{Q}$.
This contrasts with the (generalized) \Nystrom approximations (\cref{sec:lowrank-nystrom,sec:lowrank-generalized-nystrom}), which require only a single pass.

\paragraph{Sensitivity and numerical stability.}
When the matrix $\mat{A}$ is ill-conditioned, the output $\mat{Q}$ of the rangefinder is not, itself, a stable quantity.
The issue is that, when $\mat{A}$ has small singular values, the matrix $\mat{Y} = \mat{A}\mat{\Omega}$ is ill-conditioned and $\mat{Q} \gets \Call{QR}{\mat{Y}}$ is liable to change significantly under small perturbations \cite[\S19.9]{higham2002accuracy}.
Nonetheless, this ill-conditioning is localized to the \emph{small singular directions} of $\mat{A}$, so $\mat{Q}$ is stably able to capture the directions that are important in the matrix $\mat{A}$ and the projection approximation $\mat{Q}\mat{Q}^\tp\mat{A}$ is a numerically well-posed object.
\citet{connolly2022randomized} provide a detailed rounding-error analysis,
showing that the additional error incurred by roundoff is typically negligible when working in double precision, but it may become significant when working in low precision.
\citet{drineas2019low} show that the low-rank approximation \emph{error} in spectral,
Frobenius, and more generally Schatten norms is a well-posed quantity that does
not require a singular-value gap; see also~\cite{hackbusch2016new}.
As discussed above, such gaps become important when the goal is
to recover the dominant singular subspaces themselves.

\subsubsection{\Nystrom approximation} \label{sec:lowrank-nystrom}

For a \psd matrix $\mat{A}$, the preferred method of low-rank approximation is typically the \emph{\Nystrom approximation} \cite{williams2000using,drineas2005nystrom,halko2011finding,gittens2016revisiting,tropp2017fixed}, defined algebraically as 
\begin{equation*}
    \Ahat = \mat{A}\langle \mat{\Omega} \rangle \coloneqq \mat{A}\mat{\Omega}(\mat{\Omega}^\tp \mat{A} \mat{\Omega})^\dagger \mat{\Omega}^\tp \mat{A} = \mat{Y}(\mat{\Omega}^\tp\mat{Y})^\dagger \mat{Y}^\tp.
\end{equation*}
The \Nystrom approximation can be motivated as follows: we want to find a low-rank approximation $\Ahat$ to $\mat{A}$ that is spanned by the columns of $\mat{Y} = \mat{A}\mat{\Omega}$.
Since $\mat{A}$ is symmetric, we seek a symmetric approximation of the form
$\Ahat=\mat{Y}\mat{M}\mat{Y}^\tp$ with $\mat{M}=\mat{M}^\tp$.
We know that $\mat{A}\mat{\Omega} = \mat{Y}$, so it is natural to demand $\Ahat\mat{\Omega} = \mat{Y}$ as well.
Provided $\mat{\Omega}^\tp \mat{Y}$ is invertible, this condition uniquely determines $\mat{M} = (\mat{\Omega}^\tp \mat{Y})^{-1}$, leading to the \Nystrom approximation $\Ahat = \mat{Y}(\mat{\Omega}^\tp\mat{Y})^{-1} \mat{Y}^\tp$. When $\mat{\Omega}^\tp\mat{Y}$ is singular, the same construction extends naturally by replacing the inverse with the Moore--Penrose pseudoinverse $(\mat{\Omega}^\tp \mat{Y})^\pinv$.

The \Nystrom approximation has a deep connection to the general theory of Schur complements \cite{anderson1975shorted,ando2005schur}.
From this theory, it may be deduced that $\mat{A} - \mat{A}\langle \mat{\Omega}\rangle$ is \psd, so that the \Nystrom approximation is always an under-approximation to $\mat{A}$ in the Loewner order of \psd matrices.
Moreover, $\mat{A}\langle\mat{\Omega}\rangle$ is the maximal matrix $\mat{X}$ in the Loewner order such that $\mat{A} - \mat{X}$ is \psd and
$\range(\mat{X})\subseteq\range(\mat{Y})$.
Consequently, for any unitarily invariant norm $\uinorm{\cdot}$,
\begin{equation*}
    \Ahat = \argmin_{\mat{X} : \range(\mat{X}) \subseteq \range(\mat{Y}), \mat{A} - \mat{X} \succeq \mat{0}} \uinorm{\mat{A} - \mat{X}}.
\end{equation*}
See \cite[Thm.~1]{anderson1975shorted} and its proof.
\cite[Prop.~2.9]{epperly2025make} provides a summary of several useful properties of the \Nystrom approximation.
As another perspective on the \Nystrom approximation, observe that it is a projection of the matrix $\mat{A}$ onto the subspace $\range(\mat{Y})$ under the \emph{oblique} projector $\mat{Y}(\mat{\Omega}^\tp \mat{Y})^\dagger \mat{\Omega}^\tp$.

\begin{algorithm}[t]
\caption{Randomized \Nystrom approximation for a \psd matrix} \label{alg:randomized-nystrom}
\begin{algorithmic}[1]
\Require Test matrix $\mat{\Omega} \in \R^{n\times s}$, sketch $\mat{Y} \coloneqq \mat{A}\mat{\Omega} \in \R^{n\times s}$, and tolerance $\tau \in [0,1)$
\Ensure \Nystrom approximation in factored form $\Ahat = \mat{F}\mat{F}^\tp$ or $\Ahat = \mat{U}\mat{\Lambda}\mat{U}^\tp$
\State Compute the core matrix $\mat{M} \coloneqq \mat{\Omega}^\tp \mat{Y}$
\State Perform a \emph{pivoted Cholesky decomposition} $\mat{P}^\tp \mat{M} \mat{P} = \mat{R}^\tp \mat{R}$
\State Identify the set of large diagonal entries $\iset \gets \{ i : \mat{R}(i,i) \ge \tau \cdot \max_j \mat{R}(j,j) \}$
\State Form the factor matrix $\mat{F} \gets \mat{Y}\mat{P}(:,\iset)\mat{R}(\iset,\iset)^{-1}$
\State \textit{Optional:} Compute an economy-size SVD $\mat{F} = \mat{U}\mat{\Sigma}\mat{V}^\tp$ and set $\mat{\Lambda} \gets \mat{\Sigma}^2$
\end{algorithmic}
\end{algorithm}

Given $\mat{Y} = \mat{A}\mat{\Omega}$, computing the \Nystrom approximation $\Ahat$ in factored form $\Ahat = \mat{F}\mat{F}^\tp$ can be done in $\order(n s^2)$ operations; see \cref{alg:randomized-nystrom} for pseudocode.
At the cost of performing one extra SVD on $\mat{F}$, the \Nystrom approximation can also be returned as a compact eigendecomposition $\Ahat = \mat{U}\mat{\Lambda}\mat{U}^\tp$.
Note that the \Nystrom approximation is entirely determined by $\mat{\Omega}$ and $\mat{Y} = \mat{A}\mat{\Omega}$, so it can be computed in a single pass over the matrix $\mat{A}$.

Care must be taken to compute the \Nystrom approximation in a stable way.
The issue is that the core matrix $\mat{\Omega}^\top \mat{Y}$ can be numerically singular and thus requires care to invert stably.
There are two popular strategies for addressing this issue: (1) factoring the core using a pivoted matrix decomposition or (2) adding a numerical shift.
Approach (1) was proposed by \citet{bucci2025numerical}, who proposed a variant of \cref{alg:randomized-nystrom} that runs partial pivoted Cholesky until the next pivot falls below
a tolerance and applies the resulting (rectangular) Cholesky factor stably using a backward-stable least-squares solve. 
Approach (2) was proposed by \citet{li2017algorithm,tropp2017fixed}, who suggest adding a small multiple of the identity matrix to the core matrix $\mat{\Omega}^\top\mat{Y}$ to ensure Cholesky decomposition succeeds without fail; see \cite[eq.~(2.11)]{epperly2025make} for discussion and see~\cite{carson2024single} for a rounding-error analysis.

When approximating a \psd matrix, the \Nystrom approximation is typically preferred over the projection approximation; see \cref{sec:lowrank-discussion-nystrom-projection}.

\subsubsection{Generalized \Nystrom approximation} \label{sec:lowrank-generalized-nystrom}

As we saw, the standard projection approximation (\cref{alg:rsvd}) requires two passes over $\mat{A}$, while a \Nystrom approximation for \psd $\mat{A}$ requires only one pass.
It is natural to ask: Is there a single-pass algorithm for computing a low-rank approximation to a general (not necessarily \psd or square) matrix $\mat{A}\in\R^{m\times n}$?
This question is answered in the affirmative by the \emph{generalized \Nystrom approximation}.

The generalized \Nystrom approximation can be derived by combining the projection approximation (\cref{sec:lowrank-projection}) with the sketch-and-solve paradigm (\cref{sec:sas}).
Given the matrix $\mat{Y} = \mat{A}\mat{\Omega}$, observe that the projection approximation $\Ahat_{\rm proj}$ can be written as 
\begin{equation*}
    \Ahat_{\rm proj} = \mat{Y} (\mat{Y}^\dagger \mat{A}).
\end{equation*}
The matrix $\mat{Y}^\dagger \mat{A}$ solves the least-squares problem $\min_{\mat{X}} \norm{\mat{Y}\mat{X} - \mat{A}}_{\rm F}$.
To obtain a quick approximate solution, we can apply sketch-and-solve.
Drawing a second test matrix $\mat{\Psi} \in \R^{m\times t}$ with $t \ge s$, we consider the sketched problem $\min_{\mat{X}} \norm{\mat{\Psi}^\tp(\mat{Y}\mat{X} - \mat{A})}_{\rm F}$, and
obtain the approximate least-squares solution
\begin{equation*}
    (\mat{\Psi}^\tp \mat{Y})^\dagger (\mat{\Psi}^\tp \mat{A}) \approx \mat{Y}^\dagger \mat{A},
\end{equation*}
yielding the \emph{generalized \Nystrom approximation}
\begin{equation*}
    \Ahat_{\rm GN} \coloneqq \mat{Y} (\mat{\Psi}^\tp \mat{Y})^\dagger (\mat{\Psi}^\tp \mat{A}) = \mat{A}\mat{\Omega} (\mat{\Psi}^\tp \mat{A}\mat{\Omega})^\dagger \mat{\Psi}^\tp \mat{A}.
\end{equation*}
As the name ``generalized \Nystrom approximation'' suggests, the generalized \Nystrom approximation reduces to the \Nystrom approximation when $\mat{A}$ is \psd and $\mat{\Psi} = \mat{\Omega}$.
Provided that $\mat{\Psi}^\tp\mat{Y}$ has full column rank, the generalized Nystr\"om approximation can be understood as the projection of $\mat{A}$ onto $\range(\mat{Y})$ under the \emph{oblique} projector $\mat{Y} (\mat{\Psi}^\tp \mat{Y})^\dagger \mat{\Psi}^\tp$.

As will be made clear from the error analysis (\cref{cor:gen-nys}), one should typically implement the generalized \Nystrom approximation using \emph{oversampling}, choosing the left embedding dimension $t$ modestly larger than the right embedding dimension $s$, say, $t = 2s$ or $t=1.5s$.

\begin{algorithm}[t]
\caption{Generalized \Nystrom approximation} \label{alg:generalized-nystrom}
\begin{algorithmic}[1]
\Require Test matrices $\mat{\Omega} \in \R^{n\times s}$, $\mat{\Psi} \in \R^{m\times t}$ ($t \ge s$); sketches $\mat{Y} \coloneqq \mat{A}\mat{\Omega} \in \R^{m\times s}$ and $\mat{W} \coloneqq \mat{A}^\tp\mat{\Psi} \in \R^{n\times t}$; tolerance $\tau \in [0,1)$
\Ensure Generalized \Nystrom approximation $\Ahat = \mat{F}\mat{G}^\tp$
\State Compute the core matrix $\mat{C} \coloneqq \mat{W}^\tp\mat{\Omega} \;(= \mat{\Psi}^\tp\mat{Y}) \in \R^{t \times s}$
\State Perform a \emph{pivoted QR decomposition} $\mat{C}\mat{P} \eqqcolon \mat{Q}\mat{R}$, with $\mat{Q}\in\R^{t \times s}$, $\mat{R}\in\R^{s\times s}$
\State Identify the set of large diagonal entries $\iset \coloneqq \{ i : |\mat{R}(i,i)| \ge \tau \cdot \max_j |\mat{R}(j,j)| \}$
\State Set $\mat{F} \coloneqq \mat{Y}\mat{P}(:,\iset)\mat{R}(\iset,\iset)^{-1}$ and $\mat{G} \coloneqq \mat{W}\mat{Q}(:,\iset)$
\end{algorithmic}
\end{algorithm}

The two sketches $\mat{A}\mat{\Omega}$ and $\mat{\Psi}^\tp\mat{A}$ required by 
the generalized \Nystrom approximation can be
accumulated simultaneously during a single streaming pass over
$\mat{A}$.
As with the \Nystrom approximation, one should take care to implement the generalized \Nystrom approximation in a stable way, particularly in evaluating the pseudoinverse of $\mat{\Psi}^\tp \mat{A}\mat{\Omega}$ \cite{nakatsukasa2020fast}.
\Cref{alg:generalized-nystrom} presents one possible stable implementation that uses a pivoted QR decomposition (CPQR, \cref{sec:classical-pivoting}).
The paper \cite{nakatsukasa2020fast} suggests an alternative stable implementation using an $\varepsilon$-truncated pseudoinverse.\footnote{Despite requiring the truncated pseudoinverse for the paper's analysis,
the algorithm of \cite{nakatsukasa2020fast} is observed to be stable in practice without truncation.}
We refer to \cite[Alg.~5]{camao2025faster} for an implementation that outputs the approximation in compact SVD form.

\begin{remark}[History]
    Versions of the generalized \Nystrom approximation have been studied going back to the twentieth century \cite{wedderburn1934lectures,chu1995rank}.
    The first references in the modern RNLA literature are \citet{woolfe2008fast,clarkson2009numerical}.
    \citet{tropp2017practical} developed stable implementations for computing the decomposition in compact SVD form, devised extensions, and ran extensive numerical testing.
    The name ``generalized \Nystrom approximation'' was proposed by \citet{nakatsukasa2020fast}, who developed fast, stable implementations for computing the decomposition in outer product form.
\end{remark}

\begin{remark}[Streaming]
    An advantage of the Nystr\"om and generalized Nystr\"om methods is that they can be implemented in a streaming setting where the matrix is provided as a sum of increments $\mat{A} = \mat{A}_1 + \mat{A}_2 + \cdots + \mat{A}_q$, with each increment arriving one by one in sequence.
    To form, say, the generalized Nystr\"om approximation from streaming data, initialize $\mat{Y} \coloneqq \mat{0}$ and $\mat{W} \coloneqq \mat{0}$ and update the sketches as each increment arrives:
    \begin{equation*}
        \mat{Y} \gets \mat{Y} + \mat{A}_i\mat{\Omega} \quad \text{and} \quad \mat{W} \gets \mat{W} + \mat{A}_i^\top\mat{\Psi}.
    \end{equation*}
    Once all increments have arrived, a low-rank approximation to $\mat{A} = \sum_i \mat{A}_i$ may be computed by \cref{alg:generalized-nystrom}.
    Critically, the algorithm processes each increment as it arrives and need not store an increment after it is processed.
    The literature on streaming algorithms for linear algebra is vast; see \cite{clarkson2009numerical,ghashami2016frequent,tropp2017fixed,henriksen2019adaoja,tropp2019streaming, braverman2020near} and the references therein.
\end{remark}

\begin{remark}[A better algorithm for the low-memory setting]
    To conclude this subsection, let us briefly mention an alternate single-pass approach due to \citet{tropp2019streaming}.
    This approach requires three sketches with the matrix ($\mat{A}\mat{\Omega}$, $\mat{\Psi}^\tp \mat{A}$, and $\mat{\Phi}^\tp \mat{A}\mat{\Upsilon}$), making it more expensive to compute than the generalized \Nystrom approximation.
    Compensating for this, the \citet{tropp2019streaming} method can achieve better accuracy than the generalized \Nystrom approximation given a fixed \emph{storage} budget.
    We recommend the generalized \Nystrom approximation for most use cases and the approach of \citet{tropp2019streaming} for memory-constrained settings.
\end{remark}

\subsection{Choice of test matrix for low-rank approximation} \label{sec:lowrank-test-matrix}

In the previous section, we saw three different approximation formats based on a
right test matrix $\mat{\Omega}$ and, for generalized \Nystrom, an
additional left test matrix $\mat{\Psi}$. By implementing these formats with different types of test matrices, we can obtain several different algorithms for low-rank approximation, each with different properties.
As a general principle, investing more effort into the construction of test matrices can improve approximation quality, at the price of increasing the cost of the resulting algorithms.

In \cref{sec:lowrank-test-matrix-random-embedding}, we begin with the simplest choice of test matrix: a random embedding.
\Cref{sec:lowrank-subspace-iteration,sec:lowrank-block-krylov} discuss ways of obtaining an improved test matrix by beginning from a random embedding $\mat{\Omega}_{\rm init}$ and improving it by powering, $\mat{\Omega} = (\mat{A}^\tp \mat{A})^q \mat{\Omega}_{\rm init}$, or by enlarging to a block Krylov subspace $\mat{\Omega} = [\mat{\Omega}_{\rm init}, (\mat{A}^\tp \mat{A})\mat{\Omega}_{\rm init}, \ldots, (\mat{A}^\tp \mat{A})^q \mat{\Omega}_{\rm init}]$.
These approaches lead to better-quality approximations at increased cost.
Finally, in \cref{sec:lowrank-coordinate-sampling}, we discuss \emph{coordinate test matrices} $\mat{\Omega} = \Id(:,\isettwo)$ that select a family of row or column indices $\isettwo$.
Methods to choose the index set will be discussed later in this survey (\cref{sec:selection}).

\subsubsection{Random embedding} \label{sec:lowrank-test-matrix-random-embedding}

The simplest choice for the test matrix $\mat{\Omega}$ is one of the random embeddings discussed in \cref{sec:embeddings}.
This choice is computationally cheap, as the test matrix is constructed independently of the matrix $\mat{A}$, and the resulting  approximations are particularly effective when the matrix $\mat{A}$ has rapidly decaying singular values.
For problems with slower singular value decay, using powering (\cref{sec:lowrank-subspace-iteration}) or block Krylov iteration (\cref{sec:lowrank-block-krylov}) can lead to substantial improvements in accuracy.

When a random embedding is used to construct a projection approximation, the resulting algorithm is often called the \emph{randomized SVD} algorithm \cite{halko2011finding}; see \cref{alg:rsvd}.
For the randomized SVD, the test matrix $\mat{\Omega}$ is typically chosen to be a Gaussian matrix or other dense random embedding.
There are limited benefits to implementing the randomized SVD with a structured embedding, as it requires two passes over the matrix $\mat{A}$, once to compute $\mat{Y} = \mat{A}\mat{\Omega}$ and a second time to compute $\mat{B} = \mat{A}^\tp \mat{Q}$.
The first pass can be accelerated by using a structured embedding, but the matrix $\mat{Q}$ is typically unstructured, so the second pass cannot be accelerated.

When a random embedding is used to construct a \Nystrom approximation, the resulting algorithm is often called the \emph{single-pass \Nystrom algorithm}, as it requires only a single pass over $\mat{A}$ to compute; see \cite{gittens2016revisiting,tropp2017fixed}.
Similarly, the generalized \Nystrom approximation can also be implemented in a single pass; see \cite{tropp2017practical,nakatsukasa2020fast} for discussion of the generalized \Nystrom approximation in modern RNLA.
For these single-pass algorithms, using a structured embedding can lead to substantial computational savings, so their use is generally recommended on large problems.

\subsubsection{Powering} \label{sec:lowrank-subspace-iteration}

The ideal choice of test matrix for low-rank approximation contains the dominant right singular vectors or, more generally, a basis for the dominant right singular subspace.
With this choice of $\mat{\Omega}$, both the projection approximation and the \Nystrom approximation coincide with the optimal rank-$s$ approximation $\lowrank{\mat{A}}_s$.
Of course, this choice is unrealistic since it requires knowledge of the dominant singular vectors of $\mat{A}$, which is precisely what we are trying to compute.
Nevertheless, we can improve the quality of the approximation by choosing the test matrix $\mat{\Omega}$ so that its range is more closely aligned with the leading singular vectors of $\mat{A}$.

The simplest way to do this is by powering \cite{rokhlin2010randomized,halko2011finding,gu2015subspace,saibaba2019randomized,tropp2023randomized}. 
First, consider a general (possibly rectangular) matrix $\mat{A}$, and pick a random initial matrix $\mat{\Omega}_{\rm init} \in \R^{n \times s}$ such as a Gaussian embedding.
To construct a test matrix that is aligned with the dominant right singular vectors of $\mat{A}$, we apply powering
\begin{equation} \label{eq:poweredomega}
    \mat{\Omega} = (\mat{A}^\tp \mat{A})^{q} \mat{\Omega}_{\rm init}.
\end{equation}
Here, $q\ge 0$ is a parameter that controls the amount of powering.
When $q=0$, we have $\mat{\Omega} = \mat{\Omega}_{\rm init}$, so no powering is applied.

To compute $\mat{\Omega}$, we apply $\mat{A}$ and $\mat{A}^\tp$ successively.
In particular, we do not form the matrix $\mat{A}^\tp \mat{A}$ or its powers explicitly, as this would be expensive and numerically unstable. For numerical stability, one should orthonormalize after  application of $\mat{A}$ and/or $\mat{A}^\tp$, e.g., 
\begin{equation*}
    \mat{\Omega} \gets \mat{\Omega}_{\rm init},\ \underbrace{\mat{\Omega} \gets \orth(\mat{A}^\tp (\mat{A} \orth(\mat{\Omega}))),\ldots,\mat{\Omega} \gets \orth(\mat{A}^\tp (\mat{A} \orth(\mat{\Omega})))}_{\text{$q$ times}}.
\end{equation*}
Owing to its connection to the classical subspace iteration method for computing eigenvectors, this approach is often called \emph{randomized subspace iteration} (RSI), particularly when combined with the projection approximation.
When combined with the \Nystrom approximation, such an approach has been called \emph{\Nystrom with subspace iteration} (NysSI, \citet{tropp2023randomized}). Note that it is possible and common to replace $\mat{A}^\tp \mat{A}$ by $\mat{A}$ in the algorithm above when $\mat{A}$ is psd.

\subsubsection{Block Krylov iteration} \label{sec:lowrank-block-krylov}

While powering is a simple way to improve the quality of a test matrix $\mat{\Omega}_{\rm init}$, it appears to be wasteful to construct the entire block Krylov sequence 
\begin{equation*}
    \mat{\Omega}_{\rm init}, (\mat{A}^\tp \mat{A}) \mat{\Omega}_{\rm init}, (\mat{A}^\tp \mat{A})^2 \mat{\Omega}_{\rm init}, \ldots, (\mat{A}^\tp \mat{A})^q \mat{\Omega}_{\rm init},
\end{equation*}
but only use the last block $(\mat{A}^\tp \mat{A})^q \mat{\Omega}_{\rm init}$ for the approximation.
A more effective way to use the block Krylov sequence is to use the entire block Krylov sequence as the test matrix \cite{rokhlin2010randomized,musco2015randomized,tropp2023randomized}:
\begin{equation*}
    \mat{\Omega} = [\mat{\Omega}_{\rm init}, (\mat{A}^\tp \mat{A}) \mat{\Omega}_{\rm init}, (\mat{A}^\tp \mat{A})^2 \mat{\Omega}_{\rm init}, \ldots, (\mat{A}^\tp \mat{A})^q \mat{\Omega}_{\rm init}].
\end{equation*}
This approach is often called \emph{randomized block Krylov iteration} (RBKI), particularly when combined with the projection approximation.
For a \psd matrix, we use the block Krylov sequence
\begin{equation*}
    \mat{\Omega} = [\mat{\Omega}_{\rm init}, \mat{A} \mat{\Omega}_{\rm init}, \mat{A}^2 \mat{\Omega}_{\rm init}, \ldots, \mat{A}^q \mat{\Omega}_{\rm init}],
\end{equation*}
and the resulting algorithm is called \emph{\Nystrom with block Krylov iteration} (NysBKI, \citet{tropp2023randomized}).

To implement the block Krylov approach, we generate a random embedding $\mat{\Omega}_{\rm init}$ and then apply $\mat{A}$ and $\mat{A}^\tp$ to it successively to form the block Krylov sequence.
For numerical stability, it is essential to orthogonalize the block Krylov sequence as it is generated.
In principle, this orthogonalization could be done efficiently by using the block Lanczos method \cite[Alg.~9.2]{chen2024lanczos}, but it is simpler and more numerically stable to apply full orthogonalization to the entire block Krylov sequence as it is generated.
See \cite{tropp2023randomized} for stable and efficient implementations.

\begin{remark}[How to pick the block size?] \label{rem:block-size}
    To compute a low-rank approximation using RBKI, there are two natural strategies for choosing the block size $\mathsf{b}$ of the starting block $\mat{\Omega}_{\rm init} \in \R^{n \times \mathsf{b}}$ given a fixed budget of matrix--vector products: (1) Use large $\mathsf{b}$ and a small number of rounds $q$ or (2) use small $\mathsf{b}$ and run for many rounds $q$.
    One argument for a larger block size $\mathsf{b}$ is computational.
    Multiplying $\mat{A}$ with a block of $\mathsf{b}$ vectors can be much faster in practice than $\mathsf{b}$ sequential matrix--vector products \cite[Fig.~2]{tropp2023randomized}.
    On the other hand, for matrices $\mat{A}$ that do not have clustered singular values, small block sizes can give more accurate approximations for a fixed target rank \cite[Fig.~1]{meyer2024unreasonable}.
    On the basis of these competing effects, the optimal value in practice tends to lie somewhere in the middle, with the precise optimal value depending on the machine, problem size, and storage format for $\mat{A}$.
    In our experience, when the target rank is sufficiently large, a value between 10 and 100 seems to be a good default. Theoretical analyses for RBKI with moderate block sizes can be found in \cite{meyer2024unreasonable,chen2026does,shao2025structural}.
\end{remark}

\subsubsection{Coordinate test matrices via index selection} \label{sec:lowrank-coordinate-sampling}

In many applications, it is desirable to compute a low-rank approximation to a matrix $\mat{A}$ that is built from selected columns and/or rows of $\mat{A}$.
Such an approximation could be more interpretable than a general low-rank approximation \cite{mahoney2009cur}, it may be more efficient to compute \cite{musco2017sublinear,chen2025randomly}, and it may be more efficient to use in certain downstream computations \cite{oseledets2010tt,martinsson2019fast}.
While these approximations fit into the general framework of low-rank approximation through coordinate test matrices $\mat{\Omega} = \Id(:,\isettwo)$ (see \cref{tab:lowrank-algorithms}), their importance merits a separate terminology and treatment. In particular, a projection approximation based on a selected subset of columns is called an \emph{interpolative decomposition} and discussed in~\cref{sec:selection}.
The corresponding \Nystrom approximation is called a \emph{column \Nystrom approximation} (\cref{sec:column-nystrom}), while a generalized \Nystrom approximation based on selected columns \emph{and} rows becomes \emph{CUR} or \emph{cross approximation} (\cref{sec:cur}).
As we will see, these approximations are closely connected to partial QR, Cholesky, and LU factorizations with pivoting, respectively.

\subsection{Error analysis} \label{sec:lowrank-error-analysis}

The mathematical analysis of randomized low-rank approximations is an entire industry unto itself, and we will not attempt to give a comprehensive overview here.
Instead, we will summarize a few key results that illustrate some important points, and we will provide pointers to the relevant literature.

\subsubsection{Analysis of the randomized SVD} \label{sec:analysisrandsvd}

The most fundamental randomized low-rank approximation is the randomized SVD, and several works provide analysis for it \cite[e.g.,][]{halko2011finding,witten2015randomized,tropp2023randomized,epperly2026sharp,camao2025faster,woodruff2014sketching,saibaba2025randomized}.
Here is one representative result.

\begin{theorem}[Randomized SVD] \label{thm:rsvd}
    For a matrix $\mat{A} \in \R^{m\times n}$ and a random embedding $\mat{\Omega} \in \R^{n\times s}$, let $\mat{Q}=\Orth(\mat{A}\mat{\Omega})$ and fix a target rank $r\ge 1$.
    Then for the randomized SVD approximation $\mat{Q}\mat{Q}^\tp \mat{A}$:
    \begin{enumerate}[label=(\alph*)]
        \item \textbf{Gaussian embedding \cite[Thm.~10.5]{halko2011finding}.} If $\mat{\Omega}$ is a Gaussian random embedding with $s\ge r+2$, then 
        \begin{equation} \label{eq:hmtexpectedvalue}
            \expect\left[\norm{\mat{A} - \mat{Q}\mat{Q}^\tp \mat{A}}_{\rm F}^2 \right]\le \left( 1 + \frac{r}{s-r-1} \right) \norm{\mat{A} - \lowrank{\mat{A}}_r}_{\rm F}^2.
        \end{equation}
        This bound is sharp in the limit $m,n\to\infty$ for a worst-case matrix $\mat{A}$ \cite[Thm.~3.3]{epperly2026sharp}. \label{item:rsvd-gaussian}
        \item \textbf{Subspace injection \cite[Thm.~2.2]{camao2025faster}.}
        If $\mat{\Omega}$ is an oblivious subspace injection with dimension $r$, injectivity parameter $\eta$, and success probability 3/4, then 
        \begin{equation} \label{eq:rsvd-osi}
            \Median\big(\norm{\mat{A} - \mat{Q}\mat{Q}^\tp \mat{A}}_{\rm F}^2\big) \le 5(1+\eta) \norm{\mat{A} - \lowrank{\mat{A}}_r}_{\rm F}^2.
        \end{equation} \label{item:rsvd-osi}
    \end{enumerate}
\end{theorem}

\begin{proof}
Begin by decomposing $\mat{A}$ using an SVD
\begin{equation*}
    \mat{A} = 
        \begin{tikzpicture}[node distance=0.3mm,baseline=-1mm]
        \node[draw, rectangle, solid, fill=blue!10, minimum width=\rowsize/1.3, minimum height=\colsize,inner xsep=0pt] (U) {\raisebox{-3.5pt}{\smash{$\Scale[.9]{\mat{U}_1}$}}};
        \node[draw, rectangle, solid, fill=blue!10, minimum width=\rowsize/1.3, minimum height=\rowsize/1.3,inner ysep=0pt,inner xsep=0pt,right =of U] (S) {\raisebox{-3.5pt}{\smash{$\Scale[.8]{\mat{\Sigma}_1}$}}};
        \node[draw, rectangle, solid, fill=blue!10, minimum width=\colsize, minimum height=\rowsize/1.3,inner ysep=0pt,right =of S] (V) {\raisebox{-3.5pt}{\smash{$\Scale[.8]{\mat{V}_1^{\tp}}$}}};
        \end{tikzpicture}
    +\mat{E} \quad \text{with } \mat{E} = \mat{U}_2\mat{\Sigma}_2\mat{V}_2^{\tp}.
\end{equation*}
We have separated the SVD into two parts: the dominant rank-$r$ part $\mat{U}_1\mat{\Sigma}_1\mat{V}_1^{\tp}$ and the trailing part $\mat{E}=\mat{U}_2\mat{\Sigma}_2\mat{V}_2^{\tp}$. 
The matrix $\mat{Y} \coloneqq \mat{A}\mat{\Omega}$ is then given by
\begin{equation*}
    \mat{Y} =  
        \begin{tikzpicture}[node distance=.3mm,baseline=-1mm]
        \node[draw, rectangle, solid, fill=blue!10, minimum width=\rowsize/1.3, minimum height=\colsize,inner xsep=0pt] (U) {\raisebox{-3.5pt}{\smash{$\Scale[.9]{\mat{U}_1}$}}};
        \node[draw, rectangle, solid, fill=blue!10, minimum width=\rowsize/1.3, minimum height=\rowsize/1.3,inner ysep=0pt,inner xsep=0pt,right =of U] (S) {\raisebox{-3.5pt}{\smash{$\Scale[.8]{\mat{\Sigma}_1}$}}};
        \node[draw, rectangle, solid, fill=blue!10, minimum width=\rowsize/2, minimum height=\rowsize/1.3,inner ysep=0pt,right =of S] (V) {\raisebox{-3.5pt}{\smash{$\Scale[.8]{\mat{V}_1^{\tp}\mat{\Omega}}$}}};
        \end{tikzpicture}
        + \mat{E}\mat{\Omega}.
\end{equation*}
To continue the analysis, we assume 
that $\mat{V}_1^{\tp}\mat{\Omega}$ has full row rank.
We shall verify this property later.

From \cref{eq:projection-approximation-optimality}, we know that the projection approximation is the optimal approximation to $\mat{A}$ among all approximations spanned by $\mat{Y}$.
Therefore, any low-rank approximation of the form $\mat{Y}\mat{Z}^\tp$ must achieve larger error:
\begin{equation*}
    \norm{\mat{A} - \mat{Q}\mat{Q}^\tp \mat{A}}_{\rm F}^2 \le \norm{\mat{A} - \mat{Y}\mat{Z}^\tp}_{\rm F}^2.
\end{equation*}
To instantiate this bound, we want to find a matrix $\mat{Z}$ that makes $\mat{Y}\mat{Z}^\tp$ as close as possible to the leading rank-$r$ part $\mat{U}_1\mat{\Sigma}_1\mat{V}_1^{\tp}$ of $\mat{A}$.
To this end, we choose $\mat{Z}^\tp \coloneqq (\mat{V}_1^{\tp}\mat{\Omega})^\dagger \mat{V}_1^{\tp}$, obtaining
\begin{equation*}
    \mat{Y}\mat{Z}^\tp = 
   \begin{tikzpicture}[node distance=.3mm,baseline=-1mm]
    \node[draw, rectangle, solid, fill=blue!10, minimum width=\rowsize*1.3, minimum height=\colsize] (U) {$\!\Scale[.8]{\mat{A}\mat{\Omega}}\!$};
    \node[draw, rectangle, solid, fill=blue!10, minimum width=\rowsize, minimum height=\rowsize*1.3, right=of U] (S) {$\!\!\Scale[0.6]{(\mat{V}_1^{\tp}\mat{\Omega})^\dagger}\!\!$};
    \node[draw, rectangle, solid, fill=blue!10, minimum width=\colsize, minimum height=\rowsize/1.2, inner ysep=0pt, right=of S] (V) {\raisebox{-3.5pt}{\smash{$\Scale[.8]{\mat{V}_1^{\tp}}$}}};
\end{tikzpicture}
=\mat{U}_1\mat{\Sigma}_1 \mat{V}_1^{\tp} + \tilde{\mat{E}} \quad \text{with} \quad \tilde{\mat{E}} \coloneqq \mat{E}\mat{\Omega}(\mat{V}_1^{\tp} \mat{\Omega})^\dagger \mat{V}_1^{\tp}.
\end{equation*}
Consequently, the error of the randomized SVD is at most
\begin{align*}
    \norm{\mat{A} - \mat{Q}\mat{Q}^\tp \mat{A}}_{\rm F}^2
    & \le \norm{\mat{A} - \mat{Y}\mat{Z}^\tp}_{\rm F}^2 = \norm{\smash{\mat{E} - \tilde{\mat{E}}}}_{\rm F}^2 \\
    & = \norm{\mat{U}_2\mat{\Sigma}_2\mat{V}_2^\tp - \mat{U}_2\mat{\Sigma}_2\mat{V}_2^\tp\mat{\Omega}(\mat{V}_1^{\tp} \mat{\Omega})^\dagger \mat{V}_1^{\tp}}_{\rm F}^2.
\end{align*}
In the last line, we inserted the definitions of $\mat{E}$ and $\tilde{\mat{E}}$.
Now, use the fact that $\mat{V}_2$ and $\mat{V}_1$ are orthogonal to each other and the unitary invariance of the Frobenius norm to deduce
\begin{equation} \label{eq:randomized-svd-structural-bound}
    \norm{\mat{A} - \mat{Q}\mat{Q}^\tp \mat{A}}_{\rm F}^2 \le \norm{\smash{\mat{\Sigma}_2^{\vphantom{\tp}}}}_{\rm F}^2 + \norm{\smash{\mat{\Sigma}_2^{\vphantom{\tp}}\mat{V}_2^\tp\mat{\Omega}(\mat{V}_1^{\tp} \mat{\Omega})^\dagger}}_{\rm F}^2.
\end{equation}
The bound \cref{eq:randomized-svd-structural-bound} is often referred to as a \emph{structural bound} for the error of the randomized SVD, as it holds for any choice of test matrix $\mat{\Omega}$ (under the full-rank assumption on $\mat{V}_1^{\tp}\mat{\Omega}$) and is deduced purely from linear algebraic arguments.
Observe that $\norm{\mat{\Sigma}_2}_{\rm F}^2 = \norm{\mat{A} - \lowrank{\mat{A}}_r}_{\rm F}^2$ is the error of the optimal rank-$r$ approximation $\lowrank{\mat{A}}_r$.

We now complete the analysis separately for cases \ref{item:rsvd-gaussian} and \ref{item:rsvd-osi}.
For case \ref{item:rsvd-gaussian} (Gaussian embedding), rotational invariance implies that the matrices $\mat{G}_1 \coloneqq \mat{V}_1^{\tp} \mat{\Omega} \in \R^{r\times s}$ and $\mat{G}_2 \coloneqq \mat{V}_2^\tp\mat{\Omega} \in \R^{(n-r)\times s}$ both have iid $\Normal(0,1/s)$ entries, and they are independent of each other.
Moreover, $\mat{G}_1$ has full row rank with probability 1.
The structural bound \cref{eq:randomized-svd-structural-bound} takes the form
    \begin{equation*}
        \norm{\mat{A} - \mat{Q}\mat{Q}^\tp \mat{A}}_{\rm F}^2 \le \norm{\smash{\mat{\Sigma}_2^{\vphantom{\tp}}}}_{\rm F}^2 + \norm{\smash{\mat{\Sigma}_2^{\vphantom{\tp}}\mat{G}_2^{\vphantom{\tp}}\mat{G}_1^\dagger}}_{\rm F}^2.
    \end{equation*}
    The expectation of the right-hand side is known in closed form \cite{halko2011finding}, which leads immediately to \cref{eq:hmtexpectedvalue}.

For case \ref{item:rsvd-osi}, the subspace injection property implies,
with probability at least $3/4$, that  $\mat{V}_1^{\tp} \mat{\Omega}$ has full row rank and $\norm{(\mat{V}_1^{\tp} \mat{\Omega})^\dagger}^2 = 1/\sigma_{\rm min}(\mat{\Omega}^\top\mat{V}_1)^2 \le 1+\eta$.
Thus, by the structural bound \cref{eq:randomized-svd-structural-bound} and the submultiplicativity of the Frobenius and spectral norms, we have
    \begin{equation} \label{eq:structosi}
        \norm{\mat{A} - \mat{Q}\mat{Q}^\tp \mat{A}}_{\rm F}^2 \le \norm{\mat{\Sigma}_2}_{\rm F}^2 + \norm{\mat{\Sigma}_2\mat{V}_2^\tp\mat{\Omega}(\mat{V}_1^{\tp} \mat{\Omega})^\dagger}_{\rm F}^2 \le \norm{\mat{\Sigma}_2}_{\rm F}^2 + (1+\eta)\norm{\mat{\Sigma}_2\mat{V}_2^\tp\mat{\Omega}}_{\rm F}^2.
    \end{equation}
    Now invoke isotropy and Markov's inequality to conclude that
    \begin{equation*}
        \norm{\mat{A} - \mat{Q}\mat{Q}^\tp \mat{A}}_{\rm F}^2 \le \norm{\mat{\Sigma}_2}_{\rm F}^2 + 4(1+\eta)\norm{\mat{\Sigma}_2\mat{V}_2^\tp}_{\rm F}^2 \le 5(1+\eta) \norm{\mat{A} - \lowrank{\mat{A}}_r}_{\rm F}^2
    \end{equation*}
    holds with probability at least 1/2.
    The conclusion \cref{eq:rsvd-osi} follows.
\end{proof}

As discussed in \cref{sec:lowrank-generalized-nystrom}, the generalized Nystr\"om approximation can be understood as a sketch-and-solve approximation to the output of the randomized SVD algorithm.
Therefore, combining the randomized SVD analysis (\cref{thm:rsvd}) with a bound for sketch-and-solve (\cref{thm:sketch-solve-gaussian} or \ref{thm:sketch-and-solve-ose-osi}) immediately yields an error bound for the generalized Nystr\"om approximation.
Here is one example \cite[Thm.~4.3]{tropp2017practical}.

\begin{corollary}[Generalized Nystr\"om approximation] \label{cor:gen-nys}
Let $\Ahat$ denote the generalized Nystr\"om approximation to a matrix $\mat{A} \in \R^{m\times n}$ generated by Gaussian embeddings $\mat{\Omega} \in \R^{n\times s}$ and $\mat{\Psi} \in \R^{m\times t}$.
Fix a target rank $r$, and assume $t \ge s+2 \ge r+4$.
Then
\begin{equation*}
    \expect \norm{\mat{A} - \Ahat}_{\rm F}^2 \le \left( 1 + \frac{s}{t-s-1} \right) \left( 1 + \frac{r}{s-r-1} \right) \norm{\mat{A} - \lowrank{\mat{A}}_r}_{\rm F}^2.
\end{equation*}
\end{corollary}

\subsubsection{Analysis of randomized subspace iteration and block Krylov iteration}

Subspace iteration and block Krylov iteration substantially improve the quality of a randomized low-rank approximation.
Again, there is a large body of work on the analysis of these methods \cite[e.g.,][]{rokhlin2010randomized,halko2011finding,gu2015subspace,musco2015randomized,bakshi2022low,tropp2022randomized,tropp2023randomized,chen2026does}.
\Cref{thm:subspacekrylov} below is one representative result, which is a simplified version of \cite[Thm.~9.2]{tropp2023randomized}.

\begin{theorem}[Subspace iteration and block Krylov iteration] \label{thm:subspacekrylov}
    Fix integers satisfying $r < k < s$, and fix $p\in [2,\infty]$.
    For a matrix $\mat{A} \in \R^{m\times n}$ of rank at least $r$, define the $r$-to-$k$ singular value gap 
    \begin{equation*}
        \mathrm{gap} \coloneqq \frac{\sigma_r(\mat{A}) - \sigma_k(\mat{A})}{\sigma_r(\mat{A}) + \sigma_k(\mat{A})} \in [0,1],
    \end{equation*}
    and let $\mat{\Omega}_{\rm init} \in \R^{n\times s}$ be a Gaussian embedding.
    Then
    \begin{enumerate}[label=(\alph*)]
        \item \textbf{Subspace iteration.}
        If $\mat{\Omega} = (\mat{A}^\tp \mat{A})^q\mat{\Omega}_{\rm init}$ and $\Ahat = (\mat{A}\mat{\Omega})(\mat{A}\mat{\Omega})^\pinv\mat{A}$ is the associated projection approximation, then
        \begin{equation*}
            \expect \norm{\mat{A} - \Ahat}_{(p)}^2 \le \norm{\mat{A} - \lowrank{\mat{A}}_r}_{(p)}^2 + \e^{-8q\cdot \mathrm{gap}} \cdot \frac{k-1}{s-k} \norm{\mat{A} - \lowrank{\mat{A}}_{k-1}}_{\rm F}^2
        \end{equation*}
        \item \textbf{Block Krylov iteration.}
        If $\mat{\Omega} = [\mat{\Omega}_{\rm init},(\mat{A}^\tp\mat{A})\mat{\Omega}_{\rm init},\ldots,(\mat{A}^\tp\mat{A})^q\mat{\Omega}_{\rm init}]$ and $\Ahat = (\mat{A}\mat{\Omega})(\mat{A}\mat{\Omega})^\pinv\mat{A}$ is the associated projection approximation, then
        \begin{equation*}
            \expect \norm{\mat{A} - \Ahat}_{(p)}^2 \le \norm{\mat{A} - \lowrank{\mat{A}}_r}_{(p)}^2 + 4\e^{-8q\cdot \sqrt{\mathrm{gap}}} \cdot \frac{k-1}{s-k} \norm{\mat{A} - \lowrank{\mat{A}}_{k-1}}_{\rm F}^2
        \end{equation*}
    \end{enumerate}
\end{theorem}

There are several things worth commenting on.
First, we observe that these bounds hold in any Schatten $p$-norm for $p\ge 2$.
The two most important cases are $p=2$ (Frobenius norm) and $p=\infty$ (spectral norm).
Second, notice that this result is only meaningful when $\mathrm{gap} > 0$.
Under this condition, the bound shows that RSI and RBKI produce approximations that come exponentially close in the number of iterations $q$ to achieving the error of the best rank-$r$ approximation.
The convergence rate is proportional to $\mathrm{gap}$ for RSI, and proportional to 
the \emph{square root} of $\mathrm{gap}$ for RBKI.
This demonstrates a fundamental advantage for RBKI over RSI; see \cref{fig:rsi-vs-rbki} for an illustration.
Gapless results, which make no assumptions on singular value gaps, also show a fundamental advantage of RBKI over RSI~\cite[Thm.~9.1]{tropp2023randomized}.
Last, note that \cref{thm:subspacekrylov} is a \emph{large block-size result}, in that the size $s$ of the seed matrix $\mat{\Omega}_{\rm init}$ must be larger than the target rank $r$.
As explained in \cref{rem:block-size}, it is usually possible to choose a fixed, moderate block size.
Theoretical results on such small-block Krylov methods for low-rank approximation can be found in~\cite{meyer2024unreasonable,chen2026does}.

\subsection{The Gram correspondence: Connecting projection and \Nystrom approximations} \label{sec:lowrank-gram-correspondence}

Projection and \Nystrom approximations are closely related to each other by the \emph{Gram correspondence}.
To derive this connection, let $\mat{A} \in \R^{n \times n}$ be a \psd matrix, and let $\mat{A} = \mat{B}^\tp \mat{B}$ be any factorization of $\mat{A}$ as a Gram matrix.
We call any such $\mat{B} \in \R^{m \times n}$ a \emph{Gram square root} for $\mat{A}$.
Unlike the principal matrix square root $\mat{A}^{1/2}$, a Gram square
root need not be square and is not unique.
Introduce a test matrix $\mat{\Omega}$ and, for simplicity, assume that $\mat{B}\mat{\Omega}$ has full column rank.
Then, consider the projection approximation $\Bhat \coloneqq \mat{B}\mat{\Omega}(\mat{B}\mat{\Omega})^\dagger \mat{B}$ to the Gram square root $\mat{B}$, and reassemble its Gram matrix:
\begin{equation*}
    \Ahat \coloneqq \Bhat^\tp \Bhat = \mat{B}^\tp [(\mat{B}\mat{\Omega})^{\dagger}]^\tp \mat{\Omega}^\tp \mat{B}^\tp \mat{B} \mat{\Omega} (\mat{B}\mat{\Omega})^\dagger \mat{B}.
\end{equation*}
Since $\mat{B}\mat{\Omega}$ has full column rank, its pseudoinverse satisfies the identity $(\mat{B}\mat{\Omega})^\dagger = (\mat{\Omega}^\tp \mat{B}^\tp \mat{B} \mat{\Omega})^{-1} \mat{\Omega}^\tp \mat{B}^\tp$, so we can simplify $\Ahat$ as
\begin{align*}
    \Ahat &= \mat{B}^\tp \mat{B} \mat{\Omega} (\mat{\Omega}^\tp \mat{B}^\tp \mat{B} \mat{\Omega})^{-1} (\mat{\Omega}^\tp \mat{B}^\tp \mat{B} \mat{\Omega}) (\mat{\Omega}^\tp \mat{B}^\tp \mat{B} \mat{\Omega})^{-1} \mat{\Omega}^\tp \mat{B}^\tp\mat{B} \\
    &= \mat{A} \mat{\Omega} (\mat{\Omega}^\tp \mat{A} \mat{\Omega})^{-1} \mat{\Omega}^\tp \mat{A} = \mat{A} \langle \mat{\Omega} \rangle.
\end{align*}
The \Nystrom approximation $\Ahat = \mat{A} \langle \mat{\Omega} \rangle$ is the Gram matrix of the projection approximation $\Bhat$ to the Gram square root $\mat{B}$!
While we have derived this connection under the assumption that $\mat{B}\mat{\Omega}$ has full column rank, the connection holds in general.

This connection between projection and \Nystrom approximations is known as the \emph{Gram correspondence}, and it is widely used to transfer results about projection approximations to \Nystrom approximations and vice versa.
Letting
$
    \mat{\Pi}_{\mat B \mat\Omega}
$
denote the orthogonal projector onto $\range(\mat B\mat\Omega)$ and 
choosing the particular Gram square root $\mat{B} = \mat{A}^{1/2}$ leads to the \emph{projection formula}
\begin{equation*}
    \mat{A}\langle \mat{\Omega} \rangle = \mat{A}^{1/2} \mat{\Pi}_{\mat{A}^{1/2} \mat{\Omega}} \mat{A}^{1/2}
\end{equation*}
for the \Nystrom approximation, which is also widely used.
Original references in RNLA for the Gram correspondence include \cite{belabbas2009landmark,gittens2011spectral}; see \cite[\S2.6]{epperly2025make} for more history.

Note that the above expression directly proves 
the \Nystrom error is psd, as 
$\mat{A}-\mat{A}\langle\mat{\Omega}\rangle =    \mat{A}^{1/2} (\mat{I}-\mat{\Pi}_{\mat{A}^{1/2} \mat{\Omega}}) \mat{A}^{1/2}\succeq 0$ since $(\mat{I}-\mat{\Pi}_{\mat{A}^{1/2} \mat{\Omega}})$ is an orthogonal projector.
The Gram correspondence also gives an exact relation between the
error $\mat{A}-\mat{A}\langle\mat{\Omega}\rangle$ of the \Nystrom approximation and the error of the projection approximation $\mat{B}-\Bhat$ for a general Gram square root $\mat{B}$:
\begin{multline*}
    (\mat{B}-\Bhat)^{\tp}(\mat{B}-\Bhat)  =
    \bigl((\Id-\mat{\Pi}_{\mat{B}\mat{\Omega}})\mat{B}\bigr)^{\tp}
    \bigl((\Id-\mat{\Pi}_{\mat{B}\mat{\Omega}})\mat{B}\bigr) = 
    \mat{B}^{\tp}(\Id-\mat{\Pi}_{\mat{B}\mat{\Omega}})^2\mat{B} 
     \\ = \mat{B}^{\tp}(\Id-\mat{\Pi}_{\mat{B}\mat{\Omega}} )\mat{B} 
      = {\mat{B}^\tp \mat{B} - \mat{B}^\tp \mat{\Pi}_{\mat{B}\mat{\Omega}} \mat{B}} 
     = {\mat{B}^\tp \mat{B} - \Bhat^\tp \Bhat} = {\mat{A} - \mat{A}\langle \mat{\Omega} \rangle}.
\end{multline*}
The first and fifth equalities use the definition of $\Bhat = \mat{\Pi}_{\mat{B}\mat{\Omega}}\mat{B}$, the third and fifth equalities use the idempotence of projectors, and the final equality is the Gram correspondence.
Because the squared Schatten $2p$-norm of a matrix equals the 
Schatten $p$-norm of its Gram matrix, it follows that 
\[ 
    \norm{\mat{B} - \Bhat}_{(2p)}^2 = \norm{\mat{A} - \mat{A}\langle \mat{\Omega} \rangle}_{(p)}, \qquad 1\leq p\leq\infty.
\]
In particular, using that $\mat{A} - \mat{A}\langle \mat{\Omega} \rangle$ is psd, we have
\begin{align*}
    \norm{\mat{A} - \mat{A}\langle \mat{\Omega} \rangle}_* &= \tr(\mat{A} - \mat{A}\langle \mat{\Omega} \rangle) = \norm{\mat{B} - \Bhat}_{\rm F}^2, \\
    \norm{\mat{A} - \mat{A}\langle \mat{\Omega} \rangle}_{\hphantom{*}} &= \hphantom{\tr(\mat{A} - \mat{A}\langle \mat{\Omega} \rangle) = } \,\,\norm{\mat{B} - \Bhat}^2.
\end{align*}
Using these identities, we may immediately transfer the randomized SVD analysis (\cref{thm:rsvd}) to analogous results for the single-pass \Nystrom approximation.
In particular, we have the following result.

\begin{corollary}[Single-pass Nystr\"om approximation] \label{cor:nys}
Let $\mat{A}\in\R^{n\times n}$ be a \psd matrix, and let $\mat{\Omega} \in \R^{n\times s}$ be a Gaussian embedding.
Fix a target rank $r$, and assume $s\ge r+2$.
Then
\begin{equation*}
    \expect \tr(\mat{A} - \mat{A}\langle\mat{\Omega} \rangle) \le  \left( 1 + \frac{r}{s-r-1} \right) \tr\!\big(\mat{A} - \lowrank{\mat{A}}_r\big).
\end{equation*}
\end{corollary}

\subsection{Extension: Low-rank approximations to a matrix function} \label{sec:lowrank-functions}

We now consider low-rank approximation of a matrix function
$f(\mat A)$ for a psd matrix $\mat A\in\R^{n\times n}$ and a continuous function $f:[0,\infty)\to [0,\infty)$ satisfying $f(0) = 0$. The latter condition ensures that the rank of $f(\mat A)$ is not larger than the rank of $\mat{A}$. We therefore hope that $f(\mat{A})$ admits a good low-rank approximation when $\mat{A}$ itself admits one.
If this is indeed the case, it is natural to use a randomized algorithm to form a psd low-rank approximation $\Ahat \approx \mat{A}$ and use $f(\Ahat)$ as an approximation to $f(\mat{A})$. 
Such an approach allows one to completely bypass the potentially
expensive task of forming a sketch $f(\mat A)\mat\Omega$, which
requires approximating the action of $f(\mat A)$ on
vectors using, e.g., a Krylov subspace method (\cref{sec:f(A)b}).

If $f$ is \emph{operator monotone},\footnote{A function $f : [0,\infty)\to[0,\infty)$ is said to be
\emph{operator monotone} if $f(\mat{A}_1)\preceq f(\mat{A}_2)$ in the Loewner order of \psd matrices whenever $\mat{0} \preceq \mat{A}_1\preceq\mat{A}_2$.
Examples include
$f(x)=x^\alpha$ for $\alpha\in[0,1]$,
$f(x)=\log(1+x)$, and
$f(x)=\tfrac{x}{x+\mu}$ for $\mu>0$.} then error bounds for $\Ahat$ can be transferred to $f(\Ahat)$. Specifically, \citet[Thm.~1.1]{persson2025algorithm} show the following:
If $\Ahat$ has rank $r$ and $\mat A - \Ahat$ is psd, then 
$
    \norm{\mat A-\Ahat}_*
    \leq
    (1+\varepsilon)
    \norm{\smash{\mat A-\lowrank{\mat A}_r}}_*
$
implies \[
    \norm{f(\mat A)-f(\Ahat)}_*
    \leq
    (1+\varepsilon)
    \norm{\smash{f(\mat A)-\lowrank{f(\mat A)}_r}}_*.
\]
In particular, this result applies to the truncated \Nystrom approximation $\Ahat = \lowrank{\mat{A} \langle \mat{\Omega}\rangle}_r$, in which case $f(\Ahat)$ is called the \emph{fun\Nystrom approximation} to $f(\mat{A})$ \cite{persson2023randomized}. \citet{persson2026randomized} consider randomized block Krylov
methods for more general functions.

\subsection{Rank estimation} \label{sec:rankest}

In our treatment of low-rank approximation so far, we have generally
assumed that a suitable target rank $r$, and hence a sketch size
$s>r$, are known in advance.
There are two natural approaches to choosing these parameters:
(1) One can estimate the numerical rank from a relatively inexpensive
random sketch.  (2) Alternatively, one can construct the approximation
adaptively, increasing its rank until an a posteriori error certificate is met.
As the second strategy is discussed in \cref{sec:error-estimation},
we focus here on the first approach.

To estimate the numerical rank of a matrix, we can form a right sketch $\mat{A}\mat{\Omega} \in \R^{m\times s}$.
Under appropriate conditions \cite{andoni2013eigenvalues,swartworth2023optimal,meier2024fast}, the dominant $r < s$ singular values of $\mat{A}\mat{\Omega}$ will be within a modest factor of the corresponding singular values of $\mat{A}$, allowing one to estimate the rank.
To accelerate this procedure, we can draw a second random embedding $\mat{\Psi} \in \R^{m\times t}$ with $t > s$ and use the two-sided sketch $\mat{\Psi}^\top \mat{A}\mat{\Omega}$.
The two-sided sketch yields singular value estimates of accuracy comparable to the one-sided approach:
Applying a subspace embedding with distortion $\eta$ to a tall-thin matrix preserves each singular value within the multiplicative factors $(1+\eta)^{-1/2}$ and $(1+\eta)^{1/2}$~\cite{gilbert2012sketched,park2023fast}.
Alternatively, one can estimate the 
numerical rank by approximating the trace of a suitably chosen spectral
projector~\cite{ubaru2017fast}; see also \cref{sec:sde}.

\subsection{Discussion and comparison} \label{sec:lowrank-discussion}

This section has discussed a whole family of randomized low-rank
approximation methods, obtained from combining three approximation
formats---projection, \Nystrom, and generalized \Nystrom---with
different choices of test matrices.
In this section, we shall provide discussion and comparison to help address the question of which algorithm to use. We focus primarily on the combinations listed in
\cref{tab:lowrank-algorithms}, with the last row (index selection) deferred until \cref{sec:selection}.

\subsubsection{\Nystrom approximation vs.\ projection approximation for PSD matrices} \label{sec:lowrank-discussion-nystrom-projection}

For approximating a \psd matrix, which approximation format is better: the \Nystrom approximation or the projection approximation?
The following result gives a clear indication.
\begin{proposition}[\Nystrom vs.\ projection]
    Let $\mat{A} \in \R^{n\times n}$ be a \psd matrix, and let $\mat{\Omega} \in \R^{n\times s}$ be any test matrix.
    Consider the projection approximation $\Ahat \coloneqq (\mat{A}\mat{\Omega})(\mat{A}\mat{\Omega})^\pinv \mat{A}$, 
    the single-pass \Nystrom approximation $\mat{A}\langle\mat{\Omega} \rangle$,
    and the \emph{two-pass} \Nystrom approximation $\mat{A}\langle \mat{A}\mat{\Omega} \rangle$.
    Then, for any unitarily invariant norm $\uinorm{\cdot}$,
    \begin{equation*}
        \uinorm{\mat{A} - \mat{A}\langle \mat{A}\mat{\Omega}\rangle} \le \uinorm{\mat{A} - \Ahat} \le \uinorm{\mat{A} - \mat{A}\langle \mat{\Omega}\rangle}
    \end{equation*}
\end{proposition}
See \citet[Lem.~5.2]{tropp2023randomized} for the first inequality; the second inequality follows from~\cref{prop:projection}. Interestingly, when both approximations in the first inequality are truncated to rank
$r<s$, the inequality remains valid for the nuclear and Frobenius norms,
but can fail for the spectral norm~\cite[Sec.~3]{persson2025algorithm}.

This result should be interpreted carefully.
Recall that, once the test matrix $\mat{\Omega}$ is fixed, the projection approximation requires two passes to compute: one to compute $\mat{Q} = \orth(\mat{A}\mat{\Omega})$ and the second to compute $\mat{Q}^\tp \mat{A}$.
This result states that, if we have the budget to perform two passes, then we obtain lower error by using these two passes to build the two-pass \Nystrom approximation $\mat{A}\langle \mat{A}\mat{\Omega}\rangle$.
This is precisely the standard \Nystrom approximation with one step of powering.

For a \emph{fixed rank} $k$, the projection approximation $\Ahat$ is never worse and typically more accurate than the single-pass \Nystrom approximation. But the latter can be significantly cheaper, especially when a fast random embedding is used. Consequently, for a \emph{fixed computational budget}, single-pass \Nystrom approximations can often be more accurate because they can afford a larger sketch. Whether this advantage materializes depends on the eigenvalue decay 
of $\mat A$ and on the matrix access model.

\begin{figure}
    \centering

    Step-wise eigenvalue decay
    \includegraphics[width=0.99\linewidth]{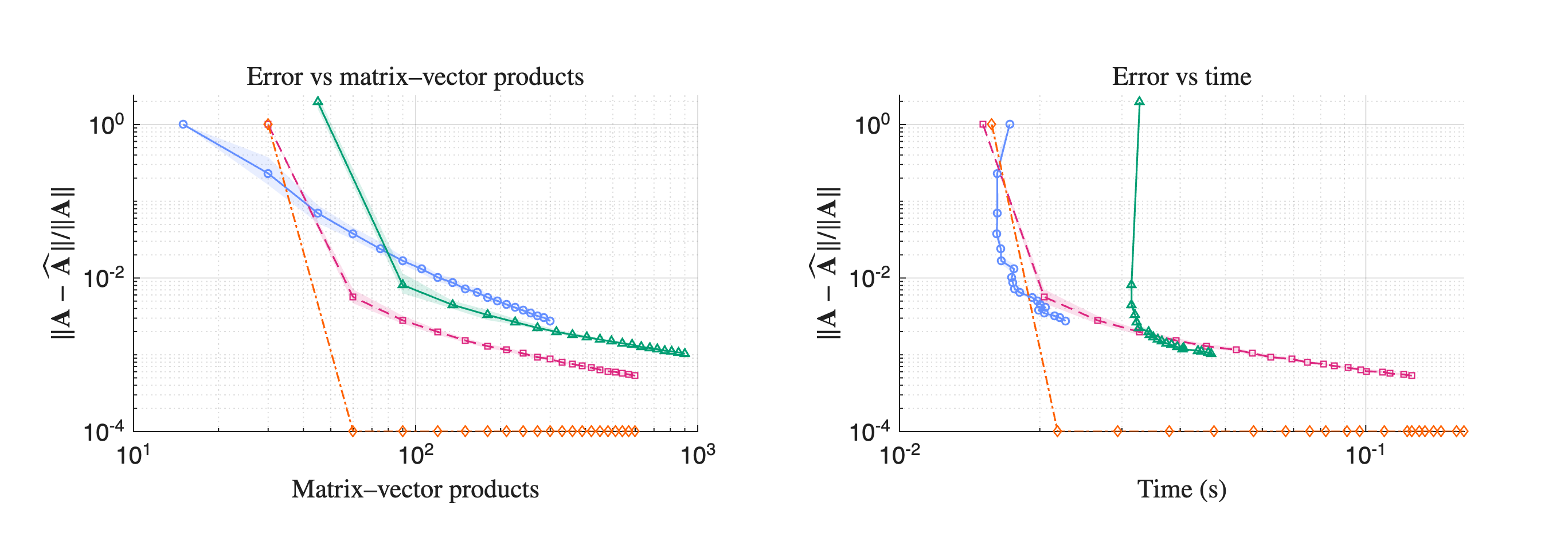}

    Inverse-square eigenvalue decay
    \includegraphics[width=0.99\linewidth]{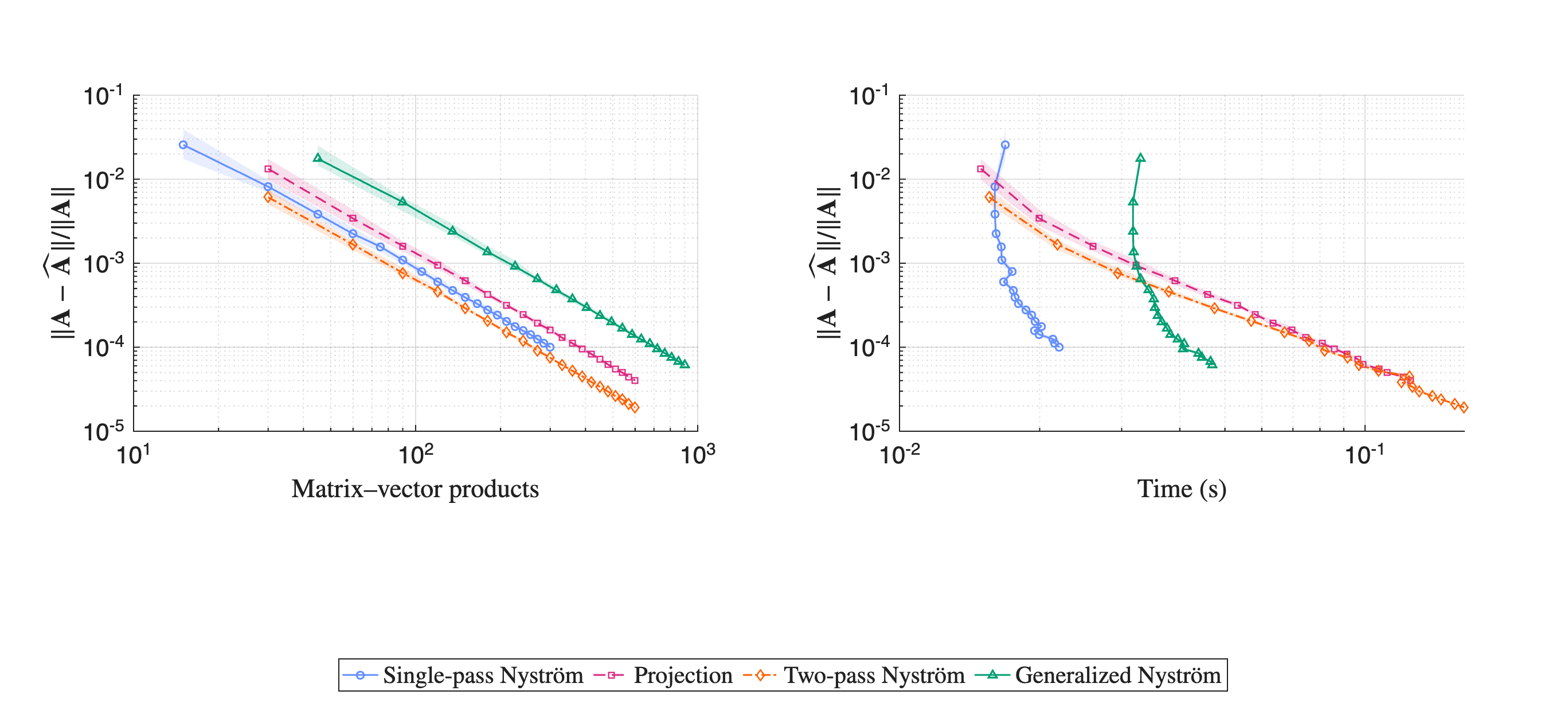}
    
    \caption{Comparison of single-pass and two-pass \Nystrom, projection
approximation, and generalized \Nystrom
on two test matrices with different eigenvalue decays. The panels show 
relative spectral norm error versus number of
matrix-vector products (left) and versus runtime (right).  }
    \label{fig:nystrom-vs-projection}
\end{figure}

These effects are demonstrated in \cref{fig:nystrom-vs-projection}, which compares the accuracy of different randomized low-rank approximations on two \psd test matrices of size $3000\times 3000$.
We include the generalized \Nystrom approximation with oversampling $t=2s$ for good measure.
The single-pass algorithms are implemented using a SparseStack embedding (\cref{sec:sparse-embeddings}), and two-pass algorithms are implemented using a Gaussian embedding.
In the top panels, we see how these algorithms perform on a matrix with step-wise eigenvalue decay
\begin{equation*}
    \lambda_i(\mat{A}) = 1 \text{ for } i=1,\ldots,20 \quad \text{and} \quad \lambda_i(\mat{A}) = 10^{-4} \text{ for } i = 21,\ldots,3000.
\end{equation*}
For a matrix with this long tail of non-trivial eigenvalues, the two-pass algorithms (both projection approximation and two-pass \Nystrom) achieve significantly better accuracy than single-pass \Nystrom for a fixed budget of matrix-vector multiplications, and the two-pass \Nystrom approximation is more accurate than all other methods by a significant margin.
The second test matrix has eigenvalues decaying at a steady rate $\lambda_i(\mat{A}) = \Theta(i^{-2})$. 
Here, we see more significant benefits of the \Nystrom approximation over the other methods, with the single- and two-pass \Nystrom approximations both being more accurate than the projection approximation for a fixed budget of matrix-vector multiplications.
This example also demonstrates the power of using fast random embeddings.
The two single-pass methods (single-pass \Nystrom and generalized \Nystrom) are meaningfully faster than both two-pass methods. 

Summarizing, there is compelling evidence to prefer \Nystrom over projection approximation for \psd matrices.
When the eigenvalues decay rapidly, single-pass \Nystrom works excellently.
For problems with slower decay, one should combine \Nystrom with powering or block Krylov approximation.

\subsubsection{No powering vs.\ powering vs.\ block Krylov} \label{sec:lowrank-discussion-powering-block-krylov}

We have three choices of how to construct test matrices in randomized low-rank approximation.
We can implement the approximation ``straight up'' using an unpowered random embedding, we can improve the quality using powering, or we can improve the quality further using block Krylov iteration. 
How do these strategies compare?

We have already seen one reason to prefer an unpowered approximation; in a single-pass approximation one can obtain substantial speedups by using a fast random embedding. For problems with rapidly decaying singular values, they can deliver excellent accuracy.

\begin{figure}
    \centering
    \includegraphics[width=0.99\linewidth]{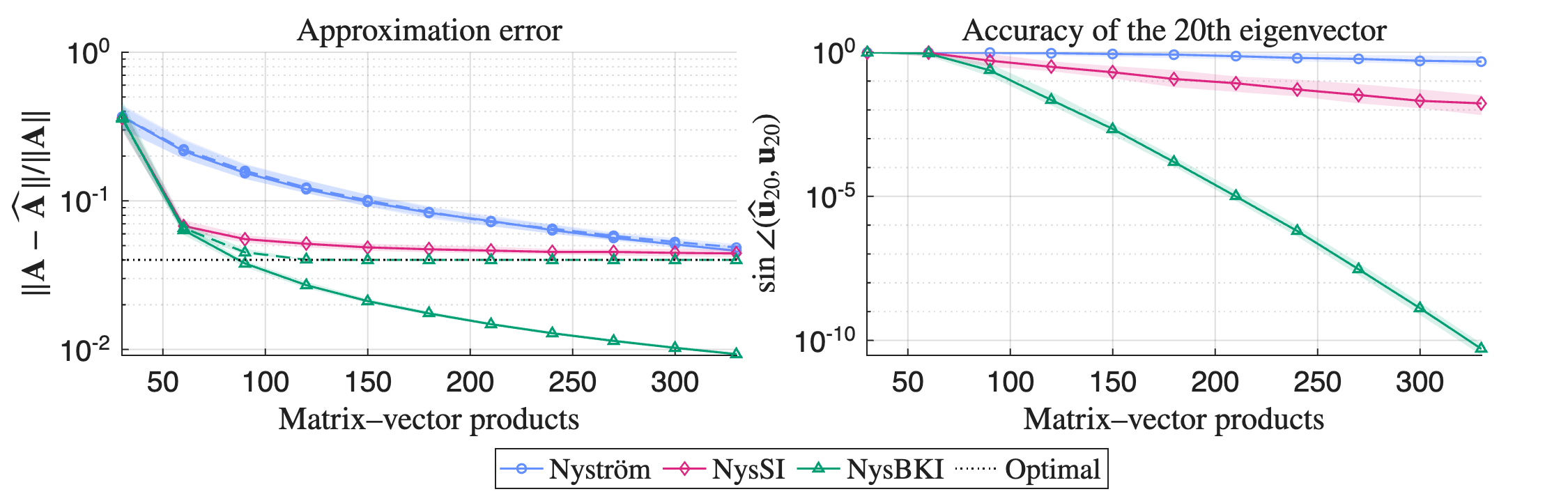}
    \caption{Comparison of single-pass \Nystrom (blue), \Nystrom with subspace iteration (NysSI, pink), and \Nystrom with block Krylov iteration (NysBKI, green).
    The panels show relative spectral norm error (left) and accuracy of the computed approximation to the 20th eigenvector (right).
    The block size for NysSI and NysBKI is 30, and the dashed lines show the error for the approximations truncated to this rank.
    The black dotted reference line is the error of the optimal rank-30 approximation.
    Lines show the median, and shaded regions show 10\% and 90\% quantiles.}
    \label{fig:rsi-vs-rbki}
\end{figure}

But when the eigenvalues or singular values decay more slowly---say, in a step-wise fashion or with a low-degree polynomial rate of decay---there can be great benefits to using block Krylov iteration.
An example is shown in \cref{fig:rsi-vs-rbki}, which evaluates the \Nystrom approximation when implemented in a single pass, with subspace iteration (NysSI), and with block Krylov iteration (NysBKI).
The matrix $\mat{A}$ is the inverse of the 2D discrete Laplacian matrix on a regular $200\times 200$ grid.
We see that NysBKI leads to substantial improvements in the low-rank approximation accuracy over the other approaches.
The difference in accuracy is even more striking when looking at the accuracy of the eigenvectors of the \Nystrom approximation as approximations to the eigenvectors of $\mat{A}$.
When implemented with 330 matrix--vector products, NysBKI computes the twentieth eigenvector to ten digits of accuracy, whereas single-pass \Nystrom and NysSI achieve zero and two digits, respectively.
The behavior of the randomized SVD, randomized subspace iteration, and RBKI is similar; see
\cite{tropp2023randomized} for many more numerical examples.

On the basis of experiments such as \cref{fig:rsi-vs-rbki}, we recommend block Krylov-based low-rank approximations as the method of choice for problems with slow or modest singular value decay if one can afford to make multiple passes over the matrix.
If memory is limited, then powering is a natural alternative.

\section{Orthogonalization and approximate orthogonalization} 
\label{sec:orthognalization}

Orthogonalization plays a central role in scientific computing; it promotes numerical stability by providing orthonormal, and hence perfectly conditioned, bases for representing subspaces.
Sketching can be used to accelerate the orthogonalization of a tall-skinny matrix. It is also able to \emph{approximately} orthonormalize a matrix (a process we refer to as \emph{whitening}) even more efficiently. 

\subsection{Whitening: Randomized approximate orthogonalization} \label{sec:whitening}

Let $\mat{A}\in\R^{m\times n}$ be a tall-skinny matrix, that is, $m\gg n$.
Assuming that $\mat{A}$ has full column rank, the standard (thin) QR factorization $\mat{A} = \mat{Q}\mat{R}$ produces a triangular matrix $\mat{R}$ with the property that $\mat{Q} = \mat{A}\mat{R}^{-1}$ is orthonormal.
If we are happy with a transformation that merely \emph{approximately orthonormalizes} or \emph{whitens} the matrix $\mat{A}$, we can do so quickly with sketching.
Draw $\mat{\Omega} \in \R^{m\times s}$, 
$s\ge n$, that serves as a subspace embedding 
for $\range(\mat{A})$. Form the sketch $\Ahat \coloneqq \mat{\Omega}^\tp \mat{A}$, and compute its thin QR factorization:
\begin{equation*}
    \Ahat = \mat{\Omega}^\tp \mat{A} = \mat{Q}_{\mat{\Omega}} \mat{R}_{\mat{\Omega}}.
\end{equation*}
The matrix $\mat{R}_{\mat{\Omega}}$ serves as a whitener.
We have the following result:

\begin{proposition}[Whitening]
    Let $\mat{\Omega}$ be a subspace embedding for $\range(\mat{A})$ with distortion $\eta$, such that 
    $\mat{\Omega}^\tp \mat{A}$ has full column rank.
    Then the whitened matrix $\mat{A}\mat{R}_{\mat{\Omega}}^{-1}$ is nearly orthonormal:
    \begin{equation*}
        (1+\eta)^{-1} \le \sigma^2_{\rm min}(\mat{A}\mat{R}_{\mat{\Omega}}^{-1}) \le \sigma^2_{\rm max}(\mat{A}\mat{R}_{\mat{\Omega}}^{-1}) \le 1+\eta
    \end{equation*}
    In particular, its condition number satisfies
$\kappa(\mat{A}\mat{R}_{\mat{\Omega}}^{-1})\le 1+\eta$.
\end{proposition}
The whitening idea and this analysis are due to \citet{rokhlin2008fast}.\footnote{We note in passing that any factorization $\Ahat = \mat{U} \cdot \mat{M}$ of the sketched matrix as a product of an orthonormal matrix $\mat{U}$ and a nonsingular matrix $\mat{M}$ gives a whitened matrix $\mat{A}\mat{M}^{-1}$ satisfying the same bounds; there are occasionally reasons to use an SVD \cite{meng2014lsrn,epperly2026fasta,epperly2026fastb} or a polar decomposition \cite{palitta2026s} instead.}
The proof is a direct invocation of the subspace embedding property.
Another view on whitening, developed by \citet{balabanov2022randomizeda}, is that it produces a matrix $\mat{A}\mat{R}_{\mat{\Omega}}^{-1}$ whose columns are orthonormal under the \emph{sketched (semi)inner product}
\begin{equation} \label{eq:sketched-inner-product}
    \langle \vec{x},\vec{y}\rangle_{\mat{\Omega}\mat{\Omega}^\tp} \coloneqq \langle \mat{\Omega}^\tp \vec{x},\mat{\Omega}^\tp \vec{y}\rangle_{\ell_2}.
\end{equation}
That is,
\begin{equation*}
    (\mat{A}\mat{R}_{\mat{\Omega}}^{-1})^\tp\mat{\Omega}\mat{\Omega}^\tp(\mat{A}\mat{R}_{\mat{\Omega}}^{-1}) = \mat{R}_{\mat{\Omega}}^{-\tp}\Ahat^\tp\Ahat \mat{R}_{\mat{\Omega}}^{-1} = \mat{Q}_{\mat{\Omega}}^\tp \mat{Q}_{\mat{\Omega}}^{\vphantom{\tp}} = \Id.
\end{equation*}
When using a constant-sparsity SparseStack (\cref{sec:sparse-embeddings}) or a standard SRFT (\cref{sec:srfts}) with $s = \order(n)$, the cost of computing the whitening matrix $\mat{R}_{\mat{\Omega}}$ is $\order(mn + n^3)$ and $\order(mn\log m + n^3)$, respectively.

There are several uses of the whitening primitive.
In the next subsection, we discuss the randomized CholeskyQR algorithm that uses whitening as a subroutine to rapidly compute a thin QR factorization.
Whitening can also be used to accelerate Krylov subspace methods (\cref{sec:approximate-orthogonalization-Krylov}) and to find an excellent preconditioner for a tall-skinny least-squares problem (\cref{sec:sketch-and-precondition}).

\subsection{Randomized CholeskyQR: Randomized orthogonalization} \label{sec:randcholqr}

CholeskyQR \cite[p.~253]{stewart1998matrix} is a cheap alternative to the thin QR factorization.
For a given input $\mat{A}$ of full column rank, the method computes the Gram matrix $\mat{G} \coloneqq \mat{A}^\tp \mat{A}$ and its Cholesky factorization $\mat{G} = \mat{R}^\tp \mat{R}$, then sets $\mat{Q} \coloneqq \mat{A}\mat{R}^{-1}$.
In exact arithmetic, the matrix $\mat{Q}$ is guaranteed to be orthonormal.
However, in finite precision, the computed factor $\widehat{\mat{Q}}$
can suffer a significant loss of orthogonality \cite{yamamoto2015roundoff}:
\begin{equation*}
    \norm{ \, \Qhat^\tp\Qhat-\Id\, } = \order\big(\kappa(\mat{A})^2\varepsilon_{\rm mach}\big).
\end{equation*}
Moreover, when $\kappa(\mat{A})$ approaches
$\varepsilon_{\rm mach}^{-1/2}$, forming the Gram matrix can destroy
positive definiteness numerically and the Cholesky factorization may
break down \cite{fukaya2020shifted}.
Because of these numerical problems, CholeskyQR has traditionally been discouraged for general use, reserved for situations where the matrix is known to be well-conditioned.

By combining CholeskyQR with whitening, we obtain a numerically stable orthogonalization scheme that maintains many of the benefits of CholeskyQR.
Randomized CholeskyQR proceeds by 
applying whitening and forming $\mat{B} \coloneqq \mat{A} \mat{R}_{\mat{\Omega}}^{-1}$.
Since $\mat{B}$ is well-conditioned, we can now safely apply CholeskyQR:
\begin{equation*}
    \mat{G} \coloneqq \mat{B}^\tp \mat{B} \quad \text{and} \quad \mat{G} = \mat{R}_{\mat{G}}^\tp \mat{R}_{\mat{G}}^{\vphantom{\tp}}.
\end{equation*}
Setting $\mat{Q} \coloneqq \mat{B}\mat{R}_{\mat{G}}^{-1}$ and $\mat{R} \coloneqq \mat{R}_{\mat{G}}\mat{R}_{\mat{\Omega}}$, we obtain a thin QR decomposition of $\mat{A}$; see \cref{alg:randomized-cholesky-qr} for pseudocode.
A prototype of randomized CholeskyQR was proposed by
\citet{fan2021novel}.  The approach was subsequently developed,
extended, and analyzed in
\cite{balabanov2022randomizedb,garrison2024randomized,
higgins2025analysis,melnichenko2025choleskyqr}.

\begin{algorithm}[t]
\caption{Randomized CholeskyQR for orthogonalizing a matrix} \label{alg:randomized-cholesky-qr}
\begin{algorithmic}[1]
\Require Tall-thin matrix $\mat{A} \in \R^{m\times n}$; sketch size $s \ge n$
\Ensure Thin QR factors $\mat{Q} \in \R^{m\times n}$ and $\mat{R} \in \R^{n\times n}$
\State Generate a subspace embedding $\mat{\Omega} \in \R^{m\times s}$
\State Form the sketch $\Ahat \gets \mat{\Omega}^\tp \mat{A}$
\State Compute a \emph{Householder QR} factorization $(\mat{Q}_{\mat{\Omega}},\mat{R}_{\mat{\Omega}}) \gets \Call{QR}{\Ahat,\texttt{'econ'}}$
\State Form the whitened matrix $\mat{B} \gets \mat{A}\mat{R}_{\mat{\Omega}}^{-1}$
\State Compute the Gram matrix $\mat{G} \gets \mat{B}^\tp \mat{B}$
\State Compute the Cholesky factor $\mat{R}_{\mat{G}}^{\vphantom{\tp}} \gets \Call{Chol}{\mat{G}}$
\State Set $\mat{Q} \gets \mat{B} \mat{R}_{\mat{G}}^{-1}$ and $\mat{R} \gets \mat{R}_{\mat{G}}\mat{R}_{\mat{\Omega}}$
\end{algorithmic}
\end{algorithm}

When $\mat{Q}$ is formed explicitly, randomized CholeskyQR (even when implemented with a fast random embedding) requires approximately $3mn^2$ floating point operations, compared with $4mn^2$
for Householder-based thin QR.
However, this often understates the practical speedup, as randomized CholeskyQR reorganizes the computation so that the dominant costs come from matrix multiplication and block triangular solves, operations that are often much faster than Householder QR on modern computing environments.
Randomized CholeskyQR is also well-suited to parallel computing \cite{balabanov2022randomizedb,higgins2025analysis}.

\begin{figure}
    \centering
    \includegraphics[width=0.99\linewidth]{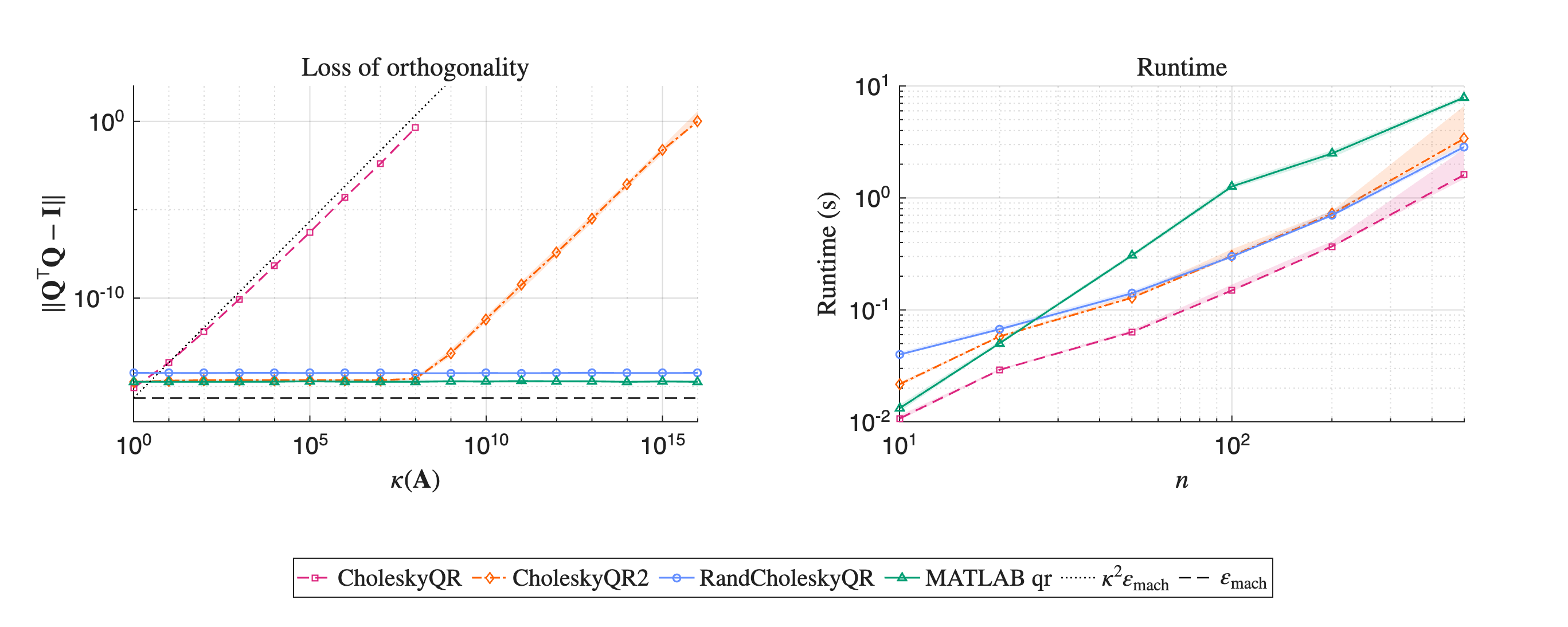}
    \caption{Comparison of MATLAB's \texttt{qr}, randomized CholeskyQR, plain CholeskyQR, and CholeskyQR2 for orthogonalizing a tall-thin matrix.
    The left panel tests the orthonormality of the $\mat{Q}$ matrix as a function of the matrix condition number $\kappa(\mat{A})$, and the right panel shows the runtime as a function of the number of columns $n$.
    Lines show the median of 30 trials, and shaded regions show 90\% and 10\% quantiles.}
    \label{fig:cholqr}
\end{figure}

\Cref{fig:cholqr} compares randomized CholeskyQR to MATLAB's \texttt{qr} command, in addition to ordinary CholeskyQR and CholeskyQR2 \cite{yamamoto2015roundoff}, another CholeskyQR variant.
We test on random matrices with geometrically spaced singular values of size $5\times 10^5$ by $n \in [10,500]$ with condition numbers $\kappa(\mat{A}) \in [10^0,10^{16}]$.
For sufficiently large matrices, randomized CholeskyQR is roughly $4\times$ faster than MATLAB's \texttt{qr}, and it reliably attains orthonormality up to working precision.
The plain CholeskyQR and CholeskyQR2 methods, by contrast, are fast but suffer from loss of orthogonality when the matrix $\mat{A}$ becomes ill-conditioned.\footnote{In order to avoid breakdown of the first Cholesky factorization in CholeskyQR2 when $\kappa(\mat{A}) > 10^8$, we implement the method by adding a tiny shift $\mat{G} \gets \mat{G} + \tr(\mat{G})\varepsilon_{\rm mach}\Id$.}

\paragraph{Extension: Oblique inner product space.}
Randomized CholeskyQR can be extended to an inner product defined by a \pd matrix $\mat{M}$, provided that whitening yields a well-conditioned basis in the corresponding $\mat{M}$-norm. 
One then forms the Gram matrix
$\mat{G}=\mat{B}^{\tp}\mat{M}\mat{B}$ and proceeds as before, obtaining
$\mat{A}=\mat{Q}_{\mat{M}}\mat{R}_{\mat{M}}$ with
$\mat{Q}_{\mat{M}}^{\tp}\mat{M}\mat{Q}_{\mat{M}}=\Id$.

\paragraph{Randomized Householder QR.}
A recent alternative is the randomized Householder QR (RHQR) factorization
of \citet{grigori2026randomized}. RHQR constructs the so-called compact WY representation of the involved Householder transformations on the sketch and then lifts it back to the original matrix.

\section{Sketching for Krylov subspace methods}
\label{sec:krylov}

As we have seen, we can use sketching to rapidly find an approximate solution to a tall-skinny least-squares problem (\cref{sec:sas}) or approximately orthogonalize a tall-skinny matrix (\cref{sec:whitening}).
The speedups are most pronounced when the aspect ratio of the matrix is significantly lopsided, with many more rows than columns.
A natural setting where this is the case is provided by Krylov subspace methods.

The role of randomization in this section is markedly different from the randomized block Krylov
iteration considered in \cref{sec:lowrank-block-krylov}: there, randomness is used to initialize a
Krylov space, whereas here sketching is used to reduce the cost of constructing and working with a Krylov basis.

\subsection{Review: Krylov subspace methods} \label{sec:krylovreview}

Krylov subspace methods are a powerful and popular class of methods for solving a plethora of linear algebra problems, including linear systems $\mat{A}\vec{x} = \vec{b}$ \cite{greenbaum1997iterative,saad2003iterative}, eigenvalue problems $\mat{A}\vec{x} =\lambda\, \vec{x}$ \cite{saad2011numerical}, and computing a matrix function times a vector $f(\mat{A})\vec{b}$ \cite[\S13.2]{higham2008functions}.
These methods find an approximate solution to these problems within the Krylov subspace
\begin{equation} \label{eq:krylov-subspace}
    \set{K}_k(\mat{A},\vec{b}) \coloneqq \spn \{ \vec{b}, \mat{A}\vec{b}, \ldots, \mat{A}^{k-1}\vec{b} \}.
\end{equation}
To simplify the discussion, we assume throughout this section that Krylov subspaces have maximal dimension $k$ for all iterations $k \le n$ under consideration.
For an eigenvalue problem, no vector $\vec{b}$ is inherent in the problem, and it is usually chosen randomly (e.g., $\vec{b} \sim \Normal(\vec{0},\Id)$) to ensure overlap with an eigenspace of interest; see also \cref{sec:smallball}. 

To frame our discussion, we turn our attention to one Krylov subspace method in particular, the generalized minimal residual method (GMRES).
We seek to find an approximate solution to a square linear system $\mat{A}\vec{x} = \vec{b}$ associated with a nonsymmetric matrix $\mat{A} \in \R^{n\times n}$.
Assuming zero initial guess, the GMRES method outputs as approximate solution the vector $\vec{x}_k \in \set{K}_k(\mat{A},\vec{b})$ within the Krylov subspace that minimizes the residual $\norm{\mat{A}\vec{x}_k - \vec{b}}$.
Given any basis matrix $\mat{B}$ for $\set{K}_k(\mat{A},\vec{b})$, this solution may be formulated as a linear least-squares problem
\begin{equation} \label{eq:gmres-as-least-squares}
    \vec{x}_k = \mat{B}\vec{y} \quad \text{with } \vec{y} \coloneqq \argmin_{\vec{z} \in \R^k} \norm{(\mat{A}\mat{B})\vec{z} - \vec{b}}.
\end{equation} 

The conventional GMRES method solves the least-squares problem \cref{eq:gmres-as-least-squares} using the Arnoldi process \cite[\S6.3]{saad2003iterative}, which iteratively builds an orthonormal basis $\mat{Q}_k \in \R^{n\times k}$, with first column $\vec{q}_1=\vec{b}/\norm{\vec{b}}$, for the Krylov subspace and a Hessenberg matrix $\mat{H}_{k+1,k} \in \R^{(k+1)\times k}$
satisfying
\begin{equation} \label{eq:arnoldi}
    \mat{A}\mat{Q}_k = \mat{Q}_{k+1} \mat{H}_{k+1,k}.
\end{equation}
From this relation, the GMRES solution is easily obtained by solving a small $(k+1)\times k$ least-squares problem
\begin{equation} \label{eq:gmres-solved}
    \vec{x}_k = \mat{Q}_k\vec{y} \quad \text{with } \vec{y} \coloneqq \argmin_{\vec{z} \in \R^k} \norm{\mat{H}_{k+1,k}\vec{z} - \norm{\vec{b}}\evec_1}.
\end{equation}
The same Arnoldi process is also used in other Krylov subspace methods beyond linear systems.

For nonsymmetric problems, the principal costs of $k$ steps of a Krylov
subspace method like GMRES are $k$ matrix--vector products with $\mat{A}$ and
$\order(nk^2)$ arithmetic operations required for orthogonalization to build $\mat{Q}_k$. When $\mat{A}$ is
sparse, possessing only $\order(n)$ entries, and $k$ is large, the cost of orthogonalization frequently dominates.
For large-scale problems, floating point operations may not be the main
bottleneck; memory consumption and the global communication required by
orthogonalization can be more important. Restarted GMRES limits the
dimension of the stored Krylov basis and hence controls these costs,
but it may also degrade performance. 
Randomization offers tools that may allow us to reduce resource usage while maintaining accuracy.

\subsection{Approximate orthogonalization for the Arnoldi process} \label{sec:approximate-orthogonalization-Krylov}

A first approach to speeding up Krylov subspace methods is to replace the traditional Arnoldi process, which fully orthogonalizes the Krylov subspace, with a sketched version that only approximately orthogonalizes (whitens) it.
The approach is due to \citet{balabanov2022randomizeda}.
We begin by drawing a random embedding $\mat{\Omega}$.
Then, we use a modified version of the Arnoldi process to orthogonalize the Krylov basis with respect to the sketched inner product \cref{eq:sketched-inner-product}.
This results in an Arnoldi-like decomposition
\begin{equation} \label{eq:approximate-arnoldi}
    \mat{A}\Qhat_k = \Qhat_{k+1} \Hhat_{k+1,k} \quad \text{with } \Qhat_{k+1}^\tp \mat{\Omega}\mat{\Omega}^\tp \Qhat_{k+1} = \Id_{k+1}.
\end{equation}
This approach is known as randomized Gram--Schmidt (RGS).

By using the Arnoldi-like decomposition \cref{eq:approximate-arnoldi} in place of the traditional Arnoldi decomposition, we obtain fast, approximate versions of Krylov subspace methods.
For instance, the RGS-GMRES method \cite[\S4.2]{balabanov2022randomizeda} produces an approximate solution to the linear system $\mat{A}\vec{x} = \vec{b}$ as follows. 
First, substitute the Arnoldi-like decomposition \cref{eq:approximate-arnoldi} into \cref{eq:gmres-as-least-squares} to see that the residual-minimizing solution is 
\begin{equation*}
    \vec{x}_k = \Qhat_k \vec{y} \quad \text{with } \vec{y} \coloneqq \argmin_{\vec{z} \in \R^k} \norm{ \Qhat_{k+1}(\Hhat_{k+1,k} \vec{z} - \beta_\Omega \evec_1) }.
\end{equation*}
with $\beta_\Omega\coloneqq\norm{\mat\Omega^\tp\vec b}$.
Instead of exactly minimizing this residual, RGS-GMRES minimizes the 
sketched residual. Using that $\mat{\Omega}^\tp\Qhat_{k+1}$ has orthonormal columns, this leads to the cheaper surrogate
\begin{equation*}
    \xhat_k = \Qhat_k \yhat \quad \text{with } \yhat \coloneqq \argmin_{\vec{z} \in \R^k} \norm{ \Hhat_{k+1,k} \vec{z} - \beta_\Omega  \evec_1 };
\end{equation*}
compare with~\cref{eq:gmres-solved}.

The RGS-GMRES method produces an approximate solution that is nearly as good as the traditional GMRES method.
Indeed, assuming that $\mat{\Omega}$ is a subspace embedding with distortion $\eta$ for $\set{K}_{k+1}(\mat{A},\vec{b}) 
    =\range(\Qhat_{k+1})$, \cref{thm:sketch-and-solve-ose-osi} (a) implies that
    the RGS-GMRES solution $\xhat_k$ achieves a residual that is within a modest factor of the solution $\vec{x}_k$ produced by the standard GMRES:
\begin{equation} \label{eq:sgmres-guarantee}
    \norm{\mat{A}\xhat_k - \vec{b}} \le (1+\eta) \norm{\mat{A}\vec{x}_k - \vec{b}}.
\end{equation}
Similarly, RGS can be combined with Arnoldi methods for eigenvalue computations and matrix functions $f(\mat{A})\vec{b}$;
see \cite{balabanov2025randomized,dedamas2025randomized,dedamas2026randomized} and, in particular, the survey by \citet{dedamas2026randomized}.

RGS-Arnoldi costs roughly $nk^2$ operations to compute the Arnoldi-like decomposition \cref{eq:approximate-arnoldi}, which is
about half the orthogonalization cost in conventional
modified Gram--Schmidt-based implementations of the Arnoldi process (MGS-Arnoldi).
More strikingly, RGS-Arnoldi requires only one global synchronization per Arnoldi step, resulting in significantly less parallel communication compared to conventional
MGS-Arnoldi.

\subsection{Sketching a non-orthogonal Krylov basis}

An alternative approach to speeding up  Krylov subspace methods is to avoid constructing an orthonormal, or even
approximately orthonormal, Krylov basis altogether and to use sketching primarily for solving subproblems.

\subsubsection{Sketched GMRES for linear systems} \label{sec:sgmres}

Standard GMRES uses Arnoldi to yield an orthonormalize basis of the Krylov subspace as we generate it.
But what if we did not do this, instead generating some non-orthogonal basis $\mat{B}$ for the Krylov subspace?
On its face, it does not seem faster, as we would still need to solve the least-squares problem \cref{eq:gmres-as-least-squares}. We can reduce this cost if we accept an approximate solution computed by sketch-and-solve (\cref{sec:sas}).
That is, we generate a random embedding $\mat{\Omega} \in \R^{n\times s}$, and return the sketched GMRES solution
\begin{equation} \label{eq:ssgmres}
    \xhat_k = \mat{B}\yhat \quad \text{for } \yhat = \argmin_{\vec{z} \in \R^k} \norm{\mat{\Omega}^\tp (\mat{A}\mat{B}) \vec z  - \mat{\Omega}^\tp \vec{b}}.
\end{equation}
With an appropriate fast random embedding (\cref{sec:fast-embedding}) and embedding dimension $s = \order(k)$, say $s = 2k+1$, the cost is $\order(nk \log k + k^3)$ operations.
This \emph{sketched GMRES} algorithm (sGMRES, \citet{nakatsukasa2024fast}) accelerates the GMRES method from requiring $\order(nk^2)$ operations to nearly $\order(nk)$ operations, when excluding the cost of matrix--vector products with $\mat{A}$.
\citet{nakatsukasa2024fast} report speedups of sGMRES over standard GMRES of up to $100\times$.

The sGMRES and RGS-GMRES methods are derived in different ways and have different computational costs ($\order(nk \log k+k^3)$ operations vs $\order(nk^2)$ operations).
Nonetheless, when using the same $\mat{\Omega}$, they produce the same approximate solution $\xhat_k$ in \emph{exact arithmetic}; see \cite[\S4.2]{dedamas2026randomized}.
In particular, if $\mat{\Omega}$ is a subspace embedding for $\set K_{k+1}(\mat A,\vec b)$ with distortion $\eta$, the sGMRES solution $\xhat_k$ satisfies the guarantee \cref{eq:sgmres-guarantee}.

To implement sGMRES, we need a cheap way of producing the non-orthogonal basis $\mat{B}$.
In exact arithmetic, this would be a triviality and we could simply build the Krylov sequence $\mat{B} = [\begin{matrix} \vec{b},\mat{A}\vec{b},\ldots,\mat{A}^{k-1}\vec{b} \end{matrix}]$.
But in finite precision, this approach typically results in a catastrophically ill-conditioned matrix $\mat{B}$, compromising the reliability of
sGMRES. 
To obtain a better-conditioned basis, a classical option is truncated Arnoldi
\cite{saad1981krylov}: for a small truncation parameter $p$, each new
Krylov vector is orthogonalized only against the $p$ most recently
generated basis vectors. \citet[\S4]{nakatsukasa2024fast} use this inexpensive $p$-partial orthogonalization as a
basis generator for sGMRES.
Even with this strategy,
$\kappa(\mat B)$ can grow rapidly with $k$,
eventually leading to numerical rank deficiency.
When this happens, we can convert $\mat{B}$ to a well-conditioned basis by applying approximate orthogonalization, but this requires the same $\order(nk^2)$ operations---leading to at most a factor-two speedup over standard GMRES. 

\begin{figure}
    \centering
    \includegraphics[width=0.99\linewidth]{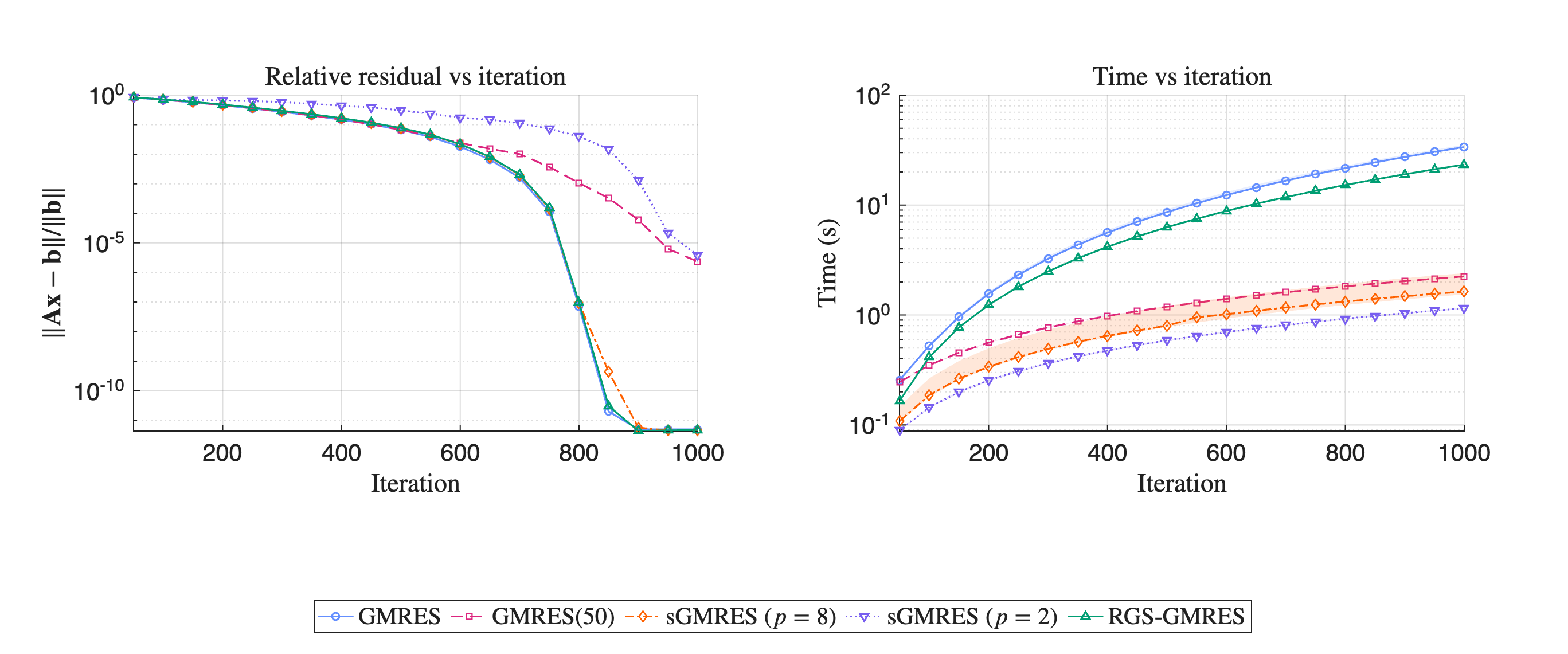}
    \caption{Comparison of GMRES, restarted GMRES(50), sGMRES with
$p=2$ and $p=8$ partial orthogonalization steps, and RGS-GMRES for an
advection--diffusion test problem. Relative residual (left) and cumulative runtime (right) versus iterations.}
    \label{fig:sgmres}
\end{figure}

\paragraph{Experimental evaluation.}

We compare sGMRES, RGS-GMRES, and traditional GMRES with and without restarting.
We test on a nonsymmetric linear system of size $n=400^2$ obtained from 
a $400\times 400$ five-point finite-difference discretization of an
advection--diffusion equation on a square domain with mesh
P\'eclet number $0.1$. Both sGMRES and RGS-GMRES use a SparseStack embedding with our recommended parameter settings $\zeta = 4$ and $s=2k$.
\Cref{fig:sgmres} shows the results.
With appropriate parameter settings, RGS-GMRES and sGMRES with $p=8$ are nearly as accurate as plain GMRES, and much more accurate than GMRES with restarting every 50 iterations.
As expected, RGS-GMRES is modestly faster ($1.5\times$ with our implementation) and sGMRES with $p = 8$ is much faster ($23\times$).
But caution is merited.
When an insufficient number of partial orthogonalization steps is used ($p=2$ in this example), the accuracy of sGMRES is substantially harmed compared to plain GMRES.

The difference between $p=2$ and $p=8$ is the conditioning of the basis.
With $p = 8$, the basis $\mat{B}$ remains numerically full rank until iteration 900 and beyond,
whereas with $p=2$ numerical rank deficiency occurs
before iteration $100$.
The conclusion is that, for serial computation, RGS-GMRES is numerically robust and leads to a modest speedup over plain GMRES.
sGMRES can lead to much larger speedups, but only if the condition number of the basis remains under control.

\paragraph{Restarting and inner-outer iterations.}
Designing inexpensive Krylov recurrences that remain reliably
well-conditioned for general nonsymmetric problems remains a
challenge.
Restarting limits the growth and, hence, the potential ill-conditioning
of the basis~\cite{burke2025numerical}, while deflation can mitigate
loss of convergence due to restarting~\cite{burke2025gmres}; see also
\citet{jang2025randomized} for randomized flexible GMRES with deflated
restarting. A complementary approach by
\citet{guttel2025stabilizing} uses sGMRES as a variable inner
preconditioner for an outer flexible GMRES iteration~\cite{saad1993flexible}.
The outer residual minimization substantially improves robustness and can
make short inner sGMRES iterations, even \emph{without} any orthogonalization ($p=0$), effective in practice.

\paragraph{Beyond Krylov subspaces.}
In sGMRES, $\mat{B}$ is chosen as a basis for a Krylov subspace. The same sketch-and-solve construction~\cref{eq:ssgmres}
applies to \emph{any} trial basis.
This viewpoint suggests investigating other ways of generating trial bases, for example, from
sketched snapshots of a dynamical system~\cite{guido2024subspace} or a neural network; see also~\Cref{sec:subappsnap}.

\subsubsection{Eigenvalue problems and sketched Rayleigh--Ritz}

Similar principles can be applied for eigenvalue problems.
Suppose that $\mat{B} \in \R^{n\times k}$ is a basis for a subspace that contains good approximations to eigenvectors of $\mat A \in \R^{n\times n}$.
To extract such eigenvector approximations, the classical Rayleigh--Ritz (RR) procedure
forms an orthonormal basis $\mat{Q} \gets \Orth(\mat{B})$ and solves
the (small) eigenvalue problem
\[
(\mat{Q}^\tp\mat{A}\mat{Q})\vec{y}_i
   =\hat{\lambda}_i\vec{y}_i, \qquad i=1,\ldots, k.
\]
The so-called Ritz pairs $(\hat{\lambda}_i,\mat{Q}\vec{y}_i)$
constitute approximate eigenpairs for $\mat{A}$, at least for some $i$ (typically the exterior ones).
As alternatives that eschew the matrix $\mat{Q}$, we can instead solve the equivalent generalized eigenvalue problem $(\mat{B}^\tp \mat{A}\mat{B})\vec{y}_i = \hat\lambda_i \mat{B}^\tp\mat{B} \vec{y}_i$ or the standard eigenvalue problem $\mat{B}^\pinv \mat{A}\mat{B}\vec{y}_i = \hat\lambda_i \vec{y}_i$. Similarity follows by noting that there is an invertible matrix $\mat{R}$ such that $\mat{B} = \mat Q \mat R$:
$
 \mat{B}^\pinv \mat{A}\mat{B} = \mat{R}^{-1} \big( \mat{Q}^\tp\mat{A}\mat{Q} \big) \mat R.
$

To involve sketching, we first note that $\mat{M}_\star = \mat{B}^\pinv \mat{A}\mat{B}$ solves the least-squares problem
\[
\min_{\mat{M}\in\R^{k\times k}}
\norm{\mat{A}\mat{B}-\mat{B}\mat{M}}_{\rm F}.
\]
Given a subspace embedding $\mat{\Omega}\in\R^{n\times s}$ for 
$\range[\mat{A}\mat{B}, \mat{B}]$, we let
\begin{equation} \label{eqn:pre-sRR-intro}
 \widehat{\mat{M}}
\coloneqq 
\argmin_{\mat{M}\in\R^{k\times k}}
\norm{\mat{\Omega}^\tp
      (\mat{A}\mat{B}-\mat{B}\mat{M})}_{\rm F}
=
(\mat{\Omega}^\tp\mat{B})^\pinv
 \mat{\Omega}^\tp\mat{A}\mat{B}.
\end{equation}
We obtain approximate eigenpairs by solving the sketched RR (sRR) eigenvalue problem 
$\hat{\mat{M}} \vec{y}_i 
	= \hat{\lambda}_i \vec{y}_i$ and returning the pairs $(\hat{\lambda}_i,\mat{B}\vec{y}_i)$.

This sRR method was proposed by \citet{nakatsukasa2024fast}.
As with sGMRES and RGS-GMRES, there is an equivalent formulation, due to \citet{balabanov2025randomized}, that uses approximate orthogonalization in place of directly solving the sketched least-squares problem \cref{eqn:pre-sRR-intro}.
The comparison of sRR to RGS-RR is much the same as the comparison of sGMRES to RGS-GMRES.

As with sGMRES, implementing sRR requires some method to obtain a well-conditioned basis $\mat{B}$ for the subspace of interest.
As another drawback, sRR introduces additional non-normality; in particular, it is not guaranteed to produce real eigenvalues and orthogonal eigenvectors when applied to symmetric $\mat{A}$.
This issue can be addressed by restoring similarity to the
standard Rayleigh--Ritz matrix, as done efficiently for
RGS-Arnoldi in~\cite{grigori2026restoring}.

\begin{remark}[Another approach to adding randomization to Rayleigh--Ritz]
 A rather different use of randomization was recently proposed by
\citet{shao2026stabilizing}: Given an $n\times k$ orthonormal basis $\mat{Q}$, an $n\times k$ Gaussian embedding  $\mat{\Omega}$ is drawn \emph{without oversampling} and
randomized Ritz pairs are extracted from
$
(\mat{\Omega}^\tp \mat{A}\mat{Q})\vec{y}
=
\widehat{\lambda}(\mat{\Omega}^\tp \mat{Q})\vec{y}.
$
Here, randomization is not used for dimension reduction but for stabilizing
eigenpair extraction, yielding probabilistic convergence guarantees for
recovering eigenvector approximations contained in $\range(\mat{Q})$. 
\end{remark}

\subsubsection{Computing matrix function times a vector $f(\mat{A})\vec{b}$} \label{sec:f(A)b}

The idea of sketching a (Krylov) subspace to speed up computation has been extended to computing $f(\mat{A})\vec{b}$ for a function $f$ defined on a domain containing the numerical range of $\mat{A}$. 
Let us consider, once more, the Arnoldi decomposition $\mat{A}\mat{Q}_k=\mat{Q}_{k+1}\mat{H}_{k+1,k}$, which reduces to the Lanczos decomposition \cite{chen2024lanczos} when $\mat{A} = \mat{A}^\tp$.
One then builds the Krylov subspace approximation~\cite{druskin1989two,saad1992analysis}
\begin{equation}    \label{eq:fAbkrylov}
f(\mat{A})\vec{b}\approx \widehat{\vec{f}} \coloneqq \norm{\vec{b}} \cdot \mat{Q}_kf(\mat{H}_k)\evec_1,
\end{equation}
where $\mat{H}_k = \mat{Q}_k^\tp \mat{A} \mat{Q}_k \in \R^{k\times k}$ consists of the first $k$ rows of $\mat{H}_{k+1,k}$. The following result establishes polynomial exactness \cite[Lemma 3.1]{saad1992analysis}. 

\begin{proposition}[Krylov approximation to $f(\mat{A})\vec{b}$]
    The approximation \cref{eq:fAbkrylov} is exact when $f$ is a polynomial of degree $\le k-1$.
\end{proposition}
In practice, we are rarely interested in using the approximation \cref{eq:fAbkrylov} for a polynomial function $f$.
Nonetheless, nice, well-behaved functions can be approximated to high accuracy by polynomials \cite{trefethen2019approximation}, so the approximation \cref{eq:fAbkrylov} is sensible for any ``nice'' function $f$; see, e.g.,~\cite{beckermann2009error} for concrete error bounds.

Randomization can be used to reduce the cost of the Arnoldi
approximation
\cite{guttel2023randomized,cortinovis2024speeding,
      palitta2025sketched,dedamas2026randomized}.
As before, there are two related strategies. One may build
an approximately orthogonal Krylov basis using RGS, or one may build a
cheaper non-orthogonal basis, using partial orthogonalization and
then apply sketching.
Both approaches result in an Arnoldi-like decomposition 
\begin{equation} \label{eq:arnoldi-like}
    \mat{A}\mat{V}_k=\mat{V}_{k+1}\widetilde{\mat{H}}_{k+1,k},
\end{equation}
from which the Krylov approximation \cref{eq:fAbkrylov} has the exact representation
\begin{equation*}
    \widehat{\vec{f}} = \mat{V}_k^{\vphantom{\pinv}}f(\mat{V}_k^\pinv \mat{A}\mat{V}_k^{\vphantom{\pinv}})\mat{V}_k^\pinv \vec{b} = \norm{\vec{b}} \cdot \mat{V}_k^{\vphantom{\pinv}}f(\mat{V}_k^\pinv \mat{A}\mat{V}_k^{\vphantom{\pinv}})\evec_1.
\end{equation*}
We have $\mat{V}_k^\pinv \vec{b}=\|\vec{b}\|\evec_1$ because the first vector of $\mat{V}_k$ is $\vec{b}/\|\vec{b}\|$.
To compute the core matrix $\mat{V}_k^\pinv \mat{A}\mat{V}_k^{\vphantom{\pinv}}$, observe that the Arnoldi-like relation \cref{eq:arnoldi-like} implies that
\begin{equation*}
    \mat{V}_k^\pinv \mat{A}\mat{V}_k^{\vphantom{\pinv}} = \mat{V}_k^\pinv \mat{V}_{k+1}^{\vphantom{\pinv}}\widetilde{\mat{H}}_{k+1,k} = \mat{V}_k^\pinv\onebytwo{\mat{V}_k}{{\vec{v}}_{k+1}}\twobyone{\widetilde{\mat{H}}_k}{\widetilde{\vec{h}}_{k+1}^\tp} = \widetilde{\mat{H}}_k^{\vphantom{\tp}} + (\mat{V}_k^\pinv \vec{v}_{k+1}^{\vphantom{\pinv}})\widetilde{\vec{h}}_{k+1}^\tp.
\end{equation*}
Thus, the core $\mat{V}_k^\pinv \mat{A}\mat{V}_k^{\vphantom{\pinv}}$ can be computed by solving a single least-squares problem $\mat{V}_k^\pinv \vec{v}_{k+1}^{\vphantom{\pinv}}$, which can be done rapidly using the sketch-and-solve method (or sketch-and-precondition---discussed below in \cref{sec:sketch-and-precondition}---for higher accuracy).
It would, of course, be cheaper to simply replace this matrix by the Hessenberg matrix
$\widetilde{\mat H}_k$ already available from the Arnoldi-like decomposition~\cref{eq:arnoldi-like}. While this often works out well in practice and polynomial exactness is preserved,
the numerical range of this matrix is difficult to 
control. This is a principal difficulty in obtaining convergence bounds for such sketched matrix-function
approximations; see \cite{guttel2023randomized,cortinovis2024speeding,
      palitta2025sketched}.

\section{Preconditioning for least squares and linear systems} 
\label{sec:precond}

Preconditioning is one of the most fundamental ideas in numerical computation.
Given a problem that is difficult for a given numerical method to solve, the goal of preconditioning is to transform the problem to an equivalent form that is easier to solve. We refer to~\cite{saad2003iterative,benzi2002preconditioning,wathen2015preconditioning} for classical treatments. 

\subsection{Sketch-and-precondition: Preconditioning tall-skinny least squares} \label{sec:sketch-and-precondition}

Perhaps the most effective instance of randomized preconditioning is the sketch-and-precondition paradigm for linear least-squares problems. 
We already saw one approach to using randomization for least squares earlier with the sketch-and-solve method (\cref{sec:sas}).
That approach gave us a cheap way to get a coarse approximation to the least-squares solution.
Using sketch-and-precondition, we expend some additional effort but, in return, we reliably obtain a high-accuracy solution, up to the level of accuracy permitted by machine precision.
While the basic approach is due to \citet{rokhlin2008fast}, \citet{avron2010blendenpik} refined the idea, provided careful numerical testing, and coined the name \emph{Blendenpik}. 

We begin with a large, highly overdetermined least-squares problem
\begin{equation} \label{eq:spre-least-squares}
    \vec{x}_\star\coloneqq \argmin_{\vec{x}} \norm{\mat{A}\vec{x} - \vec{b}} \quad \text{with } \mat{A} \in \R^{m\times n} \text{ and } \vec{b} \in \R^m,
\end{equation}
where $\mat A$ is assumed to have full column rank.
To solve such a large problem, we may employ a Krylov subspace iterative method, such as LSQR \cite{paige1982algorithm}, which only requires matrix--vector multiplications with $\mat{A}$ and its transpose.
But the convergence of LSQR and other Krylov subspace methods may be slow if the matrix $\mat{A}$ is poorly conditioned.

Fortunately, we already know a fast way to improve the conditioning of a tall-skinny matrix: whitening (\cref{sec:whitening}).
By now, the steps are familiar: we generate a random embedding $\mat{\Omega} \in \R^{m\times s}$, sketch $\Ahat \coloneqq \mat{\Omega}^\tp \mat{A}$, and factorize $\Ahat = \Qhat\Rhat$, obtaining a matrix $\Rhat$ for which $\mat{A}\Rhat^{-1}$ is well-conditioned.
Concretely, if $\mat{\Omega}$ is a subspace embedding for $\range(\mat{A})$ with distortion $\eta = 1$, then we have $\kappa(\mat{A}\Rhat^{-1}) \le 2$.

Instead of solving the sketched problem directly, we use the whitening transformation as a preconditioner:
Given the original least-squares problem \cref{eq:spre-least-squares}, we perform the change of variables $\vec{y} \coloneqq \Rhat \vec{x}$ and apply LSQR (or another iterative method) to the \emph{preconditioned} least-squares problem
\begin{equation} \label{eq:spre-transformed}
    \argmin_{\vec{y}} \norm{\mat{B}\vec{y} - \vec{b}} \quad \text{with } \mat{B} \coloneqq \mat{A}\Rhat^{-1}.
\end{equation}
After \cref{eq:spre-transformed} has been solved to sufficient accuracy, we return $\vec{x} \coloneqq \Rhat^{-1}\vec{y}$.
This is the sketch-and-precondition (Blendenpik) algorithm.
See \cref{alg:sketch-and-precondition}.

\begin{algorithm}[t]
\caption{Sketch-and-precondition for least-squares regression} \label{alg:sketch-and-precondition}
\begin{algorithmic}[1]
\Require Matrix $\mat{A} \in \R^{m\times n}$, right-hand side $\vec{b} \in \R^m$
\Ensure Least-squares solution $\vec{x} \in \R^n$
\State Generate a random embedding $\mat{\Omega} \in \R^{m\times s}$
\State Compute the sketches $\Ahat \gets \mat{\Omega}^\tp\mat{A}$ and $\bhat \gets \mat{\Omega}^\tp\vec{b}$
\State Factorize the sketch $(\mat{Q},\mat{R}) \gets \Call{QR}{\Ahat,\texttt{'econ'}}$
\State Form the initial guess $\vec{y} \gets \mat{Q}^\tp\bhat$ \Comment{Sketch-and-solve initialization} \label{line:sketchsolveinit}
\State Define LSQR primitives:
\begin{equation*}
    \texttt{Apply} : \vec{z} \mapsto \mat{A}(\mat{R}^{-1}\vec{z}) \quad \text{and} \quad \texttt{ApplyT} : \vec{z} \mapsto \mat{R}^{-\tp}(\mat{A}^\tp\vec{z}).
\end{equation*}
\State $\vec{y} \gets \Call{LSQR}{\texttt{Apply},\texttt{ApplyT},\texttt{rhs}=\vec{b},\texttt{init}=\vec{y}}$ \Comment{Preconditioned LSQR} 
\label{line:firstlsqr}
\State $\vec{x} \gets \mat{R}^{-1}\vec{y}$, $\vec{r}\gets\vec{b}-\mat{A}\vec{x}$ \Comment{Undo change of variables and compute residual}
\State $\vec{y} \gets \Call{LSQR}{\texttt{Apply},\texttt{ApplyT},\texttt{rhs}=\vec{r},\texttt{init}=\vec{0}}$ \Comment{Refine for stability}
\State Apply correction $\vec{x} \gets \vec{x} + \mat{R}^{-1}\vec{y}$
\end{algorithmic}
\end{algorithm}

We emphasize that, to achieve appealing computational complexity in this procedure, one must not form the matrix $\mat{B} = \mat{A}\Rhat^{-1}$ explicitly.
Rather, one computes matrix--vector multiplications $\vec{v} = \mat{B}\vec{z}$ by first performing the triangular solve $\vec{w} \gets \Rhat^{-1}\vec{z}$, then the matrix--vector product $\vec{v} \gets \mat{A}\vec{w}$.
Matrix--vector multiplications with $\mat{B}^\tp$ are performed similarly.

\paragraph{Analysis.}

Sketch-and-precondition satisfies the following bound.
\begin{theorem}[Sketch-and-precondition] \label{thm:analysissketchprecond}
    Let $\mat A \in \R^{m\times n}$ have full column rank, and let $\mat{\Omega} \in \R^{m\times s}$ be a subspace embedding with distortion $\eta = 1$ for $\range([\mat{A}\,,\vec b])$. Then the approximation $\vec{x}_k$ produced by sketch-and-precondition with $k$ LSQR steps (in \Cref{line:firstlsqr} of \cref{alg:sketch-and-precondition} ) 
    satisfies, in exact arithmetic, 
    \begin{equation*}
        \norm{\mat{A}\vec{x}_k - \vec{b}} \le (1 + 3^{2-k}) \norm{\mat{A}\vec{x}_\star - \vec{b}}.
    \end{equation*}
    In particular, $k = \order(\log(1/\varepsilon))$ iterations suffice to produce an approximate solution satisfying $\norm{\mat{A}\xhat - \vec{b}} \le (1+\varepsilon) \norm{\mat{A}\vec{x}_\star  - \vec{b}}$.
\end{theorem}
\begin{proof}
Using the embedding property,~\Cref{thm:sketch-and-solve-ose-osi}~\ref{item:sas-ose} gives
$\norm{\mat{A}\vec{x}_0-\vec{b}}
\le2\norm{\mat{A}\vec{x}_\star-\vec{b}}$ for the sketch-and-solve solution
$\vec{x}_0 = \mat{R}^{-1}\vec{y}$ from \Cref{line:sketchsolveinit}.
Combined with the Pythagorean decomposition~\cref{eq:sketch-solve-universal1}, this yields
$
  \norm{\mat{A}(\vec{x}_0-\vec{x}_\star)}^2
  \le3\norm{\mat{A}\vec{x}_\star-\vec{b}}^2.
$
Since LSQR is equivalent to CG applied to the normal equations,
the standard CG error bound~\cite[Thm.~38.5]{trefethen1997numerical}, combined with 
$\kappa(\mat{A}\Rhat^{-1})\le2$, implies
$
 \norm{\mat{A}(\vec{x}_k-\vec{x}_\star)}
 \le2\cdot3^{-k}\norm{\mat{A}(\vec{x}_0-\vec{x}_\star)} \le 2\cdot3^{1/2-k}\norm{\mat{A}\vec{x}_\star-\vec{b}}.
$
A second application of Pythagoras proves the result.
\end{proof}

The bound of~\Cref{thm:analysissketchprecond} is a streamlined version of the analysis of \citet{rokhlin2008fast}. It shows that the sketch-and-precondition method solves a least-squares problem to an accuracy that is exponentially small in the amount of work performed.
Compare this to the sketch-and-solve method (\cref{sec:sas}), which achieves a fixed level of accuracy for a given embedding dimension $s$ and converges only algebraically if $s$ is increased.
The proof above assumes $\mat{\Omega}$ is a subspace embedding for $\range([\mat{A}\,,\vec b])$, but for the preconditioner to be effective, $\mat{\Omega}$ need only be a subspace embedding for $\range(\mat{A})$.
In this case, using \cref{thm:sketch-and-solve-ose-osi}~\ref{item:sas-osi-uniform} in place of \cref{thm:sketch-and-solve-ose-osi}~\ref{item:sas-ose} yields $\norm{\mat{A}\vec{x}_k - \vec{b}} \le (1+ \norm{\mat{\Omega}}\cdot 3^{2-k})\norm{\mat{A}\vec{x}_\star - \vec{b}}$.
This result is important because it shows we can reuse the same preconditioner to solve the system with an arbitrary right-hand side $\vec{b}$. 

\paragraph{Runtime.}
With a fast embedding and an embedding dimension $s = \order(n)$, the total operation count of~\Cref{alg:sketch-and-precondition} is $\order(mn \log m + n^3)$ operations for dense $\mat A$ and a fixed accuracy level $\varepsilon$. This presents a fundamental asymptotic improvement over traditional least-squares solvers based on the normal equations or QR factorization, which require roughly $mn^2$ or $2mn^2$ operations, respectively.
However, meaningful speedups only manifest for large, highly overdetermined least-squares problems ($m\gg n\gg 1$).
Indeed, each step of sketch-and-precondition requires about $4mn$ floating point operations, and tens of steps may be necessary to achieve convergence to machine accuracy.
Moreover, each sketch-and-precondition iteration uses matrix--vector arithmetic, which is less hardware-efficient than the matrix--matrix arithmetic used in direct solvers.
As a rule of thumb, one should start to consider using sketch-and-precondition over a direct method when 
$m$ is in the thousands or larger, and $m\gg n$. 
A recent alternative, APLICUR \cite{huh2026adaptive}, uses a small initial sketch and adaptively updates a CUR-based preconditioner (\Cref{sec:cur}) during LSQR.

\paragraph{Initialization and numerical stability.}

Since sketch-and-precondition converges so rapidly, it can be used to produce solutions that are accurate up to numerical precision.
But care must be taken: certain versions of the sketch-and-precondition procedure are numerically unstable \cite{meier2024are} \cite[see also][]{epperly2024fast}.
To achieve accurate results, one can initialize the sketch-and-precondition iteration with $\vec{y}_0 \coloneqq \Qhat^\top (\mat{\Omega}^\top\vec{b})$, so that $\vec{x}_0 = \Rhat^{-1}\vec{y}_0 = (\mat{\Omega}^\tp\mat{A})^\pinv (\mat{\Omega}^\tp \vec{b})$ is the sketch-and-solve solution.
With this initialization, sketch-and-precondition converges to a \emph{forward stable} solution to the least-squares problem \cref{eq:spre-least-squares}. 
That is, the numerically computed solution $\xhat_1$ by sketch-and-precondition will achieve an error $\norm{\,\xhat_1 - \vec{x}_\star}$ comparable to the level of error achieved by a standard direct method based on QR factorization.
However, even with the improved initialization, sketch-and-precondition still need not produce a \emph{backward stable} solution.
To achieve this stronger stability property, we perform a single step of iterative refinement using sketch-and-precondition \emph{with zero initialization} as the solver.
These stability improvements (sketch-and-solve initialization, iterative refinement) are included in \cref{alg:sketch-and-precondition}.
These stabilizing fixes for sketch-and-precondition were identified by \citet{epperly2026fasta}, who prove rigorous forward and backward numerical stability guarantees for a slightly modified version of the algorithm using a different Krylov subspace method in place of LSQR.
See \cite{epperly2026fasta} also for additional implementation guidance, including discussion of error estimation and automatic stopping for sketch-and-precondition.

\begin{figure}
    \centering
    \includegraphics[width=0.99\linewidth]{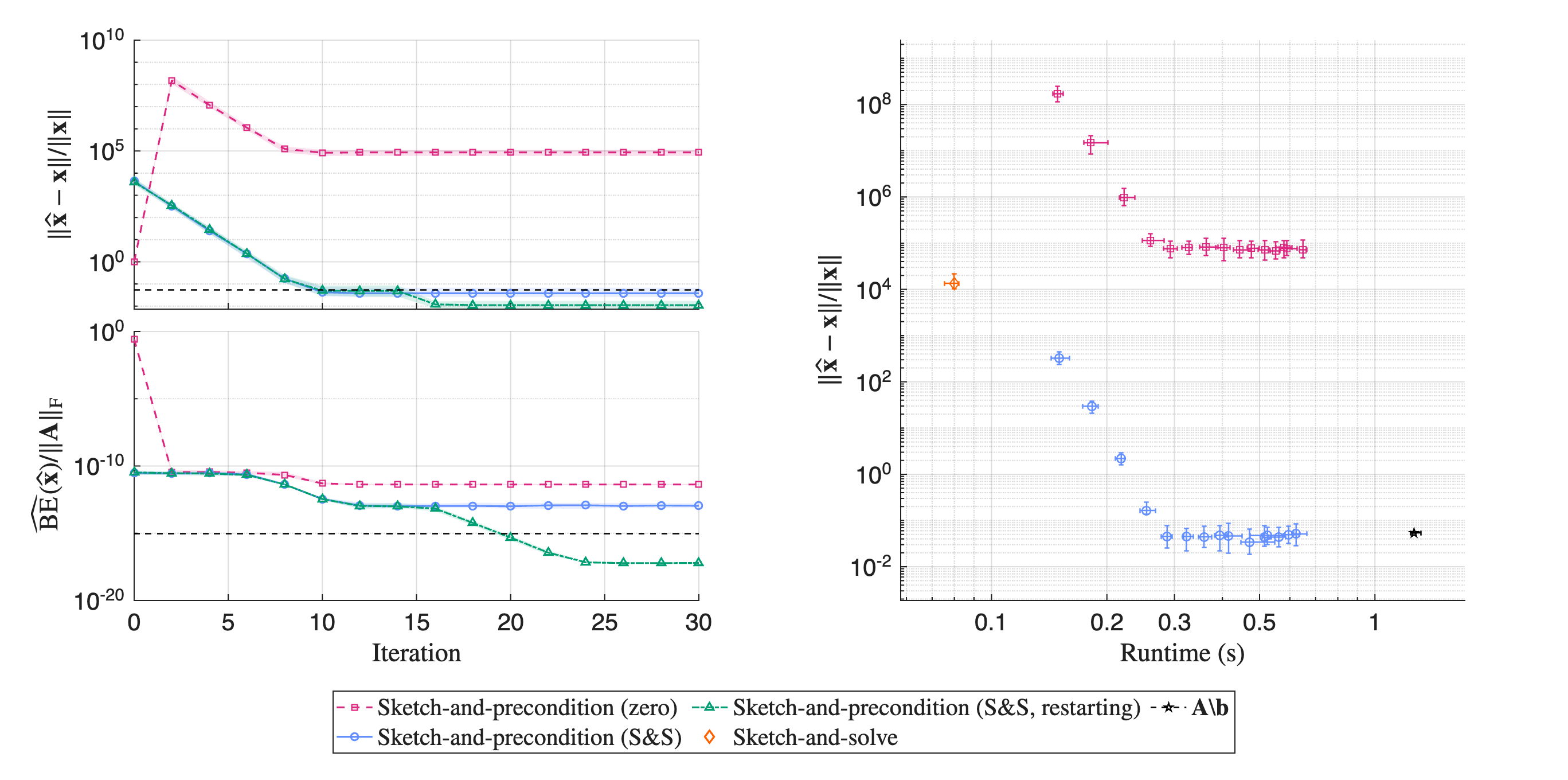}
    \caption{Sketch-and-precondition method with zero initialization (pink squares), sketch-and-solve initialization without restarting (blue circles), and sketch-and-solve initialization with restarting at iteration 15 (green triangles).
    Sketch-and-solve and MATLAB's \texttt{mldivide} are shown for reference.
    Left panels: relative forward error (top) and backward error estimate (\citet{karlson1997estimation}, bottom) versus iteration.
    Right panel: Relative forward error versus runtime.
    Data points for sketch-and-precondition are recorded every two iterations.
    Error bars and shaded regions show 90\% and 10\% quantiles computed over 30 trials.}
    \label{fig:sketch-and-precondition}
\end{figure}

\paragraph{Numerical evaluation.}
\Cref{fig:sketch-and-precondition} evaluates the sketch-and-precondition algorithm on an ill-conditioned, dense random least-squares problem with matrix $\mat{A} \in \R^{(5\times 10^5)\times 200}$.
We generate the problem as described in \cite[\S1.1.1]{epperly2024fast} and set condition number $\kappa(\mat{A}) \coloneqq 10^{11}$ and optimal residual $\norm{\mat{A}\vec{x}_\star - \vec{b}} \coloneqq 10^{-6} \norm{\mat{A}}\norm{\vec{x}_\star}$.
The left panel shows the relative forward error $\norm{\vec{x_\star} - \xhat}/\norm{\vec{x_\star}}$ and the Karlson--Wald{\'e}n estimate \citeyearpar{karlson1997estimation} of the backward error, and the right panel shows the relative forward error versus runtime.

As predicted by theory (\cref{thm:analysissketchprecond}), the sketch-and-precondition converges at an exponential rate until hitting an accuracy level governed by the numerical precision.
With the sketch-and-solve initialization, this maximum attainable accuracy matches the level of accuracy attained by MATLAB's \texttt{mldivide}, a QR factorization-based direct solver.
With restarting, the backward error converges to the level of the direct method as well.

On this problem, using sketch-and-precondition results in a healthy speedup.
On this example, achieving the maximum attainable forward error requires 10 steps, resulting in a $4\times$ speedup over MATLAB's \texttt{mldivide}.
If lower accuracy is required, the speedups can be even greater.
A baseline solution computed using the sketch-and-solve solution leads to an even faster $15\times$ speedup over the direct method, but it is five orders of magnitude less accurate.
Using sketch-and-precondition allows us to maintain some of the speed of sketch-and-solve while achieving the same numerical accuracy as the direct method.

\paragraph{Extension: Nullspace problems.}
Whitening can also be combined with preconditioned eigensolvers, such as LOBPCG~\cite{knyazev2001toward},   
for computing the smallest singular values and vectors of a tall-thin matrix.
See \cite{epperly2026fastb} for details.

\subsection{\Nystrom preconditioning: Preconditioning a \pd matrix} \label{sec:nystrom-precond}

Across different applications, we frequently encounter matrices that possess a ``signal-plus-noise'' or ``data-plus-regularization'' structure.
Such matrices have only a few large eigenvalues or singular values (the signal/data) with the remaining eigenvalues being more tightly clustered, often at a level set by observational noise or regularization.
To solve linear systems involving such a matrix, we can form a low-rank approximation that captures the dominant part of the matrix and use it to form a preconditioner.

This approach is the most developed in the setting of \pd matrices.
Consider a linear system of the form
\begin{equation} \label{eq:reg-lin-system}
    (\mat{A} + \mu \Id)\vec{x} = \vec{b} \quad \text{where } \mat{A} \text{ is \psd, $\mu\ge 0$, and $\mat{A} + \mu\Id$ is pd}.
\end{equation}
To construct a preconditioner, we first build a \Nystrom approximation $\Ahat$ to $\mat{A}$, and define the preconditioner
\begin{equation} \label{eq:nystrom-preconditioner}
    \mat{P} \coloneqq \Ahat + \hat{\mu} \, \Id \quad \text{for a shift parameter } \hat{\mu} > 0.
\end{equation}
The choice of approximation $\Ahat$ should reflect how $\mat{A}$ can be accessed.
For entry-wise access, we recommend the \RPCholesky method \cite{chen2025randomly} (\cref{sec:rpcholesky}) and, for matrix--vector access, a single-pass \Nystrom approximation using a fast random embedding to compute the approximation $\Ahat$.

To finish specifying the method, we must choose the regularization parameter $\hat{\mu}$ for the preconditioner.
\citet{diaz2026robust} consider kernel ridge regression problems with a positive regularization term $\mu > 0$ and advocate $\hat{\mu} \coloneqq \mu$. When
$\mat{A} - \Ahat$ is \psd, as is the case for the approximations recommended above, this choice ensures that
the condition number featuring prominently in the classical bound for the 
preconditioned conjugate gradient method~\cite[Thm.~38.5]{trefethen1997numerical} is controlled by
\[
 \kappa\!\big(
     \mat{P}^{-1/2}(\mat{A}+\mu\Id)\mat{P}^{-1/2}
 \big)
 \le
 1+\frac{\norm{\mat{A}-\Ahat}}{\mu}
 \le
 1+\frac{\tr(\mat{A}-\Ahat)}{\mu}.
\]
In particular, a low-rank approximation error on the level of $\mu$ ensures fast convergence.

When the problem is unregularized ($\mu = 0$) or when $\mu$ is much smaller than the smallest positive eigenvalue of $\Ahat$, the above choice can lead to a singular or nearly singular preconditioner $\mat{P}$.
Adapting \citet{frangella2023randomized}, we suggest the choice $\hat{\mu} \coloneqq \mu + \lambda_{\rm min}^+(\Ahat)$ as a safe way of setting the shift $\hat{\mu}$ in such cases.
Here, $\lambda_{\rm min}^+(\cdot)$ denotes the smallest eigenvalue that is strictly positive.
An alternative approach \cite{cutajar2016preconditioning,kaminetz2026everything} eschews the shift parameter entirely and simply restores the diagonal of the approximation
\begin{equation*}
    \mat{P} \coloneqq \Ahat + \mathcal{D}(\mat{A} + \mu \Id - \Ahat),
\end{equation*}
where $\mathcal{D}(\mat{B}) = \diag(\mat{B}(1,1),\ldots,\mat{B}(n,n))$ denotes restriction to the diagonal.

To compute the action of the inverse-preconditioner $\mat{P}^{-1}$, we represent the \Nystrom approximation as a compact eigendecomposition $\Ahat = \mat{U}\mat{\Lambda}\mat{U}^\tp$ and use the formula \cite[eq.~(5.3)]{frangella2023randomized}
\begin{equation*}
    \mat{P}^{-1}\vec{z} = \mat{U} [ (\mat{\Lambda} + \hat{\mu} \, \Id)^{-1} (\mat{U}^\tp\vec{z})] + \hat{\mu}^{\,-1}\vec{z} - \hat{\mu}^{\,-1} \mat{U}(\mat{U}^\tp \vec{z}).
\end{equation*}
As a faster---but less numerically stable---alternative, we can compute the \Nystrom approximation in outer product form $\Ahat = \mat{F}\mat{F}^\tp$ and use the Woodbury identity to derive the update (cf.\ \citet[Alg.~2]{guzmancordero2025improving}):
\begin{equation*}
    \mat{P}^{-1}\vec{z} = \hat{\mu}^{\,-1} \vec{z} - \hat{\mu}^{\,-1} \mat{F} [( \hat{\mu} \,\Id + \mat{F}^\tp \mat{F})^{-1} (\mat{F}^\tp \vec{z})].
\end{equation*}

\paragraph{Alternative: Implicit preconditioning.} \label{page:implicitprecond}
The goal of \Nystrom preconditioning is to rapidly solve a linear system that is (reasonably) well-conditioned except for a small number of outlying large eigenvalues.
A different family of methods seeks to attain this goal without the explicit use of a preconditioner.

Indeed, even the standard conjugate gradient method has such an implicit preconditioning effect, with its $k$th iterate $\vec{x}_k$ in exact arithmetic satisfying the following bound 
\cite[eqn.~(3.11)]{greenbaum1997iterative}:
\begin{equation} \label{eq:cg-improved-bound}
    \frac{ \norm{\vec{x}_k - \vec{x}}_{\mat{A} + \mu \Id}}{ \norm{\vec{x}_0 - \vec{x}}_{\mat{A} + \mu \Id} } \le 2
    \exp\!\bigg(- \frac{2(k-r)}{\sqrt{\kappa_{r+1}}} \bigg), \quad \text{where } \kappa_{r+1} \coloneqq \frac{\lambda_{r+1}(\mat{A}) + \mu}{\lambda_{\rm min}(\mat{A}) + \mu}.
\end{equation}
Here, $\vec{x}$ is the solution to \cref{eq:reg-lin-system} and $\norm{\vec{z}}_{\mat{M}} \coloneqq \sqrt{\vec{z}^\tp \mat{M}\vec{z}}$ denotes the norm associated with a \pd matrix $\mat{M}$.
We see that, after a burn-in period of $r$ steps, conjugate gradient converges at a rate insensitive to the top $r$ eigenvalues.
This desirable property can (and often does) fail catastrophically in finite-precision arithmetic, but it can be restored by using full reorthogonalization.
\citet{chen2026preconditioning} provide a detailed comparison of \Nystrom preconditioning to conjugate gradient and suggest alternative implicit preconditioning strategies based on conjugate gradient with partial reorthogonalization and block conjugate gradient \cite{oleary1980block} applied to a linear system augmented with additional random right-hand sides.

Another implicit preconditioning strategy is provided by the Kaczmarz\plusplus algorithm \cite{derezinski2025randomized}:
After a randomized preprocessing step, running a carefully tuned variant of the randomized block Kaczmarz or coordinate descent iteration \cite{elfving1980block,needell2014paved} with a block size of $k = \orderish(r)$ removes, after an initial
phase, the effect of the top $r$ eigenvalues on convergence, analogous to the improved conjugate gradient bound \cref{eq:cg-improved-bound}.

\paragraph{History.}
Versions of the \Nystrom preconditioning idea were first explored in the
context of Gaussian process regression and kernel ridge regression
\cite{cutajar2016preconditioning,gardner2018gpytorch}, and the general
idea was laid out in \cite[\S17]{martinsson2020randomized}.
The procedure was isolated and analyzed by
\citet{frangella2023randomized}, who used a slight variant of
\cref{eq:nystrom-preconditioner} in which the zero eigenvalues of
$\Ahat$ receive an additional positive shift.
\citet{diaz2026robust} propose the use of \RPCholesky (\cref{sec:rpcholesky}) to form the approximation and compare the convergence behavior of preconditioners based on different methods for constructing the \Nystrom approximation.
\citet{dereziski2025faster} propose and analyze an accelerated version
of \Nystrom preconditioning that uses multiple levels of sketching
and preconditioning to rapidly approximate the action of the
inverse-preconditioner.
For matrices with slowly decaying eigenvalues, \Nystrom preconditioning approaches struggle, and Vecchia approximation-based preconditioners can provide a powerful alternative; see \cite{kaminetz2026everything} and the references therein.

\subsection{Laplacian and diagonally dominant linear systems}

Another important instance of randomized preconditioning arises in the
solution of a consistent linear system
$\mat{A}\vec{x} = \vec{b}$, where $\mat{A}$ is the Laplacian matrix for a graph.
By an algebraic reduction \cite{gremban1996combinatorial}, the same techniques can be used to solve any symmetric diagonally dominant linear system.
To date, the empirically most effective \emph{randomized} preconditioners are based on sparsified versions of the Cholesky decomposition, pioneered by \citet{kyng2016approximate}; see \cite{livne2026lamg} and the references therein for a discussion of deterministic algebraic multigrid preconditioners.
We do not cover this material in detail, as it is already surveyed
carefully by \citet[\S18]{martinsson2020randomized}.
For subsequent numerical implementations and refinements, see
\cite{chen2021rchol,gao2026ac,baumann2026vac,chaitanya2026cast}.

\section{Trace estimation and beyond} \label{sec:trace}

In many applications, we encounter an implicit matrix $\mat{A} \in \R^{n\times n}$ for which accessing matrix entries $\mat{A}(i,j)$ is expensive but where we have (relatively) cheap access to a matrix--vector multiplication (matvec) $\vec{\omega} \mapsto \mat{A}\vec{\omega}$.
Examples include the following:
\begin{enumerate}
    \item \textbf{Products of matrices.} If $\mat{A} = \mat{A}_1 \cdots \mat{A}_t$, then we can compute $\mat{A}\vec{\omega} = \mat{A}_1(\mat{A}_2(\cdots (\mat{A}_t\vec{\omega})\cdots))$ by iterated matrix--vector multiplication.
    \item \textbf{Solutions to linear systems.}
    If $\mat{A} = \mat{M}^{-1}$, then we can compute $\mat{A}\vec{\omega} = \mat{M}^{-1}\vec{\omega}$ using a fast linear solver for $\mat{M}$ \cite{saad2003iterative,davis2006direct}.
    \item \textbf{Functions of matrices.} If $\mat{A} = f(\mat{M})$, then we can compute $\mat{A}\vec{\omega} = f(\mat{M})\vec{\omega}$ using a fast matrix-function-times-vector algorithm \cite[ch.~13, see also \cref{sec:f(A)b}]{higham2008functions}.
    \item \textbf{Automatic differentiation.}
    If $\mat{A}$ is the Jacobian or Hessian of a function, then we can compute $\mat{A}\vec{\omega}$ rapidly using automatic differentiation \cite{blondel2024elements}.
\end{enumerate}
See \cite{ubaru2018applications} for concrete examples of scientific applications.
This section will review algorithms for estimating the trace of an implicit matrix $\mat{A}$.
It will also discuss algorithms for estimating other quantities, including the Frobenius norm, diagonal, and row norms.

\subsection{The Girard--Hutchinson estimator} \label{sec:girardhutchinson}

The simplest approach to trace estimation employs a basic Monte Carlo averaging approach.
A random vector $\vec{\omega}$ is \emph{isotropic} if it satisfies the identity
\begin{equation} \label{eq:isotropic}
    \expect[\vec{\omega}\vec{\omega}^\top] = \Id;
\end{equation}
see also~\Cref{def:isotropic}.
Popular examples of isotropic vectors include a Rademacher vector $\vec{\omega} \sim \Unif \{ \pm 1 \}^n$ of random signs,  or a uniformly random vector from the sphere $\{ \vec{x} \in \R^n : \norm{\vec{x}}^2 = n \}$ of radius $\sqrt{n}$.
To estimate the trace, we compute the matrix--vector product $\mat{A}\vec{\omega}$ and form the quadratic form $\vec{\omega}^\tp (\mat{A}{\vec{\omega}})$.
This quantity is an unbiased estimate for the trace:
\begin{equation*}
    \expect [\vec{\omega}^\tp (\mat{A}\vec{\omega})] = \expect[\tr(\vec{\omega}^\tp\mat{A}\vec{\omega})] = \expect[\tr(\mat{A}\vec{\omega}\vec{\omega}^\tp)] = \tr(\mat{A}\expect[\vec{\omega}\vec{\omega}^\tp]) = \tr(\mat{A}\cdot\Id) = \tr(\mat{A}).
\end{equation*}
A single quadratic form often does not yield a sufficiently accurate estimate for the trace.
To improve it, we form independent copies of the estimator and average:
\begin{equation} \label{eq:gh}
    \hat{\tr}_{\rm GH} \coloneqq \frac{1}{s} \sum_{i=1}^s \vec{\omega}_i^\tp (\mat{A}\vec{\omega}_i) \quad \text{for } \vec{\omega}_1,\ldots,\vec{\omega}_s \simiid \vec{\omega}.
\end{equation}
The resulting estimator $\hat{\tr}_{\rm GH}$ is known as the Girard--Hutchinson trace estimator.
In its modern form, this idea was first proposed by \citet{girard1987un} and developed by \citet{girard1989fast,hutchinson1989stochastic}.
Versions of the stochastic trace estimation idea date back much earlier to work in quantum physics; see~\cite{alben1975exact} for an example.

\paragraph{Analysis.}

Without loss of generality, we may assume that $\mat{A}$ is symmetric. Indeed, for a general matrix $\mat{A}$, its symmetric part 
$\mat{A}_{\rm s}=(\mat{A}+\mat{A}^\tp)/2$ gives
$\tr(\mat{A})=\tr(\mat{A}_{\rm s})$ and 
$\vec{\omega}^\tp\mat{A}\vec{\omega}
=
\vec{\omega}^\tp\mat{A}_{\rm s}\vec{\omega}$.
Hence, many results in this section extend to nonsymmetric matrices by replacing $\mat{A}$ with
its symmetric part $\mat{A}_{\rm s}$.

To understand the typical size of the error in the Girard--Hutchinson estimator, we must compute its variance.
We have the following result.

\begin{theorem}[Girard--Hutchinson estimator] \label{thm:variancegh}
    Let $\mat{A} \in \R^{n\times n}$ be symmetric
    with eigenvalues $\lambda_1,\ldots,\lambda_n$.
    Consider the Girard--Hutchinson estimator with a base vector $\vec{\omega}$ drawn from the standard Gaussian ($\Normal(\vec{0},\Id)$), Rademacher ($\Unif \{ \pm 1\}^n$), or random sphere ($\Unif \{ \vec{x} \in \R^n : \norm{\vec{x}}^2 = n \}$) distributions.
    Then
    \begin{align*}
    \Var(\hat{\tr}_{\rm GH}) &= \frac{2}{s} \| \offdiag(\mat{A}) \|_{\rm F}^2;
         \tag*{(Rademacher)} \\
    \Var(\hat{\tr}_{\rm GH}) &= \frac{2n}{s(n+2)} \Big\| \mat{A} - \frac{\tr(\mat{A})}{n} \Id \Big\|_{\rm F}^2 
        = \frac{2n}{s(n+2)} \sum_{i=1}^n (\lambda_i - \overline{\lambda})^2;
        \tag*{(sphere)} \\
    \Var(\hat{\tr}_{\rm GH}) &= \frac{2}{s}\norm{\smash{{\mat{A}}}}_{\rm F}^2
        = \frac{2}{s}\sum_{i=1}^n \lambda_i^2.
        \tag*{(Gaussian)}
\end{align*}
    Here, $\overline{\lambda} \coloneqq n^{-1}\sum_{i=1}^n \lambda_i$ is the average of the eigenvalues
    and $\offdiag(\mat{A})$ is obtained from
    $\mat{A}$ by setting its diagonal entries to zero.
\end{theorem}
These formulas date back to \citet{girard1987un,girard1989fast} and \citet{hutchinson1989stochastic}, and all of them have simple proofs \cite{epperly2023stochastic}.
High-probability bounds on the error may be proven using Hanson--Wright-type inequalities \cite[Thm.~6.2.2]{vershynin2026high}. See~\citet{avron2011randomized,cortinovis2022randomized} for more analysis.

When $\mat{A}$ is additionally \psd, it is natural to consider the relative variance $\Var(\hat{\tr}_{\rm GH}) /\tr(\mat{A})^2$, which, according to \Cref{thm:variancegh}, admits the bound
\begin{equation}
    \frac{\Var(\hat{\tr}_{\rm GH})}{\tr(\mat{A})^2}
    \leq \frac{2}{s}\frac{\norm{\mat{A}}_{\rm F}^2}{\tr(\mat{A})^2}
    \label{eq:gh-relative-error}
\end{equation}
for each of the three probe distributions.
The ratio $\norm{\mat{A}}_{\rm F}^2/\tr(\mat{A})^2$ is small for slowly
decaying spectra, reaching its minimum value $1/n$ when all eigenvalues are
equal. Thus, Girard--Hutchinson is particularly effective in this regime.

The most important message of~\Cref{thm:variancegh} is that the variance decreases at an $\order(1/s)$ rate in all cases, so the root-mean-square error decreases like $\order(1/\sqrt{s})$.
This is the characteristic rate for a Monte Carlo method, and it suggests we must typically increase the computational work one hundred-fold to improve the accuracy of the estimator by just a single decimal digit.
The Girard--Hutchinson estimator is thus best thought of as a tool for getting a crude estimate of the trace that is accurate to a few digits.

While the $1/s$ variance scaling is the same for all three base distributions, these bounds demonstrate that the choice of probe vectors $\vec{\omega}_i$ matters.
With Rademacher vectors, the variance is independent of the diagonal entries, making these random  vectors an excellent choice for matrices $\mat{A}$ with a heavy diagonal.
This bound also reveals an important distinction between the Gaussian and sphere distributions: the variance for Gaussian probes depends on the magnitude of the eigenvalues, whereas the variance for the sphere distribution depends only on the amount of deviation of the eigenvalues from their mean.
For problems where the eigenvalues are clustered at a value far from zero, spherical probes thus result in a much lower variance than Gaussian probes.
Since the variance for Gaussian probes is larger than the variance for either Rademacher or spherical probes, Gaussian probes are not the preferred choice for standalone trace estimation.

\begin{remark}[Bad probe vectors]
    While every isotropic probe vector gives rise to an unbiased Girard--Hutchinson estimator, the variance for some choices can be unacceptably large.
    In particular, we advise against using random coordinate vectors $\vec{\omega} \sim \Unif \{ \sqrt{n} \evec_i \}$ or sparse random sign vectors unless this is indicated by the structure of $\mat{A}$; see~\cite{hallman2023monte} and~\cref{sec:stochastic-probing}.
\end{remark}

\subsection{Variance reduction by low-rank approximation: Hutch\plusplus} \label{sec:variance-reduction}

We can improve the Girard--Hutchinson estimator by employing a variance-reduction approach.
By~\Cref{thm:variancegh}, the variance of the estimator critically depends on the norm of the matrix $\mat{A}$. If we construct an approximation $\Ahat \approx \mat{A}$ with known trace, we can fashion a reduced-variance estimator for $\tr(\mat{A})$ by decomposing
\begin{equation*}
    \tr(\mat{A}) = \tr(\Ahat) + \tr(\mat{A} - \Ahat)
\end{equation*}
and applying Girard--Hutchinson only to estimate the second term.
The basic variance reduction idea dates back to the original work of \citet{girard1987un}.

A natural way of instantiating the variance reduction strategy is to use a randomized low-rank approximation algorithm to furnish the approximation $\Ahat$, an idea that was explored by \citet{saibaba2017randomized,lin2017randomized,gambhir2017deflation}.
This idea was developed most fully by~\citet{meyer2021hutch}, resulting in the Hutch\plusplus algorithm.
First, generate a test matrix $\mat{\Omega} \in \R^{n\times t}$ and apply the randomized rangefinder $\mat{Q} = \orth(\mat{A}\mat{\Omega})$, which---according to~\cref{sec:lowrank-projection}---defines an approximation $\mat{A} \approx \mat{Q}\mat{Q}^\tp \mat{A}$.
The trace then admits the decomposition
\begin{equation*}
    \tr(\mat{A}) = \tr(\mat{Q}\mat{Q}^\tp \mat{A}) + \tr((\Id - \mat{Q}\mat{Q}^\tp) \mat{A}) = \tr(\mat{Q}^\tp (\mat{A}\mat{Q})) + \tr((\Id - \mat{Q}\mat{Q}^\tp) \mat{A}(\Id - \mat{Q}\mat{Q}^\tp)).
\end{equation*}
Applying Girard--Hutchinson to the second term, we obtain the Hutch\plusplus estimator:
\begin{equation} \label{eq:hutchpp}
\hat{\tr}_{\rm H\plusplus} \coloneqq \tr(\mat{Q}^\tp(\mat{A}\mat{Q})) + \frac{1}{s} \sum_{i=1}^s \vec{\gamma}_i^\tp (\Id - \mat{Q}\mat{Q}^\tp) \mat{A}(\Id - \mat{Q}\mat{Q}^\tp) \vec{\gamma}_i^{\vphantom{\tp}}.
\end{equation}
Here, $\vec{\gamma}_i$ are iid isotropic probe vectors, generated independently of the test matrix $\mat{\Omega}$. Evaluating $\hat{\tr}_{\rm H\plusplus}$ requires a total of $2t+s$
matrix-vector products with $\mat{A}$: $t$ to construct $\mat{Q}$,
another $t$ to form $\mat{A}\mat{Q}$, and $s$ to compute
$\mat{A}\vec{\gamma}_i$, $i=1,\ldots,s$.
See \cref{alg:hutchpp} for optimized pseudocode.

\begin{algorithm}[t]
\caption{Hutch\plusplus for trace estimation} \label{alg:hutchpp}
\begin{algorithmic}[1]
\Require Symmetric matrix $\mat{A} \in \R^{n\times n}$, rank parameter $t$, number of probes $s$
\Ensure Estimate $\hat{\tr}_{\rm H\plusplus} \approx \tr(\mat{A})$
\State Draw independent test matrices $\mat{\Omega} \in \R^{n\times t}$ and $\mat{\Gamma} \in \R^{n\times s}$
with iid isotropic columns
\State Compute an orthonormal basis $\mat{Q} \gets \orth(\mat{A}\mat{\Omega})$ \Comment{$t$ matvecs}
\State Form the sketch $\mat{Y} \gets \mat{A}\mat{Q}$ \Comment{$t$ matvecs}
\State Deflate the probes $\mat{G} \gets \mat{\Gamma} - \mat{Q}(\mat{Q}^\tp\mat{\Gamma})$ \Comment{$\mat{G} = (\Id - \mat{Q}\mat{Q}^\tp)\mat{\Gamma}$}
\Comment{$s$ matvecs}
\State \Return $\widehat{\mathrm{tr}}_{\mathrm{H++}}
\gets \langle\mat{Q},\mat{Y}\rangle_{\mathrm F}
+ s^{-1}\langle\mat{G},\mat{A}\mat{G}\rangle_{\mathrm F}$
\Comment{$s$ matvecs}
\end{algorithmic}
\end{algorithm}

\paragraph{Analysis.}
Low-rank approximation is most effective for matrices with rapidly decaying spectra, whereas \cref{eq:gh-relative-error} shows that Girard--Hutchinson benefits from a slowly decaying spectrum. Hutch\plusplus combines these complementary strengths. 

\begin{theorem}[Hutch\plusplus] \label{thm:hutchpp}
    Let $\mat{A} \in \R^{n\times n}$ be a square matrix, and let $\hat{\tr}_{\rm H\plusplus}$ denote the Hutch\plusplus estimator with parameters $t = 2r+1$ and $s$, implemented with Gaussian or spherical probes.
    Then
    \begin{equation} \label{eq:hutchpp-spectral}
        \expect\bigl[ |\hat{\tr}_{\rm H\plusplus} - \tr(\mat{A})|^2 \bigr] \le \frac{4}{s} \norm{\mat{A} - \lowrank{\mat{A}}_r}_{\rm F}^2.
    \end{equation}
    Moreover, if $\mat{A}$ is \psd, then Hutch\plusplus satisfies the relative-error guarantee
    \begin{equation}
        \frac{
            \expect\bigl[
                |\hat{\tr}_{\rm H\plusplus}-\tr(\mat{A})|^2
            \bigr]
        }{
            \tr(\mat{A})^2
        }
        \leq
        \frac{1}{sr}.
        \label{eq:hutchpp-rmse}
    \end{equation}
    In particular, by setting $s = r = \lceil 1 /\varepsilon \rceil$, Hutch\plusplus achieves a root-mean-square relative error $\varepsilon$ in at most $5/\varepsilon + 7$ matvecs.
\end{theorem}

This bound has two parts.
The first guarantee \cref{eq:hutchpp-spectral} shows that the Hutch\plusplus algorithm is spectrally accurate: If the singular values decay exponentially, then the Hutch\plusplus error converges exponentially as well.
The second guarantee \cref{eq:hutchpp-rmse} shows that the relative error of Hutch\plusplus is always at most $1/\sqrt{sr}$.
By setting $s = r = \order(1/\varepsilon)$, we see that relative root-mean-square accuracy $\varepsilon$ is achieved
using only
$
    2t+s = \order(\varepsilon^{-1})
$
matrix--vector products with $\mat{A}$.
This compares favorably with the $\order(\varepsilon^{-2})$ dependence of the standalone
Girard--Hutchinson estimator.
For a worst-case psd matrix, the $\order(\varepsilon^{-1})$ complexity of Hutch\plusplus is optimal for any implicit trace estimation algorithm \cite{meyer2021hutch,meyer2024towards}.
The same $\order(\varepsilon^{-1})$ complexity follows whenever
the rangefinder achieves a constant-factor Frobenius approximation with
$t=\order(r)$. In particular, this holds whp for the
Rademacher test matrices used in the original Hutch\plusplus algorithm~\cite{meyer2021hutch} and when the Hutch\plusplus algorithm is implemented for large tensor problems with tensor-train random vectors~\cite{camano2026linear}; see~\Cref{sec:random-tt}.

\begin{proof}[Proof of \cref{thm:hutchpp}]
    Conditioned on $\mat{Q}$, the Hutch\plusplus estimator is unbiased, and
    \Cref{thm:variancegh} gives
    \begin{equation}
        \Var\bigl(\hat{\tr}_{\rm H\plusplus}\mid \mat{Q}\bigr)
        \leq
        \frac{2}{s}
        \norm{(\Id-\mat{Q}\mat{Q}^\tp)\mat{A}(\Id-\mat{Q}\mat{Q}^\tp)}_{\rm F}^2
        \leq
        \frac{2}{s}
        \norm{(\Id-\mat{Q}\mat{Q}^\tp)\mat{A}}_{\rm F}^2 .
        \label{eq:hutchpp-condvar}
    \end{equation}
    This bound is controlled by the randomized rangefinder. For a matrix $\mat{\Omega}$ with Gaussian or spherical columns, the rangefinder result (\cref{thm:rangefinder}) gives
    \begin{equation}
            \expect_{\mat{\Omega}}
     \norm{(\Id-\mat{Q}\mat{Q}^\tp)\mat{A}}_{\rm F}^2
        \leq
        2
        \norm{\mat{A}-\lowrank{\mat{A}}_r}_{\rm F}^2.
        \label{eq:hutchpp-rangefinder}
    \end{equation}
    Combining \cref{eq:hutchpp-condvar,eq:hutchpp-rangefinder} yields the first guarantee \cref{eq:hutchpp-spectral}. 

    Now assume $\mat{A}$ is \psd.
    Because the eigenvalues $\lambda_j$ of $\mat{A}$ are nonnegative, it follows that
    \begin{equation}
        \norm{\mat{A}-\lowrank{\mat{A}}_r}_{\rm F}^2
        =
        \sum_{j>r}\lambda_j^2
        \leq
        \lambda_{r+1}\sum_{j>r}\lambda_j
        \le 
        \frac{1}{r} \Big(\sum_{j \le r} \lambda_j\Big)\Big(\sum_{j > r} \lambda_j\Big)
        \leq
        \frac{\tr(\mat{A})^2}{4r}.
        \label{eq:psd-tail-trace}
    \end{equation}
    The last inequality is $\alpha\beta \le (\tfrac{\alpha+\beta}{2})^2$.
    Combining \cref{eq:psd-tail-trace} with \cref{eq:hutchpp-spectral} yields the second conclusion \cref{eq:hutchpp-rmse}.
\end{proof}

\paragraph{Resphering.}
Many matrices in applications have a ``signal plus noise'' structure, where there are a few relatively large eigenvalues and the remaining eigenvalues are relatively clustered.
In such cases, the variance of Hutch\plusplus and similar algorithms can be improved by a technique known as \emph{resphering} \cite[\S2.3]{epperly2024xtrace}; see also \cite[\S\S13.1, 13.3, \& 14.8]{epperly2025make}.
The idea is to choose the random test vectors $\vec{\gamma}_i$ not as isotropic vectors on all of $\R^n$ but as uniform vectors on the sphere of radius $\sqrt{n - t}$ in $\range(\mat{Q})^\perp$.
Such vectors are isotropic \emph{on the subspace} $\range(\mat{Q})^\perp$ in the sense that
\begin{equation*}
    \expect[\vec{\gamma}\vec{\gamma}^\top] = \Id - \mat{Q}\mat{Q}^\top.
\end{equation*}
It is easy to see that using such vectors in the Hutch\plusplus algorithm still results in an unbiased trace estimate, and \cite[Fig.~13.1]{epperly2025make} presents an example where resphering improves the accuracy of Hutch\plusplus by an order of magnitude.
To generate the vectors $\vec{\gamma}$, simply draw a Gaussian vector, orthogonalize, and rescale:
\begin{equation} \label{eq:resphering}
    \vec{\gamma} \coloneqq \sqrt{n - t}\cdot \frac{(\Id - \mat{Q}\mat{Q}^\top)\vec{\nu}}{\norm{\smash{(\Id - \mat{Q}\mat{Q}^\top)\vec{\nu}}}} \quad \text{for }\vec{\nu} \sim \Normal(\vec{0},\Id).
\end{equation}

\paragraph{Variants.}
Using other types of low-rank approximation gives rise to useful variants of
Hutch\plusplus.
\emph{NA-Hutch\plusplus}~\cite{meyer2021hutch} uses the generalized
\Nystrom approximation from~\Cref{sec:lowrank-generalized-nystrom}.
For \psd matrices, the use of a \Nystrom approximation for variance reduction
was already proposed by \citet{lin2017randomized}; this idea was developed and
analyzed as \emph{Nystr\"om\plusplus} by
\citet{persson2022improved}, using the \Nystrom approximation
from~\Cref{sec:lowrank-nystrom}.
For fixed failure probability, both methods retain the
$\order(\varepsilon^{-1})$ matrix--vector complexity of Hutch\plusplus.
Moreover, they are nonadaptive: all test vectors can be chosen in advance, so
the required matrix--vector products can be evaluated in a single block
matrix--vector product, requiring only one pass over $\mat{A}$.
\citet{jiang2021optimal} subsequently sharpened the high-probability analysis
of NA-Hutch\plusplus{}. 

\subsection{Further variance reduction by exchangeability: \XTrace} \label{sec:xtrace}

Given a fixed budget of matrix--vector products, Hutch\plusplus and its variants require a choice of 
how to divide computational resources between constructing the low-rank approximation (controlled by the rank $t$) and estimating the trace of the residual (controlled by $s$). 
The adaptive Hutch\plusplus algorithm~\cite{persson2022improved} addresses this issue by choosing the split adaptively when a target accuracy is prescribed.
A different approach, proposed by~\citet{epperly2024xtrace}, is to eliminate the split altogether by enforcing \emph{exchangeability}~\cite[Chapter 12]{vandervaart1998asymptotic}.

Hutch\plusplus draws a pool of iid random vectors but assigns different algorithmic roles to different subsets. \XTrace exploits the exchangeability of these vectors by averaging over the possible role assignments. Specifically, generate a test matrix $\mat{\Omega} = \flatmatrix{\vec{\omega}_1,\ldots,\vec{\omega}_s} \in \R^{n\times s}$ with iid isotropic columns $\vec{\omega}_i$.
By leaving out any column $\vec{\omega}_i$, we obtain a random matrix \begin{equation} \label{eq:column-deletion} \mat{\Omega}_{-i} \coloneqq \flatmatrix{\vec{\omega}_1,\ldots,\vec{\omega}_{i-1},\vec{\omega}_{i+1},\ldots,\vec{\omega}_s}\end{equation} that is independent of the left-out vector $\vec{\omega}_i$.
Using $\mat{\Omega}_{-i}$ as the test matrix and $\vec{\omega}_i$, we may fashion a lopsided Hutch\plusplus estimator
\begin{subequations} \label{eq:xtrace}
\begin{equation} \label{eq:xtrace-replicates}
    \hat{\tr}_i \coloneqq \tr(\mat{Q}_{(i)}^\tp(\mat{A}\mat{Q}_{(i)}^{\vphantom{\tp}})) + \vec{\omega}_i^{\tp}(\Id - \mat{Q}_{(i)}^{\vphantom{\tp}}\mat{Q}_{(i)}^\tp)\mat{A}(\Id - \mat{Q}_{(i)}^{\vphantom{\tp}}\mat{Q}_{(i)}^\tp)\vec{\omega}_i
\end{equation}
with $\mat{Q}_{(i)} \coloneqq \orth(\mat{A}\mat{\Omega}_{-i})$.
As we only use a single probe vector to estimate the residual trace, each individual estimator $\hat{\tr}_i$ can have high variance  and provide a poor estimate of the trace.
The \XTrace estimator averages over all possible choices of the left-out vector:
\begin{equation}
    \hat{\tr}_{\rm X} \coloneqq \frac{1}{s} \big( \hat{\tr}_1 + \cdots + \hat{\tr}_s \big).
\end{equation}
\end{subequations}
This estimator is unbiased and a symmetric function of $\vec{\omega}_1, \ldots, \vec{\omega}_s$.

At first sight, \cref{eq:xtrace} appears to require $s$ separate randomized rangefinder computations.
\XTrace avoids this overhead. It uses exactly $2s$ matrix--vector products with $\mat{A}$, organized in two sequential batches, to compute $\mat{A}\mat{\Omega}$ and $\mat{A}\mat{Q}$, where 
$\mat{Q}$ is obtained from a thin QR factorization of 
$\mat{A}\mat{\Omega}$.
The leave-one-out projectors
$\mat{Q}_{(i)}\mat{Q}_{(i)}^\tp$ can then be obtained implicitly from
$\mat{Q}\mat{Q}^\tp$ by rank-one downdates derived from this QR factorization,
without any additional matrix--vector products with $\mat{A}$. 
The entire post-processing step requires only 
$\order(ns^2)$ operations; see \cite[Alg.~SM3.1]{epperly2024xtrace} for a stable, optimized implementation.
The resphering technique \cref{eq:resphering} can be incorporated into XTrace; see \cite{epperly2024xtrace,epperly2025make}.
The paper \cite{hallman2025two} develops several useful variants of the \XTrace algorithm.

The analysis of \citet[Thm.~1.1]{epperly2024xtrace} shows that, for \psd matrices, \XTrace retains the favorable $\order(\varepsilon^{-1})$ matrix--vector complexity of Hutch\plusplus, while offering substantially improved convergence for rapidly decaying spectra.

For a \psd matrix, the exchangeability principle can be combined with 
\Nystrom approximation instead of randomized rangefinding. Using the notation from \Cref{sec:lowrank-nystrom},
define the leave-one-out \Nystrom approximation
$
    \mat{A}_{(i)}^{\rm N}
    \coloneqq
    \mat{A}\langle\mat{\Omega}_{-i}\rangle
$
and the corresponding trace estimate
\begin{equation}
    \smash{\hat{\tr}}^{\rm N}_i
    \coloneqq
    \tr\!\big(\mat{A}_{(i)}^{\rm N}\big)
    +
    \vec{\omega}_i^\tp
    \big(\mat{A}-\mat{A}_{(i)}^{\rm N}\big)
    \vec{\omega}_i.
    \label{eq:xnystrace-basic}
\end{equation}
\XNysTrace averages these estimates as in~\cref{eq:xtrace}.
Unlike \XTrace, all required quantities can be recovered from the single block
product $\mat{A}\mat{\Omega}$, so \XNysTrace requires only one pass over
$\mat{A}$; see~\cite{epperly2024xtrace} for details.
\citet{epperly2025make} develops an improved implementation.

We will meet the leave-one-out principle behind \XTrace again in~\Cref{sec:leaveoneout}, in the context of a posteriori error estimation.

\subsubsection{Numerical comparison}

\begin{figure}[t]
    \centering
    \includegraphics[width=0.99\linewidth]{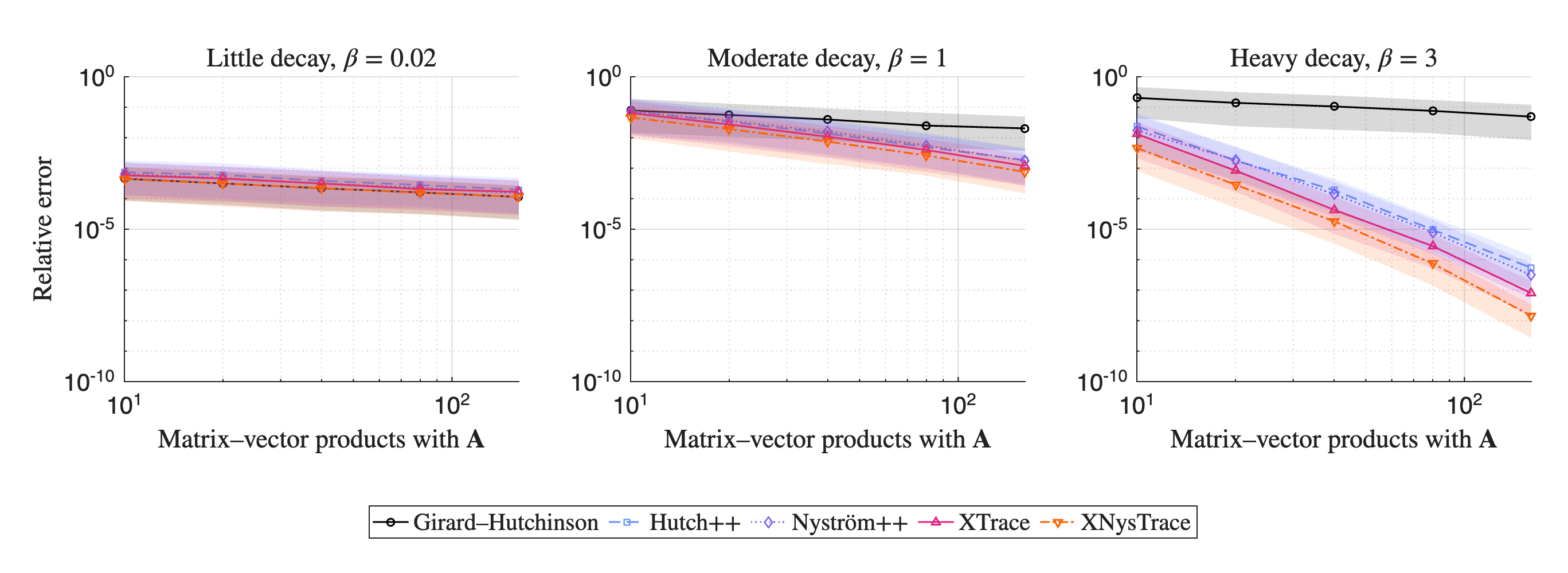}
    \caption{Relative error for Girard--Hutchinson, Hutch\plusplus, Nystr\"om\plusplus, \XTrace, and \XNysTrace estimators applied to estimate the trace of the matrix $\mat{A} = \e^{-\beta \mat{H}}$ where $\mat{H}$ is the transverse-field Ising model on 12 sites and the inverse-temperature is $\beta \in \{0.02,1,3\}$.
    Lines show the median of 1000 trials, and shaded regions show 10\% and 90\% quantiles.}
    \label{fig:trace}
\end{figure}

\Cref{fig:trace} presents a numerical comparison of the trace estimators discussed thusfar: the Girard--Hutchinson estimator, Hutch\plusplus, \XTrace, and \XNysTrace and Nystr\"om\plusplus (Nystr\"om variants of Hutch\plusplus and \XTrace).
The Girard--Hutchinson estimator is implemented with uniform vectors from the sphere, and the variance-reduced estimators use the resphering technique.
We apply these estimators to an example application from quantum statistical physics, where we seek to estimate the partition function $\tr \e^{-\beta \mat{H}}$ associated with the transverse field Ising model Hamiltonian matrix $\mat{H} \in \R^{2^{12}\times 2^{12}}$ \cite{pfeuty1970one} for a 1D chain of 12 quantum spins.\footnote{Concretely, this matrix is $\mat{H} = -\sum_{i=1}^{11} \mat{X}_i\mat{X}_{i+1} - h\sum_{i=1}^{12} \mat{Z}_i$, where $\mat{X}_i$ and $\mat{Z}_i$ denote the Pauli X and Z operators for the $i$th site.
We set $h \coloneqq 1$, which places the model at its critical point.}
To compute matrix--vector products with $\mat{A} = \e^{-\beta \mat{H}}$, we employ 60 steps of the Lanczos method \cite[\S6.1]{chen2024lanczos}, and we test three values of the inverse-temperature $\beta \in \{0.02,1,3\}$, resulting in matrices with slow, medium, and fast rates of spectral decay.

The pattern of results is instructive.
For the problem with slow decay, the relative error is relatively insensitive to the algorithm, and the algorithm achieves three digits of accuracy even with no variance reduction.
This is consistent with what we would expect from the bound \cref{eq:gh-relative-error}, which predicts the root-mean-square relative error for the Girard--Hutchinson estimator applied to a \psd matrix with a flat spectrum should have size at most $\approx s^{-1/2} \norm{\mat{A}}_{\rm F} / \tr(\mat{A}) \approx (sn)^{-1/2}\approx 0.005$ when $n = 2^{12}$ and $s = 10$.
As the rate of spectral decay increases, the accuracy of the Girard--Hutchinson estimator deteriorates, and the benefits of variance reduction become apparent.
At $\beta = 1$, all of the variance-reduced estimators behave similarly, but meaningful differences emerge at $\beta =3$.
The order of methods, from best to worst, is \XNysTrace, \XTrace, Nystr\"om\plusplus, and Hutch\plusplus.
The conclusion is that the exchangeable estimators \XNysTrace and \XTrace are preferable over Hutch\plusplus and Nystr\"om\plusplus.
For \psd matrices, the \Nystrom-based estimators (\XNysTrace, \Nystrom\plusplus) show meaningful benefits over the randomized SVD-based estimators (Hutch\plusplus, \XTrace).
See \cite{epperly2024xtrace} for more experiments, including some exceptional examples where \XTrace and Hutch\plusplus outperform \XNysTrace and Nystr\"om\plusplus, respectively.
For estimating the trace of a function of a matrix, black-box stochastic trace estimators like Hutch\plusplus and \XTrace are often not the best choice; see \cref{sec:sde} for many alternatives.

\subsection{Variance reduction by stochastic probing}
\label{sec:stochastic-probing}

Recall from~\Cref{thm:variancegh} that the variance of Girard--Hutchinson with Rademacher random vectors 
is blind to the diagonal of $\mat{A}$. \emph{Stochastic probing} extends this principle  to patterns beyond the diagonal when, for example, it is known that $\mat{A}$ is approximately banded.

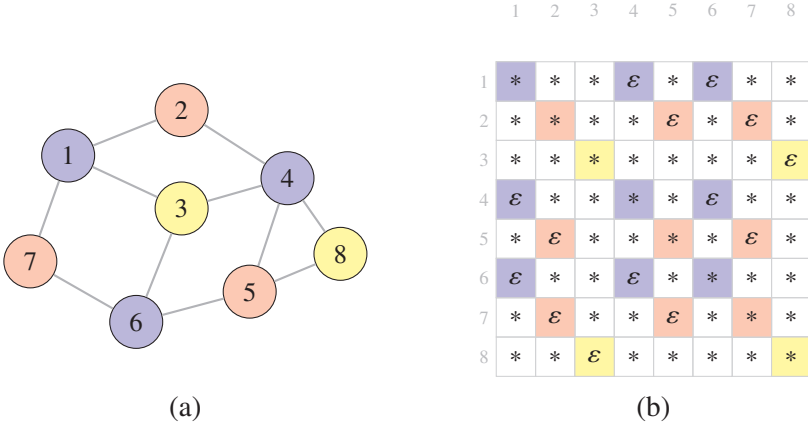
\begin{figure}
    \centering
    \begin{tikzpicture}[
        vertex/.style = {draw, circle, minimum size = 7mm, inner sep = 0pt, font = \small},
        classone/.style   = {vertex, fill = blue!25},
        classtwo/.style   = {vertex, fill = red!25},
        classthree/.style = {vertex, fill = yellow!45},
        link/.style = {thick, gray!70},
        cell/.style = {draw = gray!40, minimum size = 5.2mm, inner sep = 0pt, font = \small},
        idx/.style  = {font = \scriptsize, text = gray!60, inner sep = 1pt},
    ]
        \begin{scope}[yshift = -1.6cm]
            \foreach \i/\c/\x/\y in {1/classone/0/1.6, 2/classtwo/1.5/2.2,
                                     3/classthree/1.5/0.9, 4/classone/2.9/1.3,
                                     5/classtwo/2.4/-0.2, 6/classone/0.9/-0.6,
                                     7/classtwo/-0.5/0.2, 8/classthree/3.6/0.3}
                {\node[\c] (v\i) at (\x,\y) {$\i$};}
            \foreach \i/\j in {1/2, 2/4, 4/5, 5/6, 6/7, 1/7, 1/3, 3/4, 3/6, 4/8, 5/8}
                {\draw[link] (v\i) -- (v\j);}
        \end{scope}
        \begin{scope}[xshift = 5.4cm, yshift = 1.5cm]
            \foreach \j in {1,...,8} {\node[idx] at ({0.52*\j},0.42) {\j};}
            \foreach \i in {1,...,8} {\node[idx] at (0.1,{-0.52*\i}) {\i};}
            \foreach \i in {1,...,8}
                {\foreach \j in {1,...,8}
                    {\node[cell] at ({0.52*\j},{-0.52*\i}) {$\ast$};}}
            \foreach \i in {1,4,6}
                {\foreach \j in {1,4,6}
                    {\node[cell, fill = blue!25] at ({0.52*\j},{-0.52*\i}) {$\varepsilon$};}}
            \foreach \i in {2,5,7}
                {\foreach \j in {2,5,7}
                    {\node[cell, fill = red!25] at ({0.52*\j},{-0.52*\i}) {$\varepsilon$};}}
            \foreach \i in {3,8}
                {\foreach \j in {3,8}
                    {\node[cell, fill = yellow!45] at ({0.52*\j},{-0.52*\i}) {$\varepsilon$};}}
            \foreach \i in {1,4,6}
                {\node[cell, fill = blue!25] at ({0.52*\i},{-0.52*\i}) {$\ast$};}
            \foreach \i in {2,5,7}
                {\node[cell, fill = red!25] at ({0.52*\i},{-0.52*\i}) {$\ast$};}
            \foreach \i in {3,8}
                {\node[cell, fill = yellow!45] at ({0.52*\i},{-0.52*\i}) {$\ast$};}
        \end{scope}
        \node at (1.55,-3.35) {(a)};
        \node at (7.74,-3.35) {(b)};
    \end{tikzpicture}
    \caption{An illustration of matrix probing for a matrix $\mat{A}$ on $n=8$ indices.
    (a) The graph $(\set{V},\set{E})$ and a coloring with $q=3$ colors.
    (b) The matrix $\mat{A}$, with each diagonal block $\mat{A}(\set{C}_\ell,\set{C}_\ell)$ colored.
    The probing estimator will be accurate when all off-diagonal entries \emph{within each block} are small.}
    \label{fig:probing-coloring}
\end{figure}

\subsubsection{Deterministic probing}

Classical probing constructs deterministic vectors whose supports reflect a known pattern; see~\cite{frommer2021analysis} and the references
therein. The idea of probing is best described in terms of graph coloring~\cite{tang2012probing}.
Let $\mat{A}$ be symmetric and let $(\set{V},\set{E})$ be the undirected, unweighted graph on $\set{V}=\{1,\ldots,n\}$ such that an edge $\{i,j\}\in \set{E}$ indicates that $|\mat{A}(i,j)|$ is potentially large.
A coloring partitions the vertices into color classes
\[
    \set{V}=\set{C}_1\cup\cdots\cup\set{C}_q,
    \qquad
    \set{C}_\ell\cap\set{C}_j=\varnothing \quad (\ell\neq j),
\]
such that no two adjacent vertices have the same color; see \cref{fig:probing-coloring}.
Given color classes $\{\set{C}_\ell\}_{\ell=1}^q$, the \emph{deterministic probing estimator} is 
\begin{equation} \label{eq:vell}
     \hat{\tr}_{\rm DP} \coloneqq \sum_{\ell=1}^q
    \vec{v}_\ell^\tp\mat{A}\vec{v}_\ell^{\vphantom{\top}} \quad \text{where } \vec{v}_\ell \coloneqq \sum_{i\in \set{C}_\ell}\evec_i.
\end{equation}
The error of the deterministic probing estimator is given by the sum of the entries that are off-diagonal with respect to the block partitioning $\{ \set{C}_\ell \}$.
In particular, we have the worst-case bound
\[
    |\hat{\tr}_{\rm DP} - \tr(\mat{A})|
    \le
    \sum_{\ell=1}^q
    \norm{\offdiag(\mat{A}(\set{C}_\ell,\set{C}_\ell))}_{\ell_1},
\]
with equality obtained when all of the off-diagonal entries of the blocks $\mat{A}(\set{C}_\ell,\set{C}_\ell)$ have the same sign.
This error bound demonstrates two basic weaknesses of deterministic probing.
First, we pay a steep price for any off-diagonal entries within each block $\mat{A}(\set{C}_\ell,\set{C}_\ell)$.
In the worst case, the errors will accumulate, rather than canceling.
Second, once one fixes a graph coloring, the estimator is not systematically improvable.
The accuracy is what it is, and one needs to devise a new graph coloring to improve the error.

Notwithstanding this limitation, the deterministic probing approach can be quite effective.
To see this, consider a $(q-1)$-banded matrix $\mat{A}$:
\begin{equation*}
    \mat{A}(i,j) = 0 \quad \text{when } |i - j| > q-1.
\end{equation*}
The periodic coloring
\begin{equation}\label{eq:periodic}
     \set{C}_\ell
    =
    \{\,i: i\equiv \ell \pmod{q}\,\},
    \qquad
    \ell=1,\ldots,q,
\end{equation}
places any two indices belonging to the same color at least $q$ positions apart, making each block $\mat{A}(\set{C}_\ell,\set{C}_\ell)$ a diagonal matrix.
Consequently, deterministic probing recovers $\tr(\mat{A})$ exactly using only $q$ matrix--vector products.  More importantly, the same
construction gives a small error when $\mat{A}$ is approximately $(q-1)$-banded.

An important source of such approximate structure arises for
\emph{matrix functions}~\cite{higham2008functions}. 
Suppose that $\mat{A}=f(\mat{B})$ for symmetric sparse $\mat{B}$, where 
$f$ is sufficiently regular in a neighborhood of the spectrum of $\mat{B}$.
The inverse $\mat{B}^{-1}$ is an important
special case.  Although $f(\mat{B})$ is typically dense, its entries often
decay rapidly with graph distance in the sparsity graph of
$\mat{B}$~\cite{benzi2007decay}.  This a priori information can be incorporated
into $(\set{V},\set{E})$ by connecting every pair of vertices whose distance in the
sparsity graph of $\mat{B}$ is small.

\subsubsection{Stochastic probing}

We can improve probing by combining it with randomization.
Stochastic probing reduces the error caused by entries within the
color classes.
Unlike the deterministic query vectors used in the classical probing estimator~\cref{eq:vell}, we define for each color $\set{C}_\ell$
a Rademacher probing vector
\[
    \vec{\omega}_\ell
    \coloneqq
    \sum_{i\in \set{C}_\ell}\xi_i\evec_i,
    \qquad
    \xi_i\simiid \Unif\{\pm1\}.
\]
Letting $s$ denote a number of probes per color, we define the stochastic probing estimator as
\begin{equation}
    \hat{\tr}_{\rm SP}
    \coloneqq
    \frac{1}{s} \sum_{\ell=1}^q
    \sum_{j=1}^s
    \bigl(\vec{\omega}_\ell^{(j)}\bigr)^\tp
    \mat{A}\vec{\omega}_\ell^{(j)} \quad \text{with } \vec{\omega}_\ell^{(1)},\ldots,\vec{\omega}_\ell^{(s)} \sim \vec{\omega}_\ell \text{ iid}.
    \label{eq:stochastic-probing}
\end{equation}
It is easily seen that $\hat{\tr}_{\rm SP}$ is an unbiased estimator of $\tr(\mat{A})$, and
\citet{frommer2025analysis} show that the variance is 
\begin{equation}
    \Var(\hat{\tr}_{\rm SP})
    =
    \frac{2}{s}\sum_{\ell=1}^q
    \norm{
      \offdiag\bigl(\mat{A}(\set{C}_\ell,\set{C}_\ell)\bigr)
    }_{\rm F}^2.
    \label{eq:stochastic-probing-variance}
\end{equation}
Observe that \cref{eq:stochastic-probing-variance} reduces to
\cref{thm:variancegh} when there is a single color.
The total number of matvecs needed by the stochastic probing algorithm is $sq$.

The variance formula \cref{eq:stochastic-probing-variance} shows how randomized probing improves on the basic weaknesses of the deterministic probing estimator \cref{eq:vell}.
First, the error depends on the Frobenius norm of the off-diagonal entries, rather than the larger entrywise $\ell_1$ norm.
This difference reflects the possibility of random cancellations when the random probe vectors are randomly signed.
Second, the estimator is systematically improvable by simply taking a larger number of probes $s$ per color.

Returning to the banded example, take one random vector per color
($s=1$) and the periodic coloring~\cref{eq:periodic}.  If $\mat{A}$ is exactly
$(q-1)$-banded, \cref{eq:stochastic-probing-variance} vanishes, recovering the exactness of deterministic probing.  If, more generally, the entries satisfy
a decay of the form
\[
    |\mat{A}(i,j)|\lesssim \rho^{|i-j|},
    \qquad 0<\rho<1,
\]
then distinct indices of the same color are separated by multiples of
$q$, and \cref{eq:stochastic-probing-variance} gives
\begin{equation}
    \Var(\hat{\tr}_{\rm SP})
    \lesssim 
    n
    \rho^{2q}.
\end{equation}
If $\rho$ is not too close to one, a small number of colors thus suffices to attain small error.
Also, the root-mean-square error grows only with the square root of $n$, which compares
favorably with the linear scaling arising in bounds for deterministic probing~\cite{frommer2025analysis}.

\begin{remark}[History and variants]
    The combination of probing with stochastic estimation has a long history, particularly in lattice QCD, where related ideas are known as \emph{partitioning} or \emph{dilution} \cite{bernardson1994monte,foley2005practical,babich2012exploring}. In NLA, randomization was combined with probing in, e.g.,~\cite{stathopoulos2013hierarchical,laeuchli2020extending}.

    Several variants of probing exploit structural information in different ways. Hierarchical probing~\cite{stathopoulos2013hierarchical,laeuchli2020extending} constructs nested families of probing vectors corresponding to increasingly large graph distances. Another approach~\cite{frommer2022multilevel,hallman2022multilevel} adapts multilevel Monte Carlo techniques~\cite{giles2015multilevel}, using a hierarchy of increasingly accurate approximations to $\mat{A}$.
\end{remark}

\subsection{Small-sample estimates and small-ball probability bounds}
\label{sec:small-sample-trace}

In many applications, we require only a rough estimate of the trace that is accurate, say, to an order of magnitude.
In such cases, we are fine with an estimator of poor quality, but we do want strong assurances that the estimator does not severely under- or over-estimate the trace.
This situation arises, for example, when trace estimation approaches are used for posterior error estimation (\cref{sec:error-estimation}).
In this setting, it can be advantageous to tolerate
a rather large multiplicative error in exchange for using only one or a few
samples. This \emph{small-sample} viewpoint goes back to at least~\citet{gudmundsson1995small}. It is closely related to the notion of small-ball estimates in probability theory, which we discuss first.

\subsubsection{Small-ball probability estimates} \label{sec:smallball}

A \emph{small-ball estimate} bounds the probability that a random quantity
is unusually close to a prescribed point, typically zero~\cite{rudelson2015small}.
For a random vector $\vec{x} \in \R^s$, such a bound
often takes the form
\begin{equation} \label{eq:smallballgeneral}
    \prob\{\norm{\vec{x}} \leq \theta \}\leq \phi(\theta) \quad \text{with $\phi(\theta )\to0$ as $\theta \to0$.}
\end{equation}
In the scalar case, if the probability density function of a random variable $z$ is bounded by $K>0$, then we immediately have the small-ball bound $\prob\{ |z| \le \theta \} \le 2K\theta$.

As a simple but important example, consider a Gaussian vector $\vec{g}\sim\Normal(\vec{0},\Id)$, fix any
nonzero vector $\vec{x}\in\R^n$, and define the rescaled overlap $z \coloneqq \vec{g}^\top\vec{x} / \norm{\vec{x}}$.  The random variable $z$ follows a standard normal distribution whose pdf is bounded by $K = 1/\sqrt{2\pi}$.
Consequently,
\begin{equation}
    \prob\left\{
        |\vec{g}^\top\vec{x}|
        \leq
        \theta \norm{\vec{x}}
    \right\}
    \leq
    \sqrt{2/\pi}\,  \theta,
    \qquad \theta> 0.
    \label{eq:gaussian-small-ball}
\end{equation}
Thus, even a single Gaussian inner product $\vec{g}^\top\vec{x}$ is unlikely to underestimate the norm of $\vec{x}$
by a large factor.\footnote{Equivalently, a Gaussian vector has nontrivial overlap with any fixed vector $\vec{x}$ with high probability.
The same conclusion holds for other random vector models, such as uniform random vectors on the sphere.
This ``nontrivial overlap property'' is essential in other areas of RNLA.
In particular, this property is used to show that power and Krylov iteration converge to the maximum eigenvalue of a \psd matrix from random initialization; see \cite{kuczyski1992estimating}, \cite[\S6]{martinsson2020randomized}, and \cite[\S3.1]{kireeva2024randomized}.}

For a vector containing $s$ independent random variables with bounded densities, 
the tensorization lemma~\cite{rudelson2008littlewood} implies that~\cref{eq:smallballgeneral} holds with $\phi(\theta ) = \order(\theta^s)$ for $\theta \approx 0$, significantly reducing the probability of being close to zero. As an example, consider an $n \times s$ Gaussian embedding $\mat{\Omega}$. Then 
$\sqrt{s} \mat{\Omega}^\tp \vec{x}/\norm{\vec{x}} \sim\Normal(\vec{0},\Id_s)$ and the tensorization lemma applies. To obtain a more precise estimate, note that $s \|\mat{\Omega}^\tp \vec{x}\|^2 / \|\vec{x}\|^2 \sim \chi_s^2$ follows a $\chi^2$ distribution with $s$ degrees of freedom.
In turn,
\begin{equation}
    \prob\left\{
\|\mat{\Omega}^\tp \vec{x}\|_2 \le \theta \|\vec{x}\|_2
    \right\}
    =
    F_s(s\theta^2) \le \frac{(s/2)^{s/2}}{\Gamma(s/2+1)}\,\theta^s,
    \label{eq:gaussian-small-ball-m}
\end{equation}
where $F_s$ is the cumulative distribution function of a
$\chi_s^2$ random variable and $\Gamma(z) = \int_0^\infty t^{z-1} e^{-t}\,\mathrm{d}t$ denotes the Gamma function. Quite remarkably, we will see below that the bound~\cref{eq:gaussian-small-ball-m} also controls Frobenius norm estimation. 

\subsubsection{Small-sample trace and Frobenius norm estimation}

One of the main applications of small-sample trace estimation is to estimate the Frobenius norm of a matrix $\mat{B}$.
To do so, we can introduce the \psd matrix $\mat{A}=\mat{B}^\tp\mat{B}$, whose trace is the squared Frobenius norm $\tr(\mat{A}) = \norm{\mat{B}}_{\rm F}^2$.
To estimate, we may apply the Girard--Hutchinson estimator, generating iid isotropic vectors $\vec{\omega}_1,\ldots,\vec{\omega}_s$ and defining the squared-norm estimate
\begin{equation*}
    \psi_s(\mat{B})^2 \coloneqq \frac{1}{s} \sum_{i=1}^s \vec{\omega}_i^\top \mat{A}\vec{\omega}_i^{\vphantom{\top}} = \frac{1}{s} \sum_{i=1}^s \norm{\mat{B}\vec{\omega}_i}^2 = \norm{\mat{B}\mat{\Omega}}^2_{\rm F} \:\:\: \text{for } \mat{\Omega} = \frac{1}{\sqrt{s}} [\vec{\omega}_1,\ldots,\vec{\omega}_s].
\end{equation*}
In this section, we shall develop small-ball estimates for the estimator $\psi_s(\mat{B})$.
Since this estimator is nothing but the Girard--Hutchinson estimator applied to the \psd matrix $\mat{A} = \mat{B}^\top \mat{B}$, the estimates hold equally well for trace estimation for a \psd matrix.

We would like to use a fixed small sample size $s$ and nevertheless ensure that
$\psi_s(\mat{B})$ does not severely underestimate 
$\norm{\mat{B}}_{\rm F}$, that is, 
$\prob\{
        \psi_s(\mat{B})
        \leq
        \theta \norm{\mat{B}}_{\rm F}
\}$ should be tiny for small $\theta > 0$.
The variance and concentration bounds for trace estimation discussed in~\Cref{sec:girardhutchinson} 
do not capture this regime well; in particular, the resulting bounds do not converge to zero as $\theta\to 0$ for fixed $s$.

For a Gaussian random test matrix $\mat{\Omega}$, the SVD of $\mat{B}$
turns $\psi_s(\mat{B})^2$ into a weighted sum of independent chi-squared random variables,
with weights proportional to the squared singular values of $\mat{B}$.
Extremal results for Gaussian quadratic forms~\cite{szekely2003extremal,roostakhorasani2015assessing}
identify the singular value configurations that maximize tail probabilities. For underestimation below the mean, the worst case is attained by concentrating all mass into a single singular value, that is, a rank-one matrix $\mat{B}$. Consequently,~\cref{eq:gaussian-small-ball-m} implies
\begin{equation}
\prob\{
        \psi_s(\mat{B})
        \leq
        \theta \norm{\mat{B}}_{\rm F}
\}
    \leq
    F_s(s\theta^2), \qquad 0 < \theta < 1,
    \label{eq:small-sample-lower}
\end{equation}
for every matrix $\mat{B}$.  Thus, the probability of severe
underestimation decays like $\theta^s$ as $\theta \to 0$.
For example, we have 
    $\prob\{
        \psi_2(\mat{B})
        \leq
        \norm{\mat{B}}_{\rm F}/10
    \}
    < 0.01$.
Thus, only two Gaussian probes suffice to obtain an
order-of-magnitude upper bound with $99\%$ reliability for the Frobenius norm of any matrix $\mat{B}$.

If an actual constant-factor approximation to $\norm{\mat{B}}_{\rm F}$ is
desired, one must also control overestimation.  This direction is considerably
easier.
Introduce the \emph{stable rank}
$
    \operatorname{sr}(\mat{B})
    \coloneqq
    \norm{\mat{B}}_{\rm F}^2 / \norm{\mat{B}}^2
$
of the matrix $\mat{B}$, which is a smooth proxy for the algebraic rank \cite{ipsen2025stable}.
A Chernoff bound of~\citet{gratton2018improved} gives
\begin{equation}
    \prob\{
        \psi_s(\mat{B})
        \geq
        \tau\norm{\mat{B}}_{\rm F}
    \}
    \leq
    \exp\left(
        -\frac{s\operatorname{sr}(\mat{B})}{2}(\tau-1)^2
    \right), \qquad \tau>1.
    \label{eq:small-sample-upper}
\end{equation}
The improvement with increasing stable rank is consistent with the observation
in~\Cref{sec:girardhutchinson}. Although
$\operatorname{sr}(\mat{B})$ is generally unknown, the universal bound
obtained by replacing 
$\operatorname{sr}(\mat{B})$ with its lower bound $1$ already yields an effective estimate.
Recent work of~\citet{hallman2024extremal} also sharpens the underestimation bound~\cref{eq:small-sample-lower} by incorporating information on the stable rank.

Combining~\cref{eq:small-sample-upper} with~\cref{eq:small-sample-lower} and a union
bound immediately yields a lower bound on
$\prob\left\{
        \frac{1}{\tau}\norm{\mat{B}}_{\rm F}
        <
        \psi_s(\mat{B})
        <
        \tau\norm{\mat{B}}_{\rm F}
    \right\}$.
For illustration, replacing $\operatorname{sr}(\mat{B})$ by $1$, this gives the following
success probabilities, valid uniformly for every matrix $\mat{B}$
for a factor-ten approximation:
\begin{equation} \label{eq:small-sample-probs}
    \prob\{
            \norm{\mat{B}}_{\rm F}/10
            \leq \psi_s(\mat{B})
            \leq 10\norm{\mat{B}}_{\rm F}
        \} \ge \begin{cases}
        0.920, & s = 1, \\
        0.990, & s = 2, \\
        0.999, & s = 3.
        \end{cases}
\end{equation}

\subsection{Error estimation and automatic stopping control}

As with any randomized algorithm, it is useful to couple stochastic trace estimation algorithms with error estimation techniques, both to inform the user about the quality of the computed output and to allow the algorithm to be run adaptively, continuing to collect more information until the error meets a threshold.

The error of the Girard--Hutchinson estimator is easy to estimate empirically.
The basic trace estimates $\hat{\tr}_i = \vec{\omega}_i^\top \mat{A}\vec{\omega}_i^{\vphantom{\top}}$ are independent of one another, so an unbiased estimator for the variance of the Girard--Hutchinson estimator is given by the sample variance:
\begin{equation} \label{eq:trace-sample-stdev}
    \widehat{\Var} \coloneqq  \frac{1}{s(s-1)} \sum_{i=1}^s
      \Bigl( \hat{\tr}_i
             - \hat{\tr}_{\rm GH}(\mat{A}) \Bigr)^2 .
\end{equation}
This estimator can also be adapted for Hutch\plusplus and its variants, as the basic estimates $\hat{\tr}_i \coloneqq \tr(\mat{Q}^\tp\mat{A}\mat{Q}) + \vec{\omega}_i^\top (\Id - \mat{Q}\mat{Q}^\top)\mat{A}(\Id - \mat{Q}\mat{Q}^\top)\vec{\omega}_i$ are independent conditional on $\mat{Q}$.
An adaptive version of Hutch\plusplus that both determines the number of matvecs and apportions them between low-rank approximation and residual trace estimation is developed in \cite{persson2022improved}.
A variant of the estimator \cref{eq:trace-sample-stdev} can also be used with stochastic probing.

Error estimation for exchangeable algorithms like \XTrace is more difficult, as the algorithm produces basic trace estimates \cref{eq:xtrace-replicates} that are \emph{dependent}.
Nonetheless, the estimators seem to be weakly correlated in practice, and a standard deviation-type variance estimator seems effective in practice; see \cite[\S3]{epperly2024xtrace}.

\subsection{Estimating the diagonal and other linear quantities} \label{sec:diagonal}

The trace estimation idea extends naturally to estimate other linear quantities of a matrix.
Let $\mathcal{L} : \R^{n\times n} \to \set{V}$ denote any linear map and
let $\vec{\omega}_1,\ldots,\vec{\omega}_s$ be iid isotropic random vectors.
Replacing the identity matrix $\Id$ by its unbiased stochastic approximation $\frac1s \sum_{i=1}^s \vec{\omega}_i^{\vphantom{\tp}} \vec{\omega}_i^\tp$ gives the unbiased estimator
\begin{equation} \label{eq:linear-attribute-estimator}
    \mathcal{L}(\mat{A}) = \mathcal{L}(\mat{A}\cdot \Id) \approx \mathcal{L}\Big(\mat{A} \cdot \frac1s \sum_{i=1}^s \vec{\omega}_i^{\vphantom{\tp}} \vec{\omega}_i^\tp\Big) = \frac{1}{s} \sum_{i=1}^s \mathcal{L}\big((\mat{A}\vec{\omega}_i^{\vphantom{\tp}})\vec{\omega}_i^\tp\big).
\end{equation}
The Girard--Hutchinson estimator is generated by this construction with $\mathcal{L} = \tr$.
Many of the variance-reduction ideas discussed above can
also be adapted to other choices of $\mathcal{L}$.

\paragraph{Diagonal estimation.}
An important example is $\mathcal{L}=\diag$; applications for diagonal estimation appear in, e.g., electronic structure calculation~\cite{bekas2007estimator} and uncertainty quantification in Bayesian inversion~\cite{cho2022computationally}.
Specializing~\cref{eq:linear-attribute-estimator} to Rademacher probes gives the
classical estimator of \citet{bekas2007estimator}:
\begin{equation} \label{eq:diagrademacher}
    \hat{\diag}_{\rm BKS} \coloneqq \frac{1}{s} \sum_{i=1}^s \mat{A}\vec{\omega}_i \odot \vec{\omega}_i,
\end{equation}
where $\odot$ denotes entrywise multiplication. A direct calculation gives 
\begin{equation} \label{eq:bks-error}
\expect\bigl[\big\| \hat{\diag}_{\rm BKS} -\diag(\mat{A}) \big\|^2\bigr]
    =
    \frac{1}{s}\norm{\offdiag(\mat{A})}_{\rm F}^2,
\end{equation}
Thus, just as for Rademacher trace estimation, accuracy improves when
$\mat{A}$ is close to diagonal.
Concentration bounds and other probe distributions have been studied in~\cite{baston2022stochastic,dharangutte2023tight,hallman2023monte}.
For \psd matrices, \citet{baston2022stochastic} propose Diag\plusplus, which extends Hutch\plusplus by combining stochastic diagonal
estimation with randomized low-rank approximation. The exchangeability principle underlying \XTrace leads to the
\XDiag estimator, which applies to general
square matrices, and the \XNysDiag estimator, for \psd matrices~\cite{epperly2024xtrace,epperly2025make}.

\paragraph{Beyond the diagonal: Sparse matrix recovery.}
Diagonal estimation is the simplest instance of the more general problem of
recovering a sparse matrix from matrix-vector products. Let
$\mat{S}\in\{0,1\}^{n\times n}$ denote a prescribed sparsity pattern and consider the
linear map $\mathcal{L}_{\mat{S}}: \R^{n\times n} \to \set{V}$ defined by
\[
    \mathcal{L}_{\mat{S}}(\mat{A})
    \coloneqq \mat{S}\odot\mat{A},
\]
where $\set{V}$ is the image of $\mathcal{L}_{\mat{S}}$ and $\odot$ is again the entrywise multiplication operator.
Sparse matrix recovery has a long history, originating in methods for computing
sparse Jacobians and Hessians. In analogy to probing for trace estimation described in~\Cref{sec:stochastic-probing},
columns with non-overlapping support can be recovered from a single matrix--vector product.  Finding a small number of such groups of columns leads, again, to graph-coloring problems; see~\cite{gebremedhin2005what} for a survey. 

When $\mat{A}$ is only approximately sparse, randomization can be used to
reduce  contamination from entries outside the prescribed sparsity pattern,
in the spirit of stochastic probing in~\Cref{sec:stochastic-probing}.
Recently, \citet{amsel2026fixed} gave a general treatment of this setting,
showing that, when the prescribed sparsity pattern contains at most $\zeta$
nonzeros per row, a near-optimal approximation can be obtained
from $\order(\zeta/\varepsilon)$ Gaussian matrix--vector queries with $\mat{A}$
alone. When products with both $\mat{A}$ and $\mat{A}^\tp$ are available,
\citet{musco2026nearly} broaden the scope of sparsity patterns and show that the query complexity is primarily governed by the so-called degeneracy of $\mat{S}$.

\subsection{Row-norm estimation}
\label{sec:row-norm-estimation}

Estimating the (squared) row norms of a matrix $\mat{B}\in\R^{p\times n}$ is 
formally a diagonal estimation problem:
\begin{equation}
    \vec{r}
    \coloneqq
    \begin{bmatrix}
        \norm{\mat{B}(1,:)}^2 & \cdots & \norm{\mat{B}(p,:)}^2
    \end{bmatrix}^\tp
    = \diag(\mat{B}\mat{B}^\tp),
    \label{eq:row-norm-vector}
\end{equation}
However, there are good reasons to not estimate the squared row norms using an off-the-shelf diagonal estimator like the BKS diagonal estimator \cref{eq:diagrademacher}.
First, and less seriously, applying a diagonal estimator to $\mat{A} = \mat{B} \mat{B}^\tp$ requires matrix--vector products with both $\mat{B}$ and $\mat{B}^\tp$, and the latter may not be available.
A second issue, more serious, is corruption between rows.
As \cref{eq:bks-error} shows, the error of the BKS estimator \cref{eq:diagrademacher} depends on off-diagonal entries of the matrix $\mat{A}$, which are the pairwise row inner products $\mat{B}(i,:)\mat{B}(j,:)^\top$ for $i\ne j$.
In particular, the presence of a single large row can introduce large errors in the estimates for \emph{all} of the row norms.

The Johnson--Lindenstrauss (JL) row-norm estimator provides a simpler, more effective alternative.
Let $\vec{\omega}_1,\ldots,\vec{\omega}_s$ be iid isotropic vectors, and form the random test matrix $\mat{\Omega} = s^{-1/2}[\vec{\omega}_1,\ldots,\vec{\omega}_s]$.
We obtain an unbiased estimate for the squared row norms of $\mat{B}$ by taking the squared row norms of the sketch $\mat{Y}=\mat{B}\mat{\Omega}$:
\begin{equation}
 \widehat{\vec{r}}_{\rm JL} \coloneqq
\begin{bmatrix}
        \norm{\mat{Y}(1,:)}^2 & \cdots & \norm{\mat{Y}(p,:)}^2
    \end{bmatrix}^\tp.
    \label{eq:basicrowest}
\end{equation}

The JL estimator is unbiased, since 
\begin{equation*}
    \expect[\widehat{\vec{r}}_{\rm JL}] = \expect[\diag(\mat{B}\mat{\Omega}\mat{\Omega}^\top\mat{B}^\top)] = \diag(\mat{B} \expect[\mat{\Omega}\mat{\Omega}^\top]\mat{B}^\top) = \diag(\mat{B}\mat{B}^\top) = \vec{r}.
\end{equation*}
Moreover, if the random probes are Gaussian, Rademacher, or spherical, then \cref{thm:variancegh} implies that the root-mean-square error for each squared row-norm estimate $\hat{\vec{r}}(j)$ is at most
\begin{equation*}
    \left( \expect \left( \frac{\hat{\vec{r}}(j) - \vec{r}(j)}{\vec{r}(j)} \right)^2 \right)^{1/2} \le \sqrt{\frac{2}{s}}.
\end{equation*}
Unlike the BKS estimator applied to $\mat{B}\mat{B}^\top$, the JL estimator always computes \emph{each} squared row norm to small root-mean-square error.

As we discuss below, the JL estimator emerges naturally from the widely used Johnson--Lindenstrauss lemma (\cref{thm:jl}).
Early uses of the JL estimator and its variants in RNLA include \cite{spielman2011graph,drineas2012fast,li2013iterative}.

\subsubsection{Uniform error estimates and the Johnson--Lindenstrauss lemma} \label{sec:jl}

In many applications, we require row-norm estimates such that all norms are simultaneously estimated to small relative error.
The JL estimator can accomplish this task with small overhead.
Assume for simplicity that the base distribution is Gaussian.
Under this condition, the rescaled ratio $s\widehat{\vec{r}}(j) / \vec{r}(j) \sim \chi^2_s$ follows a $\chi^2$ distribution with $s$ degrees of freedom for each nonzero row $j$.
By exploiting standard $\chi^2$ concentration bounds \cite[Lem.~1]{laurent2000adaptive}, we may conclude the following result.

\begin{theorem}[Johnson--Lindenstrauss lemma] \label{thm:jl}
    Fix $\delta \in (0,1)$ and $\varepsilon \in (0,1/2]$.
    Let $\mat{B} \in \R^{p\times n}$ be a matrix, and let $\hat{\vec{r}} \in \R^p$ denote the JL squared row-norm estimator implemented with a Gaussian embedding $\mat{\Omega} \in \R^{n\times s}$ of size $s\ge 6\log(2p/\delta)/\varepsilon^2$.
    Then
    \begin{equation*}
        \max_{\vec{r}(j) \ne 0} \left| \frac{\hat{\vec{r}}(j) - \vec{r}(j)}{\vec{r}(j)} \right| \le \varepsilon \quad \text{with probability }\ge 1-\delta.
    \end{equation*}
\end{theorem}

This result, essentially due to \cite{johnson1984extensions}, shows that we can simultaneously estimate the (squared) norm of all $p$ rows of a matrix to relative accuracy $\varepsilon$ by sketching the matrix down to $\order(\log(p)/\varepsilon^2)$ dimensions.
Note that the sketching dimension $s$ is only \emph{logarithmic} in the number of rows $p$.

We have stated \cref{thm:jl} as a result about the JL row-norm estimation algorithm, but this is highly nontraditional.
Typically, this result is called the \emph{Johnson--Lindenstrauss lemma}, and it is rightfully presented as a deep fact about random embeddings and high-dimensional geometry.
Indeed, thinking of the rows of $\mat{B}$ as $p$ points in $\R^n$, the rows of $\mat{Y} = \mat{B}\mat{\Omega} \in \R^{p\times s}$ represent the embedding of each high-dimensional point $\mat{B}(j,:) \in \R^n$ into a low-dimensional representation $\mat{Y}(j,:) \in \R^s$.
The JL lemma states that the embedding dimension $s = \order(\log p)$ suffices to embed all $p$ points into $s$-space while preserving their lengths up to a constant factor.
Moreover, since $\log {p \choose 2} \le \log p^2 \le 2\log p$, an embedding dimension $s = \order(\log p)$ also suffices to preserve all \emph{pairwise} distances among the $p$ points.

\begin{remark}[JL vs.\ subspace embedding]
    The Johnson--Lindenstrauss lemma is a beautiful and useful fact about the capability of random embeddings.
    But its placement this far into the survey is not an accident.
    For linear algebra applications, we typically require random embeddings that preserve not just a finite set of points, but an entire linear subspace (i.e., a subspace embedding or injection, \cref{sec:embedinject}).
    The difference is important:
    To embed a finite set of $p$ points requires embedding into dimension $s = \order(\log p)$, whereas embedding an $r$-dimensional linear subspace requires a dimension $s = \order(r)$.
    When implementing an algorithm, it is essential to know which type of guarantee is needed (JL or subspace embedding) and to pick the dimension $s$ appropriately.
\end{remark}

\subsubsection{Applications}

One of the most basic applications of the JL estimator and stochastic row-norm estimation technique is the estimation of leverage scores, which we shall discuss in \cref{sec:leverage-estimation}.
Similar ideas are used to estimate effective resistances in graphs \cite{spielman2011graph}.

As we saw above, the JL estimator is capable of simultaneously estimating all squared row norms of a matrix to small relative error, even when the rows differ widely in scale.
This relative-error property is not possessed by the BKS diagonal estimator, even when applied to a \psd matrix $\mat{A}$.
This disparity suggests that it may be fruitful to go about estimating the diagonal of $\mat{A}$ through estimating the squared row norms of its square root $\mat{A}^{1/2}$.
This idea is explored by \citet{mathur2021variance}, and algorithmic applications of this \emph{square-root trick} were investigated by \citet{lindsey2026fast,fornace2026column}.
See \cite{hsu2026randomized} for more on the square-root trick and applications of row-norm estimation.

\subsubsection{Variance reduction by low-rank approximation}

As with trace and diagonal estimation, low-rank approximation can improve on the JL estimator~\cref{eq:basicrowest} by combining it with variance reduction techniques.
\citet{sobczyk2022approximate} develop and analyze a Hutch\plusplus-inspired row-norm estimator.
Exchangeable leave-one-out versions, called \textsc{XRowNorm} and, for symmetric
matrices, \textsc{XSymRowNorm}, appear in \cite{epperly2025make,hsu2026randomized}.

\section{Spectral sums and spectral density estimation}
\label{sec:sde}

An important special case of the trace estimation problem is the \emph{spectral sum problem}, where we wish to compute the trace of a function of a matrix $\tr(f(\mat{A}))$ for a symmetric matrix $\mat{A} \in \R^{n\times n}$.
Equivalently, we wish to evaluate the sum of $f$ applied to the eigenvalues of $\mat{A}$:
\begin{equation} \label{eq:spectral-sum}
    \tr(f(\mat{A})) = \sum_{i=1}^n f(\lambda_i(\mat{A})).
\end{equation}
One approach to this estimation task is to decouple the problem into two components: First, devise a scheme to (approximately) compute matrix-function-vector products $\vec{\omega} \mapsto f(\mat{A})\vec{\omega}$, for example 
using the Krylov methods discussed in~\Cref{sec:f(A)b}.
Then, apply a general-purpose stochastic trace estimator, as discussed in the last section.
As we shall see, treating these two components
jointly can be more effective, particularly when we wish to estimate the spectral sums $\tr(f_1(\mat{A})),\ldots,\tr(f_t(\mat{A}))$ associated with many different functions $f_i$.

As it happens, the spectral sum estimation problem is closely connected to the problem of estimating the \emph{distribution of eigenvalues} of a matrix.
Indeed, if we let
\begin{equation} \label{eq:global-spectral-measure}
    \mu_{\mat{A}} \coloneqq \sum_{i=1}^n \delta_{\lambda_i(\mat{A})}
\end{equation}
be the (unnormalized) spectral measure that places a unit mass on each eigenvalue, the spectral sum \cref{eq:spectral-sum} may be expressed as an integral:
\begin{equation} \label{eq:spectral-sum-and-spectral-measure}
    \tr(f(\mat{A})) = \sum_{i=1}^n f(\lambda_i(\mat{A})) = \int f \, \d\mu_{\mat{A}}.
\end{equation}
One can go the other way as well: Given estimates of spectral sums for judicious functions $f_i$, one can recover approximations to the spectral measure $\mu_{\mat{A}}$ (see \cref{sec:cheby-mom-match,sec:kpm}).

After a short discussion of background material in \cref{sec:spectral-measures}, we discuss three approaches for estimating the distribution of eigenvalues: stochastic Lanczos quadrature (\cref{sec:slq}), Chebyshev moment matching (\cref{sec:cheby-mom-match}), and the kernel polynomial method (\cref{sec:kpm}).
By integrating against any of the spectral measure estimates produced by these algorithms, we obtain estimates for $\tr(f(\mat{A}))$.
\Cref{sec:deflation-variance-reduction} discusses deflation and variance reduction techniques.

\subsection{Background: Spectral measures} \label{sec:spectral-measures}

In this section, we assume some familiarity with measure theory, mainly for notation; no deep measure-theoretic results are needed.

For our purposes, we will need two notions for the distribution of eigenvalues of a symmetric matrix $\mat{A}$.
The first is the (global) spectral measure \cref{eq:global-spectral-measure}, which places a unit mass on each eigenvalue.\footnote{This measure is often scaled to have total mass $1$; we do not do so here.}
Given the spectral measure, we can recover any spectral sum $\tr f(\mat{A})$ using \cref{eq:spectral-sum-and-spectral-measure}.
As we have already seen, a basic tool in trace estimation is to estimate a trace $\tr f(\mat{A})$ using a quadratic form $\vec{\omega}^\tp f(\mat{A})\vec{\omega}$.
The quadratic form, too, may be written as an integral
\begin{equation*}
    \vec{\omega}^\tp f(\mat{A})\vec{\omega} = \int f \, \d\mu_{\mat{A},\vec{\omega}}
\end{equation*}
where
\begin{equation*}
    \mu_{\mat{A},\vec{\omega}} \coloneqq \sum_{i=1}^n |\vec{\omega}^\tp \vec{u}_i|^2 \delta_{\lambda_i(\mat{A})}
\end{equation*}
is called the local spectral measure or eigenvector spectral measure.
It is defined using an eigendecomposition $\mat{A} = \sum_{i=1}^n \lambda_i(\mat{A}) \, \vec{u}_i^{\vphantom{\tp}}\vec{u}_i^\tp$.

For an isotropic random vector $\vec{\omega}$, we know that $\expect [\vec{\omega}^\tp f(\mat{A})\vec{\omega}] = \tr(f(\mat{A}))$.
In much the same way, it holds formally that $\expect[\mu_{\mat{A},\vec{\omega}}] = \mu_{\mat{A}}$; under an isotropic random vector, the expected value of the local spectral measure is the global spectral measure.
This suggests Monte Carlo strategies similar to those we employed in basic stochastic trace estimation: Draw different random probes $\vec{\omega}_1,\ldots,\vec{\omega}_s$ and estimate the local spectral measure associated to each.
Then, average to obtain an approximation to the global spectral measure.

\subsection{Stochastic Lanczos quadrature} \label{sec:slq}

One of the most effective methods for the spectral sum problem is the stochastic Lanczos quadrature method.
The root of the method is a powerful connection between the Lanczos method\footnote{As in~\Cref{sec:krylovreview}, we assume throughout this section that Krylov subspaces attain their maximal dimension, so that the Lanczos process does not break down.}
and Gaussian quadrature.

\begin{theorem}[Lanczos and Gaussian quadrature] \label{thm:lanczos-gaussian-quadrature}
    Let $\mat{A} \in \R^{n\times n}$ be symmetric, and let $\mat{Q} \in \R^{n\times k}$ and $\mat{T} \in \R^{k\times k}$ be the orthonormal and tridiagonal matrices obtained from $k$ steps of the Lanczos process applied to $\mat{A}$ with starting vector $\vec{\omega} \ne \vec{0}$.
    Then 
    \begin{equation*}
        \hat{\mu} \coloneqq \norm{\vec{\omega}}^2 \cdot \mu_{\mat{T},\evec_1} = \sum_{j=1}^k \norm{\vec{\omega}}^2 |\vec{u}_j(1)|^2 \cdot \delta_{\lambda_j(\mat{T})} \quad \text{with } \mat{T} = \sum_{j=1}^k \lambda_j(\mat{T})\vec{u}_j^{\vphantom{\tp}}\vec{u}_j^\tp
    \end{equation*}
    defines a $k$-point Gaussian quadrature approximation to the local spectral measure $\mu_{\mat{A},\vec{\omega}}$.
    That is,
    \begin{equation*}
        \int p \, \d\hat{\mu} = \int p \, \d\mu_{\mat{A},\vec{\omega}} = \vec{\omega}^\tp p(\mat{A})\vec{\omega} \quad \text{for every polynomial $p$ of degree $<2k$}.
    \end{equation*}
\end{theorem}

This connection between the Lanczos method and quadrature rules is one of the most beautiful pieces of mathematics in the numerical analysis canon; it is described carefully by \citet{golub2010matrices,chen2024lanczos}.
For us, the core takeaway is that the Lanczos method gives us approximations to any local spectral measure of our choosing.
These approximations are exact for low-degree polynomials, yielding the following error bound for a function $f$ defined on 
an interval $[a,b]$ containing the eigenvalues of $\mat A$:
\begin{equation}
    \left|
    \vec{\omega}^\tp f(\mat{A})\vec{\omega}
    -\norm{\vec{\omega}}^2
     \evec_1^\tp f(\mat{T})\evec_1
    \right|
    \le
    2\norm{\vec{\omega}}^2
    \min_{\deg(p) < 2k} \max_{x \in [a,b]} |f(x) - p(x)|.
    \label{eq:lanczos-quadrature-poly}
\end{equation}
This suggests excellent approximation quality for analytic functions $f$, which admit exponentially accurate polynomial approximations~\cite{trefethen2019approximation}.

To approximate the global spectral measure, we use Lanczos quadrature to approximate the local spectral measure $\hat{\mu}_i\approx \mu_{\mat{A},\vec{\omega}_i}$ associated with $s$ different random probes $\vec{\omega}_i$ and return the average $s^{-1}\sum_{i=1}^s \hat{\mu}_i$.
\Cref{alg:slq} goes back to~\citet{bai1996some}, which was later developed
and analyzed systematically as \emph{stochastic Lanczos quadrature}
by \citet{ubaru2017fast}; see also \cite{cortinovis2022randomized,chen2021analysis,chen2025randomized}.

\begin{algorithm}[t]
\caption{Stochastic Lanczos quadrature for spectral density estimation} \label{alg:slq}
\begin{algorithmic}[1]
\Require Symmetric matrix $\mat{A} \in \R^{n\times n}$, number of probes $s$, Lanczos depth $k$
\Ensure Discrete measure $\hat\mu \coloneqq s^{-1}\sum_{i=1}^s\sum_{j=1}^k w_{ij}\delta_{\lambda_{ij}} \approx \mu_{\mat{A}}$
\For{$i = 1, \dots, s$}
    \State Draw an isotropic random vector $\vec{\omega}_i$ \Comment{E.g., Rademacher or spherical}
    \State Run Lanczos: $\mat{T}\gets \Call{Lanczos}{\mat{A},\vec{\omega}_i,\texttt{steps}=k}$
    \State Compute the eigendecomposition $(\mat{U},\mat{\Lambda}) \gets \Call{Eig}{\mat{T}}$
    \State Define nodes $(\lambda_{i1},\ldots,\lambda_{ik}) \gets \diag(\mat{\Lambda})$
    \State Define weights $w_{ij} \gets  \norm{\vec{\omega}_i}^2 |\mat{U}(1,j)|^2$ for each $j=1,\ldots,k$
\EndFor
\State \Return $\hat{\mu} \gets s^{-1}\sum_{i=1}^s\sum_{j=1}^k w_{ij}\delta_{\lambda_{ij}}$
\end{algorithmic}
\end{algorithm}

\subsection{Background: Chebyshev moments}

It is a core fact of probability theory that a bounded random variable is determined by its (polynomial) moments.
That is, if $x$ is a random variable supported on $[-1,1]$ and we know the value of $\expect x^k$ for each $k\ge 0$, then we can reconstruct the distribution of $x$ to perfect accuracy \cite[{\S}VII.3]{feller1971introduction}.
Numerically, recovering a distribution from monomial moments is horribly ill-conditioned, but the conditioning issue is resolved by switching to the Chebyshev moments $\expect \mathrm{T}_k(x)$.
Here, $\mathrm{T}_k(x) = \cos(k\arccos x)$ denote the Chebyshev polynomials of the first kind.
Using the same principle, we can use estimates of the Chebyshev moments of a matrix to approximate its spectral measure.

Assume $\mat{A}$ has its spectrum in a \emph{known} interval $[a,b]$ and scale and shift to make this interval $[-1,1]$.
The Chebyshev moments of $\mat{A}$ are defined as
\begin{equation} \label{eq:cheb-moment}
    c_k(\mat{A}) \coloneqq \tr(\mathrm{T}_k(\mat{A})) \quad \text{for } k=0,1,2,\ldots.
\end{equation}
Note that $c_0(\mat A)=n$. It is straightforward to show that the 
Chebyshev moments $c_1(\mat A), \ldots, c_n(\mat A)$ uniquely determine the spectral measure $\mu_{\mat{A}}$.
As we will see, estimates of just a few Chebyshev moments can already lead to good approximations.

\begin{algorithm}[t]
\caption{Stochastic estimation of Chebyshev moments (\textsc{ChebMoments})} \label{alg:cheb-moments}
\begin{algorithmic}[1]
\Require Symmetric matrix  $\mat{A} \in \R^{n\times n}$ with spectrum in $[-1,1]$, number of probes $s$, even degree $K$
\Ensure Estimates $\hat{c}_k \approx c_k(\mat{A}) = \int \mathrm{T}_k \,\d\mu_{\mat{A}}$ for $k = 1,\ldots,K$
\State Draw $\mat{\Omega} \in \R^{n\times s}$ with columns drawn uniformly from the sphere of radius $\sqrt{n}$
\State Initialize the recurrence $\mat{V}_{0} \gets \mat{\Omega}$, $\mat{V}_{1} \gets \mat{A}\mat{\Omega}$
\State Set $\hat{c}_1\gets s^{-1}\langle \mat{V}_0, \mat{V}_1\rangle_{\rm F}$, and $\hat{c}_2\gets 2s^{-1}\norm{\mat{V}_1}_{\rm F}^2-n$
\For{$k = 2, \dots, K/2$}
    \State Advance the recurrence $\mat{V}_k \gets 2\mat{A}\mat{V}_{k-1} - \mat{V}_{k-2}$ \Comment{$\mat{V}_k = \mathrm{T}_k(\mat{A})\mat{\Omega}$}
    \State Set $\hat{c}_{2k} \gets 2s^{-1}\norm{\mat{V}_k}_{\rm F}^2 - n$ \Comment{$\mathrm{T}_{2k} = 2\mathrm{T}_k^2 - \mathrm{T}_0$}
    \State Set $\hat{c}_{2k-1} \gets s^{-1}\bigl(2\,\langle \mat{V}_{k-1},\mat{V}_{k}\rangle_{\rm F} - \langle\mat{V}_0,\mat{V}_1\rangle_{\rm F}\bigr)$ \Comment{$\mathrm{T}_{2k-1} = 2\mathrm{T}_k\mathrm{T}_{k-1} - \mathrm{T}_1$}
\EndFor
\State \Return $(\hat{c}_1, \ldots, \hat{c}_K)$
\end{algorithmic}
\end{algorithm}

We can estimate the Chebyshev moments using the Girard--Hutchinson estimator and the three-term recurrence for Chebyshev polynomials to compute $\mathrm{T}_k(\mat{A})\vec{\omega}_i$.
See \cref{alg:cheb-moments}.
This code also uses a standard ``doubling trick'' to estimate $K$ Chebyshev moments using $K/2$ matrix--vector products
with a block of $s$ vectors \cite[eqs.~(34) \& (35)]{weie2006kernel}.

\subsection{Chebyshev moment regression} \label{sec:cheby-mom-match}

The first approach to using the Chebyshev moments to estimate the distribution of eigenvalues is direct \cite{braverman2022sublinear}.
Given approximations $\hat{c}_k \approx c_k(\mat{A})$ to the leading $K$ Chebyshev moments, 
we estimate the spectral measure by approximately matching these moments.
\citet{musco2025sharper} suggest the optimization problem
\begin{equation} \label{eq:mmrs-moment-matching}
    \operatorname{minimize} \sum_{k=1}^K \frac{1}{k^2} \left( \,\hat{c}_k - \int \mathrm{T}_k \, \d\hat{\mu}\, \right)^2 \quad \text{s.t. $\hat{\mu}$ is a measure with total mass $n$}.
\end{equation}
The coefficients $1/k^2$ reflect the fact that matching the lower-order moments is more important to estimate the spectral distribution to small Wasserstein-1 error; see \cref{thm:mmrs-moment-matching}.
To make \cref{eq:mmrs-moment-matching} into a finite-dimensional problem, we can restrict $\hat{\mu}$ to be a measure supported on a 
grid of $q$ points; \citet{musco2025sharper} suggest $q \coloneqq \lceil K^{1.5}\rceil$ Chebyshev nodes.
See \cref{alg:moment-matching} for pseudocode.

\begin{algorithm}[t]
\caption{Chebyshev moment regression for spectral density estimation} \label{alg:moment-matching}
\begin{algorithmic}[1]
\Require Symmetric matrix  $\mat{A} \in \R^{n\times n}$ with spectrum in $[-1,1]$, number of probes $s$, even number of moments $K$, grid size $q$
\Ensure Discrete measure $\hat\mu \coloneqq \sum_{j=1}^q w_j \delta_{x_j} \approx \mu_{\mat{A}}$
\State Estimate $(\hat{c}_1,\ldots,\hat{c}_K) \gets \Call{ChebMoments}{\mat{A},s,K}$ \Comment{\cref{alg:cheb-moments}}
\State Set $
    x_j=\cos\!\big(\frac{(2j-1)\pi}{2q}\big)$ for $j=1,\ldots,q$
\State Form the moment matrix $\mat{M}(k,j) \gets (\mathrm{T}_k(x_j)) \in \R^{K\times q}$
\State Form the weight matrix $\mat{D} \gets \diag(1,1/2,\ldots,1/K)$
\State Solve the convex program
\Statex \qquad $\displaystyle \vec{w} \gets \argmin_{\vec{w}\in\R^q} \; \norm{\mat{D}(\hat{\vec{c}} - \mat{M}\vec{w})}^2 \quad \text{s.t. } \vec{w}\ge 0, \; \textstyle\sum_{j=1}^q w_j = n$
\State \Return $\hat\mu \gets \sum_{j=1}^q w_j \delta_{x_j}$
\end{algorithmic}
\end{algorithm}

\subsection{The kernel polynomial method} \label{sec:kpm}

The kernel polynomial method \cite{weie2006kernel} provides a different way of constructing an approximation to the spectral measure from estimated Chebyshev moments.

Let us first derive the method formally.
Begin by using the orthogonality of the Chebyshev series to deduce the relation
\begin{equation*}
    f(x) = \sum_{k=0}^\infty a_k \mathrm{T}_k(x) \int_{-1}^1 \mathrm{T}_k(y) f(y) \, w(y)\d y \quad \text{for } a_k = \begin{cases}
        \frac{1}{2} & k = 0, \\
        1 & k \ge 1.
    \end{cases}
\end{equation*}
Here, $w(y) = \tfrac{2}{\pi}(1-y^2)^{-1/2}$ is the Chebyshev density.
Consequently, integrals against the spectral measure take the form
\begin{multline*}
    \int f(x) \, \d\mu_{\mat{A}}(x) = \int \left(\sum_{k=0}^\infty a_k \mathrm{T}_k(x) \int_{-1}^1 \mathrm{T}_k(y) f(y) \, w(y)\d y \right) \, \d\mu_{\mat{A}}(x) \\
    = \int_{-1}^1 \left[\sum_{k=0}^\infty a_k \left(\int \mathrm{T}_k(x) \, \d\mu_{\mat{A}}(x)\right) \mathrm{T}_k(y) w(y)\right]  f(y) \, \d y = \int_{-1}^1 f(y) \, p(y)\d y
\end{multline*}
with
\begin{equation*}
    p(y) \coloneqq \frac{n}{2} w(y) + \sum_{k=1}^\infty c_k(\mat{A})\mathrm{T}_k(y)\,w(y).
\end{equation*}
We see that integrals against the spectral measure $\mu_{\mat{A}}$ can be computed as ordinary integrals of $f$ with respect to the ``weight function'' $p(y)$.
Moreover, this weight function is entirely determined by the Chebyshev moments $c_k(\mat{A})$.

Let us emphasize that this is a formal derivation.
The spectral measure $\mu_{\mat{A}}$ is discrete, so the relation $\tr(f(\mat{A})) = \int_{-1}^1 f(y) \, p(y)\d y$ cannot hold for any $p$ that is a function in the conventional sense of the word.
Nonetheless, this derivation is useful.
Truncating the
expansion at degree $K$ gives an absolutely continuous approximation
$p_K(y)\,\d y$ to the spectral measure, where
$p_K(y) = \tfrac{n}{2}w(y) + \sum_{k=1}^K c_k(\mat{A}) \mathrm{T}_k(y) \, w(y)$, for which $\tr(f(\mat{A})) \approx \int_{-1}^1 f(y) \, p_K(y) \,\d y$.
Directly truncating can introduce spurious oscillations in the spectral density and result in a density that is negative in some places; these issues can both be ameliorated by using a damped version of the spectral density
\begin{equation} \label{eq:kpm-estimate}
    \hat{p}_K(y) = \frac{n}{2}w(y) + \sum_{k=1}^K \beta_k c_k(\mat{A})\, \mathrm{T}_k(y)\,w(y).
\end{equation}
One natural choice for $\beta_k$ is the \emph{Jackson damping coefficients}
\begin{equation} \label{eq:jackson-damping}
    \beta_k = \frac{1}{K+2} \left[ (K+2-k) \cos \left(\frac{\pi k}{K+2}\right) + \sin\left(\frac{\pi k}{K+2}\right)\cot\left(\frac{\pi}{K+2}\right)\right],
\end{equation}
which ensures that $\hat{p}_K$ is a bona fide nonnegative density function (of total mass $n$).
The coefficients $\beta_k$ are derived from the \emph{Jackson--Rivlin kernel} used to study the best polynomial approximation for Lipschitz functions \cite{jackson1912approximation,rivlin1981introduction}.\footnote{See \citet[Foot.~2]{shah2026what} for an amusing historical footnote.}
We obtain the kernel polynomial method by using estimates $\hat{c}_k \approx c_k(\mat{A})$ for the Chebyshev moments \cref{eq:cheb-moment} in the expression \cref{eq:kpm-estimate}; see \cref{alg:kpm}.

\begin{algorithm}[t]
\caption{Kernel polynomial method for spectral density estimation} \label{alg:kpm}
\begin{algorithmic}[1]
\Require Symmetric matrix $\mat{A} \in \R^{n\times n}$ with spectrum in $[-1,1]$, number of probes $s$, even degree $K$
\Ensure Damped spectral density estimate $\hat{p}_K(y)\,\d y \approx \d\mu_{\mat{A}}(y)$
\State Estimate $(\hat{c}_1,\ldots,\hat{c}_K) \gets \Call{ChebMoments}{\mat{A},s,K}$ \Comment{\cref{alg:cheb-moments}}
\State Compute the Jackson damping coefficients $(\beta_k)_{k=1}^K$ using \cref{eq:jackson-damping}
\State \Return $\hat{p}_K(y) \gets \tfrac{n}{2}w(y) + \sum_{k=1}^K \beta_k\,\hat{c}_k\,\mathrm{T}_k(y)\,w(y)$
\end{algorithmic}
\end{algorithm}

The kernel polynomial method, as such, originated in computational physics \cite{skilling1989eigenvalues,benoit1992spectral,silver1994densities,wang1994calculating,silver1996kernel}, with intellectual roots going back decades earlier to work in numerical analysis on quadrature \cite{sack1971algorithm}.
The survey of \citet{weie2006kernel} continues to be a definitive reference.
Recently, the method and its variants have received significant attention by researchers in numerical analysis and theoretical computer science \cite{lin2016approximating,braverman2022sublinear,chen2023spectrum,chen2024lanczos,chen2025randomized}.
\citet{musco2026spectral} generalize the kernel polynomial method for spectral density estimation of normal matrices.

\subsection{Theoretical guarantees} \label{sec:sde-theory}

To understand how close the spectral measure $\mu_{\mat{A}}$ is to an approximation $\hat{\mu}$, we can integrate both measures against test functions and compare the difference:
\begin{equation*}
    \left| \int f \, \d\mu_{\mat{A}} - \int f \, \d\hat{\mu} \, \right| = \left| \, \tr(f(\mat{A})) - \hat{\tr}_f \,\right| \quad \text{for } \hat{\tr}_f \coloneqq \int f \, \d\hat{\mu}.
\end{equation*}
For spectral sums, this is already the quantity we wish to understand, as this error coincides with the error of the spectral sum approximation $\hat{\tr}_f \approx \tr(f(\mat{A}))$.
To compare measure to measure, we take the maximum difference over a class of functions $\set{F}$:
\begin{equation*}
    \mathrm{d}_{\set{F}}(\mu_{\mat{A}},\hat{\mu}) \coloneqq \sup_{f \in \set{F}} \left| \int f \, \d\mu_{\mat{A}} - \int f \, \d\hat{\mu} \, \right|.
\end{equation*}
Choosing the class $\set{F}$ to consist of all 1-Lipschitz functions gives rise to the Wasserstein-1 distance $\mathrm{W}_1(\mu_{\mat{A}},\hat{\mu})$, which measures the cost of transporting one measure onto the other, assuming a unit cost to move a unit mass over a unit distance.
We can also take another function class $\set{F}$.
For instance, \citet{chen2025randomized} consider the class $\set{F}$ of bounded functions that are analytic on an ellipse in the complex plane, which may be a more natural function class for numerical analysis and scientific computing.

Nevertheless, let us present some aspects of the Wasserstein-1 theory, as it is more complete.
The following structural bound, due to \citet[Thm.~1]{musco2025sharper}, provides a starting point for the analysis of many methods.

\begin{theorem}[Error estimates from Chebyshev moments] \label{thm:mmrs-moment-matching}
    Let $\hat{\mu}$ be a measure on $[-1,1]$ with total mass $\hat{\mu}([-1,1]) = n$, which serves as an approximation to the spectral measure $\mu_{\mat{A}}$ of a symmetric matrix $\mat{A} \in \R^{n\times n}$ with spectrum in $[-1,1]$.
    Let $\hat{c}_k \coloneqq \int \mathrm{T}_k \, \d\hat{\mu}$ denote the Chebyshev moments of $\hat{\mu}$.
    For any Lipschitz function $f:[-1,1]\to \mathbb R$ with Lipschitz constant $\Lip(f)$, define the spectral sum approximation $\hat{\tr}_f$.
    Then, for every $K\ge 1$,
    \begin{equation*}
        \left| \tr(f(\mat{A})) - \hat{\tr}_f \right| \le \left[\frac{36n}{K} + \left( \sum_{k=1}^K \frac{1}{k^2} [c_k(\mat{A}) - \hat{c}_k]^2  \right)^{1/2}\right] \cdot \Lip(f).
    \end{equation*}
    Consequently, 
    the Wasserstein 1-distance satisfies
    \begin{equation} \label{eq:wasserstein-1-bound}
       \mathrm{W}_1(\mu_{\mat{A}},\hat{\mu}) \le \frac{36n}{K} + \left( \sum_{k=1}^K \frac{1}{k^2} [c_k(\mat{A}) - \hat{c}_k]^2  \right)^{1/2}.
    \end{equation}
\end{theorem}

This theorem shows that control of the Chebyshev moments gives control of the accuracy of spectral sum approximations \(\hat{\tr}_f \approx \tr(f(\mat{A}))\) and of the Wasserstein-1 distance between the measure \(\mu_{\mat{A}}\) and its approximation \(\hat{\mu}\). The former result extends to signed approximations, but one must additionally account for the total mass of the negative part of \(\hat\mu\). The proof of \Cref{thm:mmrs-moment-matching} is pleasingly elementary, using Jackson approximation theorems and Chebyshev decay estimates for Lipschitz functions.
Recall that the maximum Wasserstein 1-distance between two measures on $[-1,1]$ with total mass $n$ is $2n$, so the factor $n$ in \cref{eq:wasserstein-1-bound} is sensible.
The constant 36 is not sharp; see \cite{shah2026what}.

As one example of how we can apply this bound, let us use it to analyze stochastic Lanczos quadrature.

\begin{corollary}[Stochastic Lanczos quadrature] \label{cor:slq}
For symmetric $\mat{A} \in \R^{n\times n}$ with spectrum in $[-1,1]$, consider the stochastic Lanczos quadrature method (\cref{alg:slq}) run for $k \ge 2$ Lanczos steps per probe with $s$ independent Rademacher or spherical probes.
    Then
    \begin{equation*}
        \expect \mathrm{W}_1(\mu_{\mat{A}},\hat{\mu}) \le \frac{18n}{k-1} + \sqrt{\frac{\pi^2n}{3s}}.
    \end{equation*}
    Thus, to achieve an expected Wasserstein error of $\varepsilon n$ requires at most $k = 38/\varepsilon$ Lanczos steps per probe with $s = \max \{ 1 , 14 / (\varepsilon^2n) \}$ random probes.
    In particular, to achieve accuracy level $\varepsilon \ge 1/\sqrt{n}$ requires a total number of $ks \le \order(1/\varepsilon)$ matrix--vector products.
\end{corollary}

\begin{proof}
    Let $\vec{\omega}_1,\ldots,\vec{\omega}_s$ denote the random probes used by the algorithm.
    By \cref{thm:lanczos-gaussian-quadrature}, the Chebyshev moments of the stochastic Lanczos quadrature approximation $\hat{\mu}$ are
    \begin{equation*}
        \hat{c}_\ell = \frac{1}{s} \sum_{i=1}^s \vec{\omega}_i^\tp \mathrm{T}_\ell(\mat{A})\vec{\omega}_i^{\vphantom{\tp}} \quad \text{for each } \ell < 2k.
    \end{equation*}
    Consequently, by \cref{thm:variancegh} and the fact that all eigenvalues of $\mat{A}$ are in $[-1,1]$, it holds that
    \begin{equation*}
        \expect[ (\hat{c}_\ell - c_\ell(\mat{A}))^2 ] \le \frac{2 \norm{\mathrm{T}_\ell(\mat{A})}^2_{\rm F}}{s} \le \frac{2n}{s}.
    \end{equation*}
    Applying \cref{thm:mmrs-moment-matching} with $K=2k-1$, we obtain
    \begin{equation*}
        \expect \mathrm{W}_1(\mu_{\mat{A}},\hat{\mu}) \le \frac{36n}{2k-1} + \left( \sum_{\ell=1}^{2k-1} \frac{1}{\ell^2} \expect [ (c_\ell(\mat{A}) - \hat{c}_\ell)^2]  \right)^{1/2} \le \frac{18n}{k-1} + \sqrt{\frac{\pi^2n}{3s}}.
    \end{equation*}
    This is the desired result.
\end{proof}

\Cref{cor:slq} has interesting consequences for the practitioner.
In particular, it shows that \emph{just one} random probe is usually sufficient to produce one or two digits of accuracy in the Wasserstein 1-distance.\footnote{In practice, it still may be worth using a small block of, say, five to ten probes to leverage hardware efficiencies for block matrix arithmetic.}
However, one should be cautious in taking this ``one probe'' guidance too far.
If the goal is to estimate spectral sums $\hat{\tr}_f \approx \tr(f(\mat{A}))$ or one uses a different error metric to evaluate the error, more stochastic probes may lead to meaningful improvements in accuracy.

For more on the analysis of spectral density estimation algorithms, see \cite{braverman2022sublinear,chen2023spectrum,chen2024lanczos,chen2025randomized,musco2025sharper,bhattacharjee2025improved}.

\subsection{Deflation and variance reduction} \label{sec:deflation-variance-reduction}

As with trace estimation, we can often improve the quality of spectral density estimation and spectral sum algorithms by using variance reduction or deflation techniques to reduce the influence of a few outlying eigenvalues.
\citet{bhattacharjee2025improved} propose versions of Chebyshev moment regression (\cref{alg:moment-matching}) that first use the block Lanczos algorithm to compute, and deflate, any large eigenvalues.
They also show that stochastic Lanczos quadrature has an \emph{implicit deflation effect}, where the algorithm gracefully handles extremal eigenvalues without additional modification, in analogy to the discussion on Page~\pageref{page:implicitprecond}.

In the context of smoothed spectral density estimation, \citet{lin2017randomized}
applies what is essentially Nystr\"om\plusplus\ to a family of localized matrix functions
$g_\sigma(t\Id-\mat A)$, using a common random sketch as $t$ sweeps across the spectrum.
\citet{matti2026stochastic} extend and analyze this idea in the more general context of
parameter-dependent trace estimation.

For spectral sum problems, the Krylov-aware stochastic trace estimation algorithm of \citet{chen2023krylov} provides a Hutch\plusplus-style variance reduced trace estimator that can be used to rapidly compute $\tr(f(\mat{A}))$ or, more generally, $\tr(f_i(\mat{A}))$ for a family of functions $f_i$.
The algorithm uses block Lanczos~\cite{golub1977block} both to construct a low-rank approximation
to $f(\mat{A})$ for variance reduction and to estimate its action on random probe vectors.
Based on the ideas discussed in \cref{sec:lowrank-functions},
an alternative variance reduction strategy uses a low-rank approximation of $\mat{A}$ itself \cite{persson2023randomized,chen2026supercharging}. 
The recently proposed FlexTrace estimator~\cite{madhavan2026flextrace} combines this approach with
the exchangeability principle behind XNysTrace (\Cref{sec:xtrace}).

\section{Index selection} \label{sec:selection}

To this point, a large part of this survey has focused on randomized methods that solve linear algebra problems that ``mix up'' an input matrix $\mat{A}$ by multiplying it against a test matrix $\mat{\Omega}$ or random vectors $\{\vec{\omega}_i\}$.
Even when we used a sparse or structured test matrix, the goal was always to design the probes so as to maximally scramble the rows or columns of $\mat{A}$.

This and the next section will focus on a different family of methods that seek to distill information about a matrix $\mat{A}$ by carefully selecting a representative subset of rows or columns.
Such \emph{index selection} methods can be vital in applications.

\subsection{Model problems: Active least squares and interpolative decomposition} \label{sec:activelsq}

To motivate index selection, let us consider two model problems in least-squares regression and low-rank approximation.

\begin{paragraph}{Active least-squares regression.}
First, suppose we are interested in solving a least-squares regression problem
\begin{equation} \label{eq:active-ls}
    \min_{\mat{X} \in \R^{n\times p}} \norm{\mat{A}\mat{X} - \mat{B}}_{\rm F}^2 \quad \text{for } \mat{A} \in \R^{m\times n} \text{ and } \mat{B} \in \R^{m\times p}
\end{equation}
with the additional stipulation that we do not wish to read the entire matrix $\mat{B}$.
For instance, we could be in a setting where the number of right-hand sides $p$ vastly exceeds the number of columns $n$ in $\mat{A}$, making computing with the entire $\mat{B}$ matrix very expensive.
Or, the entries of the matrix $\mat{B}$ may not be immediately accessible to us, requiring work to generate or collect.
In the most severe case, the entries of $\mat{B}$ could be physically observable quantities, and we must run an expensive scientific experiment or computer simulation to generate a row of $\mat{B}$. 
This situation mimics the basic setup of \emph{active machine learning}, so we may call this setup an \emph{active least-squares} problem~\cite{chen2019active}.

Index selection methods offer an immediate way to approximately solve such an active least-squares problem.
First, we inspect the matrix $\mat{A}$ to identify a representative index set $\iset \subseteq \{1,\ldots,m\}$ of $s := |\iset| \ll m$ elements.
Then, we solve a subsampled least-squares problem
\begin{equation} \label{eq:active-ls-subsampled}
    \Xtilde \coloneqq \argmin_{\mat{Z} \in \R^{n\times p}} \norm{\mat{A}(\iset,:)\mat{Z} - \mat{B}(\iset,:)}_{\rm F}^2 \quad \text{or} \quad \Xtilde \coloneqq \argmin_{\mat{Z} \in \R^{n\times p}} \norm{\mat{D}\mat{A}(\iset,:)\mat{Z} - \mat{D}\mat{B}(\iset,:)}_{\rm F}^2.
\end{equation}
For reasons that will be discussed, a diagonal weight matrix $\mat{D} \in \R^{s\times s}$ is often used.
Observe that this approach to active least-squares is an instance of the sketch-and-solve paradigm (\Cref{sec:sas}) with the choice of test matrix $\mat{\Omega} = \Id(:,\iset)$ or $\mat{\Omega} = \Id(:,\iset)\mat{D}$.
The core feature of this approach is that it requires access only to the rows of \(\mat{B}\) indexed by \(\iset\), which can lead to substantial savings when accessing \(\mat{B}\) is expensive. If \(\mat{B}\) is already fully available, index selection may still be useful in some settings, but this principal advantage disappears, and oblivious sketching is often more attractive. 
\end{paragraph}

\begin{paragraph}{CSSP and interpolative decompositions.}
A second model problem concerns low-rank approximation.
Elementary linear algebra tells us that a matrix $\mat{A}$ has rank $s$ if and only if it has a spanning set $\mat{A}(:,i_1),\ldots,\mat{A}(:,i_s)$ of $s$ linearly independent columns.
Algebraically, this fact gives rise to a matrix factorization
\[
  \begin{tikzpicture}[node distance=1mm,baseline=-1mm]       \node[draw, rectangle, solid, fill=gray!10, minimum width=\colsize, minimum height=\colsize] (A)  at (0,0) {\phantom{$\Scale[2]{A}$}};
  \foreach \x in {0.13,0.6,0.85}{
    \draw[fill=red!20,draw=none]
      ($(A.south west)+(\x*\colsize,0)$)
      rectangle
      ($(A.south west)+(\x*\colsize+\colsize/10,\colsize)$);      
  }
\node[draw, rectangle, solid, minimum width=\colsize, minimum height=\colsize] (A)  at (0,0) {$\Scale[2]{A}$}; 
   \end{tikzpicture}=
\begin{tikzpicture}[node distance=1mm,baseline=-1mm]
    \node[
        draw, rectangle, solid, fill=red!20,
        minimum width=\rowsize,
        minimum height=\colsize
    ] (C) {$\Scale[1.2]{C}$};

    \node[
        draw, rectangle, solid, fill=blue!20,
        minimum width=\colsize,
        minimum height=\rowsize,
        right=of C
    ] (W) {$\Scale[1]{\mat{W}^{\tp}}$};
\end{tikzpicture}
\]
with $\mat{C}=\mat{A}(:,\iset)$, $\iset \coloneqq \{i_1,\ldots,i_s\}$, containing the selected columns, and $\mat{W}^\tp = \mat{C}^\pinv \mat{A}$.

For a matrix that is only approximately rank \(s\), we seek an index set
\(\set{I}\) for which the selected columns
provide an accurate approximation to
\(\mat{A}\), in the sense that
\begin{equation} \label{eq:id}
 \mat{A} \approx \widehat{\mat{A}}:=\mat{A}(:,\iset)\mat{W}^\tp, \qquad \mat{W}^\tp\coloneqq \mat{A}(:,\iset)^\pinv \mat{A}.
\end{equation}
This is the \emph{column subset selection problem (CSSP)}.
The factorization~\cref{eq:id} is often called an \emph{interpolative decomposition}~\cite{cheng2005compression} because the columns in $\iset$ are ``interpolated'' exactly\footnote{The name is somewhat of a misnomer, because 
interpolation in the selected columns is a property shared by other low-rank approximations, most notably the CUR cross approximation discussed in~\cref{sec:cur} that interpolates both the selected columns \emph{and} rows. Also, it is almost always an  \emph{approximation} rather than a \emph{decomposition}. 
}:
\[
  \mat{A}(:,\iset) - \widehat{\mat{A}}(:,\iset) = \mat{A}(:,\iset) - \mat{A}(:,\iset) \mat{A}(:,\iset)^\pinv \mat{A}(:,\iset) = \mat{0}.
\]
For general $\mat{W}$, approximations of the form $\mat{A} \approx \mat{A}(:,\iset)\mat{W}^\tp$ are sometimes referred to as CX decompositions as the letters ``C'' and ``X'' are used for the factors. 

Interpolative low-rank approximations are useful in many contexts,
including fast direct solvers for PDEs \cite{martinsson2019fast} and
tensor decompositions~\cite{oseledets2010tt}.
In the language of \Cref{sec:lowrank-projection}, the interpolative decomposition~\cref{eq:id} is a projection approximation with the coordinate sampling test matrix \(\mat{\Omega} = \Id(:,\iset)\). Combining index selection with the \Nystrom and generalized \Nystrom methods leads to related approximation formats, discussed in
\cref{sec:column-nystrom,sec:cur}.
\end{paragraph}

\begin{remark}[Index selection beyond linear algebra]
Given the focus of this survey, we are largely focusing on index selection as a computational  tool to solve linear algebra problems.
In many applications, however,
identifying the subset of indices 
is itself the objective.
Given a collection of data elements $\{\vec{a}_i\}$, 
one may wish to select a small subset $\{\vec{a}_i : i \in \iset \}$  that ``represents'', ``summarizes'', or ``explains'' the whole.
This viewpoint appears in, among others, feature and variable 
selection, active learning and experimental design, recommender systems, and the selection of informative genetic markers;
see, for example,~\citet{boutsidis2008unsupervised,mahoney2009cur,kulesza2012determinantal}. Other applications include graph sparsification \cite{spielman2011graph} and approximation of functions \cite{adcock2025optimal}; see also \Cref{sec:leverage}.
\end{remark}

\subsection{Classical pivoting methods} \label{sec:classical-pivoting}

A classical strategy for index selection is by means of a pivoted matrix factorization.
Suppose we are interested in selecting $s$ columns from a matrix $\mat{A} \in \R^{m\times n}$.
To do so, we may perform a pivoted QR decomposition, selecting the column of largest norm as the next pivot, orthogonalizing the remaining columns against it, and continuing.
After initializing $\iset \gets \emptyset$, the
pivoting process can be described conceptually by the iteration
\begin{equation} \label{eq:cpqr}
    i\gets \argmax \norm{\mat{A}(:,i)}^2, \quad \iset \gets \iset \cup \{i\}, \quad \mat{A} \gets \left(\Id - \frac{\mat{A}(:,i)\mat{A}(:,i)^*}{\norm{\mat{A}(:,i)}^2}\right)\mat{A}.
\end{equation}
This greedy strategy is the classical Businger--Golub QR factorization
with column pivoting (CPQR)~\citep{businger1965linear}.

CPQR usually works well in practice, typically producing a rank-$s$
interpolative approximation~\cref{eq:id} whose error is comparable to that
of the best rank-$s$ approximation, for example in the Frobenius norm.
Nevertheless, CPQR has two notable deficiencies.
First, its worst-case behavior can be poor: a classical example due to
\citet{kahan1966numerical} shows that the approximation error produced by CPQR can be
exponentially larger in $s$ than the best rank-$s$ approximation error; see
also~\cite{gu1996efficient}.
Second, even for typical, non-pathological cases,
CPQR can devote too many pivots to isolated columns with unusually large (residual) norms, rather than selecting columns representative of the bulk of the matrix. This behavior is especially relevant when the singular values decay slowly and no small set of dominant directions clearly separates from the rest. The analogous phenomenon for greedy pivoted Cholesky has been explained by \citet{chen2025randomly}; the Gram correspondence (see also~\Cref{sec:gramcorrespondeindex}) transfers the same mechanism to CPQR.

The poor worst-case behavior of CPQR can be addressed by 
iteratively swapping selected and unselected columns, until certain rank-revealing properties are 
controlled. This strategy was proposed by~\citet{gu1996efficient} and it is known as \emph{strong rank-revealing QR} (strong RRQR).
In particular, for a prescribed parameter $f>1$, strong RRQR achieves
a worst-case approximation error bound with suboptimality factor
$\sqrt{1+f^2s(n-s)}$, in both the spectral and Frobenius norms.
This is a much stronger worst-case guarantee than for CPQR,
but it comes with non-negligible computational overhead.
The Gu--Eisenstat strategy was recently re-evaluated by \citet{damle2026estimating}.

An alternative greedy strategy selects the column that gives the largest decrease in the Frobenius-norm projection error, rather than the column of largest norm. Under the Gram correspondence, this becomes nuclear-score selection for \Nystrom approximation; we discuss this criterion and its randomized implementation in~\Cref{sec:nuclear-selection}.

\subsection{Combining sketching and pivoting} \label{sec:sketchy-pivoting}

A particularly effective way to accelerate deterministic pivoting is to perform
the pivot decisions on a randomized sketch.  For column selection, one forms
\[
    \Ahat=\mat{\Omega}^{\tp}\mat{A}\in\R^{\ell\times n}
    \quad \text{with } \ell\ll m,
\]
and applies deterministic pivoting, such as CPQR, to $\Ahat$; the
resulting pivot indices $\iset$ are then used to select the columns
$\mat{A}(:,\iset)$ of the original matrix.  More practical variants of this basic idea select
only blocks of pivots and update the sketch as the factorization proceeds,
leading to randomized CPQR and HQRRP
\cite{voronin2017efficient,duersch2017randomized,martinsson2017householder}.
This addresses a major computational bottleneck of CPQR: its
sequential pivot decisions inhibit blocking and incur substantial data movement
and synchronization. Performing these decisions on the much smaller sketch
reduces communication and enables efficient blocked implementations.
We refer to
\citet[\S\S13.4 and 16.6]{martinsson2020randomized} for a detailed exposition.
An analysis by \citet{dong2023simpler} separates the approximation error into the error of capturing the row space of \(\mat A\) with the sketch and an a posteriori computable factor measuring the quality of the selected pivots.
Recent developments replace CPQR
on the sketch by the cheaper LUPP \cite{dong2023simpler,pearce2025adaptive}, or
by strong RRQR, for which strong rank-revealing guarantees can be
obtained when the sketch is a subspace embedding for the range of $\mat{A}$~\cite{grigori2025randomized}; see also \cite{fakih2025efficient}. \citet{melnichenko2025choleskyqr} combine sketch-based CPQR with randomized CholeskyQR (\cref{sec:randcholqr}).

\subsection{Randomized pivoting} \label{sec:random-pivoting}

In contrast to the sketch-and-pivot methods above, which retain a deterministic pivot rule, randomized pivoting injects randomness directly into the pivot selection. 
Instead of selecting the column of largest norm as pivot as in CPQR, we draw 
the pivot $i_t$ in step $t$ from the squared-norm distribution
\begin{subequations} \label{eq:rpqr}
 \begin{equation} \label{eq:rpqr-a}
    \prob\{i_t=j\,|\, \mat{A}\}
    =
    \frac{\norm{\mat{A}(:,j)}^2}{\norm{\mat{A}}_{\rm F}^2},
 \end{equation}
and update
\begin{equation}
     \mat{A}
    \gets
    \left(
        \Id-
        \frac{\mat{A}(:,i_t)\mat{A}(:,i_t)^*}
             {\norm{\mat{A}(:,i_t)}^2}
    \right)\mat{A}.
\end{equation}
\end{subequations}
After $s$ steps, the selected index set is $\iset=\{i_1,\ldots,i_s\}$.
This algorithm is known as randomly pivoted QR (RPQR) or adaptive sampling.

Selecting a pivot randomly in this way hedges against pathological outcomes and produces pivot selections that are less influenced by ``outlier'' columns.
\citet{frieze2004fast} already suggested sampling columns with probabilities
proportional to their squared norms in the \emph{original} matrix.
\emph{Adaptive} sampling, that is, recomputing the probabilities from the residual as in~\cref{eq:rpqr-a},
was first studied by~\citet{deshpande2006matrix,deshpande2006adaptive} and revisited by \citet{chen2025randomly,dong2025robust}.
The iteration \cref{eq:rpqr} is only a conceptual description.
Like CPQR, it can and should be implemented as a Householder-based partial QR factorization as in~\Cref{alg:rpqr}; see \cite[p.~149]{epperly2025make} for more details. In particular, the Householder transformations ensure that  $\mat{A}(t:m,j)$ represents the component of column $j$ orthogonal to the span of the $t$ 
selected columns. Hence,
$
    \norm{\mat{A}(t:m,j)}^2
$
matches the squared residual norm in~\cref{eq:rpqr-a}.

\begin{algorithm}[t]
\caption{Randomly pivoted QR for index selection (RPQR)} \label{alg:rpqr}
\begin{algorithmic}[1]
\Require Matrix $\mat{A} \in \R^{m\times n}$ and number of pivots
    $1\le s\le \operatorname{rank}(\mat{A})$
\Ensure Column index set $\iset=\{i_1,\ldots,i_s\}$
\State Initialize $\iset\gets\emptyset$
\For{$t=1,\ldots,s$}
    \State Set $p_j\gets\norm{\mat{A}(t:m,j)}^2$
        for $j\notin\iset$
    \State Draw $i_t\notin\iset$ with probability
        $p_{i_t}/\sum_{j\notin\iset}p_j$
    \State \label{line:householder} Choose Householder reflector $\mat{H}$ such that
        $\mat{H}\mat{A}(t:m,i_t)=\alpha\vec{e}_1$
    \State \label{line:householderapply} Update
        $\mat{A}(t:m,:)\gets\mat{H}\mat{A}(t:m,:)$
    \State $\iset\gets\iset\cup\{i_t\}$
\EndFor
\end{algorithmic}
\end{algorithm}

A key advantage of RPQR over CPQR is that it admits stronger guarantees. In particular, the Gram correspondence transfers results for randomly pivoted Cholesky~\cite{chen2025randomly,epperly2026new} to~\cref{eq:rpqr}. Letting $\iset$ denote
the index set obtained after $s$ steps of RPQR,
\Cref{thm:rpcholesky} shows for a target rank $r$ that the columns $\mat{A}(:,\iset)$ selected from the original matrix $\mat{A}$ satisfy
\[
    \expect\norm{\mat{A}-\mat{A}(:,\iset)\mat{A}(:,\iset)^\pinv\mat{A}}_{\rm F}^2
    \le (1+\varepsilon)
       \norm{\mat{A}-\lowrank{\mat{A}}_r}_{\rm F}^2,
\]
provided that
$s =  r/\varepsilon + \order\!\big(r \sqrt{\log r} + r\log(1/\varepsilon)\big)$.
Thus, the exponential worst-case behavior of CPQR is avoided, at the price of a probabilistic guarantee and modest
oversampling.

The basic RPQR method~\cref{eq:rpqr} is inherently sequential: each pivot changes the residual and hence the probability distribution for the next pivot. This limits the extent to which the method can exploit parallelism on modern hardware, for example through efficient matrix-matrix operations. To address this, one can instead draw blocks of $b$ iid pivots:
\begin{subequations} \label{eq:brpqr}
\begin{equation}
    \iset' = \{i_1,\ldots,i_b\} \quad \text{with } \prob \{i_\ell = j\} = \frac{\norm{\mat{A}(:,j)}^2}{\norm{\mat{A}}_{\rm F}^2}
\end{equation}
and perform the block update
\begin{equation}
    \iset \gets \iset \cup \iset', \quad \mat{A} \gets \left( \Id - \mat{A}(:,\iset')\mat{A}(:,\iset')^\pinv \right) \mat{A}.
\end{equation}
\end{subequations}
This algorithm is known as block randomly pivoted QR \cite{deshpande2006matrix,deshpande2006adaptive,chen2025randomly,epperly2025embrace}.
The basic approach~\cref{eq:brpqr} can produce index sets of poor quality because all 
indices \emph{within} a block $\iset'$ are sampled from the \emph{same} distribution.
Consequently, the 
set $\iset'$ may contain many redundant columns that contain essentially the same information.
Two recent approaches address this difficulty.  Robust blockwise random
pivoting (RBRP) uses CPQR within each block to filter redundant
information~\cite{dong2025robust}.
Experimentally, the accuracy of RBRP is close to sequential random pivoting while maintaining 
good hardware efficiency. Accelerated randomly
pivoted QR uses rejection sampling instead to retain the pivoting behavior of
the sequential method while carrying out most of the work in blocks~\cite{epperly2025embrace}.

Random pivoting algorithms---both unblocked and blocked---require care to implement in a numerically stable way.
Stable implementations maintain an orthonormal basis using reorthogonalized Gram--Schmidt or Householder transformations; see
\cite[\S3.4 and App.~B]{chen2025randomly}, \cite[\S5]{epperly2025embrace}, and \cite[\S9.3]{epperly2025make}.
The stability of block algorithms depends both on the
orthogonalization against previously selected blocks and on the local
(intrablock) QR factorization; see~\cite{carson2022block} for a systematic
discussion of these choices and the associated stability--efficiency
tradeoffs.

\subsection{Leverage score sampling} \label{sec:leverage}

So far, our discussion has emphasized column selection, motivated by the connection to CPQR.
For active least squares, however, the natural task is to select rows. Of course, row selection for $\mat{A}$ is
equivalent to column selection for $\mat{A}^\tp$.

While the pivoting-based strategies discussed above are very effective at selecting a small subset of representative
rows, it is less obvious how to apply them directly to active least squares. 
Indeed, when selecting rows of a tall-thin matrix $\mat{A} \in \R^{m\times n}$, these methods necessarily terminate after identifying $|\iset| = n$ rows. To obtain a near-optimal least-squares residual, however, we often need to oversample.
For that purpose, sampling-based strategies can be effective.

Let us first discuss a few basic approaches that do \emph{not} work well in general.
The most basic sampling approach is to simply select rows uniformly at random~\cite{drineas2006fasta}.
This approach naturally will fail when the rows of $\mat{A}$ differ wildly in scale and ``importance''; this relates to the notion of coherence discussed in~\Cref{sec:sparse-embeddings}; see also~\cite[\S9.6]{martinsson2020randomized}.
To remedy this deficiency, early RNLA work of~\citet{frieze2004fast} explored index selection by (non-adaptive) \emph{squared-norm sampling}, drawing each index $i_t$ iid with distribution $\prob \{i_t = j\} = \norm{\mat{A}(j,:)}^2 / \norm{\mat{A}}_{\rm F}^2$.
This strategy often works better than uniform sampling, but it too can fail when, for example, rows of large norm are nearly linearly dependent.

There is a conceptual reason why norm is a bad measure for the importance of a row.
Indeed, the least-squares problem $\min_{\mat{X}} \norm{\mat{A}\mat{X} - \mat{B}}_{\rm F}^2$ can be reparametrized by any change of variables $\mat{X}' \gets \mat{R}\mat{X}$ and $\mat{A}' \gets \mat{A}\mat{R}^{-1}$.
The optimal least-squares residual is invariant under such changes, but the row norms are not.
This suggests measuring row importance only after removing the arbitrary parametrization of the column space, for instance by 
replacing $\mat{A}$ with an orthonormal basis for $\range(\mat{A})$.

\begin{definition}[Leverage scores] \label{def:leverage}
    Let $\mat{A} \in \R^{m\times n}$ be a tall-thin matrix of full column rank $n$. Let $\mat{Q}$ be any orthonormal basis of $\range(\mat{A})$.
    Then the squared row norms of $\mat{Q}$ are called the \emph{leverage scores} of $\mat{A}$:
    \begin{equation*}
        \ell_i(\mat{A}) \coloneqq \norm{\mat{Q}(i,:)}^2, \qquad i = 1,\ldots, m.
    \end{equation*}
\end{definition}

Equivalently, the leverage scores are the diagonal entries of the orthogonal
projector $\mat{A}\mat{A}^\pinv$, known in classical regression diagnostics as
the \emph{hat matrix}~\cite{hoaglin1978hat}.
It follows that the leverage scores are well-defined and admit the variational characterization
\begin{equation*}
    \ell_i(\mat{A}) = \max \big\{\vec{v}(i)^2/\norm{\vec{v}}^2 \colon \vec{v} \in \range(\mat{A}), \vec{v} \not=0 \}.
\end{equation*}
In particular, we see that the leverage scores depend only on the subspace $\set{V} = \range(\mat{A})$, so it is sensible to talk about the leverage scores \emph{for a subspace}.

\begin{remark}[Christoffel functions and leverage scores] \label{remark:Christoffel}
Leverage scores have a direct analogue in function approximation.
Given an orthonormal basis of functions $\phi_1,\ldots,\phi_n$ in an $L^2$ space,
the \emph{reciprocal Christoffel} function of the subspace $\set{P}$ spanned by this basis is defined as
\[
    K_{\set{P}}(x)=\sum_{j=1}^n |\phi_j(x)|^2.
\]
As emphasized by \citet{adcock2025optimal}, this is the continuous analogue of
the leverage scores in \Cref{def:leverage}. Sampling proportionally to
$K_{\set{P}}$ and solving the corresponding weighted least-squares problem leads to a stable near-best approximation from
$\order(n\log n)$ point evaluations~\cite{cohen2017optimal}. 
\end{remark}

\subsubsection{Active least squares by leverage score sampling}\label{sec:leverage-ls}

We now apply leverage score sampling to the active least squares problem from~\Cref{sec:activelsq}.
For this purpose, we first choose a row set $\iset = \{i_1,\ldots,i_s\}$ by sampling the indices iid from the leverage score distribution $\prob \{i_t = j\} = \ell_j(\mat{A}) / n$.
Observe that we are sampling \emph{with replacement}; indices can be selected multiple times and, hence, $\iset$ may be a multiset (a subtlety which we do not focus on for the sake of simplicity).
As we saw in \cref{sec:isotropy}, it is desirable for sketching techniques---index selection methods included---to be isotropic, preserving the size of a matrix on average.
To promote this property, we introduce the diagonal weight matrix
\begin{equation} \label{eq:leverage-score-weights}
    \mat{D} \coloneqq \sqrt{\tfrac{n}{s}}\diag\!\Big( \tfrac{1}{\sqrt{\ell_{i_1}(\mat{A})}},\ldots,\tfrac{1}{\sqrt{\ell_{i_s}(\mat{A})}}\Big).
\end{equation}
Observe that, under this approach, rows with lower leverage scores are sampled less frequently but, to compensate, are weighted more highly when they are sampled.
From here, the leverage score sketch-and-solve solution is 
\begin{equation} \label{eq:ls-least-squares}
    \Xtilde \coloneqq \argmin_{\mat{Z} \in \R^{n\times p}} \norm{\mat{D}\mat{A}(\iset,:)\mat{Z} - \mat{D}\mat{B}(\iset,:)}_{\rm F}^2.
\end{equation}

The following result is standard; see \cite[Thm.~16]{larsen2022practical} and \cite{meyer2026active}.
\begin{theorem}[Active least squares: Leverage scores] \label{thm:leverage}
    Assume $\mat{A} \in \R^{m\times n}$ has full column rank, and consider the active least-squares problem \cref{eq:active-ls}.
    Then the leverage score approximate least-squares solution \cref{eq:ls-least-squares} satisfies
    \begin{equation} \label{eq:lev-guarantee}
        \Median(\norm{\mat{A}\Xtilde - \mat{B}}_{\rm F}^2) \le (1+\varepsilon) \min_{\mat{X} \in \R^{n\times p}} \norm{\mat{A}\mat{X} - \mat{B}}_{\rm F}^2
    \end{equation}
    with $s = \order(n \log n + n /\varepsilon)$.
\end{theorem}

\Cref{thm:leverage} suggests that, unlike some of the oblivious sketching methods discussed in~\Cref{sec:embeddings}, leverage score sampling requires $s \propto n \log n$ to solve a least-squares problem to constant-factor accuracy.
This is a genuine phenomenon, not an artifact of the analysis. Assuming that $n$ divides $m$, this may be seen by inspecting the example $\mat{A} = \operatorname{blockdiag}(\onevec_{m/n},\ldots,\onevec_{m/n})$ and $\vec{b} = \onevec_m$ where $\onevec_q$ denotes a vector of ones of length $q$.
For the sampled problem to be consistent, one must sample at least one row from every block. Achieving this
with non-vanishing probability requires $\Theta(n\log n)$ samples,
an instance of the famous ``coupon collector'' phenomenon~\cite[\S3.6]{motwani1995randomized}.

\begin{remark}[Optimal-complexity active least squares] \label{rem:optimal-active-ls}
    The $\order(n\log n)$ complexity is a property of iid leverage score sampling and not inherent to active linear regression.
    Indeed, \citet{chen2019active} derived an algorithm that achieves the guarantee \cref{eq:lev-guarantee} with an index set of size $|\iset| = \order(n/\varepsilon)$, which is the best-possible rate.
    This algorithm is adaptive in the sense that later choices depend on previously selected indices,
    just as in random pivoting and adaptive random pivoting; see \Cref{sec:random-pivoting,sec:index-rangefinding}.
    However, the Chen--Price algorithm is quite complicated, and we are unaware of efforts to implement the algorithm and compare it to alternatives. Recent work has explored variants of leverage-score sampling that introduce
dependence between the sampled indices to improve coverage and reduce redundant
samples~\cite{shimizu2024improved}.
\end{remark}

\begin{remark}[Role of randomization] 
Randomization can protect against unfavorable alignment between
the selected indices and an unknown right-hand side.
Consider a \emph{deterministic} sketch-and-solve method $$\xtilde \coloneqq \argmin_{\vec{x}}\norm{\mat{\Omega}^\top(\mat{A}\vec{x}-\vec{b})}$$ in the active least-squares setting, where the sketching matrix $\mat{\Omega}$ depends on $\mat{A}$ but not on $\vec{b}$.
Since the algorithm is deterministic, we might be unfortunate and encounter a vector $\mat{b}$ that is particularly bad for the sketching matrix $\mat{\Omega}$.

To see this, take
$\vec{b}=\mat{A}\vec{x}_\star +\vec{e}$, where $\vec{x}_\star$ is the desired solution, and $\vec e$ is orthogonal to $\range(\mat A)$. Letting $\mat Q=\Orth(\mat A)$,
the computed solution satisfies
$\|\mat A\xtilde-\vec b\| = \|(\Id-\mat P)\vec e\|$, where 
$\mat P = \mat Q (\mat\Omega^\tp\mat Q)^\dagger\mat\Omega^\tp$
is the oblique projector from the proof of \cref{thm:sketch-and-solve-ose-osi}.
For the worst-case choice of $\vec e$, the ratio between
the computed and optimal residual norms thus equals 
$\|\Id - \mat P\| = \|\mat P\| = \|(\mat\Omega^\tp\mat Q)^\dagger\mat\Omega^\tp\|$.
For a deterministic selection of $s$ rows,
this ratio can grow at least like $\sqrt{m/s}$,
even with optimally chosen indices and weights.
(This is already the case for $\mat A=\vec 1 \in \R^{m\times 1}$.)
In the context of function approximation by least squares, the dimension $m$ is infinite, and this ratio can be arbitrarily large.
Randomization, by contrast, can avoid such unfavorable alignment with high probability, 
thereby allowing us to derive bounds that do not grow with the dimension, 
as described above. 
\end{remark}

\subsubsection{Computing or approximating leverage scores} \label{sec:leverage-estimation}

In many settings, we are lucky, and the leverage scores are available to us at little cost.
For instance, the leverage scores are known in closed form for the design matrix in Shapley value estimation \cite{musco2025provably}.
In the continuous setting (see~\Cref{remark:Christoffel}), the Christoffel
function is known explicitly for polynomial approximation on a number of
standard domains~\cite{shustin2022semi,adcock2025optimal}.
In addition, efficient samplers from the leverage-score distribution are available for 
certain tensor-structured matrices \cite{bharadwaj2023fast,bharadwaj2024efficient}.
Finally, in some contexts \cite{peherstorfer2020stability}, the matrix $\mat{A}$ of interest is already orthonormal. In such cases, we can obtain the leverage scores cheaply by computing the row norms.

Outside of these favorable settings, we must put in the work to compute or approximate the leverage scores ourselves.
To compute the leverage scores exactly, we may orthogonalize the matrix $\mat{A}$ using a method of our choice, requiring usually
$\order(mn^2)$ operations.

To approximate leverage scores, we combine whitening (\cref{sec:whitening})
with the JL row-norm estimator from~\Cref{sec:row-norm-estimation}, following \citet{drineas2012fast}.
Begin by whitening: Draw an oblivious sketching matrix $\mat{\Omega} \in \R^{m\times s}$, compute the sketch $\Ahat \coloneqq \mat{\Omega}^\tp \mat{A}$ as well as a thin QR factorization $\Ahat = \mat{Q}_{\mat{\Omega}}\mat{R}_{\mat{\Omega}}$.
When $\mat\Omega$ is a subspace embedding with bounded distortion for $\range(\mat A)$,
the whitened matrix $\mat A\mat R_{\mat\Omega}^{-1}$
is nearly orthonormal, so it is reasonable to suspect its row norms are nearly equal to the leverage scores of $\mat{A}$.
But computing the whitened matrix outright would require $\order(mn^2)$ operations, leading to at most a constant-factor speedup over simply computing the leverage scores directly.
So we instead estimate its row norms using the JL estimator.

Specifically, use a Rademacher embedding  $\mat{\Phi} \in \R^{n \times q}$ with entries $\pm 1 / \sqrt{q}$, and form $\mat{Z} \coloneqq \mat{A}(\mat{R}_{\mat{\Omega}}^{-1}\mat{\Phi})$. The leverage score estimates are then
\begin{equation*}
    \hat{\ell}_i \coloneqq \norm{\mat{Z}(i,:)}^2 \quad \text{for } i =1,2,\ldots,m.
\end{equation*}
The JL guarantee in \cref{thm:jl} also holds for Rademacher embeddings
\cite[Thm.~1.1]{achlioptas2003database}. Consequently, 
$q=\order(\log m)$ suffices to estimate, with high probability, all leverage scores
simultaneously to constant-factor accuracy.

Using a fast embedding with small output dimension $s = \order(n)$, the runtime is $\order(mn \log m + n^3)$ operations.
Fortunately, in practice, a crude estimate of the leverage scores is sufficient, so one can get away with small values for the parameters $s$ and $q$.
In our experiments, we implement leverage-score estimation using a SparseStack embedding (parameters $s=2n$ and $\zeta=4$) with $q \coloneqq \lceil \log m \rceil$.
For sparse matrices, leverage scores can be approximated even faster, in
input-sparsity time up to lower-order polynomial terms in $n$, using sparse
random embeddings~\cite{clarkson2017low}.

\subsubsection{When to use leverage score sampling} \label{sec:leverage-when-to-use}

Leverage score sampling is a mathematically elegant approach to index selection.
Once the leverage scores are in hand, leverage score sampling is trivially easy to perform, merely requiring iid sampling from a given list of probability weights.
Nonetheless, iid leverage score sampling also has limitations: It can require $\Theta(n\log n)$ samples to reliably produce a good approximate solution to an $m\times n$ least-squares problem and, except in special cases, it requires an upfront preprocessing step to compute or estimate the leverage scores (though crude estimates are fine).
Given these benefits and drawbacks, it is our view that leverage score sampling is a useful method, but only in certain cases.
Let us summarize our view by considering three possible use cases.

\begin{itemize}
    \item \textbf{Good use case: Active least-squares problem where the leverage scores are easy to compute.}
    This is the best case for leverage score sampling.
    For an active least-squares problem, we want to solve the problem by subsampling and index selection, and the leverage scores offer a natural choice.
    This is especially true if the leverage scores are easy to compute in a given problem.
    \item \textbf{Okay use case: Active least-squares problem, but the leverage scores need to be computed.}
    To solve an active least-squares problem where the leverage scores are not immediately accessible, plain leverage score sampling can still be a good approach.
    But it is worth testing alternatives.
    If one is already going to have to expend effort to compute or estimate the leverage scores, more sophisticated index selection methods such as (leveraged) volume-squared sampling (\cref{sec:volume-sampling}) can lead to better accuracy at similar cost.
    \item \textbf{Bad use case: General least-squares problem.}
    Leverage score sampling is usually not competitive for approximately solving a least-squares problem
    when the coefficient matrix $\mat{A}$ is an unstructured matrix and
    we have no issues reading the entire right-hand side $\mat{B}$.
    In this case, computing the leverage scores exactly requires a QR factorization of $\mat{A}$, while approximating the leverage scores requires applying an oblivious embedding to $\mat{A}$.
    But such a factorization or sketch already gives us enough information to compute the exact least-squares solution or a sketch-and-solve solution; there is typically no reason to add additional complication by performing leverage score sampling.
\end{itemize}

\subsection{Volume-squared sampling} \label{sec:volume-sampling}

The notion of volume has played an important role in the theory of index selection for several decades~\cite{gu1996efficient,goreinov2001maximal}.
For a tall matrix $\mat{B} \in \R^{m\times s}$, its volume is defined as the volume of the parallelepiped spanned by the columns of $\mat{B}$.
Algebraically, it is the product of the singular values
\begin{equation} \label{eq:volume}
    \vol(\mat{B}) \coloneqq \prod_{i=1}^{\min(m,s)} \sigma_i(\mat{B}) = \sqrt{\det\!\big( \mat{B}^\tp \mat{B} \big)}.
\end{equation}

The goal of column index selection is to identify an index set
that captures important directions while avoiding redundant columns.
Geometrically, it is clear that a column subset of large volume is well suited to achieving these goals. 
Computing the subset of $s$ columns of maximal volume is NP-hard \cite{ccivril2009selecting}, but---fortunately---sampling a subset of columns using the volume is tractable.

\begin{definition}[Volume-squared distribution] \label{def:volume}
    Given a matrix $\mat{A} \in \R^{m\times n}$ of rank at least $s$, 
    a random index set $\iset$ of fixed cardinality $s$ is said to follow the \emph{volume-squared distribution} if
    \begin{equation*}
        \prob \{ \iset = \isettwo \} = \frac{\vol(\mat{A}(:,\isettwo))^2}{\sum_{|\isetthree|= s} \vol(\mat{A}(:,\isetthree))^2} \quad \text{for every $s$-element subset } \isettwo \subseteq \{1,\ldots,n\}.
    \end{equation*}
\end{definition}

The volume-squared sampling distribution was introduced to computational linear algebra by \citet{deshpande2006matrix} and is often simply referred to as the volume sampling distribution.

\subsubsection{Interpolative decomposition by volume-squared sampling}

The volume sampling distribution yields interpolative decompositions with
near-optimal expected-error guarantees; see \cite[Thm.~1.3 \& Prop.~3.3]{deshpande2006matrix} and \cite{guruswami2012optimal}.

\begin{theorem}[Interpolative decomposition: Volume-squared sampling] \label{thm:id-volume}
    Let $\mat{A} \in \R^{m\times n}$ be a matrix of rank at least $s$, and let $\iset$ be a random $s$-element subset following the volume-squared distribution.
    Let $\Ahat \coloneqq \mat{A}(:,\iset)\mat{A}(:,\iset)^\pinv\mat{A}$ be the corresponding interpolative decomposition.
    Then
    \begin{equation} \label{eq:expectedvolume}
        \expect \norm{\mat{A} - \Ahat}_{\rm F}^2 \le \min_{r\le s} \left( 1 + \frac{r}{s-r+1} \right) \norm{\mat{A} - \lowrank{\mat{A}}_r}_{\rm F}^2 \le (s+1)\norm{\mat{A} - \lowrank{\mat{A}}_s}_{\rm F}^2.
    \end{equation}
    Moreover, the factor $s+1$ is optimal for an interpolative decomposition: For every $s\ge 1$ and $0 < \varepsilon < 1$, there exists a matrix $\mat{A}$ with $\rank(\mat A)>s$ for which 
    \begin{equation*}
        \norm{\mat{A} - \mat{A}(:,\iset)\mat{A}(:,\iset)^\pinv\mat{A}}^2_{\rm F} \ge (s+1-\varepsilon) \norm{\mat{A} - \lowrank{\mat{A}}_s}^2_{\rm F} \quad \text{for every $s$-element set $\iset$}.
    \end{equation*}
\end{theorem}

In particular,~\cref{eq:expectedvolume} implies the sharp existence result of
\citet{deshpande2006matrix}: Every matrix admits a subset $\iset$ of $s$ columns such that
the corresponding projection error $\norm{\mat{A} - \Ahat}_{\rm F}$
is at most a factor $\sqrt{s+1}$ larger than the best rank-$s$
approximation error.

Interestingly, the bound for interpolative decomposition with volume-squared sampling is nearly identical to the error bound~\cref{eq:hmtexpectedvalue} for the randomized SVD.
As with the randomized SVD, only a modest degree of oversampling is needed for an
$s$-column interpolative decomposition to be comparable to the best rank-$r$
approximation. \Cref{thm:id-volume} shows that volume-squared sampling achieves essentially sharp worst-case
Frobenius-norm guarantees, but its practical use has historically been
limited by the cost of sampling; see \Cref{sec:volume-sampling-algs}.

\subsubsection{Active least squares by volume-squared sampling} \label{sec:volume-sampling-id}

One can also use volume-squared sampling to solve active least-squares problems.
Unlike other index selection methods like leverage score sampling, volume-squared sampling guarantees reasonable accuracy with no oversampling, and it produces unbiased approximations to the optimal least-squares solution. The following result is due to~\citet[Thms.~5 \& 8; Prop.~9]{dereziski2018reverse}.

\begin{theorem}[Active least squares: Volume-squared sampling] \label{thm:volume-active-ls}
    Consider the active least-squares problem \cref{eq:active-ls}, assume $\mat{A}$ has full column rank, and let $\iset$ be a subset of size $s= n$ drawn from the volume-squared distribution on $\mat{A}^\tp$.
    Then $\Xtilde \coloneqq \mat{A}(\iset,:)^{-1} \mat{B}(\iset,:)$ is an unbiased approximation to the optimal least-squares solution $\mat{A}^\pinv \mat{B}$ in that $\expect \Xtilde = \mat{A}^\pinv \mat{B}$.
    Moreover,
    \begin{equation*}
        \expect \norm{\mat{A}\Xtilde - \mat{B}}_{\rm F}^2 \le (n+1) \min_{\mat{X} \in \R^{n\times p}} \norm{\mat{A}\mat{X}-\mat{B}}_{\rm F}^2.
    \end{equation*}
    Equality holds when $\mat{A}$ is in general position, that is, every $n \times n$ submatrix of $\mat{A}$ is nonsingular.
    The factor $n+1$ is optimal in the sense that there exists an inconsistent least-squares problem with 
    $\mat{A} \in \R^{(n+1)\times n}$
    such that 
    $$\norm{\mat{A} (\mat{A}(\iset,:)^{-1}\mat{B}(\iset,:)) - \mat{B}}_{\rm F}^2 = (n+1) \min_{\mat{X} \in \R^{n\times p}} \norm{\mat{A}\mat{X}-\mat{B}}_{\rm F}^2,$$
    for every $n$-element subset $\iset \subseteq \{1,\ldots,n+1\}$.
\end{theorem}

\begin{remark}[Volume-squared sampling with $s>n$]
There are two natural extensions of volume sampling to sample sizes $s>n$.
One may sample an $s$-element subset without replacement with probability
proportional to
$
    \det\!\bigl(\mat A(\iset,:)^\tp\mat A(\iset,:)\bigr),
$~\cite{dereziski2017unbiased,nikolov2022proportional},
or one may allow repeated indices, leading to \emph{volume-rescaled sampling};
see \cite{dereziski2022unbiased}. In the latter case, a sample of size
$s$ matches in distribution the (multiset) union of a
size-$n$ volume sample and $s-n$ iid uniformly sampled indices
\cite[Thm.~2.4]{dereziski2022unbiased}.
As we detail next, we can do better by using leverage score sampling in place of uniform sampling.
\end{remark}

\subsubsection{Leveraged volume sampling} \label{sec:leveraged-volume-sampling}

Volume sampling can be combined with leverage-score sampling to permit
oversampling while retaining some of the favorable properties of both
distributions. The idea is simple: To select a representative subset of $s > n$ rows of a matrix $\mat{A} \in \R^{m\times n}$ having full column rank, we draw $n$ elements from the volume-squared distribution on $\mat{A}^\tp$ and augment them with $s-n$ iid samples from the leverage score distribution of $\mat{A}$. As usual with leverage score sampling, duplicates are allowed.
To form the least-squares estimator, we weight all sampled rows by their
reciprocal square-root leverage scores, as in~\cref{eq:leverage-score-weights}.
The resulting weighted estimator is unbiased and achieves an expected
least-squares loss that is optimal within a factor $1+\varepsilon$ when using
$\order(n\log n+n/\varepsilon)$ sampled rows
\cite{dereziski2018leveraged,dereziski2022unbiased}.
Thus, leveraged volume sampling provides an attractive alternative to plain
leverage-score sampling, although it requires a more involved sampling procedure.
As an alternative to leveraged volume sampling, one can also use a union of multiple volume-squared samples, together with leverage score weighting \cite{nouy2025weighted}.

\subsubsection{Aside: Determinantal point processes}

Before going further, it will be useful to realize a connection between volume-squared sampling and another random distribution on subsets.
Given a \psd matrix $\mat{M} \in \R^{n\times n}$ of rank at least $s$, a \emph{determinantal point process of fixed size $s$} ($s$-DPP) is a random subset $\iset \subseteq \{1,2,\ldots,n\}$ with distribution
\begin{equation*}
    \prob \{ \iset = \isettwo \} = \frac{\det \mat{M}(\isettwo,\isettwo)}{\sum_{|\isetthree| = s} \det \mat{M}(\isetthree,\isetthree)} \quad \text{for every $s$-element subset } \isettwo \subseteq \{1,\ldots,n\};
\end{equation*}
see \cite{kulesza2011k}.
Determinantal point processes have a long history in probability theory~\cite{hough2006determinantal},
and they have recently received significant attention in the machine learning community as a way of selecting diverse, representative subsets from data \cite{kulesza2012determinantal}.
When $s\le \rank(\mat{A})$, volume-squared sampling on $\mat{A} \in \R^{m\times n}$ coincides with $s$-DPP sampling on $\mat{M} = \mat{A}^\tp\mat{A}$.
Volume-squared sampling and $s$-DPP sampling are different computational problems (they have different inputs, $\mat{A}$ vs.\ $\mat{M}$), but nonetheless there is significant overlap.

\subsubsection{Sampling from the volume-squared sampling distribution} \label{sec:volume-sampling-algs}

We have seen that volume sampling produces index sets of excellent quality for our two model problems.
But can we actually sample from this distribution in practice?
Are there efficient algorithms?
Efficient algorithms for general volume sampling were first developed by~\citet{deshpande2010efficient}.
By now, there is a wide literature on efficient samplers for the volume-squared and $s$-DPP distributions.
In the following, we provide a truncated summary, focusing on approaches most relevant to computational linear algebra.

\paragraph{Volume-squared sampling on a matrix with orthonormal rows.}

Let us start our discussion with an important special case.
Suppose we wish to perform volume-squared sampling to select $s = m$ columns of a matrix $\mat{A} \in \R^{m\times n}$ \emph{with orthonormal rows}. 
This special case will be used in~\Cref{sec:index-rangefinding} and as a
building block for volume-squared sampling with $s\le m$ on a general matrix.

It turns out that one can perform volume sampling for $\mat{A}$ by simply applying randomly pivoted QR~\cref{eq:rpqr}.
The correctness of this procedure is certified by the following result \cite[Prop.~19]{hough2006determinantal}.
\begin{proposition}[Volume-squared sampling $=$ RPQR]
    Let $\mat{A} \in \R^{m\times n}$ have orthonormal rows.
    Then the index set $\iset$ produced by applying $m$ steps of RPQR to $\mat{A}$ is a sample from the volume-squared sampling distribution on $\mat{A}$.
\end{proposition}
See \cite[\S2]{barthelme2023faster} for another proof.
Since the accelerated RPQR algorithm \cite[Alg.~5.2]{epperly2025embrace} produces the same distribution over pivot sequences as (plain) RPQR, it can also be used to perform volume-squared sampling.
Both variants of RPQR sample from the volume-squared sampling distribution in $\order(m^2n)$ operations.

There is a faster approach based on rejection sampling with $\order(mn + m^3 \log m)$ expected runtime.
It was proposed by \citet{dereziski2019minimax} and revisited by \citet{barthelme2023faster}.
\citet[\S4]{epperly2026adaptivea} discusses fast block implementations.
We recommend this fast block implementation, as it can be much faster than even accelerated RPQR in practice.

\paragraph{A general matrix.}

Now suppose we want to sample $s\le \rank(\mat{A})$ indices from a general matrix $\mat{A} \in \R^{m\times n}$.
To do so, we can reduce to sampling on a matrix with orthonormal rows.

The easiest case is when $s = \rank(\mat{A})$, where we can simply compute an orthonormal row basis $\mat{V}^\tp \in \R^{s\times n}$ for $\mat{A}$ and apply the orthonormal volume-squared sampler discussed above to it.
This observation will reappear in \Cref{sec:index-rangefinding} in the context of adaptive
randomized pivoting.

When $s < r \coloneqq \rank(\mat{A})$, we use a three-step procedure.
First, compute a (compact) SVD $\mat{A} = \mat{U}\diag(\sigma_1,\ldots,\sigma_r)\mat{V}^\tp$.
Next, sample a set $\isettwo \subseteq \{1,\ldots,r\}$ of cardinality $s$ with the distribution
\begin{equation} \label{eq:kdpp-subsample}
    \prob \{ \isettwo = \isetthree \} = \frac{\prod_{j \in \isetthree} \sigma_j^2}{\sum_{|\isetthree'| = s} \prod_{j \in \isetthree'} \sigma_j^2}.
\end{equation}
Finally, apply the orthonormal volume-squared sampler to
$\mat{V}(:,\isettwo)^\tp \in \R^{s\times n}$ to 
sample $s$ indices $\iset$.
This procedure was proposed in \cite{kulesza2011k} for $s$-DPP sampling, which also explains
how to sample from the distribution \cref{eq:kdpp-subsample}.
Due to the SVD, this volume sampling procedure costs $\order(mn \min \{m,n\})$ operations, making the procedure efficient only for matrices that are small in at least one dimension.
In principle, the cost of exact volume sampling can be reduced to $\order(mn \poly(s))$ operations \cite{derezinski2019exact,calandriello2020sampling}, but the polynomial factor and implicit constants of these asymptotically faster methods are large and, in our experience, the available implementations are brittle \cite[App.~A]{chen2025randomly}.

\paragraph{Markov chain Monte Carlo samplers.}
Volume-squared and dual volume sampling also admit Markov chain Monte
Carlo methods that repeatedly propose exchanges between selected and
unselected columns, with acceptance probabilities obtained from
determinant ratios \cite{anari2016monte,li2017polynomial,anari2025optimal}.
These methods come with rigorous mixing-time guarantees, but produce only approximate samples
after a finite number of iterations.

\subsection{Index selection by rangefinding} \label{sec:index-rangefinding}

So far, we have considered index selection methods that work directly with a matrix $\mat{A} \in \R^{m\times n}$ or a sketch of it. An alternative approach, originating in~\cite{golub1976rank} and called GKS in the following,
selects columns of $\mat{A}$ by performing index selection on an auxiliary matrix $\mat{V}^\tp \in \R^{s\times n}$ with orthonormal rows, chosen so that it captures the row space  of $\mat{A}$ in the sense that $\mat{A} \approx \mat{A}\mat{V}\mat{V}^\tp$. To execute this strategy, a natural choice of $\mat{V}$ is to let it contain the dominant $s$ right singular vectors of $\mat{A}$.
As a cheaper alternative, one can apply the randomized rangefinder (\cref{sec:rangefinder}) to
$\mat A^\tp$, possibly with powering (\cref{sec:lowrank-subspace-iteration}), to build $\mat{V}^\tp$.
Related developments include~\cite{boutsidis2009improved,boutsidis2014near}.

Once we have $\mat{V}^\tp$, what index selection strategies can we apply to it?
The GKS method uses CPQR to perform index selection on $\mat{V}^\top$.
An alternative pivoting strategy due to~\citet{bischof1990incremental,stewart1990incremental} provides stronger worst-case guarantees in this setting.
The effectiveness of the Bischof--Stewart strategy was demonstrated recently by \citet{damle2026computing}, who also developed a faster randomized version of the scheme.

\subsubsection{Adaptive randomized pivoting} \label{sec:arp}

Another effective choice is to apply randomly pivoted QR (\cref{alg:rpqr}) to $\mat V^\tp$, corresponding to an algorithm called adaptive randomized pivoting (ARP) by~\citet{cortinovis2026adaptive}.
Because $\mat{V}^\tp$ has orthonormal rows, this produces a sample from the volume-squared distribution. Such a sample can be generated very efficiently using the rejection sampling approach mentioned in~\cref{sec:volume-sampling-algs}; see~\cite{epperly2026adaptivea}.

\begin{theorem}[Adaptive randomized pivoting] \label{thm:arp}
Let $\mat A\in\R^{m\times n}$ and let $\mat{V}^\tp\in \R^{s\times n}$ be any matrix with orthonormal rows.
Let $\iset\subseteq\{1,\ldots,n\}$, $|\iset|=s$, be drawn from the
volume-squared distribution on $\mat V^\tp$, and consider the interpolative decomposition $\Ahat \coloneqq \mat{A}(:,\iset)\mat{A}(:,\iset)^\pinv \mat{A}$.
    Then
    \begin{equation*}
        \expect \norm{\mat{A} - \Ahat}_{\rm F}^2 \le (s+1) \norm{\mat{A}(\Id - \mat{V}\mat{V}^\tp)}_{\rm F}^2.
    \end{equation*}
    In particular, if $\mat{V}$ consists of the dominant $s$ right singular vectors of $\mat{A}$, then $\Ahat$ satisfies the following sharp worst-case bound (cf.~\cref{thm:id-volume}):
    \begin{equation*}
        \expect \norm{\mat{A} - \Ahat}_{\rm F}^2 \le (s+1) \norm{\mat{A} - \lowrank{\mat{A}}_s}_{\rm F}^2.
    \end{equation*}
\end{theorem}

\begin{proof}
    We prove this result by drawing a connection to active least squares and volume-squared sampling (\cref{sec:volume-sampling-id}).
    Given a good row-space approximation $\mat{V}^\tp$, we can think of low-rank approximation as a regression problem, where we try to find a matrix $\mat{F} \in \R^{m\times s}$ such that $\mat{F}\mat{V}^\tp \approx \mat{A}$.
    Once $\mat{V}$ is available, we can construct an interpolatory approximation by reading only \(s\) columns of \(\mat A\). To see this, we consider the least-squares problem
    \begin{equation*}
        \min_{\mat{F} \in \R^{m\times s}} \norm{\mat{A} - \mat{F}\mat{V}^\tp}_{\rm F} =
        \min_{\mat{F} \in \R^{m\times s}} \norm{\mat{A}^\tp - \mat{V}\mat{F}^\tp}_{\rm F}.
    \end{equation*}
    Drawing a volume-squared sample $\iset$ and solving the sketched problem gives
    \begin{equation} \label{eq:definterparp}
     \hatbold{\mat{F}} = \mat{A}(:,\iset)  \mat{V}(\iset,:)^{-\tp}, \qquad 
     \tildebold{\mat{A}} := \hatbold{\mat{F}} \mat{V}^\tp = 
     \mat{A}(:,\iset)\mat{V}(\iset,:)^{-\tp}\mat{V}^\tp.
    \end{equation}
    By \cref{thm:volume-active-ls}, this interpolative decomposition $\tildebold{\mat{A}}$ satisfies 
    \begin{equation*}
        \expect \norm{\smash{\mat{A} - \tildebold{\mat{A}} }}^2_{\rm F}  \le (s+1) \norm{\mat{A}(\Id - \mat{V}\mat{V}^\tp)}_{\rm F}^2.
    \end{equation*}
    Now, we can use the fact that every column of $\tildebold{\mat{A}}$ is in the range of $\mat{A}(:,\iset)$.
    By optimality of the projection approximation \cref{eq:projection-approximation-optimality}, we know that $\Ahat = \mat{A}(:,\iset)\mat{A}(:,\iset)^\pinv \mat{A}$ is the closest approximation to $\mat{A}$ whose columns are in this subspace.
    We conclude that
    \begin{equation*}
        \expect \norm{\mat{A} - \Ahat}_{\rm F}^2 \le \expect \norm{\smash{\mat{A} - \tildebold{\mat{A}} }}^2_{\rm F} \le (s+1) \norm{\mat{A}(\Id - \mat{V}\mat{V}^\tp)}_{\rm F}^2.
    \end{equation*}
    This is the desired conclusion.
    See also \cite[Thm.~2.4]{cortinovis2026adaptive} for a direct proof using the randomly pivoted QR iteration.
\end{proof}

\subsubsection{Osinsky's method} \label{sec:osinsky}

\citet{osinsky2025close} developed a deterministic method that can be
recovered by derandomizing ARP; see \cite[\S2.4]{cortinovis2026adaptive}.
To connect Osinsky's method to the proof of~\Cref{thm:arp}, let us impose
the interpolation conditions one at a time.
After selecting $t$ indices $\iset_t=\{i_1,\ldots,i_t\}$ for $1 \le t \le s$, let \(\mat F_t\) minimize \(\|\mat  A-\mat  F\mat  V^\tp\|_{\mathrm F}\) subject to the interpolation conditions \(\mat F\mat V(\iset_t,:)^\tp=\mat A(:,\iset_t)\), and define the residual \(\mat E_t:=\mat A-\mat F_t\mat V^\tp\). Initially, \(\mat F_0=\mat A\mat V\).
Provided that $\mat V(\iset_t,:)$ has full row rank, the
constrained minimizer admits the explicit formula
\begin{equation} \label{eq:fk}
    \mat F_t
    = \mat F_0
    + \mat E_0(:,\iset_t)
      (\mat V(\iset_t,:)\mat V(\iset_t,:)^\tp)^{-1}\mat V(\iset_t,:).
\end{equation}
After \(t = s\) steps,this gives
$\mat F_s=\mat A(:,\iset_s)\mat V(\iset_s,:)^{-\tp}
=\hatbold{\mat F}$ from~\cref{eq:definterparp},
and hence $\mat E_s = \smash{\mat{A} - \tildebold{\mat{A}} }$.

We now consider the first index selection in more detail.
For a candidate first pivot $j$ with
$\ell_j=\norm{\mat{V}(j,:)}^2>0$, the expression~\cref{eq:fk} gives
\[
        \mat{F}_1
        =\mat{F}_0+\mat{E}_0(:,j)
            \frac{\mat{V}(j,:)}{\ell_j},\qquad 
        \mat{E}_1^{(j)}
        =\mat{A}-\mat{F}_1 \mat{V}^\tp
         =\mat{E}_0-\mat{E}_0(:,j)
            \frac{\mat{V}(j,:)\mat{V}^\tp}{\ell_j}.
\]
Since $\mat{E}_0\mat{V}=0$ and $\mat{V}^\tp\mat{V}=\Id$,
Pythagoras gives
\[
    \norm{\mat{E}_1^{(j)}}_{\rm F}^2
    =\norm{\mat{E}_0}_{\rm F}^2
     +\frac{\norm{\mat{E}_0(:,j)}^2}{\ell_j}.
\]
ARP draws the first pivot with probability $\ell_j/s$.
Osinsky's method instead chooses a pivot
minimizing the next squared residual norm, equivalently minimizing
$\norm{\mat{E}_0(:,j)}^2/\ell_j$ over $\ell_j>0$.
Clearly, the minimum is no larger than the expected value under the ARP
distribution:
\[
    \begin{aligned}
        \min_{j:\,\ell_j>0}\norm{\mat{E}_1^{(j)}}_{\rm F}^2
        \le \sum_{j:\,\ell_j>0}\frac{\ell_j}{s}
              \norm{\mat{E}_1^{(j)}}_{\rm F}^2
        =\norm{\mat{E}_0}_{\rm F}^2
          +\frac{1}{s}\sum_{j:\,\ell_j>0}
              \norm{\mat{E}_0(:,j)}^2
        \le \frac{s+1}{s}\norm{\mat{E}_0}_{\rm F}^2.
    \end{aligned}
\]
Thus, an expected bound for ARP becomes a deterministic bound by
choosing a minimizing pivot. At subsequent steps, the same comparison
applies conditionally on the previous pivots, using the remaining
$s-t+1$ basis directions; see
\cite[proof of Thm.~2.4 and \S2.4]{cortinovis2026adaptive}.
Choosing a minimizing pivot at each step therefore gives
\[
    \norm{\mat{E}_t}_{\rm F}^2
    \le \frac{s-t+2}{s-t+1}\norm{\mat{E}_{t-1}}_{\rm F}^2.
\]
The factors telescope, and projection optimality, as in the 
proof of~\Cref{thm:arp}, yields the deterministic guarantee
\begin{equation} \label{eq:osinskybound}
          \norm{\mat{A}-\Ahat}_{\rm F}^2 \le \| \mat{A}-\tildebold{\mat{A}} \|_{\rm F}^2 
    = \norm{\mat{E}_s}_{\rm F}^2
    \le (s+1)
       \norm{\mat{A}(\Id-\mat{V}\mat{V}^\tp)}_{\rm F}^2,
\end{equation}
for $\tildebold{\mat{A}} = \mat{A}(:,\iset)\mat{V}(\iset,:)^{-\tp}\mat{V}^\tp$
and $\Ahat=\mat{A}(:,\iset)\mat{A}(:,\iset)^\pinv\mat{A}$, where $\iset = \iset_s$ is the index set returned by Osinsky's method; see \cref{alg:osinsky}.
\begin{algorithm}[t]
\caption{Osinsky's method for CSSP}
\label{alg:osinsky}
\begin{algorithmic}[1]
\Require Matrix $\mat{A}\in\R^{m\times n}$ and matrix $\mat{V}\in\R^{n\times s}$ with orthonormal columns
\Ensure Column index set $\iset=\{i_1,\ldots,i_s\}$
\State Initialize $\mat{W}\gets\mat{V}^\tp$ and $\iset\gets\emptyset$
\State \label{line:initres} Form $\mat{E}\gets\mat{A}-(\mat{A}\mat{V})\mat{V}^\tp$
\For{$t=1,\ldots,s$}
    \State Set $\ell_j\gets\norm{\mat{W}(t:s,j)}^2$ for $j\notin\iset$
    \State Choose $i_t\in\argmin_{j\notin\iset:\,\ell_j>0}
        \norm{\mat{E}(:,j)}^2/\ell_j$
    \State Choose Householder reflector $\mat{H}$ such that
        $\mat{H}\mat{W}(t:s,i_t)=\alpha\vec{e}_1$
    \State Update $\mat{W}(t:s,:)\gets\mat{H}\mat{W}(t:s,:)$
    \State \label{line:resupdate} Update $\mat{E}\gets\mat{E}
        -\mat{E}(:,i_t)\mat{W}(t,:)/\mat{W}(t,i_t)$
    \State Set $\iset\gets\iset\cup\{i_t\}$
\EndFor
\end{algorithmic}
\end{algorithm}The Householder reflectors used in \cref{alg:osinsky} provide a numerically stable way of maintaining an orthonormal basis for the complement of the span of the rows selected so far, in complete analogy to~\cref{alg:rpqr}.

It is instructive to compare Osinsky's method with ARP, that is, \cref{alg:rpqr} applied to $\mat{V}^\tp$. The crucial difference from ARP is that the pivot criterion in Osinsky's method
requires information from the residual, not just from
$\mat{V}$.
In particular, forming the initial residual in
\cref{line:initres} requires an additional access to $\mat A$,
which can be a major disadvantage when $\mat V$ has been obtained
cheaply, for example from a sketch, while another pass over
$\mat A$ or additional products with $\mat A$ are expensive.
Forming and subsequently updating the residual requires
$\order(mns)$ additional operations for dense $\mat A$.

\subsection{Applications of active least squares in scientific computing} 

Having described many methods for index selection, we provide a few comments on how they can be used, focusing on active least squares use cases in particular.
\Cref{sec:DEIM} describes applications in model reduction, and \cref{sec:subappsnap} describes an application to parameter-dependent linear systems.

\subsubsection{Model reduction and the DEIM method} \label{sec:DEIM}

Model reduction provides an important application of active least squares in scientific computing and, interestingly, many developments in index selection for model reduction parallel the methods discussed above.
A particularly prominent example is the discrete empirical interpolation method (DEIM) of \citet{chaturantabut2010nonlinear}, introduced as a finite-dimensional variant of the empirical interpolation method of \citet{barrault2004empirical}. 

The setting of DEIM is as follows.
Suppose we have a vector-valued function $\vec{f} = (f_1,\ldots,f_m) : \set{S} \to \R^m$ that we wish to evaluate, up to some error, at a large number of locations $\vec{z}_1,\vec{z}_2,\ldots$.
We assume that each component function $f_i : \set{S} \to \R$ is expensive to evaluate, and each measurement may be corrupted by noise $\varepsilon_i$.
We assume access to an orthonormal matrix $\mat{U} \in \R^{m\times n}$ that we expect to capture the range of $\vec{f}$.\footnote{In practice, such a basis $\mat{U}$ is typically obtained as the dominant $n$ left singular vectors of a matrix of snapshots, whose columns are the evaluations $\vec{f}(\vec{t}_j)$ at a large number of query locations $\vec{t}_j$.
The columns of $\mat{U} \in \R^{m\times n}$ are known as \emph{proper orthogonal decomposition} (POD) modes.
The computation of these modes presents another opportunity for
randomization, as they can be approximated efficiently via a randomized
rangefinder (\cref{sec:rangefinder}); see, e.g.,
\cite{alla2019randomized}.}

Given these ingredients, the idea of DEIM is simple and natural.
To approximate $\vec{f}(\vec{z})$, we can interpolate or regress through the noisy function values $f_i(\vec{z}) + \varepsilon_i$ for an appropriate index set $\set{I} \subseteq \{1,\ldots,m\}$.
This is an active least-squares problem of the form~\cref{eq:active-ls-subsampled} with
$\mat{A}=\mat{U}$ and
$\mat{B}=\vec{f}(\vec{z})+\vec{\varepsilon}$.
DEIM approximates $\vec{f}(\vec{z})$ in
$\range(\mat{U})$ using only a small selected subset of its entries, which can greatly reduce the cost of function evaluation.
When $|\set{I}| = n$ and $\mat{U}(\iset,:)$ is nonsingular, the resulting coefficient vector can be obtained by interpolation,  which gives
\begin{equation} \label{eq:deim-interpolate}
    \widehat{\vec{f}}(\vec z) \coloneqq
    \mat{U} \cdot \mat{U}(\iset,:)^{-1}
    \Id(:,\iset)^\tp ( \vec{f}(\vec{z})+\vec{\varepsilon} ) = 
    \mat{U} \cdot \mat{U}(\iset,:)^{-1} (\vec{f}_{\iset}(\vec{z}) + \vec{\varepsilon}(\iset)).
\end{equation}
Here, we write $\vec{f}_{\iset} = (f_{i_1},\ldots,f_{i_s})$ when $\iset = \{i_1,\ldots,i_s\}$, and we introduce the vector of noise components $\vec{\varepsilon} = (\varepsilon_1,\ldots,\varepsilon_m)$.
By applying the first inequality in~\cref{eq:sketch-solve-universal2}
with the sketching matrix $\mat{\Omega} = \Id(:,\iset)$,
it follows that the DEIM approximation~\cref{eq:deim-interpolate} satisfies
\begin{equation} \label{eq:deimerror}
 \norm{\vec{f}(\vec z)-\widehat{\vec{f}}(\vec z)}
\le
\norm{\smash{\mat{U}(\iset,:)^{-1}}}
\left(
\norm{(\Id-\mat{U}\mat{U}^{\tp})\vec{f}(\vec z)}
+
\norm{\vec{\varepsilon}(\iset)}
\right).
\end{equation}
Thus, $\|\mat{U}(\iset,:)^{-1}\|$ controls both the loss relative to the best approximation
from $\range(\mat{U})$ and the amplification of measurement noise.

\paragraph{Choosing interpolation points.}

The error bound~\cref{eq:deimerror} raises the question of whether one can always select
$\iset$, $|\iset|=n$, such that the norm of $\mat U(\iset,:)^{-1}$ is not too large.
A result by~\citet{goreinov1997theory}
states that one can always choose $\iset$ such that
\begin{equation}
    \label{eq:deim-maxvol}
    \|\mat U(\iset,:)^{-1}\|
    \leq \sqrt{1+n(m-n)}.
\end{equation}
For example, this bound holds for the submatrix $\mat U(\iset,:)$ of maximum volume, which---as already mentioned---is NP-hard to find. The original DEIM method of \cite{chaturantabut2010nonlinear} selects $\iset$ by greedily maximizing the volume, which is equivalent to Gaussian elimination with partial pivoting applied to $\mat{U}$~\cite{sorensen2016deim}. In contrast, Q-DEIM~\citep{drmavc2016new} applies CPQR to
$\mat{U}^{\tp}$. 

Interestingly, Osinsky's method (\Cref{alg:osinsky}) applied to $\mat{A}=\Id_m$
with $\mat V=\mat U$ results in an index set $\iset$ that satisfies
the same guarantee~\cref{eq:deim-maxvol} by performing only
$\order(mn^2)$ operations! As noted in~\cite{cortinovis2026adaptive}, applying randomly pivoted QR (\Cref{alg:rpqr})
to $\mat{U}^\tp$ yields the same bound in expectation. More remarkably, for every fixed
$\vec{f}\in\R^m$, the corresponding DEIM approximation satisfies 
\begin{equation}
    \label{eq:deim-arp}
    \expect\left[
      \norm{\vec{f}
      -\mat{U}\mat{U}(\iset,:)^{-1}\vec{f}(\iset)}^2
    \right]
    \leq
    (n+1)
    \norm{(\Id-\mat{U}\mat{U}^\tp)\vec{f}}^2
\end{equation}
in the noiseless case. Thus, in root-mean-square error, the loss relative to the best
approximation from $\range(\mat U)$ is at most $\sqrt{n+1}$,
independently of $m$. In fact, \cref{eq:deim-arp} is simply an instance of the active
least-squares result in~\Cref{thm:volume-active-ls} with  $\mat{A}=\mat{U}$ and
$\mat{B}=\vec{f}$. 

\paragraph{Oversampling.}
Particularly when the noise is nonzero $\vec{\varepsilon} \ne \vec{0}$, it can be helpful to oversample, taking $|\iset| > n$, and fit the measurements using (weighted) least-squares regression.
This gives the oversampled DEIM approximation
\begin{equation} \label{eq:deim-regression}
    \widehat{\vec{f}} \coloneqq \mat{U} \cdot (\mat{D}\mat{U}(\iset,:))^\pinv \mat{D}(\vec{f}_{\iset}(\vec{z}) + \vec{\varepsilon}(\iset)),
\end{equation}
commonly referred to as GappyPOD for $\mat D = \Id$.
Index selection strategies for GappyPOD include the accelerated missing-point-estimation procedure of \citet{zimmermann2016accelerated}, the related GappyPOD+E method of~\citet{peherstorfer2020stability}, and a randomized 
greedy variant by~\citet{zimmermann2024randomized}.
We allow a diagonal weight matrix $\mat{D}$, as in the weighted active
least-squares formulation~\cref{eq:active-ls-subsampled}.
Such weighting is natural for leverage-score sampling and, along these lines,~\citet{saibaba2020randomized}
uses oversampled leverage-score sampling with the corresponding rescaling for
DEIM, paralleling the construction~\cref{eq:ls-least-squares}.
The first inequality in~\cref{eq:sketch-solve-universal2} also applies to
this oversampled formulation by taking
$\mat{\Omega}^\tp=\mat{D}\Id(\iset,:)$.

\subsubsection{SubApSnap: Approximately solving parameter-dependent linear systems} \label{sec:subappsnap}

As another application, we can use active least squares to solve parameter-dependent linear systems.
The idea is due to \citet{jones2025subapsnap}, and a variant for eigenvalue problems appears in \cite{hrobat2026solving}.
Suppose we are interested in quickly solving a parameter-dependent system of linear equations 
\begin{equation} \label{eq:param-dependent}
    \mat{A}(\vec{p})\vec{x}(\vec{p})=\vec{b}(\vec{p})
\end{equation}
for many values of the parameter $\vec{p}$.
To begin, suppose that we have already solved the system \cref{eq:param-dependent} to high accuracy at a series of parameter values $\vec{p}_1,\ldots,\vec{p}_r$, which we may assemble into a snapshot matrix
\begin{equation*}
    \mat{X} \coloneqq \begin{bmatrix}
        \vec{x}(\vec{p}_1) & \cdots & \vec{x}(\vec{p}_r)
    \end{bmatrix} \in \R^{n\times r}.
\end{equation*}
If the solution changes ``slowly'' as the parameter changes, it is reasonable to expect that the solution $\vec{x}(\vec{p})$ should nearly lie in the span of these previous solutions, which motivates us to consider the least-squares problem
\begin{equation} \label{eq:subapsnap}
    \min_{\vec{c} \in \R^r} \norm{(\mat{A}(\vec{p})\mat{X})\vec{c} - \vec{b}(\vec{p})},
\end{equation}
which yields the best approximate solution $\xtilde = \mat{X}\vec{c} \approx \vec{x}(\vec{p})$ to the parameter-dependent system \cref{eq:param-dependent} in the subspace $\range(\mat{X})$.
This parallels the construction used in~\Cref{sec:sgmres}, with snapshots replacing a Krylov basis. Snapshots are a crucial idea in model order reduction~\cite{quarteroni2015reduced}, and in the context of linear systems, have been used and analyzed in, e.g.,~\cite{markovinovic2006accelerating,guido2024subspace}.

Solving the least-squares problem \cref{eq:subapsnap} is already cheaper than solving the original system \cref{eq:param-dependent}, but we can go even faster using index selection.
Given a good index set $\iset$ and, if necessary, diagonal weights $\mat{D}$, the subsampled system
\begin{equation} \label{eq:subapsnap-for-real}
    \min_{\vec{c} \in \R^r} \norm{\mat{D}\Id(\iset,:)[(\mat{A}(\vec{p})\mat{X})\vec{c} - \vec{b}(\vec{p})]}
\end{equation}
is very cheap to solve, requiring just 
$\order(nrs)$ operations to set up and 
$\order(sr^2)$ operations to solve, where $s\coloneqq |\iset|$ is the cardinality of $\iset$.
Moreover, it only requires us to read $s$ rows of the matrix $\mat{A}(\vec{p})$ at each parameter value.

In principle, a fresh index set $\iset$ and weight matrix should be chosen for every value of the parameter $\vec{p}$, obviating most of the computational benefits of subsampling.
But SubApSnap makes one more optimization to make the solve process faster.
Namely, it uses an index selection \emph{once} on \cref{eq:subapsnap} for a particular parameter value $\vec{p}_\star$ and then \emph{reuses} the same index set and weights for all future parameter values $\vec{p}$ (cf.~\cref{sec:continuouscur}).
For a general parameter-dependent family, reuse of the index set in this way could fail badly.
But, for some use cases in scientific computing, where the system \cref{eq:param-dependent} varies ``slowly'' with the parameter $\vec{p}$, this procedure is very effective, resulting in significant speedups over direct solution of \cref{eq:param-dependent} at each parameter value while maintaining high accuracy.

\subsection{Experimental evaluation for active least squares: The DEIM method} \label{sec:modred}

In this section, we provide an experimental evaluation of different index selection methods for DEIM, as described in \cref{sec:DEIM}.
We use a variant of the experiment from \cite[Fig.~3 (b)]{peherstorfer2020stability}, which tests the noise resilience of different index selection methods for DEIM.
We set the function to be
\begin{equation*}
    \vec{f}(z) = (f_1(z),\ldots,f_m(z)) \coloneqq \left( \e^{-(x_i - z)^2 / 0.005} : i=1,\ldots,m \right) \quad \text{for } m = 8192
\end{equation*}
where $x_1,\ldots,x_m$ are equidistant points on $[-2\pi,2\pi]$, and choose the basis matrix $\mat{U}$ to consist of the dominant 500 left singular vectors of the snapshot matrix $[\vec{f}(z_1),\ldots,\vec{f}(z_{2500})]$ for equidistant points $z_1,\ldots,z_{2500} \in [1,3]$.
The value 500 is in the middle of the range of values tested in \cite{peherstorfer2020stability} and is chosen to evaluate performance of the methods at a moderately large problem.
We test at a single uniformly random location $z_\star \in [1,3]$, and we assume the measurement noise is Gaussian $\vec{\varepsilon} \sim \Normal(\vec{0},10^{-8}\Id)$.

We compare the true function values $\vec{f} \coloneqq \vec{f}(z_\star)$ to the DEIM approximations \cref{eq:deim-interpolate,eq:deim-regression}.
We test four methods discussed above: CPQR (corresponding to Q-DEIM), volume sampling (equivalent to RPQR), leverage-score sampling, and leveraged volume sampling.
As another deterministic baseline, we also implemented the GappyPOD+E from~\cite[\S5.3]{peherstorfer2020stability}.
We also compare with the approximation $\mat{U}\mat{U}^\tp(\vec{f} + \vec{\varepsilon})\approx \vec{f}$ obtained by projecting all $m$ noisy evaluations $f_i(z_\star) + \varepsilon_i$ onto the approximation space.

\begin{figure}
    \centering
    \includegraphics[width=0.99\linewidth]{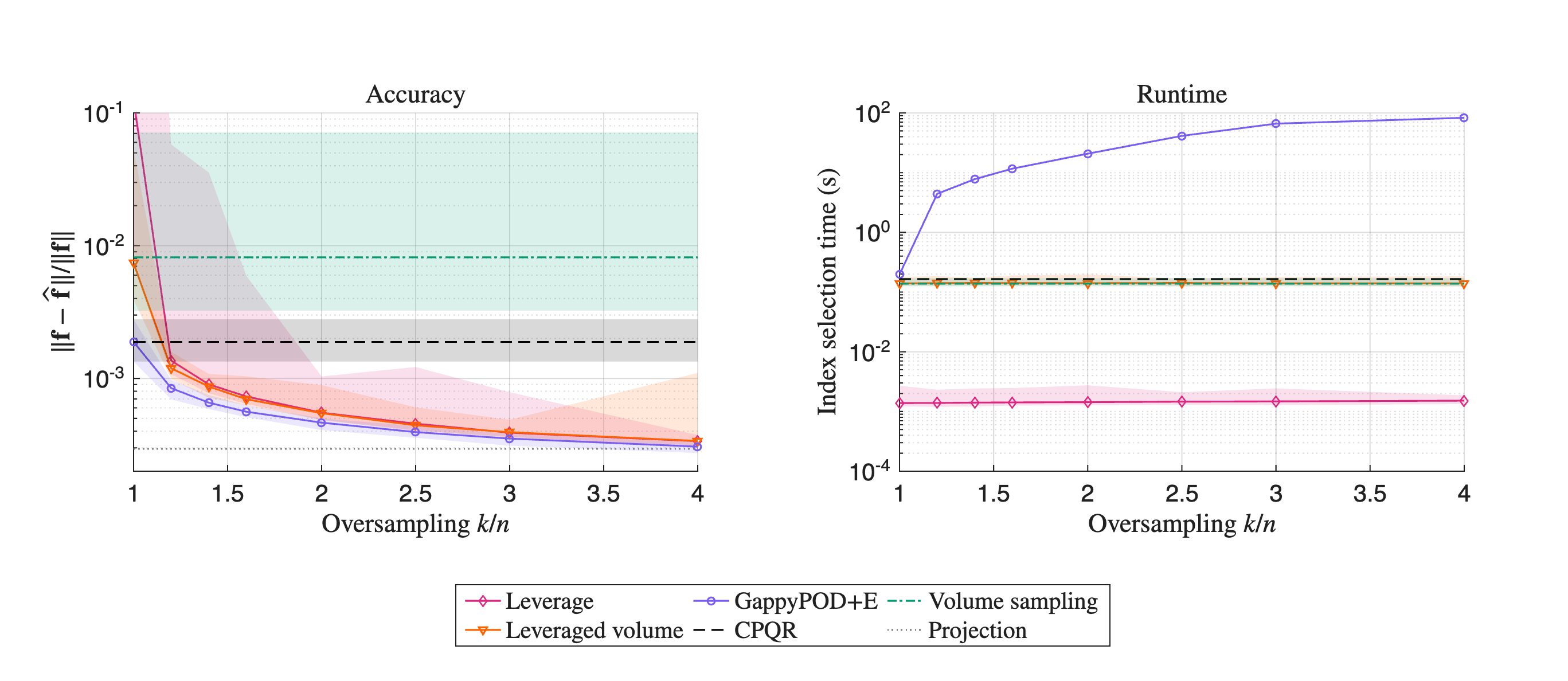}
    \caption{DEIM reconstruction error and time needed for index selection for different
sampling and pivoting strategies as a function of the oversampling ratio.}
    \label{fig:deim}
\end{figure}

\Cref{fig:deim} shows the results.
On this example, volume sampling is clearly suboptimal. 
With sufficient oversampling, leverage score-based methods achieve lower reconstruction errors than CPQR.
GappyPOD+E gives the smallest reconstruction errors among the tested methods, but at substantially higher index-selection costs as the oversampling increases.
By contrast, with the fast volume sampler described in \cref{sec:volume-sampling-algs}, volume sampling and leveraged volume sampling are fast, even marginally beating MATLAB's built-in CPQR implementation.
Leverage score sampling is faster still by over an order of magnitude.

\subsection{Experimental evaluation for interpolative decomposition: Kernel matrices}\label{sec:id-application-kernel}

Interpolative decompositions and CUR decompositions (\cref{sec:cur}) are frequently used to approximate kernel
matrices $\mat A$ \cite{march2017far,xing2020interpolative}, whose entries $\mat{A}(i,j) = \kappa(\vec{x}_i,\vec{y}_j)$ are obtained by evaluating a kernel function $\kappa$ at pairs of points $(\vec{x}_i,\vec{y}_j)$ from two given point sets.
Note that, unlike the positive-definite kernels typically used in machine learning and discussed in~\cref{sec:column-nystrom-applications}, the kernel functions and matrices considered here need \emph{not} be positive definite.
For suitable kernels and point configurations, such a kernel matrix (or off-diagonal blocks associated with well-separated point sets) admits an accurate low-rank approximation. An interpolative decomposition replaces the original sources by a small subset of representative sources, which can accelerate subsequent kernel summations.

To compare methods, we test on an example from \cite[\S7.3]{epperly2026adaptivea}.
Here, the kernel function is $\kappa(\vec{x},\vec{y}) = 1/\norm{\vec{x} - \vec{y}}$, and $\{\vec{x}_i\}$ and $\{\vec{y}_j\}$ are each $10^2\times 10^2$ Cartesian grids on the squares $[0,1)^2$ and $[0,1)\times [1,2)$, respectively.

\begin{figure}
    \centering
    \includegraphics[width=0.99\linewidth]{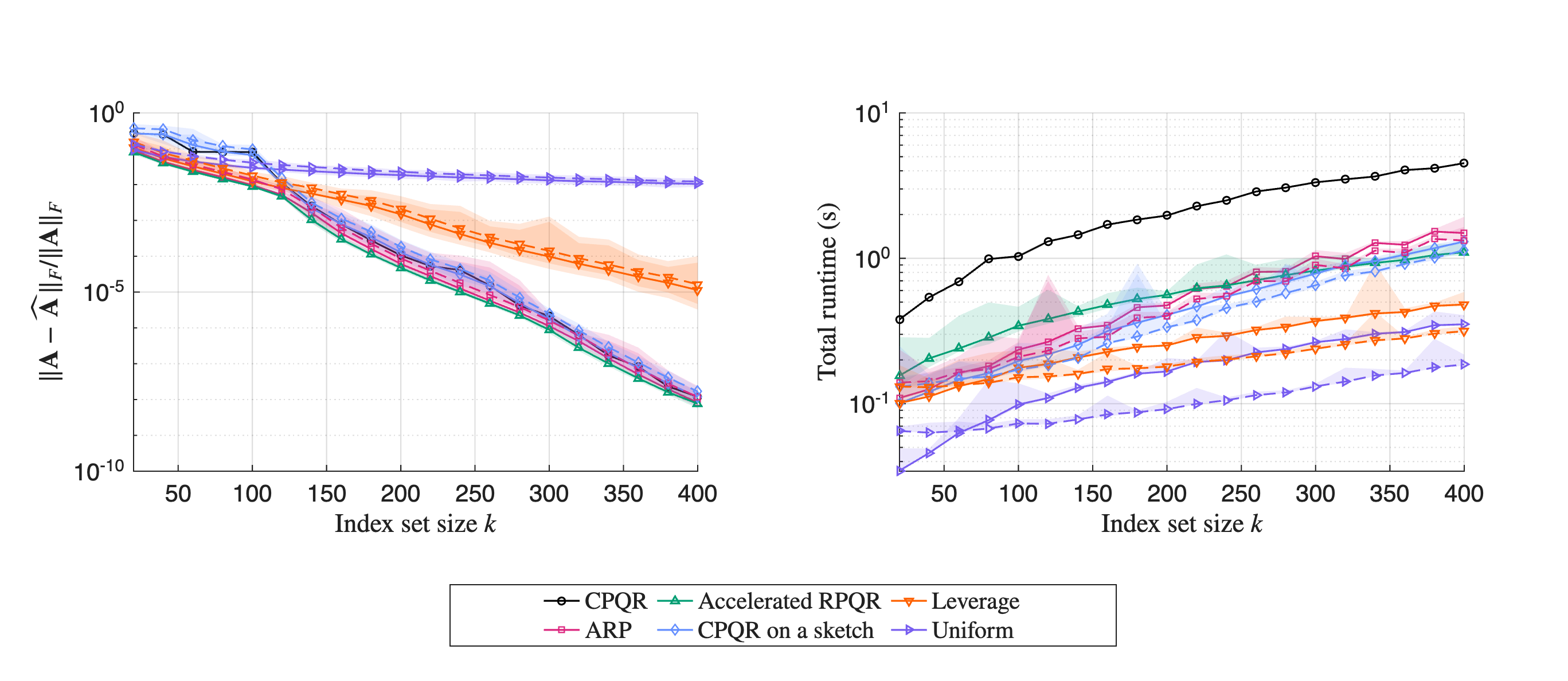}
    \caption{Accuracy and runtime of interpolative decompositions constructed
with different index selection strategies. Dashed curves use OSID to
construct the interpolation matrix.}
    \label{fig:id-comparison}
\end{figure}

\Cref{fig:id-comparison} presents a comparison of interpolative decompositions constructed using different index selection methods discussed in this section: CPQR (\cref{sec:classical-pivoting}), CPQR on a sketch (\cref{sec:sketchy-pivoting}), accelerated randomly pivoted QR (\cref{sec:random-pivoting}), ARP (\cref{sec:arp}), leverage score sampling combined with rangefinding (\cref{sec:leverage,sec:index-rangefinding}), and uniform random sampling.
For methods using sketching (CPQR on a sketch, ARP, and leverage), we use a SparseStack with embedding dimension $2s$ and sparsity level $\zeta = 4$.
We report the total runtime required both to choose the index set and to construct the interpolative decomposition in factored form; these are the solid lines in \cref{fig:id-comparison}.

On this example, ARP, (accelerated) RPQR, and CPQR (both on the matrix and on a sketch) all achieve similar levels of accuracy, and ARP, accelerated RPQR, and CPQR on a sketch all achieve similar levels of performance.
We recommend any of these methods (ARP, accelerated RPQR, and CPQR on a sketch) as a good default for practical use.
Accelerated randomly pivoted QR and robust blockwise random pivoting \cite{dong2025robust} are especially convenient when the rank is to be determined adaptively from a residual tolerance. Sketch-based methods can also be made adaptive, but may require additional sketching \cite{pearce2025adaptive}.

Lastly, we note a useful modification to make computing an interpolative decomposition faster.
When index selection does not already produce a partial QR factorization of the original matrix, computing the interpolation matrix $\mat{W}^\tp = \mat{A}(:,\iset)^\pinv \mat{A}$ can comprise a meaningful fraction of the computational cost.
To reduce cost, we recognize this computation as a tall-skinny least-squares problem, which is amenable to fast, approximate solution by the sketch-and-solve method (\cref{sec:sas}).
This results in an approximate interpolation matrix $\hatbold{W}^\tp = (\mat{\Omega}^\tp\mat{A}(:,\iset))^\pinv \mat{\Omega}^\tp\mat{A}$ that achieves a comparable approximation error to the optimal interpolation matrix:
\begin{equation*}
    \norm{\smash{\mat{A} - \mat{A}(:,\iset)\hatbold{W}^\tp}}_{\rm F} \le (1+\eta) \norm{\smash{\mat{A} - \mat{A}(:,\iset)\mat{W}^\tp}}_{\rm F},
\end{equation*}
provided $\mat{\Omega}$ is a subspace embedding (e.g., for $\operatorname{range}(\mat{A})$) with distortion $\eta$.
\citet[\S5.2]{dong2025robust} study this oversampled construction under the name \emph{oversampled sketchy interpolative decomposition} (OSID).
The dashed curves in \cref{fig:id-comparison} use OSID for methods that do not already provide a partial QR factorization of the original matrix.
We see that OSID leads to a small loss of accuracy, and, even on this small example, a meaningful reduction in runtime.

\subsection{Practical recommendations}

As we have highlighted, there are two primary settings for index selection in RNLA.
In one setting (the ``interpolative decomposition setting''), we have a matrix $\mat{A} \in \R^{m\times n}$ that is large in both dimensions and we wish to identify a row or column set of small cardinality $s \ll m,n$ that ``represent'' the matrix, spanning a near-optimal low-rank approximation and maybe possessing other nice properties such as large volume.
In the ``active least-squares setting'', we are presented with, say, a tall-thin matrix $\mat{A} \in \R^{m\times n}$ and want to find a complete or overcomplete set of rows of size $s\ge n$ that we can use to solve a subsampled fitting problem.

For the ``interpolative decomposition setting'', pivoting on a sketch (\cref{sec:sketchy-pivoting}), accelerated RPQR (\cref{sec:random-pivoting}), and ARP (\cref{sec:arp}) all consistently give high-quality results and are efficient on hardware.
The traditional CPQR method (\cref{sec:classical-pivoting}), despite unappealing behavior in the worst case and its less efficient ``one column at a time'' approach, benefits from existing high-quality implementations in LAPACK and remains an effective choice as well.
For the highest accuracy, Osinsky's method (\cref{alg:osinsky}) and strong rank-revealing QR (\cref{sec:classical-pivoting}) are among the most effective algorithms.

For the ``active least-squares setting'', we generally recommend oversampling $s > n$ when possible.
For selecting $s=n$ samples, CPQR (\cref{sec:classical-pivoting}) and volume-squared sampling (\cref{sec:volume-sampling}) both yield subsets of reasonably high quality and are fast in practice.
When $s \gg n$, leverage score sampling (\cref{sec:leverage}) performs quite well.
In the intermediate range, say $n < s \le n\log n$, leveraged volume-squared sampling (\cref{sec:leveraged-volume-sampling}) and dependent variants of leverage score sampling \cite{shimizu2024improved} can yield noticeable benefits in accuracy at somewhat higher cost.
In our view, it remains possible that more effective algorithms for index selection with oversampling remain to be discovered, as we do not know a simple, practical algorithm that achieves the theoretically optimal error scaling (\cref{rem:optimal-active-ls}).

\section{CUR approximation} \label{sec:cur}

\Cref{sec:selection} considered approximations from a selected set of columns of a matrix $\mat A\in\R^{m\times n}$.
In this section, we discuss approximations that simultaneously involve a selected set of columns $\mat{C}=\mat{A}(:,\isettwo)$ \emph{and} a selected set of rows $\mat{R}=\mat{A}(\iset,:)$.
Such approximations are often called CUR decompositions and they take the form 
\begin{equation} \label{eq:curpicture}
   \begin{tikzpicture}[baseline=-1mm]
  \node[draw, fill=gray!10,
        minimum width=\colsize+(\colsize/50),
        minimum height=\colsize,
        anchor=west] (A) at (0,0) {\phantom{$\Scale[2]{A}$}};
  \def\stripe{\colsize/10}
  \def\barh{\colsize/10}
\draw[fill=blue!20,draw=none]
  ($(A.south west)+(0,4mm)$) rectangle ($(A.south west)+(\colsize+\stripe-1.7mm,4mm-\barh)$);
\draw[fill=blue!20,draw=none]
  ($(A.south west)+(0,8mm)$) rectangle ($(A.south west)+(\colsize+\stripe-1.7mm,8mm-\barh)$);
\draw[fill=blue!20,draw=none]
  ($(A.south west)+(0,16mm)$) rectangle ($(A.south west)+(\colsize+\stripe-1.7mm,16mm-\barh)$);
  \foreach \x in {0.13,0.6,0.85}{
    \draw[fill=red!20,draw=none]
      ($(A.south west)+(\x*\colsize,0)$)
      rectangle
      ($(A.south west)+(\x*\colsize+\stripe,\colsize)$);      
  }
  \node[draw, 
        minimum width=\colsize+(\colsize/50),
        minimum height=\colsize,
        anchor=west] (A) at (0,0) {$\Scale[2]{A}$};  
\end{tikzpicture}
\approx 
   \begin{tikzpicture}[node distance=1mm,baseline=-1mm]
       \node[draw, rectangle, solid, fill=red!20, minimum width=\rowsize, minimum height=\colsize] (U) {$\Scale[1.2]{\mat{C}}$};
       \node[draw, rectangle, solid, fill=yellow!30, minimum width=\rowsize, minimum height=\rowsize,right =of U] (S) {$\Scale[1]{\mat{U}}$};
       \node[draw, rectangle, solid, fill=blue!20, minimum width=\colsize, minimum height=\rowsize,right =of S] (V) {$\Scale[1]{\mat{R}}$};
   \end{tikzpicture}\end{equation}for some choice of middle factor $\mat{U}$.\footnote{A word on notation: Following standard usage in the CUR literature, we denote the middle factor by $\mat{U}$, although this conflicts with the usual notation for the matrix of left singular vectors and the right factor in an LU decomposition.}

While the term CUR emerged in the early 2000s~\cite{drineas2003pass},
closely related constructions had already been studied as
(pseudo)skeleton approximations~\cite{goreinov1995pseudo}
and cross approximations~\cite{tyrtyshnikov2000incomplete}.
In this section, we take a broad view, and call any approximation of the form $\Ahat = \mat{C}\mat{U}\mat{R}$ a CUR \emph{approximation} to $\mat{A}$ for column and row submatrices $\mat{C}$ and $\mat{R}$, even if the matrix $\Ahat$ is presented differently (say, as we shall see, as a partial pivoted LU factorization, \cref{sec:CURandLU}).

\subsection{Choice of $\mat{U}$}

For fixed $\mat{C}$ and $\mat{R}$, the choice
$\mat{U}_{\mathrm{opt}}=\mat{C}^\dagger\mat{A}\mat{R}^\dagger$
minimizes $\|\mat{A}-\mat{C}\mat{U}\mat{R}\|_{\mathrm{F}}$ \cite[\S4]{stewart1999four}.
Moreover, since $\mat{C}\mat{C}^\dagger$ and $\mat{R}^\dagger\mat{R}$
are orthogonal projectors,
\begin{equation}
\begin{split} \label{eq:curcolsel}
\|\mat{A}-\mat{C}\mat{U}_{\mathrm{opt}}\mat{R}\|_{\mathrm{F}}^2
&= \|\mat{A}-\mat{C}\mat{C}^\dagger\mat{A}\|_{\mathrm{F}}^2
 + \|\mat{C}\mat{C}^\dagger
      (\mat{A}-\mat{A}\mat{R}^\dagger\mat{R})\|_{\mathrm{F}}^2  \\
&\leq \|\mat{A}-\mat{C}\mat{C}^\dagger\mat{A}\|_{\mathrm{F}}^2
 + \|\mat{A}-\mat{A}\mat{R}^\dagger\mat{R}\|_{\mathrm{F}}^2. 
 \end{split}
\end{equation}
For the equality we used Pythagoras, as the column spaces of the two terms are orthogonal. 
Thus, algorithms and error bounds for $\mat{C}\mat{U}_{\mathrm{opt}}\mat{R}$ 
follow by
applying column subset selection to $\mat{A}$ and $\mat{A}^\tp$,
using the techniques of \Cref{sec:selection}, including oversampling.
Examples include leverage score sampling~\cite{mahoney2009cur}, DEIM-based selection~\cite{sorensen2016deim},
derandomized volume sampling~\cite{cortinovis2020low}, nuclear-score selection~\cite{fornace2026column}, and random pivoting~\cite{epperly2025make}.

Computing $\mat{U}_{\mathrm{opt}}$ requires additional access to $\mat A$
beyond the selected rows and columns.
We therefore focus much of our discussion on 
\begin{equation} \label{eq:Ucross}
\mat{C}\mat{U}\mat{R} \quad \text{with} \quad \mat{C} = \mat{A}(:,\isettwo), \quad \mat{U}=\mat{A}(\iset,\isettwo)^\dagger, \quad 
\mat{R} = \mat{A}(\iset,:).
\end{equation}
This is sometimes called \emph{cross approximation}; the matrix $\mat{A}(\iset,\isettwo)$ is the intersection of the selected rows and columns highlighted in~\cref{eq:curpicture}.
Except when mentioned otherwise, all CUR approximations in this section will be unweighted cross approximations, $\mat{U}=\mat{A}(\iset,\isettwo)^\dagger$.
What makes this choice amazing is that, 
once suitable index sets $\iset,\isettwo$ are known, there is {\bf no need to even look at the whole matrix $\mat A$!}
All we need is the column submatrix $\mat C$ and row submatrix $\mat R$. This can be a decisive advantage in applications where entries of $\mat A$ are expensive to obtain~\cite{bebendorf2000approximation,park2025low,martinsson2019fast}.

Another natural question is how to pick the subset sizes $|\iset|$ and $|\isettwo|$: Should they be equal or should we oversample by choosing one larger than the other?
There are good reasons to use oversampling.
Indeed, the CUR approximation is a special case of the generalized \Nystrom approximation
$\mat A\approx \mat{A}\mat{\Omega} (\mat{\Psi}^\tp \mat{A}\mat{\Omega})^\dagger \mat{\Psi}^\tp \mat{A}$
(\Cref{sec:lowrank-generalized-nystrom}) with $\mat{\Omega} = \Id_n(:,\isettwo)$ and $\mat{\Psi} = \Id_m(:,\iset)$, and we know from that algorithm that oversampling is generally necessary to obtain the most accurate results.
Using oversampling can lead to CUR approximations that are more accurate \cite{drineas2008relative} and numerically stable \cite{park2025accuracy}.
Moreover, if we use oversampling, we are free to use a cheaper method to select the larger index set, such as leverage score sampling (\cref{sec:leverage}).
However, there are also reasons to prefer equal sizes $|\iset| = |\isettwo|$.
First, choosing equal sizes results in an approximation that is simpler and more easily to manipulate, as the core matrix is simply a matrix inverse $\mat{U} = \mat{A}(\iset,\isettwo)^{-1}$.
Second, many algorithms---such as those based on pivoted LU decompositions (\cref{sec:CURandLU})---naturally return sets of equal size $|\iset| = |\isettwo|$.
Last, while allowing oversampling can lead to improvements in accuracy on challenging problems (see, e.g., \cite[Thm.~10.8]{epperly2025make}), many matrices $\mat{A}$ that we encounter in scientific computing do not exhibit worst-case behavior, and well-chosen subsets $\iset$ and $\isettwo$ of equal size often produce approximations of satisfactory quality.
For concreteness, this section will largely focus on the case $|\iset| = |\isettwo|$.

\subsection{Why CUR?}

CUR approximations are widely used in scientific computing, for instance in boundary element methods \cite{bebendorf2000approximation}, seismic modeling \cite{fomel2013seismic}, X-ray spectro-microscopy \cite{meier2026reducing}, image and video analysis \cite{cai2021robust}, genetics \cite{mahoney2009cur,bodor2012rcur}, and chemical machine learning \cite{imbalzano2018automatic}.
In these applications, the CUR approximation is typically preferred over alternatives for one of two reasons: \emph{computational efficiency} and \emph{interpretability}.

The core efficiency benefit of the CUR approximation was stated above:
\actionbox{Once good index sets $\iset$ and $\isettwo$ are known, the CUR approximation can be formed  without reading the entire matrix $\mat{A}$.}
\noindent 
This property makes the CUR approximation ideal for applications where each entry of $\mat{A}$ is expensive to obtain, requiring, say, an expensive formula to be evaluated or a physical experiment to be run in order to access each matrix entry \cite{bebendorf2000approximation,fomel2013seismic,meier2026reducing}.

Now there is a catch: 
The CUR approximation can be generated without reading the whole matrix \emph{given the index sets $\iset,\isettwo$}.
For a worst-case matrix with a single large entry in an unknown location,  we must examine nearly the entire matrix to find good sets $\iset,\isettwo$.
Fortunately, most matrices do not exhibit such worst-case behavior, and there are effective heuristics for identifying the index sets that read only a small number of matrix entries \cite{bebendorf2003adaptive,chiu2013sublinear,cortinovis2025sublinear}; see \cref{subsec:CURconnections}.
Additionally, we can often use a more computationally expensive method to identify index sets $\iset,\isettwo$, then reuse these index sets as the matrix undergoes changes such as addition of rows and columns or smooth changes under variation of a parameter; see \cref{sec:continuouscur}.

The second core advantage of the CUR approximation is interpretability \cite{mahoney2009cur}.
In statistics and machine learning, we frequently deal with tabular data where the $i$th row lists a collection of numerical attributes (``features'') for the $i$th data element.
The CUR decomposition identifies a small set of rows (data elements) and columns (features) from which one can approximately reconstruct the full data set.
Such a representation can be easier to interpret than another representation such as the singular value decomposition that uses row and column bases that are not submatrices.
For instance, compare 
\begin{align*}
    \text{weight in kg} &\approx 0.45 \cdot \text{height in cm} + 0.05 \cdot \text{age in years}, \\
    \text{weight in kg} &\approx 340 \cdot \text{first singular vector} - 110 \cdot \text{second singular vector}.
\end{align*}
An illustration of this effect is given in~\cite[Fig.~1]{sorensen2016deim}.
The interpretability of the CUR approximation plays a fundamental role in genetics applications \cite{mahoney2009cur,bodor2012rcur}. 
If one needs only ``column interpretability'' but not ``row interpretability'' in a given application, then the interpolative decomposition (\cref{sec:activelsq}) provides another option. 

\subsection{Existence of good CUR approximations}
\label{sec:cur-quality} 

As we saw in \cref{thm:id-volume}, any matrix $\mat{A}$ possesses a collection $\mat{C}$ of $r$ columns such that the interpolative approximation $\Ahat \coloneqq \mat{C}\mat{C}^\dagger\mat{A}$ satisfies
\begin{equation} \label{eq:interp-best}
    \norm{\mat{A} - \Ahat}_{\rm F} \le C_r \norm{\mat{A} - \lowrank{\mat{A}}_r}_{\rm F} \quad \text{with } C_r = \sqrt{r+1}.
\end{equation}
The factor $\sqrt{r+1}$ is sharp for a worst-case instance.
Thus, a rank-$r$ interpolative approximation can get within a factor $\sqrt{r+1}$ of the best rank-$r$ approximation to $\mat{A}$.
Since the interpolative approximation is the closest matrix to $\mat{A}$ in the span of $\mat{C}$ (see \cref{prop:projection}), it stands to reason that the best-possible CUR approximation will satisfy a guarantee like \cref{eq:interp-best} with a prefactor larger than $\sqrt{r+1}$.

The fundamental theorem of the CUR approximation (FTCUR) shows that CUR (cross) approximations achieve the guarantee \cref{eq:interp-best} with $C_r = r+1$.

\begin{theorem}[FTCUR: Fundamental theorem of CUR] \label{thm:ftcur}
Given $\mat{A} \in \R^{m \times n}$ of rank at least $r$, there exist index sets with 
$|\iset|=|\isettwo|=r$, such that the resulting CUR~\eqref{eq:Ucross} satisfies
    \begin{equation}
        \label{eq:CURboundbest}
        \|\mat A-\mat C\mat U\mat R\|_{\rm F}
        \leq
        (r+1)\|\mat A-\lowrank{\mat A}_r\|_{\rm F}.
    \end{equation}
\end{theorem}

This theorem was originally proven by \citet{zamarashkin2018existence}.
The sharpness of the factor $r+1$ follows from the \pd example discussed after~\Cref{thm:fundamental-nystrom} below. An earlier result by \citet{goreinov2001maximal} controls the entrywise error
$\|\mat A-\mat C\mat U\mat R\|_{\rm max} \le (r+1) \|\mat A-\lowrank{\mat A}_r\|$
when $\mat A(\iset,\isettwo)$ has maximum volume, with the maximum norm $\norm{\mat{B}}_{\rm max} \coloneqq \max_{i,j} |\mat{B}(i,j)|$.
A probabilistic proof of the FTCUR is easily obtained by choosing $\isettwo$ using volume sampling on $\mat{A}$, choosing $\iset$ by volume sampling on $\mat{C}^\top$, and applying \cref{thm:id-volume,thm:volume-active-ls}.
It is quite remarkable that the suboptimality factor depends only on $r$, not on $m$ or $n$.
In contrast, existing worst-case bounds in the spectral norm involve factors growing with $m,n$; see, e.g.,~\cite{miranian2003strong}. In~\cref{sec:cur-rowspace}, we will provide a constructive proof for FTCUR.

If one is willing to use the optimal CUR approximation $\mat{C}\mat{U}_{\rm opt}\mat{R}$, then one obtains a better prefactor.
Indeed, combining the interpolative decomposition result \cref{eq:interp-best} and the error bound \cref{eq:curcolsel} yields index sets $\iset$ and $\isettwo$ of size $r$ such that
\begin{align} \label{eq:curopt}
\|\mat{A}-\mat{C}\mat{U}_{\mathrm{opt}} \mat{R}\|_{\rm F} & \le  \sqrt{2(r+1)}\|\mat{A}-\lowrank{\mat A}_r\|_{\rm F}
\end{align}
By allowing for $\order(r/\varepsilon)$ selected rows and columns, and
using a middle factor of rank at most $r$, \citet{boutsidis2017optimal} obtain a
$(1+\varepsilon)$ relative-error guarantee in the Frobenius norm. 

\subsection{Equivalence between CUR and LU with row--column pivoting}\label{sec:CURandLU}
As it happens, a rank-$r$ CUR approximation coincides with the approximation obtained by performing $r$ steps of a LU factorization with row--column pivoting (provided that all pivots are nonzero).  
To see the connection to LU, we examine the CUR residual
\[
 \mat{E} = \mat{A}-\mat{C}\mat{U}\mat{R} = \mat{A} - \mat{A}(:,\isettwo) \mat{A}(\iset,\isettwo)^{-1} \mat{A}(\iset,:).
\]
We have that
\[
 \mat{E}(:,\isettwo) = \mat{A}(:,\isettwo)-\mat{A}(:,\isettwo) \mat{A}(\iset,\isettwo)^{-1} \mat{A}(\iset,\isettwo) = \mat{0},
\]
and, similarly, $\mat{E}(\set{I},:)=\mat{0}$. Hence, the \emph{CUR approximation recovers the columns in $\isettwo$ and the rows in $\iset$ exactly!}
By permuting these rows and columns to the top-left corner, we may assume without loss of generality that 
$\iset = \isettwo = \{1,\ldots,r\}$. By partitioning $\mat{A}=\begin{bmatrix}
\mat{A}_{11} & \mat{A}_{12}\\
\mat{A}_{21} & \mat{A}_{22}
\end{bmatrix}$ with $\mat{A}_{11} \in \mathbb R^{r\times r}$, it follows that $\mat{C} = 
\begin{bmatrix}
\mat{A}_{11} \\
\mat{A}_{21} 
\end{bmatrix}$, 
$\mat{U} = \mat{A}_{11}^{-1}$, 
and
 $\mat{R} = 
\begin{bmatrix}
\mat{A}_{11} & \mat{A}_{12}
\end{bmatrix}$. 
Thus, we have the interpolation property of CUR in the chosen indices
\begin{equation}\label{eq:CUR0}
\mat{A}-\mat{C}\mat{U}\mat{R} = \begin{bmatrix}
0_{r\times r} & 0_{r\times (n-r)} \\ 0_{(m-r)\times r} & \mat{F}    
\end{bmatrix}. 
\end{equation}
 Here, $\mat{F} = \mat{A}_{22} - \mat{A}_{21}\mat{A}_{11}^{-1}\mat{A}_{12}$ is  the \emph{Schur complement}, a central object in linear algebra~\cite{unknown2006schur}. 

To make the connection with LU, we first note~\cite{golub2013matrix} that applying
$r$ steps of GE also yields a decomposition of the form
\[
 \mat{A} - \mat{L} \tildebold{R} = 
 \begin{bmatrix}
0_{r\times r} & 0_{r\times (n-r)} \\ 0_{(m-r)\times r} & \mat{F}    
\end{bmatrix}, 
\]
where $\mat{L}\in\R^{m\times r}$ is lower trapezoidal, and $\tildebold{R}\in\R^{r\times n}$ is upper trapezoidal.
We use a tilde to distinguish $\tildebold{R}$ (right factor in LU) from $\mat{R}$ (row submatrix).
Partitioning $\tildebold{R} = 
\begin{bmatrix}
\tildebold{R}_{11} & \tildebold{R}_{12}
\end{bmatrix}$ and 
considering the first $r$ columns of this relation, it follows that $\mat{C} = \mat{L} \tildebold{R}_{11}$. Analogously, $\mat{R} = \mat{L}_{11} \tildebold{R}$. But $\mat{U} = \mat{A}_{11}^{-1} = ( \mat{L}_{11} \tildebold{R}_{11})^{-1}$. Thus, $\mat{C}\mat{U}\mat{R} = \mat{L} \tildebold{R}$. 

The equivalence between CUR and row--column pivoted LU means that an algorithm for choosing index sets $\iset$ and $\isettwo$ for CUR yields a pivoting order for pivoted LU and vice versa.\footnote{Take care to note that equivalence applies only to the CUR cross approximation \cref{eq:Ucross} with equal sizes $|\iset| = |\isettwo|$.}
The use of GE with row--column pivoting to compute a ``rank-revealing'' decomposition is classical \cite{hwang1992rank}, and this equivalence underlies early work on CUR approximations
\cite{bebendorf2000approximation,tyrtyshnikov2000incomplete,bebendorf2003adaptive}.
The observation that GE can be viewed as an iterative method for low-rank approximation is also highlighted by \cite{townsend2013gaussian}; see also~\cite{damle2026estimating}.

\subsection{Numerical stability} \label{sec:CUR-stability}

For the CUR approximation to produce reliable results, care must be taken to store the low-rank matrix decomposition in a stable way.
Take warning:
\actionbox{The natural way of representing a CUR decomposition where one forms and stores the matrix entries of $\mat{U}$ is numerically unstable.}
\noindent The issue, highlighted by \citet[\S13.1]{martinsson2020randomized}, is that when the matrix $\mat{A}$ is numerically low-rank, the core matrix $\mat{U}=\mat{A}(\iset,\isettwo)^{-1}$ is the inverse of a matrix that can be very ill-conditioned, resulting in a matrix $\mat{U}$ with large norm.
Computing and storing this inverse directly can be disastrous.
The instabilities of the direct CUR representation are not cured by using the optimal core matrix $\mat{U}_{\rm opt}$, and they are not cured by oversampling.
See \cite[Fig.~10.2]{epperly2025make} and \cite[Fig.~2]{anderson2015spectral} for examples.

The solution is to use an alternative representation of the CUR approximation with improved stability properties.
For instance, as we saw last section, a pivoted LU decomposition computes the same approximation as the CUR decomposition, but uses the more numerically stable representation $\mat{L}\tildebold{R}$.
Other (empirically) stable factorizations for representing a CUR approximation are discussed in \cite{anderson2015spectral,park2025accuracy}.

A second strategy is to maintain the CUR structure but represent the $\mat{U}$ matrix \emph{implicitly} \cite[\S10.3]{epperly2025make}.
For instance, we can compute a QR factorization $(\mat{Q},\mat{T})\gets \Call{QR}{\mat{A}(\iset,\isettwo)}$ and store the matrices $\mat{Q}$ and $\mat{T}$.
The core is now represented as $\mat{U} = \mat{T}^{-1}\mat{Q}^\top$.
By using an implicit representation, we maintain most of the structure of the CUR decomposition while improving its stability.
See \cite[\S10.3]{epperly2025make} for more details.

Last, let us note that detailed numerical stability analysis for the CUR approximation is available only in special cases, such as for Gaussian elimination with complete pivoting (\cref{subsec:CURconnections}).
We regard the development of a more complete theory for the numerical stability of the CUR approximation to be a natural subject for future research.

\subsection{Algorithms for CUR}

The computation of a good CUR approximation is a nontrivial task, as we must simultaneously identify a good row set and a good column set.
This subsection will describe algorithms that meet this task.
Some algorithms will meet the optimal guarantee provided by the FTCUR (\cref{thm:ftcur}), and others will not attain this guarantee but nonetheless are effective. 

\subsubsection{A general strategy: Columns then rows} \label{sec:curfromcssp}

The accuracy of the CUR approximation is always no better than the projection approximation $\mat{C}\mat{C}^\dagger \mat{A}$.
For this reason, it is natural to begin searching for a good CUR approximation by first solving the CSSP to select a good set of columns $\mat{C}$.
Naturally, we already have seen many strategies for doing so in \cref{sec:selection}.

Once the columns have been selected, choosing a good set of rows can be seen as an active least-squares problem (\cref{sec:activelsq}).
Indeed, the ``UR'' part of a CUR approximation can be seen as an approximation $\mat{C}^\dagger \mat{A}\approx\mat{C}(\iset,:)^\dagger \mat{A}(\iset,:) = \mat{U}\mat{R}$ that is constructed using only a subset of the rows of the ``right-hand side'' $\mat{A}$.
Thus, to pick a good set of rows, we can apply an index selection method to pick \emph{rows} of $\mat{C}$.
Again, \cref{sec:selection} provides a buffet of options.\footnote{Note that many index selection methods for active least squares are most natural with both weighting and oversampling.
In such cases, it is natural to use a weighted version of the core matrix $\mat{U} = (\mat{D}\mat{A}(\iset,\isettwo))^\dagger\mat{D}$, where $\mat{D}$ is a diagonal weight matrix populated with, say, inverse square-roots of the leverage scores (\cref{sec:leverage-ls}).}

Combining these elements yields the columns-then-rows strategy.
Informally,
\actionbox{
Apply an index selection method to $\mat{A}$ to pick good columns $\isettwo$, then apply an index selection method to $\mat{C} = \mat{A}(:,\isettwo)$ to pick good rows $\iset$.
}
\noindent Naturally, one can also do rows-then-columns.

\begin{remark}[Bad strategy: Independent row and column selection]
    An important feature of the columns-then-rows strategy is that the rows should be picked by doing row selection \emph{on the matrix of columns $\mat{C}$}.
    Picking the row and column sets independently and directly on $\mat{A}$ can lead to disastrous results.
    Indeed, consider
    \begin{equation}\label{eq:curbadexample}
        \mat{A} = \twobytwo{\varepsilon}{1}{1}{0} \quad \text{for } 0 < \varepsilon \ll 1.
    \end{equation}
    Clearly, in isolation, $\iset = \{1\}$ and $\isettwo = \{1\}$ are 
    the optimal choices for selecting one column from $\mat{A}^\tp$ and $\mat{A}$, respectively. Nevertheless,
    $\mat A(1,1)=\varepsilon$, and the resulting CUR  error can be arbitrarily large:
    \[
        \left\|
        \mat A-
        \mat A(:,1)\,
        \mat A(1,1)^{-1}\,
        \mat A(1,:)
        \right\|_{\rm F}
         =\varepsilon^{-1}.
    \]
\end{remark}

\subsubsection{ARP-CUR and Osinsky-CUR: Algorithms achieving the FTCUR} \label{sec:cur-rowspace}

Using the columns-then-rows strategy, we can derive practical algorithms achieving the FTCUR error bound~\cref{eq:CURboundbest}.
As in \cref{sec:index-rangefinding}, 
we start with a matrix $\mat{V}^\tp\in\R^{r\times n}$ that has orthonormal rows and represents a 
row space approximation of $\mat{A}$, in the sense that
$\mat{A} \approx \mat{A} \mat{V}\mat{V}^\tp$.
We apply Osinsky's method (\cref{sec:osinsky}) or ARP (\cref{sec:arp}) to select an index set $\isettwo$ such that 
\[
 \mat{A} \approx \mat{A} \mat{V}\mat{V}^\tp \approx \mat{A}(:,\isettwo)\mat{V}(\isettwo,:)^{-\top} \mat{V}^\tp.
\]
We set $\mat{C} = \mat{A}(:,\isettwo)$ and consider, as above, a
thin QR factorization $\mat C=\mat Q_C\mat R_C$.
Applying 
Osinsky's method (\cref{alg:osinsky} with 
$\mat{A}\leftarrow\mat{A}^\tp$,  $\mat{V}\leftarrow\mat{Q}_{\mat{C}}$) or volume sampling yields another index set $\iset$ with 
$|\iset| = r$.
This idea leads to the two algorithm variants, Osinsky-CUR and ARP-CUR, displayed in~\cref{alg:osinsky-cur}.

\begin{algorithm}[t]
\caption{Osinsky-CUR \emph{or} ARP-CUR}
\label{alg:osinsky-cur}
\begin{algorithmic}[1]
\Require Matrix $\mat A\in\R^{m\times n}$ and
         $\mat V\in\R^{n\times r}$ with $\mat V^\tp\mat V=\Id$
\Ensure CUR approximation $\mat A\approx\mat C\mat U\mat R$ with implicitly represented $\mat{U} = \mat{T}^{-1}\mat{Q}^\top$
\State {\bf Call} $\isettwo\gets\textsc{Osinsky}(\mat A,\mat V)$ {\bf or}
$\isettwo\gets\textsc{RPQR}(\mat V^\tp,r)$
       \Comment{Algorithms~\ref{alg:osinsky},~\ref{alg:rpqr}}
\State $\mat C\gets\mat A(:,\isettwo)$
\State $\mat Q_C\gets\operatorname{Orth}(\mat C)$
\State {\bf Call} $\iset\gets\textsc{Osinsky}(\mat A^\tp,\mat Q_C)$ {\bf or} 
$\iset\gets\textsc{RPQR}(\mat Q_C^\tp, r)$
       \Comment{Algorithms~\ref{alg:osinsky},~\ref{alg:rpqr}}
\State $\mat R\gets\mat A(\iset,:)$
\State $(\mat{Q},\mat{T}) \gets \Call{QR}{\mat{A}(\iset,\isettwo)}$ \Comment{Implicit CUR representation (\cref{sec:CUR-stability})}
\State \Return $\mat C,\mat Q, \mat T,\mat R$
\end{algorithmic}
\end{algorithm}

Given $\mat{V}$, Osinsky-CUR is deterministic and requires $\order(mnr)$ operations involving both $\mat{A}$ and $\mat{V}$.
Using randomization, ARP-CUR allows us to work with $\mat{V}$ and selected rows/columns of $\mat{A}$ only, reducing the complexity to $\order((m+n)r^2)$ operations. The following result is due to \citet{osinsky2025close,cortinovis2026adaptive}. The assumption on the range of $\mat V$ ensures that $\mat{A}(\iset,\isettwo)$ is invertible.

\begin{theorem}[Osinsky-CUR and ARP-CUR] \label{thm:curosinsky}
Let $\mat A \in \mathbb R^{m\times n}$ and let $\mat V\in\R^{n\times r}$ have orthonormal columns, such that $\range(\mat V) \subseteq \range(\mat A^\tp)$.
Then \Cref{alg:osinsky-cur}, when combined with Osinsky's method, returns a CUR approximation such that
\[
 \| \mat{A} - \mat{CUR}\|^2_{\mathrm F} \le (r+1)^2  \|\mat{A}- \mat{AVV}^\tp\|_{\rm F}^2.
\]
When combined with ARP, this bound holds in expectation.
\end{theorem}
\begin{proof}
Let 
$\mat{C}=\mat{Q}_{\mat{C}}\mat{R}_{\mat{C}}$ be a thin QR factorization. 
Then $\mat{A}(\iset,\isettwo)=\mat{Q}_{\mat{C}}(\iset,: )\mat{R}_{\mat{C}}$, so 
\begin{equation}
\label{eq:curasprojector}
    \mat{C}\mat{U}\mat{R} = 
   \mat{Q}_{\mat{C}}     \mat{Q}_{\mat{C}}(\iset,:)^{-1} \mat{A}(\iset,:)
    = 
    \mat{Q}_{\mat{C}}
    \mat{Q}_{\mat{C}}(\iset,:)^{-1}
    \Id_m(:,\iset)^\tp
    \mat{A}.
\end{equation}
Then 
the result for Osinsky-CUR follows from applying the error bound~\cref{eq:osinskybound} for \cref{alg:osinsky} twice:
\begin{align*}
\|\mat{A} - \mat{CUR}\|_{\rm F}^2 &= \|\mat A - \mat{Q}_{\mat{C}}     \mat{Q}_{\mat{C}}(\iset,:)^{-1} \mat{A}(\iset,:) \|_{\rm F}^2 \\
                &\le (r+1) \|(\Id-\mat{Q}_{\mat{C}}^{\vphantom{\tp}}\mat{Q}_{\mat{C}}^\tp)\mat{A}\|_{\rm F}^2 && \text{(by~\cref{eq:osinskybound})} \\
                &\le (r+1)\|\mat{A} - \mat{A}(:,\isettwo)\mat{V}(\isettwo,:)^{-\tp}\mat{V}^\tp\|_{\rm F}^2 && \text{(by optimality of $\mat{Q}_{\mat{C}}^{\vphantom{\tp}}\mat{Q}_{\mat{C}}^\tp$)} \\ 
                &\le (r+1)^2 \|\mat{A}- \mat{AVV}^\tp\|_{\rm F}^2. && \text{(by~\cref{eq:osinskybound})}
\end{align*}
The result for ARP-CUR follows in a similar manner, using \cref{thm:volume-active-ls,thm:arp}.
\end{proof}

\Cref{thm:curosinsky} recovers the FTCUR bound~\cref{eq:CURboundbest}. 
However, to obtain such a $\mat{V}$ via the SVD would require $\order(m n^2)$ operations for $m\ge n$.
One can instead extract $\mat{V}$ from a sketch $\mat{\Omega}^\tp \mat{A}$ through, e.g., an SVD. This reduces the complexity to $\order(mnr)$ or even below. When using a Gaussian test matrix $\mat{\Omega} \in \R^{m\times s}$, the error satisfies the bound of \cref{thm:curosinsky} in expectation, with a small  additional multiplicative factor for $s = 2r$; see~\cref{eq:hmtexpectedvalue}.

\subsubsection{LU- and ACA-based methods}
\label{subsec:CURconnections}

A different class of algorithms is obtained from the close connection 
of CUR with classical Gaussian elimination (GE), described in \cref{sec:CURandLU}. 

This observation leads to~\cref{alg:cur-lu-template}, which is just GE with pivoting stopped after $r$ steps. Because each step adds a ``cross'' consisting of one selected row and column, this algorithm is often called (adaptive) \emph{cross approximation}. In contrast to the CUR algorithms discussed so far, elimination-based
methods allow for selecting one row and column index at a time, \emph{adaptively} in each step. All common pivot selection rules are based on the CUR residual. Letting
$\mat L_{t-1} \in \R^{m\times(t-1)}$, $\mat R_{t-1}\in \R^{(t-1) \times n}$ denote the factors obtained after $t-1$ steps, we define
\begin{equation*}
 \mat E_{t-1} := \mat A - \mat L_{t-1} \mat R_{t-1}.
\end{equation*}

\begin{algorithm}[t]
\caption{CUR by pivoted LU (algorithm template)}
\label{alg:cur-lu-template}
\begin{algorithmic}[1]
\Require Row and column access to $\mat A\in\R^{m\times n}$;
         target rank $r$; pivot selection rule.
\Ensure Index sets $\iset,\isettwo$ defining a rank-$r$ CUR approximation ${\mat A} \approx \mat L\tildebold R$.
\State Initialize $\iset\gets\emptyset$, $\isettwo\gets\emptyset$,
       $\mat L\gets[\,]\in\R^{m\times0}$,
       $\tildebold R \gets[\,]\in\R^{0\times n}$
\For{$t=1,\ldots,r$}
    \State Use pivot selection rule to select pivots $i_t \not\in \iset$, $j_t \not\in \isettwo$ 
    \State $\vec v\gets\mat A(:,j_t)-\mat L\tildebold R(:,j_t)$
    \State $\vec w^\tp\gets\mat A(i_t,:)-\mat L(i_t,:)\tildebold R$
    \State $\mat L\gets[\mat L,\,\vec v/ \vec v(i_t)]$,
           $\tildebold R\gets\begin{bmatrix}\tildebold R\\\vec w^\tp\end{bmatrix}$
    \State $\iset \gets \iset \cup \{i_t\}$, $\isettwo \gets \isettwo \cup \{j_t\}$
\EndFor
\end{algorithmic}
\end{algorithm}

\paragraph{Complete pivoting and volume maximization.} Gaussian elimination with complete pivoting (GECP) chooses 
\[
    (i_t,j_t)\in\operatorname*{arg\,max}_{i,j}
                      |\mat E_{t-1}(i,j)|.
\]
This choice greedily maximizes the volume of the selected submatrix in step $t$, considering the previous $t-1$ steps fixed~\cite{bebendorf2000approximation,tyrtyshnikov2000incomplete}. This interpretation connects GECP to the existence result 
by~\citet{goreinov2001maximal}, discussed briefly in~\Cref{sec:cur-quality}. While GECP is observed to work very well in practice, the worst-case error bounds~\cite{cortinovis2020maximum} are fairly pessimistic:
$
 \|\mat A-\mat L \tildebold R\|_{\max}
 \le 4^r\rho_r\,\sigma_{r+1}(\mat A),
$
where $\rho_r$ is the classical worst-case pivot growth factor for GECP \cite{wilkinson1961error,shah2026entry}. 
Volume-based methods can improve a selected submatrix by exchanging
selected and unselected indices; see \cite{tyrtyshnikov2000incomplete,goreinov2010how,damle2026estimating}.
One major disadvantage of GECP is that it requires searching and, hence, accessing the whole matrix $\mat{A}$.

\paragraph{Adaptive cross approximation (ACA).} The adaptive cross approximation (ACA) algorithm by \citet{bebendorf2003adaptive} searches rows and columns of the residual: Given a row index $i_{t}$, the column index $j_{t}$ is determined by the maximal entry in row $i_{t}$ of $|\mat{E}_{t-1}|$,
and the next row index $i_{t+1} \not= i_t$
by the maximal entry in column $j_{t}$ of $|\mat{E}_{t-1}|$, etc.
ACA has the major advantage that
only one row and one column of $\mat{A}$ and $\mat L_{t-1} \mat R_{t-1}$ need to be evaluated in each step.
As such, the algorithm runs in $\order((m+n)r^2)$ operations,  \emph{linear} in the matrix dimensions.

Because the search is local, the ACA algorithm can miss important parts of $\mat{A}$; practical implementations therefore employ suitable restart or
safeguard mechanisms.
Rook pivoting~\cite{neal1992geometric} yields a variation of ACA that alternates between a column search and a row search until the candidate pivot is maximal in magnitude in both its row and its column.
The complexity of rook pivoting depends on the number of scans needed to identify a suitable pivot.
In the worst case, $2\cdot \min \{m,n\}+1$ scans may be necessary, resulting in an overall cost for CUR approximation by rook pivoting that grows with proportional to $m\cdot n$.
However, in practice, the number of scans needed is often a small constant 
to identify each pivot, resulting in the same $\order((m+n)r^2)$ runtime as standard ACA.

\paragraph{Random pivoting.}
Randomly pivoted LU (RPLU) of
\citet{gilles2026low} draws the pivot $(i_t,j_t)$ in step $t$ randomly, according to the distribution 
$
p_{ij} = |\mat E_{t-1}(i,j)|^2 / \|\mat E_{t-1}\|_{\rm F}^2.
$
The Frobenius norm error is bounded by $2^r \|  \mat{A}-\lowrank{\mat A}_r \|_{\mathrm F}$ in expectation, and \citet{divan2026convergence} prove convergence results for matrices with steadily decaying singular values.
RPLU is an attractive option when sampling pivots from this
distribution (or a suitable approximation) is substantially
cheaper than locating the largest residual entry.

\subsubsection{Rank-adaptive algorithm with a small sketch: strong a posteriori guarantee}

When the input rank $r$ is unknown and a prescribed accuracy $\tau$ is desired, IterativeCUR~\cite{pritchard2025iterativecur} chooses $r$ adaptively. The algorithm only requires a single small sketch of $\mat{A}$, and adds a few indices to $\iset,\isettwo$ at each iteration by pivoting a sketch of the residual $\mat{A}-\mat{CUR}$, as in~\cref{sec:sketchy-pivoting}, until the sketch-based estimate indicates that $\|\mat{A}-\mat{C}\mat{U}\mat{R}\|_{\rm F}\leq \tau$ is satisfied for a prescribed tolerance $\tau$. 
IterativeCUR is often significantly faster than randomized SVD or its rank-adaptive variant~\cite{yu2018efficient}, while usually outputting a comparable rank. Aside from the small sketch, the complexity is $\order((m+n)r^2)$. 
Computing the small sketch is also cheap, requiring only $\order(mn)$ operations when $\mat{A}$ is dense and the sketch size is $\order(1)$.
Its accuracy is estimated using the tools from~\cref{sec:error-estimation}, so the a posteriori error guarantee is very reliable. However, a priori convergence theory for IterativeCUR is largely lacking. 

\subsubsection{Linear-time heuristic algorithms}

There is an even faster, and arguably simpler, approach to CUR: choose the row and column indices uniformly at random. The clear advantage is that in $\order((m+n)r)$ operations (just reading the $r$ rows and columns), one has a CUR approximation. 
This heuristic can obviously fail~\cite{pan2017superfast}; nonetheless, under incoherence assumptions on the leading singular vectors,  oversampled and regularized variants admit
spectral-norm error bounds with high probability~\cite{chiu2013sublinear}.
Verifying such assumptions is often infeasible or inefficient.

Another popular (and more powerful) heuristic~\cite{li2015butterfly} is to start with randomly selected rows, use these to find a good subset of columns (e.g., by CPQR), together with some additional random samples, and repeat the process to both improve and increase the subsets. This comes with guarantees under milder conditions~\cite{cortinovis2025sublinear}; nevertheless, some assumption on incoherence or matrix structure is  necessary for a method based on uniform random samples to be guaranteed to succeed. Note that these guarantees concern the
optimal core $\mat U_{\mathrm{opt}}$, whose computation generally requires access to the full matrix.

\subsubsection{Illustration}

Let us provide some practical guidelines for computing CUR. 
In the entry access model, where one wishes to avoid reading the whole matrix, tools based on ACA are highly efficient, but one needs to be aware that convergence is not guaranteed.

When matrix--vector multiplications with $\mat{A}$ are feasible, we recommend IterativeCUR when the rank is unknown, and a prescribed accuracy is desired. When the rank is known, we recommend the algorithms in 
\cref{sec:cur-rowspace}, namely Osinsky's method or ARP, as they come with optimal a priori guarantees. 

\begin{figure}
  \begin{minipage}[t]{0.49\hsize}
  \includegraphics[width=0.99\textwidth]{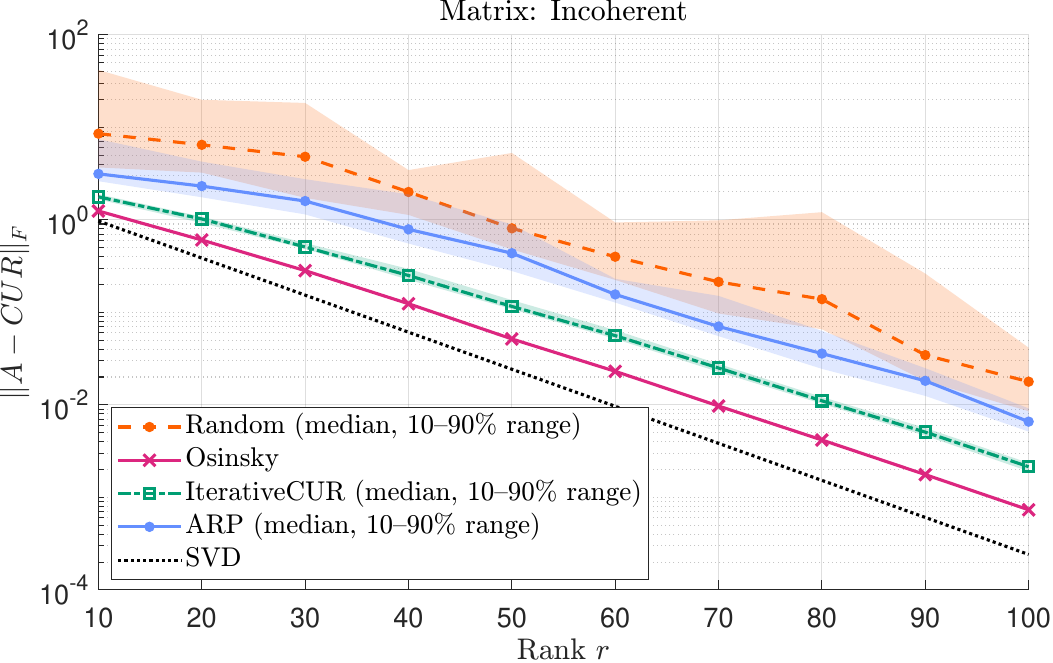}
  \end{minipage}    
\begin{minipage}[t]{0.49\hsize}    \includegraphics[width=0.99\textwidth]{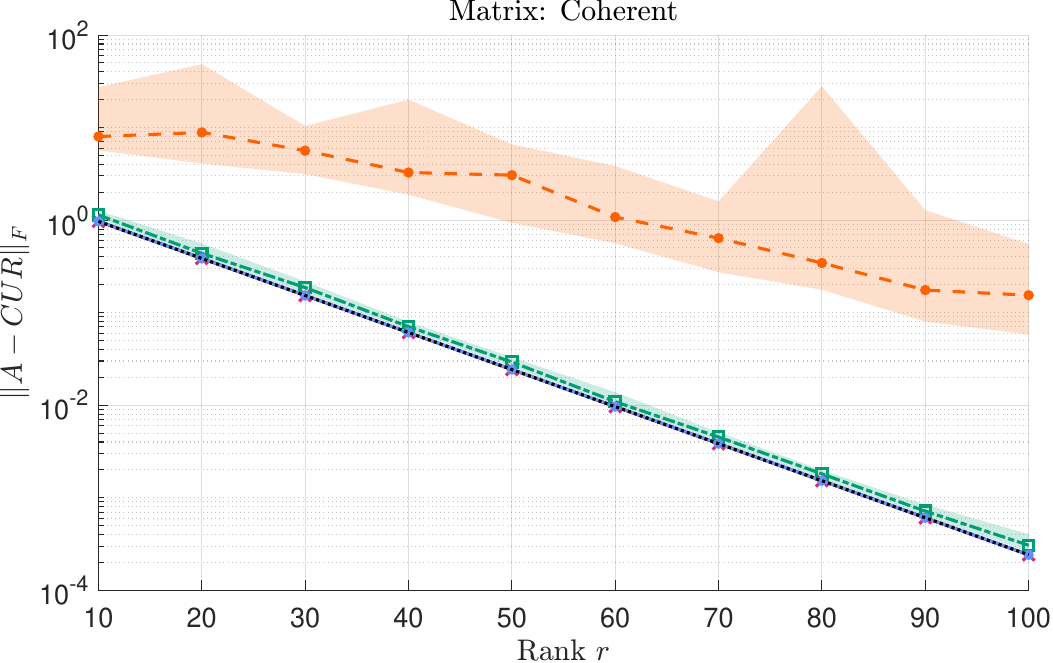}
  \end{minipage}   
    
    \caption{Frobenius norm error for CUR approximations of varying ranks for $\mat{A}\in\mathbb{R}^{1000\times 10000}$, with Osinsky-CUR, ARP-CUR, 
    and IterativeCUR. Uniform sampling (random) and the truncated SVD are included
for comparison.
    Left: incoherent example. Right: coherent example.}
    \label{fig:cur_exp}
\end{figure}

We illustrate these algorithms in \cref{fig:cur_exp}, where we take $\mat{A}=\mat{U\Sigma V}^\top$ so that the singular values are geometrically decaying. We examine two cases distinguishing coherent and incoherent problems; coherence is often associated with the difficulty of the column/row subset selection problem~\cite{mahoney2011randomized}. 
 In the incoherent example we take $\mat{U},\mat{V}$ to be Haar distributed. In the coherent example 
$\mat{U},\mat{V}$ are the orthogonal factors of $\Id+10^{-5}\mat{G}$, where 
$\mat{G}$ has iid standard Gaussian entries.
In Osinsky-CUR and ARP-CUR we take $\mat{V}$ to be the exact dominant right singular vectors; when $\mat{V}$ is obtained by a sketch, the left figure looks qualitatively similar, while the right figure will have a gap between SVD and the rest. For reference we also test choosing the indices uniformly at random; we do not show ACA as it can easily be made to fail by, e.g., taking $\mat{A}$ to be block diagonal~\cite{martinsson2019fast}. 
For the randomized methods (CUR-ARP, random, and IterativeCUR) we indicate the $10\%,90\%$ percentiles in addition to the median values. 
We observe that CUR-Osinsky and IterativeCUR are extremely robust with respect to varying runs; the variation is hardly visible. 
It is hardly surprising that choosing pivots uniformly at random does not do well, especially in the coherent example.

Let us highlight that, in the  coherent example here, all the methods except random perform exceptionally well---their errors are almost indistinguishable from that of truncated SVD! This is a common phenomenon; when there is a set of $r$ rows and columns that overwhelmingly span the dominant ranges (e.g. the dominant singular vectors have $r$ large leverage scores relative to the rest), it is extremely important to find them (and failure to do so can result in a disaster), but the reward is enormous. For incoherent matrices, choosing reasonably good $\iset,\isettwo$ is relatively easy, but choosing them (near-)optimally yields limited improvements.

\subsection{CUR with an evolving matrix}\label{sec:continuouscur}

One powerful advantage of the CUR approximation is that it is defined using the index sets $\iset$ and $\isettwo$.
Identifying good index sets may require some work, but once we have computed them we may reuse them if the matrix undergoes changes.

\paragraph{Row or column additions.}
A first type of change a matrix can undergo is addition of rows or columns.
Begin with an initial matrix $\mat{A}$, and suppose we have computed an accurate CUR approximation $\mat{A} \approx \mat{C}\mat{U}\mat{R}$ described by index sets $\iset,\isettwo$.
For simplicity, we consider the column addition case and assume columns are added to the end, i.e., $\tildebold{A} = \flatonebytwo{\mat{A}}{\mat{B}}$.
Without doing any additional work, we immediately obtain a CUR approximation to the augmented matrix
\begin{equation} \label{eq:augmented-CUR}
    \tildebold{A} \approx \mat{C}\mat{U}\tildebold{R} \quad \text{for } \tildebold{R} = \tildebold{A}(\iset,:).
\end{equation}
Pictorially, 
\[   
\begin{tikzpicture}[baseline=-1mm]
  \def\stripe{\colsize/10}
  \def\barh{\colsize/10}
 \def\grayext{\dimexpr1\colsize/5\relax}

  \node[draw, fill=gray!10,
        minimum width=\colsize+\stripe+\grayext,
        minimum height=\colsize,
        anchor=west] (A) at (0,0) {$\phantom{\Scale[2]{\tilde A}}$};

  \draw[fill=blue!20,draw=none]
    ($(A.south west)+(0,4mm)$)
    rectangle
    ($(A.south west)+(\colsize+\stripe+\grayext-0mm,4mm-\barh)$);

  \draw[fill=blue!20,draw=none]
    ($(A.south west)+(0,8mm)$)
    rectangle
    ($(A.south west)+(\colsize+\stripe+\grayext-0mm,8mm-\barh)$);

  \draw[fill=blue!20,draw=none]
    ($(A.south west)+(0,18mm)$)
    rectangle
    ($(A.south west)+(\colsize+\stripe+\grayext-0mm,18mm-\barh)$);

  \foreach \x in {0.03,0.25,0.6,0.85}{
    \draw[fill=red!20,draw=none]
      ($(A.south west)+(\x*\colsize,0)$)
      rectangle
      ($(A.south west)+(\x*\colsize+\stripe,\colsize)$);
  }

\node[draw, 
        minimum width=\colsize+\stripe+\grayext,
        minimum height=\colsize,
        anchor=west] at (0,0) {$\Scale[2]{\tilde A}$};
  
\end{tikzpicture}
\approx 
   \begin{tikzpicture}[node distance=1mm,baseline=-1mm]
  \def\stripe{\colsize/10}
  \def\barh{\colsize/10}
 \def\grayext{\dimexpr1\colsize/5\relax}
       \node[draw, rectangle, solid, fill=red!20, minimum width=\rowsize, minimum height=\colsize] (U) {$\Scale[1.2]{C}$};
       \node[draw, rectangle, solid, fill=yellow!30, minimum width=\rowsize, minimum height=\rowsize,right =of U] (S) {$\Scale[1]{U}$};
       \node[draw, rectangle, solid, fill=blue!20, minimum width=\colsize+\stripe+\grayext, minimum height=\rowsize,right =of S] (V) {$\Scale[1]{\tilde R}$};
   \end{tikzpicture}, 
\]
Observe that computing the updated CUR approximation only requires reading $|\iset|$ rows of the new columns, so we can construct the augmented approximation \emph{without even reading most of the new entries in $\tildebold{A}$!}

The augmented CUR approximation \cref{eq:augmented-CUR} will continue to be accurate provided
that the new columns $\mat{B}$ are nearly in the span of $\mat{C}$.
To see this, one may employ the DEIM error bound~\cref{eq:deimerror} to obtain
\begin{equation}\label{eq:curextrap_deim}
    \norm{\smash{\tildebold{A}-\mat{C}\mat{U}\tildebold{R}}}_{\rm F}
    \leq
    \norm{\smash{\mat{Q}_{\mat{C}}(\iset,:)^{-1}}} \cdot 
    \norm{\smash{\tildebold{A}-\mat{C}\mat{C}^\dagger\tildebold{A}}}_{\rm F}.    
\end{equation}
Here, $\mat{Q}_{\mat{C}} = \orth(\mat{C})$.
We see that, provided we have chosen a good row index set $\iset$ (so that $\norm{\smash{\mat{Q}_{\mat{C}}(\iset,:)^{-1}}}$ is small), the quality of the extended CUR approximation depends on how close the augmented matrix $\tildebold{A} = \flatonebytwo{\mat{A}}{\mat{B}}$ is to being spanned by the selected columns $\mat{C}$.

\paragraph{Parameter-dependent matrices.}
Another case where we can reuse index sets in a CUR decomposition is for a parameter-dependent matrix $\mat{A}(\vec{p})$ \cite{park2025low}.
The idea is that, we can compute the index sets $\iset,\isettwo$ once for a fixed value $\vec{p}_\star$ of the parameter $p$.
Then, provided $\mat{A}(\vec{p})$ ``changes slowly'' as a function of $\vec{p}$, we can reuse the same index sets $\iset$ and $\isettwo$ to obtain CUR approximations to any matrix $\mat{A}(\vec{p})$ in the parametric family.
Critically, building an approximation at each new parameter value only requires reading $m \cdot |\isettwo| + n \cdot |\iset|$ entries of each new matrix $\mat{A}(\vec{p})$.
Again, we avoid reading the entire matrix.
An error estimator from~\cref{sec:error-estimation} is crucial to guarantee its success.
These examples highlight the extrapolation power of CUR.

\paragraph{CUR for approximating functions.}

The ideas underlying CUR approximation have a long history
in interpolation and function approximation.
For example, consider the classical subject of (Chebyshev) polynomial interpolation~\cite[Ch.~2]{trefethen2019approximation}, wherein a function $f:[-1,1]\rightarrow \mathbb{R}$ is approximated by a polynomial $p(x) = \sum_{i=0}^{r-1} c_iT_i(x)$, such that $f(x_i) = p(x_i)$ at $x_i = \cos(i\pi/(r-1))$ for $i=0,1,\ldots,r-1$. 
This can be viewed as a 
rank-$r$
CUR approximation 
$\mat{A}:=[T_0, T_1, \ldots T_{r-1}, f]\approx \mat{CUR}$, where the columns of $\mat A$ and $\mat C$ are functions. 
Here, 
the $i$th column of 
$\mat{C}$ is the $(i-1)$th Chebyshev polynomial $T_{i-1}$, 
$\mat{R}$ is a $r$ row-subset of $\mat{A}$, namely its evaluation at $x_0,\ldots, x_{r-1}$, and 
$\mat{U}^{-1}$ is 
the $r\times r$
intersection matrix $(\mat{U}^{-1})_{ij}=T_{j-1}(x_{i-1})$. 
Here $\mat{A}$ 
and $\mat{CUR}$ are 
equal in all but the last column, which gives the polynomial $p$ that interpolates $f$ in $x_0,\ldots, x_{r-1}$. Note that this argument holds for any univariate function $f$. 

Another example is  low-rank approximation of a bivariate function $f(x,y)\approx 
\sum_{i=1}^r g_i(x)\ell_i(y)$ (note that $g_i,\ell_i$ are univariate; the right-hand side is a bivariate function of rank at most $r$). 
Following the idea of ACA~\cite{bebendorf2000approximation,townsend2013extension},
such an approximation can be constructed by initializing
the residual $h(x,y)=f(x,y)$.
At the $i$th iteration, a pivot point $(x_i,y_i)$ is chosen, usually where the current $h$ is maximized, and $h$ 
undergoes a rank-1 update:
defining a rank-1 function 
$h_i(x,y):=\frac{1}{h(x_i,y_i)} h(x,y_i)h(x_i,y)$, update
 $h(x,y):= h(x,y)-h_i(x,y)$. 
 We repeat this until $h$ is small in the domain of interest. 
 After $r$ steps this results in a rank-$r$ approximation $f\approx \sum_{i=1}^r h_i$, which 
is a pivoted LU approximation of $f$, and therefore 
equivalent to CUR. In Chebfun2~\cite{townsend2013extension}, the univariate factors $g_i$ and $\ell_i$ are approximated by Chebyshev polynomials of adaptively chosen degrees.

\section{Column \Nystrom approximation 
} \label{sec:column-nystrom}
    
When $\mat{A}$ is symmetric, it makes sense to specialize the CUR approximation to the case where the column and row sets coincide (that is, $\mat{R} = \mat{C}^\tp$ and $\iset=\isettwo$).
The resulting approximation is known as the \emph{column \Nystrom approximation}:
\begin{equation} \label{eq:column-nystrom}
    \mat{A}\langle \iset \rangle \coloneqq \mat{A}(:,\iset) \mat{A}(\iset,\iset)^\pinv \mat{A}(\iset,:).
\end{equation}
It exactly coincides with the general notion of \Nystrom approximation introduced in \cref{sec:lowrank-nystrom} with test matrix $\mat{\Omega} = \Id(:,\iset)$.
While the formula \cref{eq:column-nystrom} is sensible for any symmetric matrix $\mat{A}$, there are issues of stability in the symmetric indefinite case, which we discuss below in \cref{sec:column-nystrom-indefinite}.
Except in that subsection, we shall assume throughout this section that $\mat{A}$ is \psd.

\subsection{Why column \Nystrom approximation?}

The column \Nystrom approximation is one of the most powerful forms of matrix approximation.
As with the CUR decomposition, the \Nystrom approximation is described by a list of $s$ indices $\iset$, and it may be generated in factored form $\mat{A}\langle\iset\rangle = \mat{F}\mat{F}^\tp$ by accessing only $sn$ entries of $\mat{A}$ and performing $\order(s^2n)$ arithmetic operations.
But the column \Nystrom approximation has a property that the CUR decomposition does not: it admits a near-optimal low-rank approximation algorithm using time linear in $n$, not even requiring us to read the entire matrix $\mat{A}$!
The crucial difference here is positive semidefiniteness.
For a general matrix, large entries can always ``hide'' at an arbitrary location, and we must inspect the entire matrix if we want to be assured that we have not missed one.
By contrast, the largest entry in absolute value of a \psd matrix
always lies on the diagonal, which---as we will see---allows us
to cheaply control the error of an approximation.
Algorithms that produce near-optimal column \Nystrom approximations in linear time will be detailed in \cref{sec:colum-nystrom-linear}.

\subsection{\Nystrom, Cholesky, and Schur}

In \cref{sec:CURandLU}, we saw that a CUR (cross) approximation with $|\iset| = |\isettwo| = s$ exactly coincides with a partial completely pivoted LU factorization using the row and column pivot sequences $\iset$ and $\isettwo$.
Specializing to the case where $\mat{A}$ is \psd and $\iset = \isettwo$, it follows that the column \Nystrom approximation $\mat{A}\langle \iset\rangle$ coincides with the output of a partial Cholesky decomposition with diagonal pivoting defined by $\iset$.
Thus, choosing a good index set $\iset$ for \Nystrom approximation is equivalent to selecting a good pivot sequence for the Cholesky decomposition.

Let us be a bit more explicit.
Introduce a permutation matrix $\mat{P}$ that brings the indices in $\iset$ to the first $s \coloneqq |\iset|$ positions, so that $\mat{P}\mat{A}\mat{P}^\tp$ has the block representation
\begin{equation*}
    \mat{P}\mat{A}\mat{P}^\tp = \twobytwo{\mat{A}(\iset,\iset)}{\mat{A}(\iset,\overline{\iset})}{\mat{A}(\overline{\iset},\iset)}{\mat{A}(\overline{\iset},\overline{\iset})}
\end{equation*}
Performing $s$ steps of Cholesky decomposition and assumin a zero pivot is never encountered, we obtain a Cholesky decomposition of the upper left $\mat{A}(\set{I},\set{I}) = \mat{L}_{\set{I}}^{\vphantom{\tp}}\mat{L}_{\set{I}}^{\tp}$ and a representation of $\mat{P}\mat{A}\mat{P}^\tp$ as a sum
\begin{equation*}
    \mat{P}\mat{A}\mat{P}^\tp = \underbrace{\twobyone{\mat{L}_{\set{I}}}{\mat{A}(\overline{\iset},\iset)\mat{L}_{\set{I}}^{-\tp}}\twobyone{\mat{L}_{\set{I}}}{\mat{A}(\overline{\iset},\iset)\mat{L}_{\set{I}}^{-\tp}}^\tp}_{\text{approximation } \mat{A}\langle \iset\rangle} + \underbrace{\twobytwo{\mat{0}}{\mat{0}}{\mat{0}}{\mat{A}(\overline{\iset},\overline{\iset})-\mat{A}(\overline{\iset},\iset)\mat{A}(\iset,\iset)^{-1}\mat{A}(\iset,\overline{\iset})}}_{\text{residual } \mat{A} - \mat{A}\langle \iset\rangle}.
\end{equation*}
The first term is the column \Nystrom approximation, represented in factored form using the leading $s$ columns of a Cholesky factorization of the permuted matrix $\mat{P}\mat{A}\mat{P}^\tp$.
As a consequence, we see that we can compute the \Nystrom approximation in factored form $\mat{A}\langle \iset\rangle = \mat{F}\mat{F}^\tp$ by performing $s$ steps of the pivoted Cholesky algorithm.
The second term is the residual, known in matrix theory as the Schur complement \cite{unknown2006schur}.
One of the fundamental results of matrix analysis is that the Schur complement of a \psd matrix remains \psd \cite[Thm.~1.12]{unknown2006schur}.
As such, the residual of the \Nystrom approximation $\mat{A} - \mat{A}\langle \iset \rangle$ is a \psd matrix.

\begin{remark}[History]
    The column \Nystrom approximation has an interesting history.
    The use of pivoted Cholesky factorizations to compute ``rank-revealing'' decompositions of \psd matrices dates at least to the 1970s \cite{dongarra1979linpack,higham1990analysis}.
    This approach received renewed attention in the early 2000s in the kernel machine learning community \cite{williams2000using,fine2001efficient}.
    The name ``\Nystrom'' was attached to the approximation by \citet{williams2000using}, drawing a connection to classical work by \citet{nystrm1930uber} on integral equations.
    Accuracy guarantees for \Nystrom approximation for random index sets $\iset$ were studied by \citet{drineas2005nystrom} and later by \citet{gittens2016revisiting,musco2017recursive}.
    At some point, the notion of \Nystrom approximation was generalized to allow for an arbitrary test matrix $\mat{\Omega}$; see, e.g.,~\cite{halko2011finding}.
\end{remark}

\subsection{The Gram correspondence: Redux} \label{sec:gramcorrespondeindex}

The \emph{Gram correspondence} (\cref{sec:lowrank-gram-correspondence}) provides a useful connection between projection approximations and \Nystrom approximations.
It has particularly interesting consequences in the index selection context.

Let $\mat{A}$ be a \psd matrix, and let $\mat{B}$ be any Gram square root, so that $\mat{A} = \mat{B}^\tp \mat{B}$.
The Gram correspondence states that the column \Nystrom approximation $\mat{A}\langle \iset\rangle$ is the Gram matrix of the interpolative approximation $\Bhat \coloneqq \mat{B}(:,\iset)\mat{B}(:,\iset)^\pinv \mat{B}$:
\begin{equation*}
    \mat{A}\langle \iset\rangle = \Bhat^\tp \Bhat.
\end{equation*}
Moreover, the errors of the approximations are equivalent: 
\begin{equation*}
    \tr( \mat{A} - \mat{A}\langle \iset\rangle ) = \norm{\mat{A} - \mat{A}\langle \iset\rangle}_* = \norm{\mat{B} - \Bhat}_{\rm F}^2 \quad \text{and} \quad \norm{\mat{A} - \mat{A}\langle \iset\rangle} = \norm{\mat{B} - \Bhat}^2.
\end{equation*}
In this sense, \Nystrom approximation of a \psd matrix $\mat{A}$ and interpolative decomposition of a general matrix $\mat{B}$ are equivalent.
This is an interesting algebraic fact, but its most important consequence is \emph{algorithmic}.
\actionbox{\textbf{Transference of algorithms principle.}
Any algorithm for computing an interpolative decomposition of a general matrix $\mat{B}$ gives rise to an algorithm for producing a \Nystrom approximation to the \psd matrix $\mat{A} = \mat{B}^\tp \mat{B}$.}
\noindent See \cref{tab:gram-correspondence} for a partial list of examples, and see \cite[Cor.~9.5]{epperly2025make} for more information on the transference principle.

\begin{table}[t]
  \centering
  \caption{List of algorithm pairs for computing an interpolative decomposition or column \Nystrom approximation that are equivalent under the Gram correspondence.}
  \label{tab:gram-correspondence}
  \begin{tabular}{@{}ll@{}}
    \toprule
    Interpolative decomposition & Column \Nystrom approximation \\ \midrule
    Partial pivoted QR & Partial pivoted Cholesky \\
    Length-squared sampling & Diagonal sampling \\
    Businger--Golub QR & Greedy pivoted Cholesky \\
    Randomly pivoted QR & Randomly pivoted Cholesky \\
    Volume-squared sampling & $s$-DPP sampling 
    \\ \bottomrule
  \end{tabular}
\end{table}

\subsection{How accurate can column \Nystrom approximations be?}

With a judicious choice of the index set $\iset$, column \Nystrom approximations can be nearly as accurate as the optimal rank-$r$ approximation.

\begin{theorem}[Fundamental theorem of (column) \Nystrom approximation] \label{thm:fundamental-nystrom}
    Let $\mat{A}$ be a \psd matrix.
    Then, for every $r\ge 1$ and $\varepsilon \in (0,r]$, there exists an index set $\iset$ of size $|\iset| = \lceil r/\varepsilon + r - 1\rceil$ such that
    \begin{subequations} \label{eq:ftn}
    \begin{equation} \label{eq:ftn-1}
        \tr(\mat{A} - \mat{A}\langle \iset\rangle) \le (1+\varepsilon) \tr( \mat{A} - \lowrank{\mat{A}}_r ).
    \end{equation}
    In particular, there exists a subset of size $|\iset| = r$ for which
    \begin{equation} \label{eq:ftn-2}
        \tr(\mat{A} - \mat{A}\langle \iset\rangle) \le (r+1) \tr( \mat{A} - \lowrank{\mat{A}}_r ).
    \end{equation}
    \end{subequations}
\end{theorem}

The bound \cref{eq:ftn-2} was first established, as such, by \citet[Thm.~1]{belabbas2009spectral}.
It also may be deduced from \cite[Thm.~1.3]{deshpande2006matrix} and the Gram correspondence.
The bound \cref{eq:ftn-1} follows from \cite[Thms~1.1 \& 1.2]{guruswami2012optimal} and the Gram correspondence; see also \cite[Thm.~6]{derezinski2021determinantal}.
The prefactor $r+1$ in \cref{eq:ftn-2} cannot be improved, as can be seen from the example $\mat{A} = (1+\eta) \Id - (r+1)^{-1}\onevec\onevec^\tp \in \R^{(r+1)\times (r+1)}$ for $\eta\downarrow0$, 
and there exists an example on which at least $\lceil r/\varepsilon \rceil$ are needed to achieve the guarantee in \cref{eq:ftn-1}; see \cite[Prop.~3.3]{deshpande2006matrix} and \cite[\S6]{guruswami2012optimal}.

Computationally, how can we identify subsets satisfying the bounds \cref{eq:ftn-1} and \cref{eq:ftn-2}?
One approach is given by $s$-DPP sampling, which was discussed briefly in \cref{sec:volume-sampling-algs}.
If the subset $\iset$ is drawn from the $s$-DPP with kernel $\mat A$,
both guarantees \cref{eq:ftn-1,eq:ftn-2} hold in expectation.
See \citet{derezinski2021determinantal} for a survey on $s$-DPP sampling in RNLA.
With some extra effort, volume sampling and DPP-type algorithms can be derandomized to achieve deterministic algorithms with the same guarantees \cite{deshpande2010efficient,guruswami2012optimal,cortinovis2020low,massei2022some}.
An alternative, more computationally appealing approach to achieving the guarantee \cref{eq:ftn-2} is by means of an ARP-type algorithm \cite[Alg.~5.1]{cortinovis2026adaptive}.
This requires an orthonormal matrix $\mat V\in\R^{n\times r}$
of the dominant eigenspace of $\mat A$, which can be approximated from a randomized sketch $\mat A\mat\Omega$.

It should be stressed that most of the algorithms that achieve the sharpest error guarantees \cref{eq:ftn-1,eq:ftn-2} require processing the entire matrix at cubic cost.
In some cases, they may, in principle, be accelerated to near-linear time using a fast $s$-DPP sampling algorithm.
However, even with the fastest-known $s$-DPP samplers, all of these approaches appear to be significantly slower than the best alternative approaches discussed below in \cref{sec:colum-nystrom-linear}.

\paragraph{Numerical stability.}
Column \Nystrom also exhibits the potential pitfall with an ill-conditioned pseudoinverse $\mat{A}(\iset,\iset)^\pinv $. 
If good pivots are chosen (specifically, locally maximum-volume pivots~\cite{damle2026estimating}),
\citet{bucci2025numerical} establish stability for a carefully implemented, truncated
column \Nystrom approximation; see also~\cite{carson2024single}.
While $\varepsilon$-truncation is needed in theory for the stability analysis, the procedure is  observed to be stable in practice without truncation.
An alternative approach is to employ a shift~\cite{tropp2017fixed} to ensure Cholesky is performed on a numerically nonsingular matrix; see \cref{sec:lowrank-nystrom}.

\subsection{Linear-time algorithms} \label{sec:colum-nystrom-linear}

This section describes algorithms for computing an index set $\iset$ comprising a column \Nystrom approximation that run in \emph{linear time}, here defined as running in $\order(n)$ operations to produce an index set of constant size.\footnote{In computer science, these algorithms in this subsection are often described as \emph{sublinear}, since the runtime is a sublinear function of the total size $n^2$ of the input. We assume that an entry of $\mat A$ can be evaluated in $\order(1)$ operations.}

\subsubsection{Uniform random selection and diagonal sampling}

The simplest algorithm for computing an index set for a column \Nystrom approximation is simply to choose $\iset$ to be a uniformly random subset of $s$ indices, sampled without replacement.
This basic approach, dating back to the seminal work of \citet{williams2000using}, is surprisingly effective on many problems; see \cite{gittens2011spectral} for an analysis.
When the diagonal entries of $\mat{A}$ have vastly different scales, uniform sampling can be poor.
In this case, a simple modification is to sample indices $i$ with probabilities proportional to the diagonal entries, possibly raised to some power $p$ \cite{drineas2005nystrom}.

Uniform sampling can work surprisingly well on some problems, but it is easy to find examples in theory and practice where it fares very badly; see, e.g., \cite[Thm.~C.3(b)]{chen2025randomly}.
To perform better, we must use more sophisticated sampling distributions or adaptive selection procedures.

\subsubsection{Ridge leverage score sampling}

A more effective sampling approach to identifying an index set defining a good \Nystrom approximation is \emph{ridge leverage score sampling}.
The $\lambda$-ridge leverage scores of a \psd matrix $\mat{A}$ are the diagonal entries of $\mat{A}(\mat{A} + \lambda \Id)^{-1}$ \cite{elalaoui2015fast}.
By tuning the parameter $\lambda > 0$ appropriately and using ridge leverage score-weighted sampling to pick the index set $\iset$, we obtain, with high probability, a column \Nystrom approximation  comparable (in terms of the trace error) with the best rank-$r$ approximation, using $|\iset| = \order(r\log r)$ \cite[Thm.~17]{musco2017recursive}.
Remarkably, despite the ridge leverage scores being a nontrivial rational function of the matrix entries, several algorithms exist that approximately sample from the ridge leverage score distribution in linear time \cite{musco2017recursive,rudi2018fast}.

Ridge leverage score sampling is a mathematically elegant approach, and it can give substantial improvements over uniform and diagonal-power sampling.
As limitations, the available algorithms are somewhat more expensive than Cholesky decomposition-based approaches, and there are examples where they produce approximations of much lower accuracy than Cholesky approaches for a fixed budget of entry accesses \cite[Tab.~1 \& Fig.~1]{chen2025randomly}.
Therefore, we generally recommend against ridge leverage score sampling for general-purpose \Nystrom approximation.

\subsubsection{Greedy pivoted Cholesky}

So far, we have seen algorithms that use random sampling to select an index set.
An alternative approach utilizes the connection between \Nystrom approximation and Cholesky decomposition, building an index set by selecting a new index at each step of the pivoted partial Cholesky algorithm.
The simplest approach, known as either complete, diagonal, or greedy pivoting, selects the pivot entry at each step as the largest diagonal entry of the current residual.
This basic approach is classical, dating at least to the 1970s \cite{dongarra1979linpack}.
It is equivalent to the Businger--Golub CPQR algorithm under the Gram correspondence (\cref{sec:gramcorrespondeindex}).
Other important references are \citet{higham1990analysis,fine2001efficient,harbrecht2012lowrank}.

The properties of greedy pivoted Cholesky are much the same as Businger--Golub QR, discussed above in \cref{sec:classical-pivoting}.
The algorithm often works terrifically, but it can fail catastrophically on pathological ``Kahan-type'' examples.
More worryingly, the algorithm can produce meaningfully sub-optimal results on matrices $\mat{A}$ assembled from datasets with outliers.
See \citet[\S4]{chen2025randomly} for examples.

\subsubsection{Randomly pivoted Cholesky} \label{sec:rpcholesky}

We have now seen two extremal strategies for index selection on a \psd matrix.
The first, uniform sampling, discards information in the diagonal---exploring all possible choices but not exploiting information from the diagonal.
The second, greedy selection, exploits the diagonal information but never explores, always deterministically using the largest diagonal entry as a pivot.
As is often the case in life, good things happen when one takes a middle road between these extremes, balancing exploration and exploitation.

The randomly pivoted Cholesky algorithm (\RPCholesky)
performs a pivoted partial Cholesky decomposition of the matrix $\mat{A}$.
Begin with the trivial approximation $\Ahat \gets \mat{0}$ and initialize the residual $\mat E\gets\mat A$.
Conceptually, each step of \RPCholesky draws a random pivot
\begin{subequations}
\begin{equation}\label{eq:rpcholesky-1}
    \prob \{ i_t = j \mid \mat{E} \} = \mat{E}(j,j)/\tr(\mat{E})
\end{equation}
and updates the low-rank approximation $\Ahat$,
\begin{equation}\label{eq:rpcholesky-2}
    \Ahat \gets \Ahat + \mat{E}(:,i_t)\mat{E}(i_t,:)/\mat{E}(i_t,i_t),
\end{equation}
as well as the residual,
\begin{equation} \label{eq:rpcholesky-3}
    \mat{E} \gets \mat{E} - \mat{E}(:,i_t)\mat{E}(i_t,:)/\mat{E}(i_t,i_t).
\end{equation}
\end{subequations}
For efficiency, the approximation $\Ahat$ should be stored in factored form $\Ahat = \mat{F}\mat{F}^\top$, and columns of the residual matrix $\mat{E} = \mat{A} - \Ahat$ can be generated on the fly using the entries of $\mat{A}$ and the factor $\mat{F}$.
See \cref{alg:rpcholesky} for pseudocode.
With this optimized implementation, the \RPCholesky algorithm requires only $(s+1)n$ entry accesses, $\order(sn)$ storage, and $\order(s^2n)$ arithmetic operations to select an index set of size $|\iset| = s$.
The trace error provides a cheap stopping criterion for \cref{alg:rpcholesky}, as it cam be computed at negligible cost by summing the diagonal entries of the residual: 
\[
\vec{d}(1) + \vec{d}(2) + \cdots + \vec{d}(n) = 
\tr(\mat A-\mat F\mat F^\tp)
    =\norm{\mat A-\mat F\mat F^\tp}_*.
\]

\begin{algorithm}[t]
\caption{Randomly pivoted Cholesky (\RPCholesky)} \label{alg:rpcholesky}
\begin{algorithmic}[1]
\Require Psd matrix $\mat{A} \in \R^{n\times n}$ and rank $s$
\Ensure Pivot set $\iset = \{i_1,\ldots,i_s\}$ and factored approximation $\mat{A}\langle\iset\rangle = \mat{F}\mat{F}^\tp$
\State Initialize $\mat{F} \gets \mat{0}_{n\times k}$ and the residual diagonal $\vec{d} \gets \diag(\mat{A})$
\For{$t= 1, \dots, s$}
    \State Draw $i_t$ with probability $\vec{d}(i_t)/\sum_{j=1}^n \vec{d}(j)$ \Comment{Sample from the diagonal}
    \State Evaluate the residual column $\vec{g} \gets \mat{A}(:,i_t) - \mat{F}(:,1:t-1)\,\mat{F}(i_t,1:t-1)^\tp$
    \State Update the factor matrix $\mat{F}(:,t) \gets \vec{g} / \sqrt{g_{i_t}}$
    \State Update $\vec{d} \gets \max\{\vec{d} - |\mat{F}(:,t)|^2, 0\}$ \Comment{Entrywise maximum to clip roundoff}
\EndFor
\State \Return $\iset$ and $\mat{F}$
\end{algorithmic}
\end{algorithm}

Let us emphasize: The \RPCholesky algorithm ``just'' builds a column Nystr\"om via a partial pivoted Cholesky decomposition, just like the greedy pivoted Cholesky method discussed in the previous subsection.
However, \RPCholesky is often more effective than greedy selection, as the random pivot selection rule \cref{eq:rpcholesky-1} prioritizes the selection of larger diagonal elements without \emph{always} selecting the largest diagonal element.

\begin{remark}[History]
    The \RPCholesky algorithm may be derived from the RPQR algorithm (\cref{sec:random-pivoting}) by the Gram correspondence (\cref{sec:gramcorrespondeindex}).
    The paper \cite{musco2017sublinear} appears to be the first work to explicitly mention \RPCholesky as a linear-time algorithm for \psd low-rank approximation.
    The algorithm was revisited, numerically evaluated, and re-analyzed by \citet{chen2025randomly}, with further analysis and implementation improvements developed in \cite{epperly2023kernel,epperly2025embrace,epperly2026new,divan2026convergence}.
    The first author's thesis \cite{epperly2025make} provides a detailed treatment of \RPCholesky, variants, and connections to other algorithms.
\end{remark}

\paragraph{Accelerated \RPCholesky.}

The \RPCholesky algorithm can be made over an order of magnitude faster in practice by using a blocked variant of the algorithm instead of \cref{alg:rpcholesky}.
The method of choice is usually the \emph{accelerated \RPCholesky method}, developed in \cite{epperly2025embrace}.
Using a block rejection sampling approach, accelerated \RPCholesky produces the \emph{same random distribution of pivots} while leveraging block matrix computations.
We recommend accelerated \RPCholesky as a drop-in replacement for the standard implementation for almost all uses.
See \cite{epperly2025embrace} for discussion of the algorithm and comparison to alternatives.

\paragraph{Analysis.}
How effective is the \RPCholesky algorithm?
One way of quantifying its effectiveness is by comparison to the fundamental theorem of \Nystrom approximation (\cref{thm:fundamental-nystrom}), which states there exists a \Nystrom approximation with an index set of size $|\iset| = \lceil r/\varepsilon + r - 1\rceil$ that is within a $1+\varepsilon$ factor of the trace error of the best rank-$r$ approximation.
Here is one representative result.

\begin{theorem}[\RPCholesky] \label{thm:rpcholesky}
    Let $\mat{A}$ be a \psd matrix, fix $r\ge 1$ and $\varepsilon \in (0,1)$, and let $\iset$ denote the index set generated by $s$ steps of \RPCholesky.
    Then \RPCholesky achieves the guarantee
    \begin{equation} \label{eq:r-eps}
        \expect \tr(\mat{A} - \mat{A}\langle\iset\rangle) \le (1+\varepsilon) \tr(\mat{A} - \lowrank{\mat{A}}_r)
    \end{equation}
    in $s = r/\varepsilon + \order\!\big(r \sqrt{\log r} + r\log(1/\varepsilon)\big)$ steps.
\end{theorem}

This result was obtained by GPT 5.6, and its proof is documented in \cite{epperly2026new}.
The complexity can be improved to $s = \order(r/\varepsilon + r\log \log r)$ using the techniques of \cite{divan2026convergence}; see \cite{epperly2026adaptiveb}.

The paper \cite{chen2025randomly} was the first to prove guarantees similar to \cref{thm:rpcholesky} for \cref{alg:rpcholesky}, though results for a blocked version of \RPCholesky can be deduced from earlier works \cite{deshpande2006adaptive,deshpande2006matrix} and the Gram correspondence.
\citet{epperly2025embrace} present analysis for the accelerated \RPCholesky algorithm, \citet{epperly2026new} proves \cref{thm:rpcholesky} and provides an error bound in the spectral norm, and \citet{divan2026convergence} present convergence bounds for \RPCholesky applied to matrices with eigenvalues decaying at polynomial or exponential rates.

\subsubsection{Nuclear score selection} 
\label{sec:nuclear-selection}

The nuclear score selection strategy can be concisely described as follows:
\emph{At each step of pivoted partial Cholesky, always choose the index that causes the greatest single-step decrease in the trace error.}
Using the formula \cref{eq:rpcholesky-3}, we see that the decrease in trace for eliminating index $i$ is 
\begin{equation*}
    n_i(\mat{A}) \coloneqq \tr\left(\frac{\mat{A}(i,:)\mat{A}(:,i)}{\mat{A}(i,i)}\right) = \frac{\norm{\mat{A}(i,:)}^2}{\mat{A}(i,i)}.
\end{equation*}
This quantity, called the $i$th \emph{nuclear score} of $\mat{A}$, is the ratio of the $i$th squared row norm to the $i$th diagonal entry.
Notice that while the nuclear selection strategy is \emph{locally optimal} in the sense it always selects the index that causes the greatest single-step decrease in the trace error, it typically will not identify the \emph{globally optimal} index set that minimizes the trace error $\tr(\mat{A} - \mat{A}\langle \iset\rangle)$ over all $s$-element index sets $\iset$.

As with other column \Nystrom algorithms, there is a natural variant of nuclear score selection to select indices for computing an interpolative decomposition of a general rectangular matrix.
We have elected to focus only on the \psd variant in this work.
The nuclear score selection strategy may have first been studied by \citet{farahat2011novel}, with subsequent work by \citet{altschuler2016greedy,sood2025statistical,fornace2026column}.

As presented so far, the nuclear score selection algorithm is \emph{not} a linear-time algorithm, as computing the nuclear scores requires computing the row norms which requires $\order(n^2)$ operations.
However, as \citet{fornace2026column} noticed, if we are interested in performing index selection on a matrix $\mat{A}$ that possesses a fast matrix--vector product $\vec{\omega} \mapsto \mat{A}\vec{\omega}$, then we can use a row-norm estimator (\cref{sec:row-norm-estimation}) to efficiently estimate the nuclear scores.
In this setting, nuclear score selection can also be implemented as a linear-time algorithm, and \citet{fornace2026column} present numerical evidence that it can outperform other column \Nystrom algorithms in terms of the quality of the index set it generates.

\subsubsection{Comparison}

\Cref{tab:column-nystrom-compare} provides a comparison of the different linear-time column \Nystrom methods on twenty kernel matrices from data sets from the LIB-SVM, OpenML, and UCI repositories, reproducing the experimental setup from \cite[Tab.~1]{chen2025randomly}.
Each dataset is subsampled to $n= 10^4$ points, the data is standardized, and we use a squared-exponential kernel with bandwidth $\sqrt{d}$, where $d$ is the number of features in each dataset.
We report the median error over ten trials.

\begin{table}[t]
  \centering
  \scriptsize
  \begin{tabular}{lllllll}
    \toprule
    Dataset & {Unif} & {RLS} & {Greedy} & {RPChol} & {Nuclear} & {Opt} \\
    \midrule
    sensit\_vehicle          & 1.57e-01 & 1.40e-01 & 2.07e-01 & 1.37e-01 & 1.24e-01 & 8.77e-02 \\
    yolanda                  & 1.47e-01 & 1.39e-01 & 2.08e-01 & 1.36e-01 & 1.23e-01 & 8.42e-02 \\
    YearPredictionMSD        & 1.30e-01 & 1.20e-01 & 1.73e-01 & 1.16e-01 & 1.01e-01 & 6.79e-02 \\
    w8a                      & 1.41e-01 & 1.17e-01 & 1.91e-01 & 1.05e-01 & 8.04e-02 & 6.09e-02 \\
    MNIST                    & 1.21e-01 & 1.10e-01 & 1.67e-01 & 1.06e-01 & 9.35e-02 & 5.83e-02 \\
    jannis                   & 1.10e-01 & 1.10e-01 & 1.28e-01 & 1.09e-01 & 1.01e-01 & 5.45e-02 \\
    HIGGS                    & 6.72e-02 & 6.29e-02 & 8.51e-02 & 6.07e-02 & 5.37e-02 & 2.96e-02 \\
    connect-4                & 5.90e-02 & 5.05e-02 & 6.48e-02 & 4.79e-02 & 4.08e-02 & 2.25e-02 \\
    volkert                  & 5.32e-02 & 4.41e-02 & 5.92e-02 & 4.17e-02 & 3.63e-02 & 1.99e-02 \\
    creditcard               & 4.83e-02 & 3.70e-02 & 5.77e-02 & 3.03e-02 & 2.34e-02 & 1.31e-02 \\
    Medical-Appointment      & 1.73e-02 & 1.43e-02 & 1.92e-02 & 1.29e-02 & 1.05e-02 & 4.69e-03 \\
    sensorless               & 1.22e-02 & 7.62e-03 & 8.70e-03 & 5.78e-03 & 4.64e-03 & 2.11e-03 \\
    ACSIncome                & 9.73e-03 & 5.67e-03 & 8.35e-03 & 4.05e-03 & 3.04e-03 & 1.27e-03 \\
    Airlines\_DepDelay\_1M   & 4.16e-03 & 2.41e-03 & 2.64e-03 & 1.77e-03 & 1.38e-03 & 5.08e-04 \\
    covtype.binary           & 1.03e-02 & 2.25e-03 & 1.41e-03 & 1.05e-03 & 7.09e-04 & 2.97e-04 \\
    diamonds                 & 1.41e-03 & 2.07e-04 & 1.12e-04 & 5.92e-05 & 4.31e-05 & 1.30e-05 \\
    hls4ml\_lhc\_jets\_hlf   & 4.27e-04 & 7.11e-05 & 6.68e-05 & 4.38e-05 & 3.31e-05 & 1.04e-05 \\
    ijcnn1                   & 3.90e-05 & 3.06e-05 & 2.86e-05 & 2.14e-05 & 1.50e-05 & 4.67e-06 \\
    cod-rna                  & 4.75e-04 & 1.43e-05 & 9.19e-06 & 5.02e-06 & 3.75e-06 & 9.86e-07 \\
    COMET\_MC\_SAMPLE        & 3.58e-03 & 1.43e-07 & 1.10e-10 & 3.44e-11 & {---}    & 3.46e-12 \\
    \bottomrule
  \end{tabular}
  \caption{Relative trace error for rank-1000 column \Nystrom approximations computed by uniform random sampling, approximate ridge leverage score sampling using the RecursiveRLS algorithm \cite{musco2017recursive}, greedy pivoted Cholesky, \RPCholesky, and nuclear Cholesky applied to a test bed of 20 \psd kernel matrices of dimension $n = 10^4$.
  We report the median error over ten trials.
  The ``Opt'' column is the Eckart--Young best approximation.
  Our implementation of the nuclear Cholesky was unstable on the COMET\_MC\_SAMPLE example, so we do not report a result.}
  \label{tab:column-nystrom-compare}
\end{table}

Among all of the methods, the nuclear Cholesky method consistently achieves the best results, with \RPCholesky only slightly less accurate (a factor 1.5 at most).
We remind the reader that the nuclear Cholesky method is linear time only when we have a cheap way to access the row norms.
Therefore, on the basis of these experiments, we generally recommend the nuclear Cholesky method when row norms are available and the \RPCholesky method when they are not.
On the most extreme example (COMET\_MC\_SAMPLE), \RPCholesky produces approximations that are $10^8\times$ more accurate than uniform sampling, $4000\times$ more accurate than ridge leverage score sampling, and $3\times$ more accurate than greedy Cholesky.

\subsection{Applications} \label{sec:column-nystrom-applications}

The column \Nystrom approximation is a powerful tool, as it allows us to compute rank-$s$ approximations to a \psd matrix using only $\order(ns^2)$ work.
In particular, one can set the rank as large as $s = \order(\sqrt{n})$ and maintain a runtime proportional to the cost of merely scanning through the input matrix a single time.

One of the major applications of column \Nystrom approximation is for approximating \psd kernel matrices and covariance matrices that arise in kernel and Gaussian process algorithms for machine learning \cite{schlkopf2002learning,rasmussen2006gaussian}.
These algorithms reduce learning tasks on a dataset $\{\vec{x}_i\}_{i=1}^n$ to linear algebra and optimization problems involving a dense kernel matrix $\mat{A}$ with entries $\mat{A}(i,j) = \kappa(\vec{x}_i,\vec{x}_j)$.
The function $\kappa$ is chosen to be a positive definite kernel, which ensures that $\mat{A}$ is \psd.
Column \Nystrom methods are appealing in this context because they may accelerate methods to run in $\order(n)$ time, which compares favorably to the $\order(n^2)$ of generating the kernel matrix a single time.
Learning tasks that can be accomplished using kernels and Gaussian processes include regression \cite{rudi2015less}, clustering \cite{fowlkes2004spectral}, and nonlinear dimensionality reduction~\cite{zhang2010clustered}.
There is a truly enormous literature on the use of column \Nystrom approximation to accelerate kernel algorithms; see, e.g., \cite{williams2000using,fine2001efficient,fowlkes2004spectral,drineas2005nystrom,rudi2015less,musco2017recursive,rudi2018fast}.
The paper \cite{chen2025randomly} compares several column \Nystrom algorithms in the context of such learning algorithms.

Column \Nystrom approximations can also be useful for more purely linear algebraic tasks.
For instance, as we mentioned in \cref{sec:nystrom-precond}, column \Nystrom preconditioners constructed using, e.g., \RPCholesky can yield effective preconditioners for \pd linear systems with a few large outlying eigenvalues \cite{frangella2023randomized,diaz2026robust}. 
Column \Nystrom approximation can also be used for approximating eigenvalues.
Original work on eigenvalue approximation through index selection~\cite{williams2000using} suggested approximating the eigenvalues of $\mat{A}$ by the eigenvalues of $\frac{n}{s} \mat{A}(\iset,\iset)$.
This basic approach can fail badly.
It is often much more accurate to use the eigenvalues of the \Nystrom approximation $\mat{A}\langle\iset\rangle$ itself as approximations to the eigenvalues of $\mat{A}$.
These eigenvalues can be computed in $\order(ns^2)$ operations, and the resulting eigenvalue approximations often even outperform the classical Rayleigh--Ritz method~\cite{nakatsukasa2026finding}.

\subsection{Column \Nystrom approximation for indefinite matrices} \label{sec:column-nystrom-indefinite}
When $\mat{A}$ is symmetric but not \psd, we lose some of the advantages of the column \Nystrom approximation (e.g., the Gram correspondence).
However, it is still valuable to maintain symmetry in the approximation, maintaining $\mat{C}=\mat{R}^\tp$ in a CUR representation.
However, the CUR cross approximation with $\mat{U}=\mat{A}(\set{I},\set{I})^\pinv$ may give a poor approximation (note that FTCUR does not apply with $\mat{C}=\mat{R}^\tp$ even if $\mat{A}=\mat{A}^\tp$) or even break down; attempting a rank-1 approximation of~\cref{eq:curbadexample}
gives a simple example, wherein $\mat{A}(\iset,\iset)=\varepsilon$, which can be $0$.
Even when $\mat{A}(\iset,\iset)$ is invertible, the approximation can be of very poor quality and severe numerical issues can arise. 
For the optimal core $\mat U_{\mathrm{opt}}=\mat C^\dagger\mat A(\mat C^\tp)^\dagger$ and a suitable column selection, the bound~\cref{eq:curopt} holds with $\mat R=\mat C^\tp$, giving a suboptimality factor of $\sqrt{2(r+1)}$~\cite{osinsky2025close}.
Moreover, given a good column submatrix $\mat{C}$, one can obtain a good core matrix $\mat{U}$ without reading the whole $\mat{A}$ by solving the two-sided least-squares problem $\min_{\mat{M}}\|\mat{A}-\mat{C}\mat{M}\mat{C}^\tp\|_{\rm F}$ using an active least-squares method, keeping much of the appeal of the cross approximation.
See \cite{chen2026nystr} for details.

\section{A posteriori error estimation and quality control} \label{sec:error-estimation}

As this survey has emphasized, many randomized algorithms can be highly reliable.
For instance, randomized low-rank approximation algorithms like the randomized SVD (\cref{alg:rsvd}), randomized \Nystrom approximation (\cref{alg:randomized-nystrom}), and generalized \Nystrom approximation (\cref{alg:generalized-nystrom}) consistently produce rank-$s$ outputs that are nearly as good as the best rank-$r$ approximation, where $r < s-1$ is some target rank  (\cref{thm:rsvd,cor:nys,cor:gen-nys}).
But \emph{a priori} guarantees like these only certify that the algorithm produces an approximation \emph{nearly as good} as the best rank-$r$ approximation.
If we do not understand the error of the best rank-$r$ approximation, then we also do not understand the error level achieved by these algorithms.
Therefore, even when using highly reliable randomized algorithms, we have to use \emph{a posteriori error estimation techniques} and other diagnostics that can be used to understand the quality of the computed output at runtime.
Error estimation techniques are naturally also very useful if we employ algorithms that are less safe, for instance using uniform random sampling or random embeddings with aggressive parameter choices (\cref{fig:sketching-comparison}).\footnote{Turing Award winner Jack Dongarra refers to such approaches as ``responsibly reckless''.} 
All sins in algorithm design are forgiven if the quality of the output is certified with an a posteriori error estimate.

Before going on to discuss error \emph{estimation} techniques, let us first suggest an alternative to error estimation---namely, simply computing the error outright.
In many applications of, say, randomized low-rank approximation, we are interested in approximating a matrix that is stored in memory, and computing the trace or Frobenius-norm error is trivially easy.
The techniques described herein are designed for settings where we do not have such cheap, direct access.
For instance, we may be in a matrix-free setting where we wish to approximate a matrix $\mat{A}$ that is only accessible by matrix--vector multiplications $\vec{\omega} \mapsto \mat{A}\vec{\omega}$.
Or, the matrix $\mat{A}$ may be very large and we wish to avoid an additional pass over the matrix entries.
In these cases, the fast error estimation techniques described here may be useful.

This section will describe techniques for estimating the error of a randomized approximation $\mat{A} \approx \Ahat$.
Throughout, we shall let $\mat{E} \coloneqq \mat{A} - \Ahat$ denote the error matrix, which we will generally assume we can access by matrix--vector products $\vec{\omega} \mapsto \mat{E}\vec{\omega}$.
Some of the approaches we discuss will be agnostic to the precise nature of the approximation $\Ahat$, whereas others will be specialized to particular randomized approximation algorithms.

\subsection{Small-sample estimation} \label{sec:small-sample-errest}

One of the basic tools for error estimation is the small-sample Frobenius-norm estimation technique
described in \cref{sec:small-sample-trace}.
To estimate $\norm{\mat{E}}_{\rm F}$, we generate a Gaussian embedding $\mat{\Gamma} \in \R^{n\times \ell}$ (independently of $\mat{E}$) and form the norm estimate $\hat{\err} \coloneqq \norm{\mat{E}\mat{\Gamma}}_{\rm F}$.\footnote{Recall our definition for a Gaussian embedding (\cref{def:gaussian}), which bakes in a scaling factor of $\ell^{-1/2}$.}
The quantity $\hat{\err}^2$ is an unbiased estimator for the squared error $\norm{\mat{E}}_{\rm F}^2$.
As demonstrated in \cref{eq:small-sample-probs}, a small number of
probes, say $\ell = 10$, suffice to estimate the norm within a modest multiplicative
factor with high confidence.

A core advantage of the small-sample estimator is that the information to compute it can be collected \emph{before the approximation $\Ahat$ has been computed}.
For instance, if we are running the randomized rangefinder algorithm (\cref{alg:rangefinder}), we can generate one large test matrix $\flatonebytwo{\mat{\Omega}}{\mat{\Gamma}} \in \R^{n\times (s+\ell)}$ and begin by computing $\mat{A}\flatonebytwo{\mat{\Omega}}{\mat{\Gamma}}$. 
The matrix $\mat{Y} \coloneqq \mat{A}\mat{\Omega}$ is used to form $\mat{Q} = \orth(\mat{Y})$, and the remaining data $\mat{Z} \coloneqq \mat{A}\mat{\Gamma}$ is held out for error estimation.
For the projection approximation
$\Ahat=\mat{Q}\mat{Q}^\top\mat{A}$, once $\mat{Q}$ has been computed,
one forms the error estimator
\[
    \hat{\err}
    \coloneq \norm{(\Id-\mat{Q}\mat{Q}^\top)\mat{Z}}_{\rm F}
    = \norm{(\mat{A}-\Ahat)\mat{\Gamma}}_{\rm F}
    \approx \norm{\mat{A}-\Ahat}_{\rm F};
\]
no additional interaction with $\mat{A}$ is required!
Using similar ideas, the small-sample estimate can be computed cheaply for many randomized matrix approximation algorithms with low overhead and no additional passes over the matrix.

Suppose we compute the error estimate $\hat{\err}$ and it is unacceptably high.
What then?
In this case, we typically wish to increase the rank of the
approximation to see whether a prescribed error tolerance can
be met.
Given that we already have computed the product $\mat{A}\mat{\Gamma}$, one natural strategy is simply to incorporate $\mat{\Gamma}$ into the test matrix $\mat{\Omega} \gets \flatonebytwo{\mat{\Omega}}{\mat{\Gamma}}$ and draw a fresh $\mat{\Gamma}$ to form a new small-sample estimate; see \cite[Alg.~2]{pearce2026randomized} for an adaptive rangefinder implementation based on this principle.

So far, we have motivated the small-sample error estimate for low-rank matrix approximation, but the technique works just as well to estimate the residual norms for other types of computations such as linear systems, eigenvalue problems, and least-squares problems.
For instance, suppose that we wish to estimate the residual for a (matrix-valued) linear system solve $\mat{A}\mat{X} = \mat{B}$.
To estimate the Frobenius norm of the residual $\mat{A}\hatbold{X} - \mat{B}$, we can use a left-sided small-sample estimate 
\begin{equation*}
    \norm{\smash{\mat{A}\hatbold{X} - \mat{B}}}_{\rm F} \approx \norm{\smash{(\mat{\Gamma}^\top\mat{A})\hatbold{X} -\mat{\Gamma}^\top\mat{B}}}_{\rm F}.
\end{equation*}
Critically, since we have used a left sketch, we are able to approximate the residual without requiring another pass over $\mat{A}$ to compute the product $\mat{A}\hatbold{X}$.

\begin{remark}[Optimal-complexity estimates]
    As an estimate for $\norm{\mat{E}}_{\rm F}^2$, the root-mean-square error of the estimator $\smash{\hat{\err}^2 = \norm{\mat{E}\mat{\Omega}}_{\rm F}^2}$ decreases at an $\order(\ell^{-1/2})$ rate.
    The recently proposed \emph{Hutch\texttt{\#} estimator}
    \begin{equation*}
        \err_{\textrm{H}\texttt{\#}}^2 \coloneqq \norm{\mat{E}\mat{\Omega}}_{\rm F}^2 + \norm{\mat{\Psi}^\top\mat{E}}_{\rm F}^2 - \norm{\mat{\Psi}^\top\mat{E}\mat{\Omega}}_{\rm F}^2
    \end{equation*}
    is defined for even $\ell$ with independent Gaussian embeddings $\mat{\Omega} \in \R^{n\times \ell/2}$ and $\mat{\Psi} \in \R^{m\times \ell/2}$.
    Using $\ell$ matrix--vector products with $\mat{E}$ and $\mat{E}^\top$, Hutch\texttt{\#}
improves the error scaling to $\order(\ell^{-1})$, which is optimal among any single-pass estimator in the worst case \cite{chen2026hutch}.
    While Hutch\texttt{\#} is a natural and elegant idea, we generally recommend the basic small-sample estimator for most use cases, as the benefits of Hutch\texttt{\#} over it manifest only for larger $\ell$. 
\end{remark}

\subsection{Leave-one-out error estimation} \label{sec:leaveoneout}

The small-sample estimation technique gives us a way of estimating the error of a matrix approximation with low overhead.
Nonetheless, there is something fundamentally dissatisfying about holding out some portion of the data we collect about $\mat{A}$ for error estimation.
If we expend effort to compute $s$ matrix--vector multiplications in total, should we not want to use \emph{all $s$} to approximate the matrix?
Put more boldly, is it possible to estimate the error with \emph{no overhead}, in terms of additional matrix--vector multiplications?
The leave-one-out error estimator \cite{epperly2024efficient} provides one tool for doing exactly this.

Let us first develop the estimator in a general way.
Fix a matrix $\mat{A} \in \R^{m\times n}$, and let $\mat{\Omega} = [\vec{\omega}_1,\ldots,\vec{\omega}_s] \in \R^{n\times s}$ be an \emph{unnormalized} test matrix consisting of iid random vectors $\vec{\omega}_i$ satisfying the isotropy property \cref{eq:isotropic}.
Consider any matrix approximation $\Ahat = \Ahat(\mat{Y},\mat{\Omega})$ depending on the test matrix $\mat{\Omega}$ and the sketch $\mat{Y} \coloneqq \mat{A}\mat{\Omega}$.\footnote{For this approach, the approximation $\Ahat$ can also depend on $\mat{A}$ and other quantities such as an independent left sketch $\mat{\Psi}^\top \mat{A}$.
For notational simplicity, we only express the dependence $\Ahat(\mat{Y},\mat{\Omega})$ on $\mat{Y}$ and $\mat{\Omega}$.}
The basic intuition is that, if $s$ is large, the matrix approximation $\Ahat$ should not depend too sensitively on any particular vector $\vec{\omega}_i$, so the approximation should be nearly stable under removal of the $i$th column \cref{eq:column-deletion}:
\begin{equation*}
    \Ahat(\mat{Y},\mat{\Omega}) \approx \Ahat(\mat{Y}_{-i},\mat{\Omega}_{-i}).
\end{equation*}
Once the $i$th vector $\vec{\omega}_i$ has been left out, it is independent of the downdated approximation $\Ahat(\mat{Y}_{-i},\mat{\Omega}_{-i})$, and we can use it to form a single-sample estimate for the error
\begin{equation*}
    \norm{\mat{A} - \Ahat}_{\rm F}^2 \approx \norm{\mat{A} - \Ahat(\mat{Y}_{-i},\mat{\Omega}_{-i})}_{\rm F}^2 \approx \norm{(\mat{A} - \Ahat(\mat{Y}_{-i},\mat{\Omega}_{-i}))\vec{\omega}_i}^2.
\end{equation*}
Using $\mat{A}\vec{\omega}_i=\mat{Y}(:,i)$ and averaging over $i$,
we obtain the \emph{leave-one-out error estimator}
\begin{equation} \label{eq:leave-one-out}
    \hat{\err}^2_{\rm L} \coloneqq \frac{1}{s} \sum_{i=1}^s \norm{\mat{Y}(:,i) - \Ahat(\mat{Y}_{-i},\mat{\Omega}_{-i})\vec{\omega}_i}^2.
\end{equation}
This construction follows the same exchangeability principle as the XTrace estimator from \Cref{sec:xtrace}.

The formula \cref{eq:leave-one-out} for the leave-one-out error estimator may seem daunting.
On the positive side, the estimator requires no interaction with the matrix $\mat{A}$ beyond that used to form the approximation $\Ahat$.
However, it may seem computationally expensive to evaluate $\Ahat(\mat{Y}_{-i},\mat{\Omega}_{-i})\vec{\omega}_i$ for each $i$.
Fortunately, for many algorithms of interest, the formula \cref{eq:leave-one-out} admits a dramatic algebraic simplification using an appropriate matrix approximation downdating formula \cite{epperly2024efficient,epperly2024xtrace,epperly2025make,lazzarino2026efficient}.
For instance, for the projection approximation (randomized SVD) 
\begin{equation*}
    \Ahat = \mat{Q}\mat{Q}^\top\mat{A} \quad \text{for } (\mat{Q},\mat{R}) \gets \Call{QR}{\mat{Y},\texttt{'econ'}},
\end{equation*}
the leave-one-out estimate admits the simple formula
\begin{equation} \label{eq:rsvd-leave-one-out}
    \hat{\err}_{\rm L}^2 = \frac{1}{s} \sum_{i=1}^s \norm{\mat{R}^{-1}(i,:)}^{-2}.
\end{equation}
The ``R'' matrix in the QR decomposition $\mat{Y} = \mat{Q}\mat{R}$, a useless byproduct in many randomized SVD implementations, contains by itself enough information to estimate the error.
The cost of forming the leave-one-out error estimator is only $\order(s^3)$ operations, independent of \emph{both} dimensions of the matrix $\mat{A}$!\footnote{In fact, one can compute this error estimator to constant-factor accuracy in an even faster $\order(s^2\log s)$ by using a randomized row-norm estimator (\cref{sec:row-norm-estimation}).}
See \cite[Alg.~4.1]{epperly2024efficient} for leave-one-out error estimation for the Nystr\"om approximation (\cref{sec:lowrank-nystrom}) and \cite[\S3]{lazzarino2026efficient} for the generalized Nystr\"om approximation (\cref{sec:lowrank-generalized-nystrom}).
It is an open question whether leave-one-out error estimation can be implemented in a stable way for a low-rank approximation using powering (\cref{sec:lowrank-subspace-iteration}) or block Krylov iteration (\cref{sec:lowrank-block-krylov}).

By construction, the leave-one-out estimator is generally \emph{not} an unbiased approximation to the error $\norm{\mat{A} - \Ahat}^2_{\rm F}$.
But the estimator \emph{is} an unbiased estimator for the mean-square error for an $(s-1)$-sample approximation.
That is,
\begin{equation*}
    \expect[\hat{\err}^2_{\rm L}] = \expect \norm{\mat{A} - \Ahat(\mat{A}\mat{\Omega}',\mat{\Omega}')}_{\rm F}^2 \quad \text{where } \mat{\Omega}' = [\vec{\omega}_1,\ldots,\vec{\omega}_{s-1}] \in \R^{n\times (s-1)}.
\end{equation*}

\subsubsection{Matrix jackknife variance estimation}

In many applications, we are interested in using randomized algorithms not to estimate $\mat{A}$ itself but some derived quantity of interest $Q(\mat{A})$ such as its eigenvalues or eigenvectors, singular values or singular vectors, projectors onto invariant subspaces, etc.
In these cases, a gross estimate of $\norm{\mat{A} - \Ahat}_{\rm F}$ may say little about the accuracy of the approximation $Q(\Ahat) \approx Q(\mat{A})$.
While it may not be computationally tractable in general to approximate the error $\norm{Q(\mat{A}) - Q(\Ahat)}_{\rm F}$, one can use jackknife resampling ideas from statistics to estimate the \emph{variance} of the quantity $Q(\Ahat)$ for a quite general class of matrix approximations $\Ahat$ and quantity of interest $Q$.
While the variance is often not a good proxy for the (squared) error, it nonetheless contains useful information.
In particular, high variance is an indicator that the output of an algorithm should not be trusted.
See \cite{epperly2024efficient} for details on jackknife variance estimation for randomized matrix computations.

\subsection{Error estimation by sketching}

Recall that a projection approximation (\cref{sec:lowrank-projection}) approximates a matrix $\mat{A}$ by drawing a test matrix $\mat{\Omega}$, computing $\mat{Y} = \mat{A}\mat{\Omega}$, and orthogonally projecting:
\begin{equation*}
    \Ahat \coloneqq \mat{Y}\mat{Y}^\pinv\mat{A} = \mat{Q}\mat{Q}^\top\mat{A} \quad \text{for } \mat{Q} = \orth(\mat{Y}).
\end{equation*}
When $\mat{\Omega}$ is a random embedding, $\Ahat$ coincides with the randomized SVD approximation (\cref{alg:rsvd}).
When $\mat{\Omega} = \Id(:,\iset)$ is a selection matrix, $\Ahat$ is an interpolative low-rank approximation (\cref{sec:activelsq}).
In many situations, determining the non-orthogonal basis matrix $\mat{Y}$ is the cheapest part of computing a projection approximation, with the dominant costs being orthogonalization $\mat{Y} = \mat{Q}\mat{R}$ and computing either $\mat{Q}^\top\mat{A}$ or $\mat{Y}^\dagger\mat{A}$.
In such cases, it can be advantageous to have an error estimate that can be computed before performing these expensive steps.

To do so, we first recognize the error of the projection approximation as the minimum residual norm of a least-squares problem
\begin{equation*}
    \norm{\mat{A} - \Ahat}_{\rm F} =  \norm{\smash{\mat{Y}\mat{Y}^\pinv \mat{A} - \mat{A}}}_{\rm F} = \min_{\mat{Z} \in \R^{n\times s}} \norm{\smash{\mat{Y}\mat{Z}^\top - \mat{A}}}_{\rm F}.
\end{equation*}
When the projection approximation is accurate, the residual for this least-squares problem will be tiny.
As always when only a quasi-optimal solution for a tall-thin least-squares problem is required, we naturally turn to the sketch-and-solve method (\cref{sec:sas}).
Employing a fast random embedding $\mat{\Psi} \in \R^{m\times t}$, we obtain an approximate least-squares problem
\begin{equation} \label{eq:sketched-projection-error-est}
    \min_{\hatbold{Z}\in\R^{n\times s}}\norm{\smash{(\mat{\Psi}^\top\mat{Y})\hatbold{Z}^\top - \mat{\Psi}^\top\mat{A}}}_{\rm F} \quad \text{with solution } \hatbold{Z}^\top = (\mat{\Psi}^\top\mat{Y})^\pinv(\mat{\Psi}^\top\mat{A}).
\end{equation}
Assuming $\mat{\Psi}$ is an oblivious subspace embedding or injection with subspace dimension $s$, \cref{thm:sketch-and-solve-ose-osi}(c) ensures that
\begin{equation*}
    \norm{\smash{\mat{Y}(\mat{\Psi}^\top\mat{Y})^\pinv(\mat{\Psi}^\top\mat{A}) - \mat{A}}}_{\rm F} \le C\, \norm{\mat{Y}\mat{Y}^\dagger \mat{A} - \mat{A}}_{\rm F} = C\, \norm{\mat{A} - \Ahat}_{\rm F}
\end{equation*}
for some small prefactor $C > 1$.
Evaluating the left-hand side of this expression is still expensive, but we can cheaply approximate it by employing a small-sample estimate, as in \cref{sec:small-sample-errest}.
We then define the \emph{sketched error estimate}
\begin{equation*}
    \hat{\err}_{\rm S} \coloneqq \norm{\mat{Y} [(\mat{\Psi}^\top \mat{Y})^\pinv(\mat{\Psi}^\top(\mat{A}\mat{\Gamma}))] - \mat{A}\mat{\Gamma}}_{\rm F}.
\end{equation*}
The sketched error estimation technique is due to \citet{martinsson2026fly}.

\begin{algorithm}[t]
\caption{Projection approximation with sketched error estimate} \label{alg:rsvd-with-errest}
\begin{algorithmic}[1]
\Require Matrix $\mat{A} \in \R^{m\times n}$ and rank $s$
\Ensure Approximation $\Ahat = \mat{Q}\mat{B}^\tp$ or $\Ahat = \mat{U}\mat{\Sigma}\mat{V}^\top$ and error estimate $\hat{\err}_{\rm S} \approx \norm{\mat{A} - \Ahat}_{\rm F}$
\State Generate Gaussian embeddings $\mat{\Omega} \in \R^{n\times s}$ and $\mat{\Gamma} \in \R^{n\times 10}$
\State Generate a SparseStack $\mat{\Psi} \in \R^{m\times 2s}$ with sparsity level $\zeta = 4$
\State Sketch on the right: $\flatonebytwo{\mat{Y}}{\mat{W}} \gets \mat{A}\flatonebytwo{\mat{\Omega}}{\mat{\Gamma}}$
\State Sketch on the left: $\smash{\flatonebytwo{\hatbold{Y}}{\hatbold{W}} \gets \mat{\Psi}^\tp \flatonebytwo{\mat{Y}}{\mat{W}}}$
\State Compute a factorization of the sketch $(\hatbold{Q},\hatbold{R}) \gets \Call{QR}{\hatbold{Y},\texttt{'econ'}}$
\State Solve the sketched least-squares problem $\mat{C} \gets \hatbold{R}^{-1}(\hatbold{Q}^\top\hatbold{W})$ \Comment{$\mat{C} = \hatbold{Y}^\pinv\hatbold{W}$}
\State Compute the error estimate $\hat{\err}_{\rm S} \gets \norm{\mat{Y}\mat{C} - \mat{W}}_{\rm F}$
\State Orthonormalize by QR factorization $(\mat{Q},\sim) \gets \Call{QR}{\mat{Y},\texttt{'econ'}}$
\State Form $\mat{B} \gets \mat{A}^\tp\mat{Q}$
\State \emph{Optional:} Compute an economy-size SVD $\mat{B}^\tp = \Utilde \mat{\Sigma}\mat{V}^\tp$ and set $\mat{U} \gets \mat{Q}\Utilde$
\end{algorithmic}
\end{algorithm}

\Cref{alg:rsvd-with-errest} provides an example implementation of the randomized SVD with the sketched error estimate.
The key power of the estimator comes from the fact that it is computed before the matrix $\mat{Y}$ is factorized and the product $\mat{Q}^\top\mat{A}$ is formed.
If the estimate is unacceptably high at that point, one could iteratively expand the matrix $\mat{\Omega}$ by appending additional columns until the error estimate falls below a tolerance.

\subsection{The spectral norm and other Schatten norms}

So far, we have discussed techniques for estimating the \emph{Frobenius norm} of an error matrix $\mat{E}$.
In many applications, the \emph{spectral norm} is the more natural error metric.
To estimate the spectral norm (equivalently, the top singular value of $\mat{E}$), we can generate a random initial guess and utilize the power method 
\begin{equation*}
    \vec{x} \gets \vec{x} / \norm{\vec{x}}, \quad \vec{y} \gets \mat{E}\vec{x}, \quad \hat{\err} \gets \norm{\vec{y}}, \quad \vec{x} \gets \mat{E}^\top\vec{y}.
\end{equation*}
Starting from a random initial guess $\vec{x} \sim \Normal(\vec{0},\Id)$, roughly $\log n$ steps are needed for the power method to approximate the spectral norm within a constant factor; see \cite{kuczyski1992estimating}, \cite[\S6.2]{martinsson2020randomized}, and \cite[\S3.1]{kireeva2024randomized}.
Using the Lanczos method instead of power iteration results in asymptotically faster convergence, though it still requires $\order(\log n)$ iterations to achieve a constant-factor estimate~\cite[\S6.3]{martinsson2020randomized}.
This $\log n$ dependence is not an artifact of the analysis of these methods, either, as it is known that roughly $\log n$ matrix--vector queries are \emph{necessary} in the worst case for ``any small-sample estimate'' to produce a constant-factor approximation to the spectral norm; see \cite{simchowitz2017gap,bakshi2023krylov} for precise results.
We thus recognize that spectral-norm estimation is fundamentally harder than Frobenius-norm estimation, which can easily be done with a small-sample estimate in a single pass.

As a halfway house between Frobenius norm (easy) and spectral norm (hard), one can also try to estimate the intermediate Schatten norms $\norm{\mat{E}}_{(p)}$ for an integer $p \in (2,\infty)$.
If $p$ is divisible by four, one approach is to use the identity $\norm{\mat{E}}_{(p)}^p = \norm{\smash{(\mat{E}^\top\mat{E})^{p/4}}}_{\rm F}^2$ and deploy a small-sample estimator.
Using $\ell=10$ samples, the cost is $5p$ matrix--vector multiplications with $\mat{E}$ and its transpose computed across $p/2$ passes over the matrix.
Alternatively, if we insist on a single-pass estimator, we can use the estimator of \citet{kong2017spectrum}; see \cite[\S5.4]{martinsson2020randomized} and \cite{chu2025improved}.

\section{Regularization and simplification by randomization} \label{sec:regularization}

This section focuses on a role of randomization that has already appeared a few
times in this survey: using randomness to transform a difficult problem into one that is
more amenable to standard numerical algorithms. For example, in
\Cref{sec:orthognalization}, whitening turns a poorly conditioned basis
into one on which CholeskyQR can be applied safely.

In the following, we explore further instances of this principle, where randomization
is used to bypass exceptional situations and numerical difficulties. One can
broadly distinguish between two classes of techniques: (1) deliberately introducing random perturbations
to the problem and controlling the resulting error, and (2) randomly transforming the problem into one that preserves the quantity of interest exactly.

Achieving the accuracy and efficiency required in large-scale scientific computing
introduces additional challenges. In particular, regularization should preserve
properties that make the original problem computationally tractable, such as
sparsity and symmetry.

While closely related, \emph{smoothed analysis}~\cite{spielman2004smoothed,sankar2006smoothed} has traditionally had a different aim: Using small random perturbations as an analytical device to
explain why algorithms with poor worst-case behavior often perform well in
practice. When such perturbations are deliberately introduced by the algorithm
itself, smoothed analysis techniques become tools for algorithm design.

Another source of beneficial randomness, which we do not discuss further here,
is \emph{stochastic rounding}. It can mitigate the accumulation of rounding
errors~\cite{croci2022stochastic} and can also have a genuine regularizing effect:
under suitable conditions, stochastic rounding lifts the smallest singular value of a tall-thin matrix
and yields full column rank with high probability~\cite{dexter2025stochastic}.

\subsection{Linear systems} \label{sec:singsystems}

Randomness is not a silver bullet for curing numerical difficulties. For example, consider a highly ill-conditioned linear system $\mat{A} \vec{x} = \vec{b}$. 
Simply adding a random perturbation usually does not resolve difficulties associated with ill-conditioning: A perturbation large enough to substantially improve the condition number may also substantially change the solution.  

\paragraph{Modification of singular linear systems.}
Interestingly, there are exceptional situations when the right-hand side has little or no component in the ill-conditioned directions of $\mat A$.
A particularly favorable example is a consistent singular system
\begin{equation} \label{eq:singular}
   \mat A\vec x=\vec b, \qquad
   \mat A\in\R^{n\times n}, \qquad
   \rank(\mat A)=n-k<n, \qquad
   \vec b\in\range(\mat A).
\end{equation}
Although the solution is
not unique, finding \emph{some} solution 
does \emph{not} need to be an ill-posed problem and, in fact, randomization helps find it.
Following~\cite{pan2010additive,sifuentes2015randomized}, 
consider the rank-$k$ modification
\begin{equation} \label{eq:pertlinsystem}
    (\mat{A}+\mat{\Psi} \mat{\Omega}^\tp)\vec y=\vec b, \qquad \mat{\Psi},\mat{\Omega}\in\R^{n\times k}.
\end{equation}
Let $\mat U,\mat V\in\R^{n\times k}$ have orthonormal columns spanning
the left and right null spaces of $\mat A$. For generic choices of $\mat{\Psi}$, $\mat{\Omega}$, the matrices $\mat U^\tp\mat{\Psi}$ and $\mat \Omega^\tp\mat V$ are invertible, which implies that $\mat A+\mat{\Psi}\mat \Omega^\tp$ is invertible.\footnote{Indeed, if
$(\mat A+\mat{\Psi}\mat \Omega^\tp)\vec z=0$, multiplication by $\mat U^\tp$
first gives $\mat \Omega^\tp\vec z=0$, and hence $\mat A\vec z=0$.
Writing $\vec z=\mat V\vec c$ and using the invertibility of
$\mat \Omega^\tp\mat V$ gives $\vec z=0$.}
In turn, a standard linear system solver can
be applied to solve~\cref{eq:pertlinsystem}; its unique solution $\vec y$ is also a 
solution of the original system~\cref{eq:singular}. To see this, first note that 
$\mat{U}^\tp \mat A = 0$ implies $\mat U^\tp\vec b=0$ because 
$\vec b\in\range(\mat A)$.
Multiplying~\cref{eq:pertlinsystem} from the left by
$\mat{U}^\tp$ thus gives
$
    (\mat{U}^\tp \mat{\Psi}) \mat{\Omega}^\tp \vec y=0.
$
The invertibility of $\mat{U}^\tp \mat{\Psi}$ implies
$
    \mat{\Omega}^\tp \vec y=0
$ and, hence, it follows from~\cref{eq:pertlinsystem}
that $\mat{A} \vec y=\vec b$.
Thus, the rank-completing modification
removes the singularity without introducing any additional error when~\cref{eq:pertlinsystem} is solved exactly.

In practice, the choice of $\mat{\Psi}$, $\mat{\Omega}$ influences the conditioning of~\cref{eq:pertlinsystem} and the robustness of the approach to error caused by, e.g., roundoff or an iterative solver. The oblique projector
$
 \mat{\Pi}_{\mat A,\mat{\Psi}} \coloneqq \Id - \mat{\Psi} (\mat{U}^\tp \mat{\Psi})^{-1} \mat{U}^\tp 
$
allows us to connect the residual for an inexact solution $\widehat{\vec y}$ of the regularized system~\cref{eq:pertlinsystem}
to the original system. Indeed, since $\mat{\Pi}_{\mat A,\mat{\Psi}}\vec b=\vec b$, $\Pi_{\mat A,\mat{\Psi}}\mat A=\mat A$,
$\mat{\Pi}_{\mat A,\mat{\Psi}} \mat{\Psi} = 0$,  it follows that
\[
 \vec b - \mat{A} \widehat{\vec y} = \Pi_{\mat A,\mat{\Psi}} \big( \vec b - (\mat{A}+\mat{\Psi} \mat{\Omega}^\tp)\widehat{\vec y} \big)
\]
The situation is analogous to sketch-and-solve (see, in particular, the proof of \Cref{thm:sketch-and-solve-ose-osi}(b)) and we obtain 
\[
 \| \vec b - \mat{A} \widehat{\vec y} \| \le \frac{\|\mat{\Psi}\|}{\sigma_{\min}(\mat{\Psi}^\tp \mat{U})} \cdot \big\| \vec b - (\mat{A}+\mat{\Psi} \mat{\Omega}^\tp)\widehat{\vec y} \big\|.
\]
If, for example, $\mat{\Psi}$ is uniformly distributed on the real Stiefel manifold, this prefactor is typically of order $\sqrt{k(n-k)}$~\cite{dumitriu2012smallest}; see~\cite{shah2021smoothed,sifuentes2015randomized} for more detailed analyses. 
An extension of this approach to nonsingular but nearly rank-deficient systems can be found in~\cite{pan2012solving}.

\paragraph{Preprocessing for Gaussian elimination.}
Another application of randomization to put linear systems in a better position for numerical solution 
is the use of fast random orthogonal transformations to enable Gaussian elimination \emph{without} pivoting.
This important development is already surveyed in detail in \cite[Section~17.2]{martinsson2020randomized}. 

\subsection{Perturbing symmetric eigenvalue problems}

The remainder of this section is concerned with regularizing eigenvalue problems and we start by adding 
a random diagonal perturbation to a symmetric matrix $\mat{A} \in \R^{n\times n}$:
\begin{equation} \label{eq:diagpert}
  \hat{\mat A} = \mat{A} + \gamma \mat{D}, \quad \mat{D} = \diag\big(d_{11}, d_{22}, \ldots, d_{nn} \big),
\end{equation}
where the diagonal entries $d_{11},\ldots,d_{nn}$ are i.i.d. random variables from a
distribution with \emph{bounded} probability density function $\rho(\cdot)$. The scalar $\gamma > 0$ is meant to control the norm of the perturbation; by Weyl's inequality, each eigenvalue of $\mat A$ is displaced by at most
$\gamma\|\mat D\|_2$.
\begin{figure}
\begin{center}
 \includegraphics[width = 0.5\linewidth]{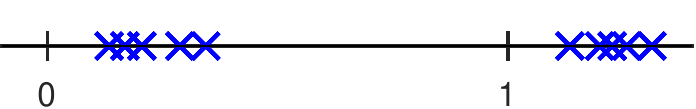}
\end{center}
\caption{\label{fig:diagpert}Eigenvalues of perturbed matrix~\cref{eq:diagpert} for a non-diagonal matrix $\mat{A}$ with eigenvalues $0$ and $1$ (of multiplicity $5$ each), $\gamma = 0.4$, and $d_{11}, \ldots, d_{nn} \sim \Unif(0,1)$.   }
\end{figure}
In the example shown in~\Cref{fig:diagpert}, one observes that the eigenvalues move away from their original locations and---more importantly---multiple eigenvalues \emph{split} into visibly separated eigenvalues.

\paragraph{Theory.}
Eigenvalue splitting and related phenomena have been studied in detail in 
mathematical physics, particularly for Schr\"odinger operators with random potentials. The classical Wegner estimate~\citeyearpar{wegner1981bounds} implies for~\cref{eq:diagpert} that
\[
 \expect N(\set{I}) \le \gamma^{-1} \|\rho\|_\infty |\set{I}|\,n,
\]
where $N(\set{I})$ counts the number of eigenvalues of $\hat{\mat A}$ in an interval $\set{I}$ of length $|\set{I}|$; see also~\cite[Eq. (2.4)]{combes2009generalized}.
By Markov's inequality,
\begin{equation} \label{eq:wegner}
  \prob\{ N(\set{I}) \ge 1 \} \le \gamma^{-1} \|\rho\|_\infty |\set{I}|\,n.
\end{equation} 
For example, when choosing a small interval around zero, this implies that $\big\| \hat{\mat{A}}^{-1} \big\| = \order\big( n \gamma^{-1} \big)$ holds with high probability, \emph{independent} of $\mat{A}$. Such results are central in the smoothed analysis mentioned above.

For more than one eigenvalue in a small interval, Minami-type estimates~\citeyearpar{minami1996local} yield stronger bounds of the form
\begin{equation} \label{eq:minami}
  \prob\{ N(\set{I}) \ge 2 \} \le  \frac12 \big( \gamma^{-1} \|\rho\|_\infty |\set{I}|\,n \big)^2;
\end{equation}
see~\cite[Corollary 2.4]{combes2009generalized}. By a union bound argument, this confirms what is observed in~\Cref{fig:diagpert}: The minimum gap among the eigenvalues of $\hat{\mat{A}}$ contained in an interval of length $\order(\gamma)$ is unlikely to be much smaller than $\gamma / n^2$; see, e.g.,~\cite[Theorem 2.2]{kressner2026randomized} for a precise statement.

\paragraph{Use in scientific computing.} A diagonal perturbation~\cref{eq:diagpert} is attractive computationally: it preserves sparsity, adding only $\order(n)$ work to a matrix--vector product.
One compelling use is to fix a shortcoming of the standard Lanczos method for computing eigenvalues and eigenvectors of $\mat A$~\cite[\S12.2.1]{parlett1998symmetric}: If $\mat A$ has a multiple eigenvalue, a Krylov subspace
$\mathcal K_k(\mat A,\vec x_0)$, see~\Cref{sec:krylovreview},
cannot, in exact arithmetic, contain or approximate more than one
independent eigenvector associated with that eigenvalue.
Such multiple eigenvalues often arise from symmetries or structural properties of the
matrix and are of particular relevance when computing nullspaces of matrices.

To a certain extent, roundoff error breaks this phenomenon by introducing new directions into the Krylov subspace. However, for eigenvalues of higher multiplicity, this benefit from finite-precision arithmetic is usually insufficient to reliably recover the full eigenspace. On the other hand, the theory above and numerical evidence reported in~\cite{meyer2024unreasonable,kressner2026randomized} suggest that 
deliberately adding a small diagonal random perturbation allows standard Lanczos to reliably and rapidly resolve clusters of (perturbed) exterior eigenvalues.

A natural alternative is block Lanczos (\Cref{sec:lowrank-block-krylov}), which can capture the full eigenspace when the block size is at least as large as the eigenvalue multiplicity.
In practice, rather small blocks often appear to be sufficient; as~\citet[\S13.12]{parlett1998symmetric}
observes, ``Researchers who have experimented with various blocksizes find the
smaller sizes, such as 2 or 3, are most efficient.''
Recent results provide further support for the effectiveness
of small blocks, especially when combined with random diagonal regularization~\cite{chen2026does,shao2025structural}.

\subsection{Perturbing non-symmetric eigenvalue problems}

For non-symmetric eigenvalue problems, the situation is more subtle than
in the symmetric case. Splitting multiple eigenvalues no longer suffices
to alleviate numerical difficulties: even a matrix with
well separated eigenvalues can have a highly ill-conditioned eigenvector basis.
Moreover, diagonal perturbations do not, in general, control this ill-conditioning.
To illustrate these two issues, consider a Jordan block
\[
\mat J_m(\lambda) \coloneqq
 \begin{bmatrix}
 \lambda & 1      &        &        \\
         & \lambda& \ddots &        \\
         &        & \ddots & 1      \\
         &        &        & \lambda
 \end{bmatrix}
 \in\C^{m\times m}.
\]
Adding a diagonal perturbation $\mat D=\diag(d_1,\ldots,d_{m})$
splits the multiple eigenvalue $\lambda \in \C$ up into $m$ (generically) different eigenvalues $\lambda + d_1,\ldots,\lambda + d_{m}$.
Even so, the condition number $\kappa(\mat V)$ of
any eigenvector matrix $\mat V$ for $J_m(\lambda) + \mat D$ 
grows at least exponentially with $m$, provided that $|d_i| < 1/2$.
Indeed, choosing $d_i$ iid  uniformly distributed in $[0,1/4]$ and $m = 10$ typically yields $\kappa(\mat V) \approx 10^{12}$ and similarly high eigenvalue condition numbers.

Throughout this and the next subsection,
we work over the complex field; not only because a real nonsymmetric matrix may
have complex eigenvalues, but also because complex random perturbations usually
provide simpler and stronger regularization guarantees.

\paragraph{Pseudospectral shattering.} Striking results of \citet{banks2021gaussian,banks2023pseudospectral} show that a full Gaussian perturbation of small norm simultaneously controls eigenvector conditioning and eigenvalue gaps: 
\[
 \kappa_V(\mat M)
 :=\inf\big\{ \kappa(\mat V) \colon \text{$\mat V^{-1} \mat M \mat V$ is diagonal}\},
 \quad
 \operatorname{gap}(\mat M)
 :=\min_{i\ne j}|\lambda_i(\mat M)-\lambda_j(\mat M)|,
\]
where eigenvalues are counted with multiplicity,
and $\kappa_V(\mat M)=\infty$ if $\mat M$ is not diagonalizable.
To describe this result, let us first recall that
the \emph{$\varepsilon$-pseudospectrum} $\Lambda_\varepsilon(\mat M)$
contains the eigenvalues of every perturbed matrix $\mat M + \mat F$
with $\|\mat F\|<\varepsilon$. Equivalently,
\begin{equation}\label{eq:shattering-pseudospectrum}
 \Lambda_\varepsilon(\mat M)
 :=\{z\in\C \colon \sigma_{\min}(z\Id-\mat M)<\varepsilon\};
\end{equation}
see~\cite{trefethen2005spectra}. 

A pseudospectrum $\Lambda_\varepsilon(\mat M)$ is said to be
\emph{shattered} with
respect to a square grid of mesh width $\omega$ on $\C$ if every grid
square contains at most one eigenvalue and
$\Lambda_\varepsilon(\mat M)$ does not intersect any grid line.
To avoid grid lines (nearly) colliding with eigenvalues, the grid
is randomly shifted by an offset
$\alpha+\mathrm{i}\beta$ with
$\alpha,\beta\simiid \Unif(0,\omega)$.

The following result is from \citet[Theorems~1.4 and~3.12]{banks2023pseudospectral}.
\begin{theorem}[Pseudospectral shattering]
\label{thm:pseudospectral-shattering}
Let $n\ge2$ and let $\mat A\in\C^{n\times n}$ satisfy
$\|\mat A\|\le1$, and set
$\mat M=\mat A+\gamma\mat G$ for $0<\gamma<1/2$, where
$G_{ij}=(X_{ij}+\mathrm{i}Y_{ij})/\sqrt{2n}$ with iid $X_{ij},Y_{ij} \sim \Normal(0,1)$.
Consider a randomly shifted square grid of mesh width
$
 \omega=\gamma^4 / (4n^5)
$.
Then, with probability at least $1-13/n$,
\[
 \|\mat M-\mat A\|_2\le4\gamma,\qquad
 \kappa_V(\mat M)\le\frac{n^2}{\gamma},\qquad
 \operatorname{gap}(\mat M)\ge\frac{\gamma^4}{n^5},
\]
and $\Lambda_{\varepsilon_\star}(\mat M)$ is shattered with respect to
this grid for
$
 \varepsilon_\star = \gamma^5 / (16n^9).
$
\end{theorem}

\begin{figure}[tp]
 \centering
 \includegraphics[width=\linewidth]{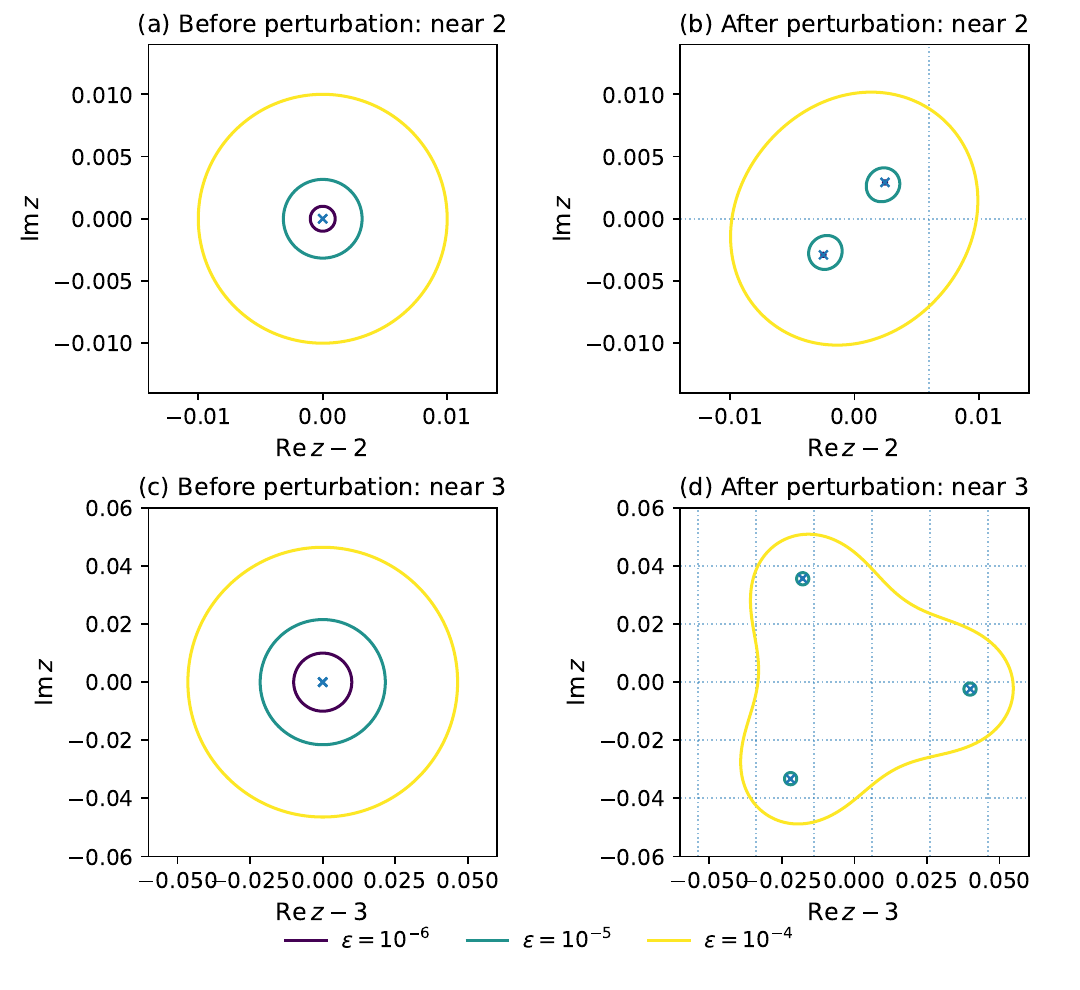}
  \caption{\label{fig:pseudospectral-shattering}
 Pseudospectral contours for $\varepsilon=10^{-6},10^{-5},10^{-4}$ of $\mat A=\diag(J_2(2),J_3(3))$ before
 (left) and after (right) the perturbation of~\cref{thm:pseudospectral-shattering}
 with $\gamma = 10^{-4}$. The upper and lower panels zoom
 around the eigenvalues $2$ and $3$, respectively. 
 The dotted lines form an illustrative square grid of mesh width $0.02$.}
\end{figure}

To illustrate~\Cref{thm:pseudospectral-shattering}, we consider 
the matrix $\mat A=\diag(J_2(2),J_3(3))$, that is $\mat A$ has eigenvalues $2$ and $3$ contained in $2\times 2$ and $3\times 3$ Jordan blocks, respectively. As seen in~\Cref{fig:pseudospectral-shattering},
the pseudospectra of $\mat A$ appear to be very large, disks of radii $\varepsilon^{1/2}$ and 
$\varepsilon^{1/3}$ around $2$ and $3$, respectively, reflecting the ill-conditioning of eigenvalues of Jordan blocks under perturbations. When $\mat A$ is randomly perturbed, 
the multiple eigenvalues split into distinct eigenvalues of $\mat M$.
Moreover, the pseudospectra of $\mat M$ become much smaller, at least for smaller $\varepsilon$, in accordance with the bound on $\kappa_V(\mat M)$: The Bauer--Fike inclusion theorem~\cite[Theorem~6.3.2]{horn2012matrix} states that
$\Lambda_{\varepsilon}(\mat M)$ is contained in discs of radius $\varepsilon\kappa(\mat V)$ around the eigenvalues of $\mat M$.

\paragraph{Diagonalization.}

The perturbed matrix $\mat M$ of~\Cref{thm:pseudospectral-shattering} can be diagonalized:
$\mat M=\mat V\mat D\mat V^{-1}$ with the condition number of $\mat V$ bounded independently of spectral properties of $\mat A$. In other words, $\mat A$ can be approximately diagonalized up to a backward error of norm $\gamma \|\mat G \|$.
For the realization in \cref{fig:pseudospectral-shattering}, $\gamma \|\mat G\| \approx 2\times 10^{-4}$ and $\kappa(\mat V)\approx626.25$. Moreover, the eigenvalue error caused by regularization is at most $4 \times 10^{-2}$.

\paragraph{Spectral division and the matrix sign function.}
Shattering supplies lines on which the resolvent is
controlled: By~\Cref{eq:shattering-pseudospectrum}, $\|(z\mat I-\mat M)^{-1}\| \le\varepsilon_\star^{-1}$ for any $z$ on a grid line. Such resolvent bounds play an important role in the design and analysis of numerical algorithms.
In particular, the vertical line $\operatorname{Re}z=s$ separates the spectrum of $\mat M$ and the matrices
\[
 \mat P_+
 =\tfrac12\big(\Id+\operatorname{sign}(\mat M-s\mat I)\big), \qquad 
  \mat P_-
 =\tfrac12\big(\mat I-\operatorname{sign}(\mat M-s\mat I)\big)
\]
are well defined; they are the spectral projectors associated with the eigenvalues to the right and to the left of the line, respectively. The matrix sign function $\operatorname{sign}(\mat M - s\mat I)$ 
(and, hence, $\mat P_{\pm}$)
can be computed efficiently by matrix iterations. A classical choice is the matrix sign function iteration~\cite{roberts1980linear},
\[
 \mat X_{k+1}=\tfrac12(\mat X_k+\mat X_k^{-1}),\qquad
 \mat X_0=\mat M - s\mat I,
\]
which is usually combined with scaling strategies to ensure fast convergence~\cite[\S5.5]{higham2008functions}. An analogous construction treats horizontal lines.
Extracting bases for the projector ranges using, e.g., a
randomized rangefinder (\Cref{sec:rangefinder}), allows one to recursively
reduce the size of the eigenvalue problem, effectively dividing its spectrum. While this spectral divide-and-conquer is classical~\cite{bai1993design}, random regularization allows one to control the numerical error throughout the recursion, using the resolvent bound stated above. \citet{banks2023pseudospectral} develop and analyze a variant of
spectral divide-and-conquer that ensures 
for
$\|\mat A\| \le1$ and $0<\delta<1$, an approximate diagonalization satisfying
\[
 \|\mat A-\mat V\mat D\mat V^{-1}\|\le\delta,
 \qquad \kappa(\mat V)\le 32n^{5/2} \delta^{-1},
\]
with probability at least $1-14/n$, in nearly matrix multiplication complexity and with 
$\order\!\left(\log^4(n/\delta)\log n\right)$
bits of precision. For Hermitian matrices, a backward stable spectral divide-and-conquer
algorithm was already developed by~\citet{nakatsukasa2013stable}.

\paragraph{Convergence of shifted QR.}
While the algorithm described above aims at approximate diagonalization, for nonsymmetric matrices it is usually preferable to avoid a potentially ill-conditioned eigenbasis and instead 
compute a Schur decomposition $\mat A=\mat Q\mat T\mat Q^*$, where $\mat Q$ is unitary
and $\mat T$ is upper triangular.
The standard approach first reduces $\mat A$ to Hessenberg form and then applies an iterative method, the shifted QR algorithm~\cite{francis1961qr,francis1962qr}.
Despite its remarkable practical success, a general theory guaranteeing global
rapid convergence of shifted QR remained elusive for decades.
In a series of recent papers,~\citet{banks2025globala,banks2025globalb,banks2025globalc} developed new shifting
strategies with such convergence guarantees, also in finite precision, for diagonalizable matrices with controlled eigenvector conditioning and eigenvalue gaps, both of which can be ensured by the random regularization described above.

\paragraph{Matrix functions.}
While specialized algorithms are available for matrix functions such as the matrix exponential
and logarithm, the Schur--Parlett algorithm
is a standard approach for general $f(\mat A)$~\cite{higham2008functions}.
\Cref{thm:pseudospectral-shattering} seems to suggest a simpler approach: Approximate $f(\mat A)$ by
$\mat V f(\mat D)\mat V^{-1}$, where
$\mat A\approx\mat V\mat D\mat V^{-1}$ is obtained by random regularization.
A difficulty is that even a small perturbation of $\mat A$ may induce a much
larger change in $f(\mat A)$. To keep perturbations on the level of working precision, \citet{higham2021multiprecision} use higher precision when applying such an approach to matrix functions of (small-sized) diagonal blocks within the Schur--Parlett algorithm.

\paragraph{Structured random perturbations and extensions.}
Extending \Cref{thm:pseudospectral-shattering} to structured random
perturbations is remarkably challenging.
\citet{jain2021real} and, independently, \citet{banks2024overlaps}
establish analogous results for \emph{real} random perturbations.
\citet{shah2024sparse} show that pseudospectral shattering can also be
achieved with sparse perturbations, adding complex Gaussian noise at only
$\order(n\log^6 n)$ randomly selected entries.

\citet{demmel2026generalized} extend pseudospectral shattering to 
matrix pencils and use it to develop a randomized spectral divide-and-conquer
algorithm for diagonalization in nearly matrix multiplication time.
\citet{demmel2026structured} further extend this approach to definite
Hermitian matrix pencils, using structure-preserving random perturbations.

\subsection{Curing singular eigenvalue problems}

A matrix pencil $\mat L(z)=\mat A_0-z\mat A_1$, with $\mat A_0,\mat A_1\in\C^{n\times n}$,
is called \emph{singular} if $\det \mat L(z)$ vanishes identically. Its
\emph{normal rank} is $r:=\nrank(\mat L)=\max_{z\in\C}\rank \mat L(z)$, and
$\lambda\in\C$ is a finite eigenvalue if $\rank \mat L(\lambda)<r$.
As a simple example, the pencil
\begin{equation}\label{eq:smallsing}
 \mat L(z)= \mat A_0-z\mat A_1 = 
 \begin{bmatrix}0&1&0\\0&0&0\\0&0&2\end{bmatrix}
 -z\begin{bmatrix}1&0&0\\0&0&0\\0&0&1\end{bmatrix}
\end{equation}
has normal rank $2$ and just one eigenvalue, $\lambda=2$.
Such singular eigenvalue problems frequently arise in a broad range of applications,
including control problems~\cite{vandooren1981generalized},
computational geometry~\cite{tu2009using},
aeroelastic flutter~\cite{pons2018multiparameter},
and wave propagation~\cite{kiefer2023computing}.

While the Kronecker canonical form~\cite{gantmacher1959theory} provides detailed theoretical insights, the numerical behavior of singular pencils has remained somewhat mysterious for a long time. Arbitrarily small perturbations can move eigenvalues to arbitrary locations in $\C$ and introduce additional, so-called spurious eigenvalues. For example, consider the perturbation
\begin{equation} \label{eq:badpert}
  \varepsilon \mat{E}_0 = 10^{-16} \times \begin{bmatrix}
  0 & 0 & 1 \\
  0 & 0 & 1 \\
  1 & 1 & 1
 \end{bmatrix}.
\end{equation}
The perturbed pencil $\mat{A}_0 +\varepsilon \mat{E}_0 - \lambda \mat{A}_1$ is regular and has eigenvalues $-1,\infty, \infty$, that is, the eigenvalue $2$ has been moved to $-1$ and two spurious, infinite eigenvalues have been created! 
Nevertheless, it was already observed by \citet{wilkinson1979kronecker}
that the finite-precision behavior of the QZ algorithm~\cite{golub2013matrix}, the standard method for computing eigenvalues of \emph{regular} pencils, is much more benign for singular pencils than such examples suggest.
A classical way to make
computations more reliable is to extract the regular part in a preprocessing step using a so-called 
staircase form~\cite{vandooren1979computation}. Randomization offers significantly simpler alternatives.

\subsubsection{Probabilistic eigenvalue sensitivity} 

A probabilistic perturbation analysis provides important insights into why the worst-case behavior caused by perturbations like~\cref{eq:badpert} is \emph{not} typical.
Consider a \emph{simple} finite eigenvalue $\lambda$ of a singular pencil $\mat L(z)$ and let $k=n-\nrank(\mat L)$.
The right and left nullspaces of $\mat L(\lambda)$ both have dimension
$k+1$. Of these, $k$ dimensions are inherited from the singular part
of the pencil: They arise as limits of the nullspaces of $\mat L(z)$
for generic $z\to\lambda$. Denote the corresponding subspaces by
$\mathcal S_R$ and $\mathcal S_L$, respectively.
Eigenvectors of a singular pencil are only determined
up to such singular directions.
We remove this ambiguity by choosing unit-norm vectors
$\vec x\in\ker\mat L(\lambda)$ and
$\vec y\in\ker\mat L(\lambda)^*$
orthogonal to $\mathcal S_R$ and $\mathcal S_L$, respectively.
These distinguished right and left eigenvectors define
\[
  \gamma_\lambda =
  \frac{|\vec y^*\mat A_1\vec x|}{\sqrt{1+|\lambda|^2}}.
\]
For a regular pencil, $\gamma_\lambda^{-1}$ is the usual (absolute) eigenvalue
condition number~\citep[\S7.7.3]{golub2013matrix}. For a singular pencil, \citet{lotz2020wilkinson}
show that the first-order amplification under a complex Gaussian random perturbation
is also unlikely to exceed a modest multiple of
$\gamma_\lambda^{-1}$, despite unbounded worst-case
sensitivity. For the $3\times 3$ example~\cref{eq:smallsing}, one has $\gamma_2=1/\sqrt5$ and thus the eigenvalue $\lambda=2$ will typically remain quite stable under small perturbations.  

This provides some intuition for the frequently observed benign behavior of the
QZ algorithm applied to singular pencils, but does not establish its
numerical reliability, because roundoff errors are not random.
The regularization methods discussed in the next subsection instead
introduce randomness deliberately to obtain regular, well-behaved pencils.

\subsubsection{Randomized regularization} \label{sec:randregular}

Suppose that the normal rank $r=n-k$ is known.\footnote{If the normal rank is not known a priori, it is commonly estimated by computing $\max_i\rank(\mat A_0-\xi_i\mat A_1)$ for a few randomly chosen
$\xi_i\in\C$.} Consider unitary matrices $\big[ \mat{\Omega}, \mat{\Omega}_\perp \big]  \in\C^{n\times n}$ and $\big[ \mat{\Psi}, \mat{\Psi}_\perp \big] \in\C^{n\times n}$
with $\mat{\Omega}, \mat{\Psi}\in\C^{n\times k}$. Then the rank-completing method of \citet{hochstenbach2019solving}
and its variants~\citeyearpar{hochstenbach2023solving} give three regularizations:
\begin{align}
 \mat L_{\rm p}(z)
 &=\mat{\Psi}_\perp^*\mat L(z)\mat{\Omega}_\perp,
 &&\text{projection},
 \label{eq:singular-compression}\\
 \mat L_{\rm rc}(z)
 &=\mat L(z)+\tau\mat{\Psi}\mat D(z)\mat{\Omega}^*,
 &&\text{rank completion},
 \label{eq:singular-rank-completion}\\
 \mat L_{\rm aug}(z)
 &=\begin{bmatrix}
      \mat L(z)&\mat{\Psi} \mat D_1(z) \\
      \mat D_2(z)\mat{\Omega}^*&0
    \end{bmatrix},
 &&\text{augmentation}.
 \label{eq:singular-border}
\end{align}
Here, $\tau\ne0$ need not be small, and
$\mat D, \mat D_1, \mat D_2$ are regular diagonal pencils, usually chosen  at random.
Generically, these constructions produce regular pencils of orders
$r$, $n$, and $n+k$, respectively, and \emph{preserve} the finite eigenvalues of $\mat L$ in exact arithmetic. Projection produces the smallest regular pencil and is the most efficient choice for dense problems. Rank completion is the analogue of the rank-$k$ modification of singular linear systems discussed in~\Cref{sec:singsystems}.
Augmentation can be attractive for large sparse problems when shifted systems with the bordered pencil can be solved efficiently.
In all three approaches, additional, so-called spurious eigenvalues may arise from the regularization of the singular part.

For rank completion, genuine eigenvalues can be distinguished from
spurious ones using the corresponding eigenvectors.
In exact arithmetic, any right and left eigenvectors
$\widehat{\vec x},\widehat{\vec y}$ of the regularized pencil
associated with a genuine eigenvalue satisfy
$\mat\Omega^*\widehat{\vec x}=0$ and
$\mat\Psi^*\widehat{\vec y}=0$.
For spurious eigenvalues, at least one of these conditions fails
generically. Thus, the norms of these quantities provide a practical
criterion for distinguishing genuine from spurious eigenvalues \cite{hochstenbach2019solving}.
Analogous tests are available for projection and augmentation~\cite{hochstenbach2023solving}.

It is common to choose independent random matrices
$\mat{\Omega},\mat{\Psi}$ by orthonormalizing complex Gaussian
matrices, which yields the uniform distribution on the complex
Stiefel manifold (cf.~\Cref{sec:dense-sketching}).
For this choice, \citet{kressner2024analysis} show that the 
eigenvalue condition number $\kappa_{\rm reg}(\lambda)$ for any of the three regularizations 
is unlikely to be significantly higher than the quantity $\gamma_\lambda^{-1}$ introduced above:
\[
  \prob\left\{
    \kappa_{\rm reg}(\lambda)>t \gamma^{-1}_\lambda
  \right\}
  \le {k^2(1+2 \log t)} t^{-2},
  \qquad t\ge 1.
\]
For the example~\Cref{eq:smallsing}, projection regularization yields a median eigenvalue condition number of about $5.2$, which is only modestly larger than $\gamma^{-1}_2 = \sqrt{5}$.

\begin{paragraph}{Extensions.}
Several recent developments extend the regularization techniques described above.
These include regularization of singular matrix polynomials and
structure-preserving variants for Hermitian pencils~\cite{hochstenbach2025solving,hochstenbach2025numerical}
as well as the analysis of eigenvalue conditioning under linearization~\cite{kressner2023singular}.
For large-scale sparse problems, \citet{meerbergen2026shift} combine a pivoting-based 
regularization with shift-and-invert Arnoldi.
\end{paragraph}

\subsection{Joint eigenvalue problems}\label{sec:joint-eigenvalues}

Given a family of $n\times n$ matrices $\mat A_1,\ldots, \mat A_d$,
joint eigenvalue problems seek common eigenvectors $\vec x \not=\vec 0$, 
$\mat A_j\vec x=\lambda_j\vec x$, and associated eigenvalue tuples
$(\lambda_1,\ldots,\lambda_d)$. For $d>1$, the existence of such common eigenvectors only holds under strong assumptions on the family, such as pairwise commuting matrices. These assumptions are easily destroyed by perturbations arising from,
for example, roundoff or modeling errors, so joint eigenvalue problems can only be solved approximately in practice.

It is common to address joint eigenvalue problems through optimization approaches~\cite{bunsegerstner1993numerical,pham2001joint}. A simpler alternative is to extract common eigenvectors from the eigenvectors of a random linear combination
\begin{equation} \label{eq:randlincomb}
  \mat A(\vec \mu)=\mu_1\mat A_1 + \mu_2\mat A_2 + \cdots + \mu_d\mat A_d, \qquad \vec \mu \sim \Normal(\vec 0, \Id_d).
\end{equation}
Such random linear combinations have been employed in numerous applications, including
the eigenvector method for multivariate polynomial root finding~\cite{corless1997reordered,cox2005using}, tensor
decompositions~\cite{kuleshov2015tensor,evert2022recursive}, learning from higher-order moments~\cite{anandkumar2012method}, and Prony methods for signal reconstruction~\cite{ehler2019randomized}.

To see why~\Cref{eq:randlincomb} makes sense,
consider pairwise commuting symmetric matrices
$\mat A_1,\ldots,\mat A_d\in\R^{n\times n}$.
If $\mat A_1$ has $n$ distinct eigenvalues, then any orthogonal
matrix $\mat Q$ diagonalizing $\mat A_1$ also diagonalizes the entire
family, yielding an orthonormal basis of common eigenvectors.
If, however, $\mat A_1$ has repeated eigenvalues, as in the extreme
case $\mat A_1=\mat I$, its eigenvectors need not form a common
eigenbasis.
Random linear combinations avoid this difficulty:
Distinct joint eigenvalue tuples are mapped, almost surely, to distinct
eigenvalues of $\mat A(\vec\mu)$, ensuring that any orthonormal
eigenvector basis of $\mat A(\vec\mu)$ jointly diagonalizes the family.
The analysis in \cite{he2024randomized} extends this observation to nearly commuting
matrices: If the given family is within distance
$\varepsilon$ of a commuting family, then an orthogonal matrix $\mat Q$
diagonalizing $\mat A(\vec\mu)$ jointly diagonalizes the given family
up to an off-diagonal error of norm $\order(\varepsilon)$ with high probability,
without any assumptions on eigenvalue gaps.
Thus, randomization enables the robust extraction of an approximate
joint eigenbasis without perturbing the given matrices.
The corresponding approximate joint eigenvalues are obtained from
Rayleigh quotients; see \citet{he2025randomized} for an extension to
a family of nonsymmetric matrices.

\section{Tensor-structured randomness} \label{sec:tensors}

Imposing tensor structure in random probe vectors and test matrices can drastically reduce computational cost in high-dimensional applications, at the price of introducing dependencies among their entries.

A (real) tensor of order $d$ is a multivariate array
$\tensor{X} \in \R^{n_1 \times \cdots \times n_d}$ with entries\[
\tensor{X}(i_1,\ldots,i_d), \quad i_{\mu} = 1,\ldots, n_\mu, \quad \mu = 1,\ldots,d.                                 
                                \]
For simplicity, throughout this section we assume that all modes have the
same size, $n_1=\cdots=n_d=n$, so that
$
    \tensor{X} \in \R^{n\times\cdots\times n}.
$
The vectorization $\vectorization(\tensor{X})$ is a vector of length $n^d$ that lists the entries of
$\tensor{X}$ with the first index varying fastest, followed by the second index, and so on.
For a matrix $\mat{A}$, this agrees with the usual definition of $\vectorization(\mat{A})$ obtained by stacking the columns of $\mat{A}$ on top of one another.

For a comprehensive introduction to tensor notation, basic operations, and decompositions, we refer to the recent book~\cite{ballard2025tensor} and the classic survey~\cite{kolda2009tensor}.

\subsection{Rank-one probe vectors} \label{sec:rankone}

Given a base distribution $\nu$ on $\R^n$, 
a random rank-one vector of length $n^d$ is obtained by vectorizing the outer product $\circ$ of 
independent random vectors $\vec{\omega}^{(1)}, \ldots, \vec{\omega}^{(d)} \sim \nu$: 
\begin{equation} \label{eq:rankonerandom}
 \vec{\omega} = \vectorization\!\big( \vec{\omega}^{(1)}\circ \vec{\omega}^{(2)} \circ \cdots\circ\vec{\omega}^{(d)} \big) = 
 \vec{\omega}^{(d)} \otimes \cdots \otimes \vec{\omega}^{(2)} \otimes \vec{\omega}^{(1)},
\end{equation}
where $\otimes$ denotes the Kronecker product. Rank-one tensors are also known as simple, pure, or decomposable tensors, and as product states in quantum physics.
In the probability literature, random tensors of the form~\cref{eq:rankonerandom} are often called
\emph{simple random tensors}.

\subsubsection{Benefits and basic guarantees}
Rank-one random vectors have two immediate computational advantages over
unstructured random vectors: they require only $dn$ instead of $n^d$ random
numbers to generate and store, and their interaction
with tensor-structured data can often be computed
without ever forming a vector of length $n^d$. Moreover, given an isotropic base
distribution $\nu$, independence of the factors implies that $\vec{\omega}$
is isotropic as well:
\begin{equation} \label{eq:sampleX}
 \expect[\vec{\omega}\vec{\omega}^\tp]
    =
    \expect\bigl[
       \vec{\omega}^{(d)}\vec{\omega}^{(d)\tp}
    \bigr]
    \otimes\cdots\otimes
    \expect\bigl[
       \vec{\omega}^{(1)}\vec{\omega}^{(1)\tp}
    \bigr]
    =
    \Id_{n^d}.
\end{equation}
Consequently, all elementary
unbiasedness properties of isotropic probe vectors remain valid:
For $\vec{x}\in\R^{n^d}$ and
$\mat{A}\in\R^{n^d\times n^d}$,
\begin{equation} \label{eq:unbiasedtensor}
    \expect\bigl[
        |\ip{\vec{x}}{\vec{\omega}}|^2
    \bigr]
    =
    \norm{\vec{x}}^2,
    \qquad
    \expect\bigl[
        \vec{\omega}^\tp \mat{A}\vec{\omega}
    \bigr]
    =
    \tr(\mat{A});
\end{equation}
see \Cref{sec:girardhutchinson} for the second equality, which implies that the Girard--Hutchinson trace estimator remains unbiased. 

\subsubsection{Small-ball probability estimates}
Consider $d=2$ and a fixed nonzero vector
$\vec{x}=\vectorization(\mat{X})$ for $\mat{X}\in\R^{n\times n}$. Note
that $\norm{\vec{x}}=\norm{\mat{X}}_{\rm F}$. Then~\cref{eq:unbiasedtensor} implies that
\[
    |\ip{\vec{x}}{\vec{\omega}}|^2
    =
    \big|(\vec{\omega}^{(1)})^\tp
    \mat{X}\vec{\omega}^{(2)}\big|^2
\]
is an unbiased estimator for $\norm{\mat{X}}_{\rm F}^2$. Evaluating this bilinear form in $\mat{X}$ can be much more attractive than probing $\mat{X}$ with an unstructured $n\times n$ random
matrix, especially when $n$ is very large and matrix-vector products with $\mat{X}$ are cheap. For independent
$\vec{\omega}^{(1)},\vec{\omega}^{(2)}\sim\normal{0}{\Id_n}$,
\citet{bujanovic2021norm} established
\begin{equation} \label{eq:bujkressner}
     \prob\left\{
       \norm{\mat{X}}_{\rm F}
       \geq \theta
       \big|(\vec{\omega}^{(1)})^\tp
       \mat{X}\vec{\omega}^{(2)}\big|
    \right\}
    \leq
    \frac{2}{\pi\theta}
    \bigl(2+\log(1+2\theta)\bigr), \qquad \theta > 1.
\end{equation}
Thus, a random bilinear form is unlikely to severely underestimate the
Frobenius norm; overestimation is easier to control, with probability at
most $\theta^{-2}$ for an overestimate by a factor $\theta$.
For comparison, the failure probability for an unstructured Gaussian
probe is $\order(1/\theta)$, so the price of imposing rank-one structure
is only a logarithmic factor in $\theta$. 
\citet{guttel2024randomized} used such random bilinear forms to construct inexpensive surrogates
for rational approximation of matrix-valued functions, with applications
to nonlinear eigenvalue problems.

The anti-concentration bound~\cref{eq:bujkressner} extends to arbitrary tensor order $d \ge 2$.
For an isotropic log-concave base distribution $\nu$,
\citet{hu2026small} established that there
are universal constants $c,C>0$ such that, for every nonzero $\vec{x}\in\R^{n^d}$,
\begin{equation}\label{eq:anticoncentration}
  \prob\left\{
    \norm{\vec{x}}
    \geq
    \theta |\ip{\vec{x}}{\vec{\omega}}|
  \right\}
  \leq
  \frac{1}{(d-1)!\,\theta}
  \bigl(C\log\theta\bigr)^{d-1},
  \qquad
  \theta>e^{cd}.
\end{equation}
The asymptotic dependence on $\theta$ is sharp when $\vec{x}$ has rank one.
Complementary concentration inequalities for rank-one random tensors were
developed, for example, in~\cite{vershynin2020concentration} and,
for quadratic forms, in~\cite{bamberger2022hanson}.

\subsubsection{The curse of dimensionality: Fourth moments}
The exponential growth of the lower bound on $\theta$ in the anti-concentration result~\cref{eq:anticoncentration} is only one manifestation of the many difficulties encountered with growing tensor order $d$. It is even more visible in the growth of higher moments. Given an isotropic random vector $\vec{\xi}\in\R^N$,
define its worst directional fourth moment by
\begin{equation} \label{eq:h4}
    h_4(\vec{\xi})
    :=
    \sup_{\norm{\vec{u}}=1}
    \expect\bigl[
        |\ip{\vec{u}}{\vec{\xi}}|^4
    \bigr].
\end{equation}
Assuming an isotropic base distribution and using that the factors of the rank-one random vector $\vec{\omega}$ are iid, it follows that
\begin{equation} \label{eq:tensor4moments}
     h_4(\vec{\omega})
    =
    h_4(\vec{\omega}^{(1)} )^d,
\end{equation}
with the worst direction $\vec{x}$, which maximizes $\expect[
        |\ip{\vec{x}}{\vec{\omega}}|^4]$, attained (once again) by a rank-one vector.
The value of $h_4(\vec{\omega}^{(1)})$ is $3$ for Gaussian and $3-2/n$ for Rademacher random vectors.
As noted by~\citet{meyer2026hutchinson}, the smallest possible value $3- 6 / (n+2)$ among all isotropic distributions on $\R^n$ is attained 
by the uniform distribution on the sphere of radius $\sqrt{n}$. 

\subsubsection{Trace estimation}
The exponential growth~\cref{eq:tensor4moments} is bad news for the reliability of
RNLA based on rank-one probes when $d$ becomes large.  To illustrate this, we consider the Girard--Hutchinson trace estimator from~\Cref{sec:girardhutchinson}.

Consider an arbitrary isotropic 
random vector $\vec{\xi}\in\R^N$ and a \psd
matrix $\mat{A} \in\R^{N\times N}$ with eigendecomposition
$\mat{A}
    =  \lambda_1 \vec{u}_1\vec{u}_1^\tp + \cdots + \lambda_N \vec{u}_N\vec{u}_N^\tp$.
By Cauchy--Schwarz,
\begin{align*}
\expect\bigl[
    (\vec{\xi}^\tp\mat{A}\vec{\xi})^2
\bigr]
&=
\sum_{i,j=1}^N
\lambda_i\lambda_j
\expect\big[
    (\vec{u}_i^\tp\vec{\xi})^2
    (\vec{u}_j^\tp\vec{\xi})^2
\big] \\
&\leq
\sum_{i,j=1}^N
\lambda_i\lambda_j
\sqrt{
    \expect[(\vec{u}_i^\tp\vec{\xi})^4]
    \expect[(\vec{u}_j^\tp\vec{\xi})^4]
} 
\leq
h_4(\vec{\xi})\,\tr(\mat{A})^2.
\end{align*}
Since isotropy gives
$\expect[\vec{\xi}^\tp\mat{A}\vec{\xi}]=\tr(\mat{A})$,
it follows that
\begin{equation}
\label{eq:h4-trace-variance}
    \Var\bigl(\vec{\xi}^\tp\mat{A}\vec{\xi}\bigr)
    \leq
    \bigl(h_4(\vec{\xi})-1\bigr)\tr(\mat{A})^2.
\end{equation}
Consequently, the Girard--Hutchinson estimator
$H_\ell(\mat{A}) = \big( \vec{\xi}_1^\tp\mat{A}\vec{\xi}_1 + \cdots + \vec{\xi}_\ell^\tp\mat{A}\vec{\xi}_\ell ) / \ell$
based on $\ell$ independent copies of $\vec{\xi}$ satisfies
\begin{equation}
\label{eq:h4-gh-variance}
    \Var(H_\ell(\mat{A}))
    \leq
    \frac{h_4(\vec{\xi})-1}{\ell}\tr(\mat{A})^2.
\end{equation}
This bound becomes an equality for a 
rank-one matrix
$\mat{A}=\vec{u}\vec{u}^\tp$, with $\vec{u}$ assuming a worst-case direction in~\cref{eq:h4}.

Because of~\cref{eq:tensor4moments}, the (tight) inequality~\cref{eq:h4-gh-variance} 
implies that a Girard--Hutchinson estimator based on 
independent copies of a rank-one random vector $\vec{\omega}$
requires up to
$\ell \sim h_4(\vec{\omega}^{(1)})^d \varepsilon^{-2}$ samples to achieve relative standard deviation $\varepsilon$; see~\cite{meyer2026hutchinson} for more details.
More generally, \citet{meyer2025understanding} show that such an exponential dependence on the tensor order is unavoidable when a matrix $\mat{A} \in \R^{n^d \times n^d}$ can be accessed only through products $\mat{A}\vec{\omega}$ with rank-one vectors $\vec{\omega}$. They establish exponential query lower bounds for several linear algebra tasks, including trace and top-eigenvalue estimation.

Despite their theoretical limitations, Frobenius norm and trace estimators based on rank-one probes have proved useful for inexpensive error estimation in tensor compression~\cite{aldaas2025adaptive} and for estimating partition functions of quantum spin systems~\cite[Section~6.2]{camao2025faster}.

\subsection{Khatri--Rao test matrices} \label{sec:khatrirao}

The Khatri--Rao product $\odot$ is the columnwise Kronecker product of matrices with the same number of columns.  It provides a convenient way to express matrices containing rank-one vectors.
Specifically, let
\[
    \mat{\Omega}^{(\mu)}
    =
    [\vec{\omega}^{(\mu)}_1,\ldots,\vec{\omega}^{(\mu)}_s]
    \in\R^{n\times s},
    \qquad \mu=1,\ldots,d,
\]
be independent random matrices with independent columns drawn from an isotropic base distribution
$\nu$.  Then
\begin{equation} \label{eq:krp}
    \mat{\Omega}
    =
    \frac{1}{\sqrt{s}}
    \bigl(
       \mat{\Omega}^{(d)}\odot\cdots\odot\mat{\Omega}^{(1)}
    \bigr)
    =
    \frac{1}{\sqrt{s}}
    [\vec{\omega}_1,\ldots,\vec{\omega}_s]
    \in\R^{n^d\times s},
\end{equation}
where $\vec{\omega}_j
=\vec{\omega}^{(d)}_j\otimes\cdots\otimes\vec{\omega}^{(1)}_j$.
This representation of $\mat{\Omega}$ requires only $dns$ random numbers.

\subsubsection{Oblivious subspace injections} \label{sec:KROSI}

The following result by~\citet{oliveira2016lower} shows that bounded
directional fourth moments imply the OSI property introduced in~\Cref{sec:embedinject}; see also~\cite{tropp2026comparisona}. 

\begin{theorem}[Fourth moment implies OSI]
\label{thm:osi-fourth-moment}
Let $\vec{\omega}\in\R^N$ be an isotropic random vector and suppose that its directional fourth moment $h_4(\vec{\omega})$ defined in~\cref{eq:h4} is bounded.
Consider 
$
    \mat{\Omega}
    =
    \frac{1}{\sqrt{s}}
    [\vec{\omega}_1,\ldots,\vec{\omega}_s]
$
for independent copies $\vec{\omega}_1,\ldots,\vec{\omega}_s$ of
$\vec{\omega}$.
Then, for $0<\eta,\delta<1$, the matrix $\mat{\Omega}$ is an oblivious
subspace injection with subspace dimension $r$, injectivity parameter $\eta$,
and failure probability $\delta$, provided that
\[
  s \geq
  C\,h_4(\vec{\omega})\,\eta^{-2}
  \big(r+\log \delta^{-1} \big)
\]
for some universal constant $C$.
\end{theorem}

Combined with~\cref{eq:tensor4moments}, \Cref{thm:osi-fourth-moment} yields OSI guarantees for Khatri--Rao
test matrices; a connection that was made explicit in~\cite{camao2025faster}.
In particular, constant injectivity is obtained with
$s=\order(3^d r)$ and $s=\order((3-2/n)^d r)$ for Khatri--Rao products of Gaussian and Rademacher random matrices, respectively.
For small tensor order $d$, the linear dependence of the sketch size on the subspace dimension makes \Cref{thm:osi-fourth-moment} a strong argument in favor of Khatri--Rao sketches; for larger $d$, however, the exponential dependence on $d$ makes the guarantee increasingly pessimistic.

\subsubsection{Oblivious subspace embeddings} \label{eq:KROSE}

The situation for OSE properties of Khatri--Rao products is somewhat less straightforward than for OSIs. \citet{ahle2020oblivious} establish a JL moment property for Khatri--Rao products whose factors have independent columns drawn
from isotropic sub-Gaussian random vector distributions.
Combined with a standard
$\varepsilon$-net argument (see~\Cref{sec:beyondsubspacemebddding}), this property implies that such an $n^d \times s$ Khatri--Rao product is an OSE for dimension $r$, failure
probability $\delta$, and \emph{fixed} distortion $\eta$, provided that
\begin{equation} \label{eq:krose}
s\gtrsim C^d(r+\log(1/\delta))^d, 
\end{equation}
where $C$ depends only on the sub-Gaussian norm of the distribution. We refer to~\cite{haselby2026fast,saibaba2025improved} for related results and details.
For $d=2$, \citet{bujanovic2026subspace} use AMM to improve the dependence of $s$ on $r$ from $r^2$ to $r^{3/2}$.
Recently, \citet{beretta2026tight} improved the dependence on the subspace dimension to nearly linear,
showing that
\[
 s \ge 2^{\order(d)} r \log( r/\delta )^{d+1}
\]
suffices for an OSE with fixed distortion.

\subsubsection{Application: Tucker compression} \label{sec:tuckercompression}
Khatri--Rao embeddings are particularly well suited for compressing
tensors of small to moderate order, where they can lead to substantial computational
savings while the exponential dependence on $d$ of the worst-case guarantees above may
still be acceptable. In the following, we highlight their use when performing a Tucker approximation. 
To simplify the discussion, we consider third-order tensors in this section.

\paragraph{Randomized Tucker approximation.}
A Tucker approximation of multilinear rank
$(r_1,r_2,r_3)$ to a tensor
$\tensor{X}\in\R^{n\times n \times n}$ takes the form
\begin{equation} \label{eq:tucker-format}
 \tensor{X}
 \approx
\tensor{G}\times_1\mat{Q}_1\times_2 \mat{Q}_2 \times_3\mat{Q}_3,
\end{equation}
with the \emph{core tensor} $\tensor{G}\in\R^{r_1\times r_2 \times r_3}$ and factor matrices
$\mat{Q}_\mu\in\R^{n\times r_\mu}$, which are often chosen to have orthonormal columns.
The symbol $\times_\mu$ denotes the $\mu$-mode matrix product, which multiplies a matrix
along mode $\mu$ of the tensor. For $r \coloneq \max\{r_1,r_2,r_3\} \ll n$, the Tucker approximation greatly reduces the memory requirements for $\tensor{X}$ from $\order(n^3)$ to $\order(nr+r^3)$.

There is a close connection between Tucker and low-rank matrix approximations. Let
\[
 \mat{X}_{(\mu)}
 \in \R^{n\times n^2}, \qquad \mu = 1,2,3,
\]
denote the \emph{mode-$\mu$ unfolding} of $\tensor{X}$, which groups the $\mu$th mode into the rows and the other two modes into the columns.
The truncated higher-order singular value decomposition (HOSVD) by~\citet{delathauwer2000multilinear}
computes an SVD of
$\mat{X}_{(\mu)}$, independently for each $\mu$, and lets the factor matrix $\mat{Q}_\mu$ contain the dominant $r_\mu$ left singular vectors. It then forms
the core 
\[
 \tensor{G}
 =
 \tensor{X}
 \times_1 \mat{Q}_1^\tp
 \times_2 \mat{Q}_2^\tp
 \times_3 \mat{Q}_3^\tp.
\]
The resulting Tucker approximation is quasi-optimal; its Frobenius norm error is within a 
factor $\sqrt{3}$ of the best approximation error among tensors of the same multilinear rank.
As the SVD of the unfolding $\mat{X}_{(\mu)}$ is fairly expensive, it is a natural idea to replace it by a randomized rangefinder (cf. \Cref{sec:rangefinder}), using a sketch $\mat{X}_{(\mu)} \mat{\Omega}_\mu$ with $\mat{\Omega}_\mu \in \R^{n^2\times s_\mu}$, $s_\mu  > r_\mu$~\cite{kressner2017recompression,minster2020randomized}, leading to \Cref{alg:kr-tucker}. In analogy to the optional 
step in \Cref{alg:rsvd}, one can optionally apply the HOSVD to the $s_1 \times s_2 \times s_3$ core computed in~\Cref{step:core} in order to reduce the multilinear rank 
of the approximation returned by~\Cref{alg:kr-tucker} to $(r_1, r_2, r_3)$.

\begin{algorithm}[t]
\caption{Randomized Tucker compression}
\label{alg:kr-tucker}
\begin{algorithmic}[1]
\Require Tensor $\tensor{X}\in\R^{n \times n \times n}$;
         sketch sizes $s_1,s_2,s_3$
\Ensure Tucker approximation
        $\widehat{\tensor{X}}
        =\tensor{G}\times_1\mat{Q}_1\times_2 \mat{Q}_2 \times_3\mat{Q}_3$
        
\State Generate random matrices $\mat \Omega_\mu \in\R^{n^2 \times s_\mu}$ for $\mu=1,2,3$
\For{$\mu=1,2,3$}
    \State Form the sketch
    $\mat{Y}_\mu\gets\mat{X}_{(\mu)}\mat{\Omega}_\mu
       \in\R^{n\times s_\mu}$
    \State Orthonormalize
    $(\mat{Q}_\mu,\sim)\gets
       \Call{QR}{\mat{Y}_\mu,\texttt{'econ'}}$
\EndFor
\State Compute the core
 $\tensor{G}
 \gets
 \tensor{X}
 \times_1 \mat{Q}_1^\tp
 \times_2 \mat{Q}_2^\tp
 \times_3 \mat{Q}_3^\tp$ \label{step:core}
\end{algorithmic}
\end{algorithm}

\begin{theorem}[Randomized Tucker compression]
\label{thm:kr-tucker}
Suppose that, for every $\mu=1,2,3$, the test matrix
$\mat{\Omega}_\mu$ in \Cref{alg:kr-tucker} is an
oblivious subspace injection for subspace dimension $r_\mu$, injectivity
parameter $\eta$, and failure probability at most $1/12$.
Then the output of \cref{alg:kr-tucker} satisfies
\begin{equation} \label{eq:kr-tucker-bound}
 \operatorname{Median}\!\big(
   \big\| {\tensor{X}-\widehat{\tensor{X}}} \big\|_{\rm F}^2
 \big)
 \leq
 5(1+\eta) \sum_{\mu=1}^3  \norm{\mat{X}_{(\mu)} - \lowrank{ \mat{X}_{(\mu)} }_{r_\mu}}_{\rm F}^2.
\end{equation}
In particular, if $\tensor{X}_\star$ denotes a best Tucker approximation
of multilinear rank $(r_1,r_2,r_3)$, then
\begin{equation} \label{eq:kr-tucker-quasioptimal}
 \operatorname{Median}\!\big(
   \big\| {\tensor{X}-\widehat{\tensor{X}}} \big\|_{\rm F}^2
 \big)
 \leq
 15(1+\eta)
 \norm{\tensor{X}-\tensor{X}_\star}_{\rm F}^2 .
\end{equation}
\end{theorem}
\begin{proof}
Setting $\mat{P}_\mu=\mat{Q}_\mu\mat{Q}_\mu^\tp$, a standard inequality in the analysis of HOSVD (see \cite[Theorem 5.1]{vannieuwenhoven2012new})
gives
\begin{equation} \label{eq:hosvdbound}
  \big\| {\tensor{X}-\widehat{\tensor{X}}} \big\|_{\rm F}^2
 \leq
 \sum_{\mu=1}^3
 \norm{(\Id-\mat{P}_\mu)\mat{X}_{(\mu)}}_{\rm F}^2.
\end{equation}
Setting $\mat{E}_\mu=\mat{X}_{(\mu)}-
\lowrank{\mat{X}_{(\mu)}}_{r_\mu}$, the bound~\cref{eq:structosi} derived in the proof of \Cref{thm:rsvd} implies that
\[
 \norm{(\Id-\mat{P}_\mu)\mat{X}_{(\mu)}}_{\rm F}^2
 \leq
 \norm{\mat{E}_\mu}_{\rm F}^2
 +(1+\eta)\norm{\mat{E}_\mu\mat{\Omega}_\mu}_{\rm F}^2
\]
holds with probability at least $11/12$.
The union bound, isotropy of $\mat{\Omega}_\mu$, and Markov's inequality give
\[
 \sum_{\mu=1}^3
 \norm{\mat{E}_\mu\mat{\Omega}_\mu}_{\rm F}^2
 \leq
 4\sum_{\mu=1}^3\norm{\mat{E}_\mu}_{\rm F}^2
\]
with probability at least $3/4$. Inserting these inequalities into~\Cref{eq:hosvdbound} concludes the proof of~\cref{eq:kr-tucker-bound}.
The second inequality~\cref{eq:kr-tucker-quasioptimal} follows from~\cref{eq:kr-tucker-bound} via a standard argument \cite[Property~10]{delathauwer2000multilinear}.
\end{proof}

Note that~\Cref{thm:kr-tucker} does not \emph{require} independence of $\mat{\Omega}_\mu\in\R^{n^2\times s_\mu}$ for different $\mu$.
For (unstructured) Gaussian embeddings, \Cref{thm:kr-tucker} applies with 
$s_\mu=\order(r_\mu)$ for fixed $\eta$. Thus,~\Cref{alg:kr-tucker} attains a constant-factor
error bound relative to the best Tucker approximation
of multilinear rank at most $(r_1,r_2,r_3)$. This uses only 
$\order(n^3r)$ operations for a dense input tensor $\tensor{X}$, 
comparing favorably with the $\order(n^4)$ operations required
by the classical HOSVD.

\paragraph{Khatri--Rao randomized Tucker (re)compression.} For a structured input tensor $\tensor{X}$, it can pay off to use 
Khatri--Rao test matrices~\cite{che2025efficient,kressner2017recompression,saibaba2025improved}.
Choose independent random matrices
$\mat{\Omega}^{(1)},\mat{\Omega}^{(2)},
\mat{\Omega}^{(3)}\in\R^{n\times s}$
with iid standard Gaussian entries, and form
\begin{equation} \label{eq:krtestmatrices}
\mat{\Omega}_1 = \frac{1}{\sqrt{s}} \mat{\Omega}^{(3)} \odot \mat{\Omega}^{(2)}, \quad 
 \mat{\Omega}_2 = \frac{1}{\sqrt{s}} \mat{\Omega}^{(3)} \odot \mat{\Omega}^{(1)}, \quad 
 \mat{\Omega}_3 = \frac{1}{\sqrt{s}} \mat{\Omega}^{(2)} \odot \mat{\Omega}^{(1)}.
\end{equation}
According to~\Cref{sec:KROSI}, the assumptions of \Cref{thm:kr-tucker} are satisfied for such Gaussian Khatri--Rao embeddings when
$s=\order(r)$.

A simple but important situation where~\cref{eq:krtestmatrices} gives an advantage is \emph{recompression}.
Consider a tensor that is already given in Tucker format:
\[
 \tensor{X}
 =\tensor{H}\times_1\mat{A}_1
            \times_2\mat{A}_2
            \times_3\mat{A}_3,
 \qquad
 \tensor{H} \in\R^{R\times R\times R},
 \quad \mat{A}_\mu\in\R^{n\times R},
\]
The multilinear rank $R$ is fairly high, the factor matrices $\mat{A}_\mu$ are not orthonormal, and one aims at recompressing 
$\tensor{X}$ to a Tucker approximation of much smaller multilinear rank $r \ll R$. The standard approach to such a recompression task is to compute thin QR decompositions of $\mat{A}_\mu$, multiply the resulting triangular factors into the core tensor $\tensor{H}$, followed by an HOSVD of the transformed $\tensor{H}$. This requires $\order(nR^2 + R^4)$ operations.

Carefully combining QR decompositions of the factor matrices with \Cref{alg:kr-tucker} using unstructured Gaussian test matrices reduces the cost to $\order(nR^2 + R^3 r)$ operations. This cost is reduced even further when applying 
\Cref{alg:kr-tucker} directly to $\tensor{X}$ using Khatri-Rao test matrices~\cref{eq:krtestmatrices}. For example, consider
\begin{equation}\label{eq:tucker-kr-fast}
 \mat{X}_{(1)}\mat{\Omega}_1
  =
 \frac{1}{\sqrt{s}}\,
 \mat{A}_1\mat{H}_{(1)}
 \bigl[
  (\mat{A}_3^\tp\mat{\Omega}^{(3)})
  \odot
  (\mat{A}_2^\tp\mat{\Omega}^{(2)})
 \bigr].
\end{equation}
The products with the factor matrices require $\order(nRr)$
operations, and the remaining contraction with the core costs
$\order(R^3 r)$ operations. In turn, the total complexity of~\Cref{alg:kr-tucker} reduces to 
$\order(nR r + R^3 r)$. Compared with unstructured Gaussian rangefinding,
this reduces the cost of processing the factor matrices
by a factor of order $R/r$.
The savings can be even more pronounced for elementwise (Hadamard) products of tensors in Tucker format, where the intermediate multilinear ranks grow
multiplicatively and randomized recompression with Khatri--Rao test matrices can exploit the
product structure without forming the enlarged core~\cite{kressner2017recompression}.

\subsection{Random tensor trains}

\label{sec:random-tt}

The theory discussed above predicts that rank-one random vectors may become ineffective when the tensor order $d$ is large. A natural way to increase the amount of randomness while retaining tensor structure is to use a (normalized) sum of $R$ independent rank-one random vectors~\cite{rakhshan2020tensorized,sun2021tensor} instead of a single rank-one vector. Existing concentration and embedding results, however, predict only limited improvements from this construction. Indeed, to control the directional fourth moment by a constant independent of $d$ would require $R$ to grow exponentially with $d$.

It turns out that tensor network formats offer a much more effective way of injecting additional randomness, avoiding exponential dependence on the tensor order for key probabilistic properties. This class includes, in particular, the widely used tensor train format~\cite{oseledets2011tensor}, known in the physics literature as matrix product states~\cite{schollwock2011density}.

\subsubsection{Tensor trains}

Let $\tensor{X}\in\R^{n\times\cdots\times n}$ be a tensor of order $d \ge 2$.
Given integers
\[
    \rho_0=\rho_d=1,
    \qquad
    \rho_1,\ldots,\rho_{d-1}\geq 1,
\]
a tensor train (TT) representation of $\tensor{X}$ consists of the so-called \emph{TT cores}
\[
    \tensor{X}_\mu
    \in
    \R^{\rho_{\mu-1}\times n\times\rho_\mu},
    \qquad \mu=1,\ldots,d.
\]
For fixed $i_\mu$, we let 
\begin{equation} \label{eq:matrixslices}
     \mat{X}_\mu(i_\mu)
    :=
    \tensor{X}_\mu(:,i_\mu,:)
    \in
    \R^{\rho_{\mu-1}\times\rho_\mu}
\end{equation}
denote the corresponding matrix slice of the $\mu$th core.
The TT representation is then compactly expressed as the matrix product
\begin{equation}
\label{eq:tt-representation}
    \tensor{X}(i_1,\ldots,i_d)
    =
    \mat{X}_1(i_1)
    \mat{X}_2(i_2)
    \cdots
    \mat{X}_d(i_d).
\end{equation}
The integers $\rho_1,\ldots,\rho_{d-1}$ are called 
\emph{bond dimensions} or, sometimes, TT representation ranks.
The componentwise minimal possible bond dimensions,
called the \emph{TT ranks} of $\tensor{X}$, are characterized by
\begin{equation} \label{eq:ttunfoldingrank}
 \rho_k = \rank\big( \mat{X}^{\le k} \big), \qquad k = 1,\ldots, d-1,
\end{equation}
where $\mat{X}^{\le k}  \in \R^{ n^k \times n^{d-k}}$
denotes the unfolding of $\tensor{X}$ that groups the first $k$ modes into the rows and the remaining modes into the columns.
This connection between matrix and TT ranks is one of the key structural properties of the TT
format and underlies many of its algorithms.

For simplicity, much of the discussion below concerns a uniform bond dimension $\rho_1=\cdots=\rho_{d-1}=\rho$. Then a TT representation of $\tensor{X}$ requires 
$\order(dn\rho^2)$ numbers.

\subsubsection{Gaussian random tensor trains} Commonly, a random TT is obtained by filling the TT cores
with independent random variables.  For example, let
all entries of
\[
    \tensor{G}_1\in\R^{1\times n\times\rho},
    \quad
    \tensor{G}_\mu\in\R^{\rho\times n\times\rho},
    \quad 2\leq\mu\leq d-1,
    \quad
    \tensor{G}_d\in\R^{\rho\times n\times1},
\]
be iid $\Normal(0,1)$, and define
$\tensor{W}$ by
\begin{equation}
\label{eq:gaussian-random-tt}
    \tensor{W}(i_1,\ldots,i_d)
    :=
    \rho^{-(d-1)/2}
    \mat{G}_1(i_1)\cdots\mat{G}_d(i_d).
\end{equation}
Then 
$\vec{\omega}:=\vectorization(\tensor{W})\in\R^{n^d}$
is called a \emph{Gaussian random TT (vector) of bond dimension $\rho$}.
One of the earliest uses of this random vector model in NLA appears in~\cite{huber2017randomized}.
As pointed out in~\cite{camano2026linear}, closely related random matrix product state models had already been studied earlier in quantum many-body physics.

The independence and normalization in~\cref{eq:gaussian-random-tt} ensure 
\begin{equation}
\label{eq:random-tt-isotropic}
    \expect[\vec{\omega}]=0,
    \qquad
    \expect[\vec{\omega}\vec{\omega}^{\tp}]
    =\Id_{n^d}.
\end{equation}
In particular, Gaussian random TT vectors are isotropic.

When $\rho=1$, \cref{eq:gaussian-random-tt} reduces exactly to the
Gaussian rank-one vector from \cref{eq:rankonerandom}.  Thus Gaussian
random TTs extend the rank-one construction from~\Cref{sec:rankone}, with
the bond dimension $\rho$ tuning the additional amount of randomness.

\subsubsection{Moments and the bond dimension}

\citet{rakhshan2020tensorized} show that the directional fourth moment introduced in~\cref{eq:h4} satisfies
\begin{equation}
\label{eq:random-tt-fourth}
    h_4(\vec{\omega})
    \leq
    3\left(1+\frac{2}{\rho}\right)^{d-1}
\end{equation}
for a Gaussian random TT of bond dimension $\rho$.
\citet{camano2026linear} obtain a finer description by computing the fourth-moment tensor.
\citet{bujanovic2026stochastic} consider moments of arbitrary order 
and derive the bound
\begin{equation}
\label{eq:random-tt-moments}
    \expect\big[
       (\vec{u}^{\tp}\vec{\omega})^{2p}
    \big]
    \leq
    (2p-1)!!
    \exp\left(\frac{p^2(d-1)}{\rho}\right)
\end{equation}
for every unit-norm vector $\vec{u}\in\R^{n^d}$ and integer $p\geq1$.
A striking feature of these bounds is that they remain bounded independently of the tensor order
once $\rho\gtrsim d$ . This already suggests that a modest bond dimension is sufficient to
remove exponential dependence on $d$ from many probabilistic guarantees.

\subsubsection{Oblivious subspace injections} 

Given independent Gaussian random TT vectors
$\vec{\omega}_1,\ldots,\vec{\omega}_s$, we form the test matrix
\begin{equation}
\label{eq:random-tt-test-matrix}
    \mat{\Omega}
    =
    \frac{1}{\sqrt{s}}
    [\vec{\omega}_1,\ldots,\vec{\omega}_s]
    \in\R^{n^d\times s}.
\end{equation}
By \cref{eq:random-tt-isotropic}, $\mat{\Omega}$ is an isotropic
embedding. The fourth-moment estimate \cref{eq:random-tt-fourth} can be combined
directly with \Cref{thm:osi-fourth-moment}.  In this way,
\citet{camano2026linear} and \citet{bujanovic2026stochastic} show that choosing
$\rho \geq d-1$ and $s \gtrsim \eta^{-2} \bigl(r + \log \delta^{-1} \bigr)$ gives an OSI
for $r$-dimensional subspaces 
with injectivity parameter $\eta$ and failure probability $\delta$. Note that such a TT test matrix is built from only
$\order(nd^3r)$ independent Gaussian random variables (for fixed $\eta$ and $\delta$).
In contrast, the corresponding guarantee for the Gaussian
Khatri--Rao construction from \Cref{sec:khatrirao} uses
$\order(nd\,3^d r)$ independent Gaussian random variables.

\subsubsection{Oblivious subspace embeddings and TTStack}

The first embedding results for Gaussian TT test matrices~\cref{eq:random-tt-test-matrix} predate the OSI theory discussed above. \citet{rakhshan2020tensorized} established a Johnson--Lindenstrauss property. Combined with a
standard $\varepsilon$-net argument, this yields an OSE, but just as for Khatri-Rao test matrices (cf. \Cref{eq:KROSE}), the resulting bound remains unfavorable for large tensor order.

A more flexible \emph{TTStack} construction was recently proposed by~\citet{cazeaux2026linear}. 
TTStack uses a block variant of Gaussian random TT vectors defined as follows. Let 
\[
    \tensor{G}_1\in\R^{\rho \times n\times\rho},
    \quad
    \tensor{G}_\mu\in\R^{\rho\times n\times\rho},
    \quad 2\leq\mu\leq d-1,
    \quad
    \tensor{G}_d\in\R^{\rho\times n\times1},
\]
have iid $\Normal(0,1)$ entries. Using the matrix slices
$\mat{G}_\mu(i_\mu)
=\tensor{G}_\mu(:,i_\mu,:)$,
define the order-$(d+1)$ tensor
$\tensor{W} \in \R^{\rho \times n \times \cdots \times n}$ by
\begin{equation}
\label{eq:gaussian-random-ttstack}
    \tensor{W}(:,i_1,\ldots,i_d)
    :=
    \rho^{-d/2}
    \mat{G}_1(i_1)\cdots\mat{G}_d(i_d).
\end{equation}
Sample $P$ independent such tensors $\tensor{W}_1, \ldots, \tensor{W}_P$ and turn them into matrices via mode-1 unfolding:
\[
 \mat{\Phi}_p := \mat{W}_p^{\le 1} \in \R^{\rho \times n^d}, \qquad p = 1,\ldots,P.
\]
The TTStack test matrix is then defined by
\begin{equation}
\label{eq:ttstack}
    \mat{\Omega}_{\rm TTS}^{\tp}
    =
    \frac{1}{\sqrt{P}}
    \begin{bmatrix}
       \mat{\Phi}_1\\
       \vdots\\
       \mat{\Phi}_P
    \end{bmatrix}
    \in\R^{P\rho \times n^d}.
\end{equation}
The normalization ensures $\expect[\mat{\Omega}_{\rm TTS}
         \mat{\Omega}_{\rm TTS}^\tp]
=\Id_{n^d}$.
The embedding dimension $s=P \rho$ depends on the two integers $P$ and $\rho$, reducing 
to a Khatri--Rao sketch for $\rho=1$. For $P = 1$,
TTStack reduces to the Gaussian TT map considered 
earlier by~\citet{ma2022cost}.
Their analysis, combined with a standard
$\varepsilon$-net argument, also yields an OSE
with an embedding dimension scaling only linearly with $d$.
However, tying the embedding dimension and the bond dimension together results in an
unfavorable dependence of the bond dimension on the distortion $\eta$.
TTStack avoids this drawback by decoupling the bond dimension from the embedding dimension through the additional parameter $P$.
\citet{cazeaux2026linear} show that an OSE for $r$-dimensional subspaces with distortion $\eta$ and
failure probability $\delta$ is obtained by TTStack with
\[
    \rho
    =
    \order\big(
       d\big(r+\log\delta^{-1}\big)
    \big),
    \ 
    P=\order(\eta^{-2}) \quad \Rightarrow \quad 
    s=P \rho 
    =
    \order\big(
       d\eta^{-2}\big(r+\log\delta^{-1}\big)
    \big).
\]
Thus, passing from the 
moment-based OSI guarantee for the Gaussian TT test matrix~\cref{eq:random-tt-test-matrix}
to an OSE increases the required bond dimension by a factor of order
$r+\log\delta^{-1}$ and the embedding dimension by a factor of order
$d$.

\subsubsection{Trace estimation}

A number of high-dimensional applications give rise to matrices
$\mat{A} \in \R^{n^d \times n^d}$ that admit highly compressed
tensor network representations, such as matrix product operators,
or can be applied efficiently using such representations.
This makes it possible to evaluate quadratic forms with TT vectors
without ever explicitly forming exponentially large matrices or vectors,
enabling efficient stochastic trace estimation with Gaussian random TT vectors.
For example, \citet{camano2026linear} estimate the partition
function $\tr(e^{-\beta\mat{H}})$ of a quantum spin system,
representing the Hamiltonian $\mat{H}$ as a matrix product operator
and approximating the action of $e^{-\beta\mat{H}}$ on Gaussian
random TT vectors using tensor-network methods.

Given \psd $\mat{A} \in \R^{n^d \times n^d}$, 
consider the Girard--Hutchinson trace estimator 
$H_\ell(\mat{A}) = \big( \vec{\omega}_1^\tp\mat{A}\vec{\omega}_1 + \cdots + \vec{\omega}_\ell^\tp\mat{A}\vec{\omega}_\ell \big) / \ell$ with $\ell$ independent 
Gaussian random TT vectors $\vec{\omega}_j$ of bond dimension $\rho$. Then~\cref{eq:h4-gh-variance}
combined with~\cref{eq:random-tt-fourth} immediately gives
\begin{equation}
\label{eq:tt-gh-simple}
    \Var\bigl(H_\ell(\mat{A})\bigr)
    \leq
    \frac{1}{\ell}
    \Big[
        3\left(1+\frac{2}{\rho}\right)^{d-1}-1
    \Big]
    \tr(\mat{A})^2.
\end{equation}
Choosing $\rho\gtrsim d$ makes the prefactor independent
of the tensor order, showing that Gaussian random TT vectors, in contrast to
rank-one probes, avoid exponential dependence on $d$.

There are two complementary ways to refine~\cref{eq:tt-gh-simple}.
First, the variance bound~\cref{eq:tt-gh-simple}, combined
with median-of-means amplification
\cite[Sec.~2.4]{vershynin2026high},
already implies that relative accuracy $\varepsilon$
and failure probability $\delta$ can be achieved
with bond dimension $\rho\geq d-1$ and
$\ell\gtrsim\varepsilon^{-2}\log(2/\delta)$ samples;
see~\cite{bujanovic2026stochastic}.
Alternatively, the higher-order moment bounds
\cref{eq:random-tt-moments} provide tail bounds
for the ordinary Girard--Hutchinson estimator,
without amplification,
at the expense of increasing the bond dimension to
$\rho=\order(d\log(1/\delta))$.
Second, by exploiting the more detailed fourth-moment analysis of Gaussian
TT vectors mentioned above, \citet{camano2026linear} show that
\Nystrom{}\plusplus achieves the optimal
$\order(\varepsilon^{-1})$ sample complexity for arbitrary
\psd\ matrices and fixed failure probability, provided that $\rho=\order(d/\varepsilon)$.

\subsubsection{Applications to sketching and tensor compression}

\label{sec:ttapplications}
Random TT tensors provide test matrices that can be
used to efficiently sketch high-dimensional data,
often without forming exponentially large arrays.
As in \Cref{sec:tuckercompression}, the computational
benefit is particularly pronounced when the input
already has a low-rank representation or arises
from operations with TT vectors and matrix product operators.

\paragraph{Recompression.}
\citet{huber2017randomized} were among the first to explore the use of random TT tensors for low-rank approximation, relating it to alternating least-squares approximation with random initialization.
\citet{aldaas2023randomized} develop randomized TT
recompression algorithms based on one-sided and two-sided
contractions with Gaussian random TT tensors.
Their \emph{randomize-then-orthogonalize} approach
follows the same principle as the use of Khatri--Rao
test matrices for Tucker recompression in
\Cref{sec:tuckercompression}: sketching avoids
an expensive initial orthogonalization.
\citet{ma2022cost} place such constructions in a broader framework
of Gaussian tensor network embeddings.
The HaTT algorithm of \citet{sun2026hatt} uses random TT contractions to recompress 
Hadamard products, paralleling the Tucker recompression approach of
\citet{kressner2017recompression}.

\paragraph{Streaming algorithms.}
\citet{kressner2023streaming} develop a streaming TT approximation,
 based on two-sided sketches in the spirit of the generalized
\Nystrom{} approximation. \citet{bucci2025sketched} use this strategy to compress linear
combinations of TT vectors in a TT variant of sketched GMRES.

\paragraph{Selected applications of random TT compression.}
\citet{yu2024reanchoring} use random TT sketches, constructed in a
weighted local basis, to recompress an ensemble of product-state
walkers for quantum Monte Carlo.
\citet[Section~6]{peng2025density} use random TT test tensors to
compress a smoothed empirical density, exploiting its representation
as a sum of separable terms before applying a deconvolution step.
\citet{strand2026adaptive} use two-sided Gaussian TT sketches to
recompress the growing sum of separable Gaussian functions defining
the bias potential in metadynamics.
\citet{verma2025lsh} use TT-structured random projections to construct
locality-sensitive hashes for approximate nearest-neighbor searches
with tensor data.
\citet{camano2026successive} develop a single-pass algorithm
that uses Khatri--Rao test matrices to construct a compressed
TT representation of the product of a matrix product operator
and a TT vector, an approach that suggests an extension to Gaussian random TT
test matrices.

\paragraph{Outlook.}
As illustrated in this section, there has been substantial
recent progress in the theory of random TT vectors and
test matrices, including OSI and OSE guarantees that
avoid exponential dependence on the tensor order.
Nevertheless, much remains to be explored in transferring
the broader range of techniques developed in this survey
to tensor-structured randomness and high-dimensional applications.

\section*{Acknowledgments}
The authors thank Chris Cama\~no, Tyler Chen, Alice Cortinovis, Jim Demmel, Stefan G\"uttel, Per-Gunnar Martinsson, Raphael Meyer, Kevin Miller, Riley Murray, Taejun Park, David Persson, Bor Plestenjak, Nathaniel Pritchard, Nian Shao, Igor Simunec, Nikhil Srivastava, Robert Webber, and Joel Tropp for helpful discussions and valuable feedback.

ENWE began work on this survey at the California Institute of Technology, where he was supported under the aegis of Joel Tropp by ONR Award N00014-24-1-2223 and the Caltech Carver Mead New
Adventures Fund.
He completed this work with support by the Miller Institute for Basic Research in Science, University of California Berkeley.
YN is supported by the EPSRC grant EP/Y030990/1.
Part of this work was conducted while the authors were visiting the Simons Institute for the Theory of Computing during the Fall 2025 program on Complexity and Linear Algebra; the support of the institute is gratefully acknowledged.

\clearpage 
\addcontentsline{toc}{section}{References}
\bibliography{bib}

@STRING(simax="SIAM J. Matrix Anal. Appl.")

@STRING(sisc="SIAM J. Sci. Comput.")

@STRING(mathpro="Math. Program.")

@STRING(apnum="Appl. Numer. Math.")

@STRING(acha="Appl. Comput. Harmon. Anal. ")

@STRING(laa="Linear Algebra Appl.")

@STRING(jmlr="J. Mach. Learn. Res.")

@STRING(annalsstats="Ann. Statist.")

@STRING(sinum="SIAM J. Numer. Anal.")

@STRING(sirev="SIAM Rev.")

@STRING(nlaa="Numer. Linear Algebra Appl.")

@STRING(nm="Numer. Math.")

@STRING(numalgs="Numer. Algorithms")

@STRING(pnas="Proc. Natl. Acad. Sci. USA")

@STRING(etna="Electron. Trans. Numer. Anal.")

@STRING(jacm="J. ACM")

@STRING(jfunct="J. Funct. Anal. ")

@STRING(fourappl="J. Fourier Anal. Appl.")

@STRING(mc="Math. Comp.")

@STRING(cpam="Comm. Pure Appl. Math.")

@STRING(imajna="IMA J. Numer. Anal.")

@STRING(acmtoms="ACM Trans. Math. Software")

@String{focm = "Found. Comput. Math."}

@String{actanum="Acta Numer."}

@String{actamath="Acta Math."}

@String{neurips="Adv. Neural Inf. Process. Syst."}

@String{tams="Trans. Amer. Math. Soc."}

@String{sicomp= "SIAM J. Comput."}

@String{siamuq= "SIAM/ASA J. Uncertain. Quantif."}

@String{simods= "SIAM J. Math. Data Sci."}

@String{bit= "BIT"}

@STRING{soda = "Proc. Annu. ACM-SIAM Symp. Discrete Algorithms"}

@STRING{stoc = "Proc. Annu. ACM Symp. Theory Comput."}

@STRING{focs = "Proc. Annu. IEEE Symp. Found. Comput. Sci."}

@String{icml        = "Proc. Int. Conf. Mach. Learn."}

@String{colt        = "Proc. Conf. Learn. Theory"}

@String{aistats     = "Proc. Int. Conf. Artif. Intell. Stat."}

@String{icalp       = "Proc. Int. Colloq. Autom. Lang. Program."}

@String{approxrandom= "Approx. Random. Comb. Optim. Algorithms Tech."}

@String{esa         = "Proc. Annu. Eur. Symp. Algorithms"}

@String{sosa        = "Proc. Symp. Simplicity Algorithms"}

@String{pods        = "Proc. ACM SIGACT-SIGMOD-SIGART Symp. Princ. Database Syst."}

@String{kdd         = "Proc. ACM SIGKDD Conf. Knowl. Discov. Data Min."}

@String{scworkshops = "Proc. SC Workshops High Perform. Comput. Netw. Storage Anal."}

@String{icm         = "Proc. Int. Congr. Math."}

@String{mathofdata  = "The Mathematics of Data"}

@String{compsens    = "Compressed Sensing and Its Applications"}

@String{relnumcomp  = "Reliable Numerical Computation"}

@String{structmat   = "Structured Matrices in Mathematics, Computer Science, and Engineering I"}

@STRING{jcss = "J. Comput. System Sci."}

@STRING{dcg = "Discrete Comput. Geom."}

@String{aada        = "Adv. Adapt. Data Anal."}

@String{advmath     = "Adv. Math."}

@String{annphys     = "Ann. Physics"}

@String{bernoulli   = "Bernoulli"}

@String{biometrika  = "Biometrika"}

@String{cams        = "Comm. Amer. Math. Soc."}

@String{cmp         = "Comm. Math. Phys."}

@String{cssc        = "Comm. Statist. Simulation Comput."}

@String{cpc         = "Comput. Phys. Comm."}

@String{ftml        = "Found. Trends Mach. Learn."}

@String{fttcs       = "Found. Trends Theor. Comput. Sci."}

@String{gmd         = "Geosci. Model Dev."}

@String{ieeetac     = "IEEE Trans. Automat. Control"}

@String{ieeetnn     = "IEEE Trans. Neural Netw."}

@String{ieeetpami   = "IEEE Trans. Pattern Anal. Mach. Intell."}

@String{ieeetsp     = "IEEE Trans. Signal Process."}

@String{infinf      = "Inf. Inference"}

@String{iclr        = "Int. Conf. Learn. Represent."}

@String{imrn        = "Int. Math. Res. Not. IMRN"}

@String{jrssb       = "J. R. Stat. Soc. Ser. B Stat. Methodol."}

@String{jsp         = "J. Stat. Phys."}

@String{noticesams  = "Notices Amer. Math. Soc."}

@String{prb         = "Phys. Rev. B"}

@String{prd         = "Phys. Rev. D"}

@String{probsurv    = "Probab. Surv."}

@String{ptrf        = "Probab. Theory Related Fields"}

@String{stspda      = "Sampl. Theory Signal Process. Data Anal."}

@String{siamnews    = "SIAM News"}

@String{tcs         = "Theoret. Comput. Sci."}

@String{toc         = "Theory Comput."}

@String{ussrcmmp    = "USSR Comput. Math. Math. Phys."}

@String{zpb         = "Z. Phys. B"}

@String{acom        = "Adv. Comput. Math."}

@article{meier2026reducing,
  author = {Meier, Maike and Lazzarino, Lorenzo and Shustin, Boris
and {Al Daas}, Hussam and Quinn, Paul},
  doi = {10.1107/S1600577526006491},
  journal = {J. Synchrotron Radiat.},
  number = {5},
  pages = {1440--1454},
  title = {Reducing acquisition time and radiation damage: data-driven subsampling for spectro-microscopy},
  volume = {33},
  year = {2026}
}

@article{halko2011finding,
  author = {Halko, Nathan and Martinsson, Per-Gunnar and
Tropp, Joel A.},
  doi = {10.1137/090771806},
  journal = sirev,
  month = jan,
  number = {2},
  pages = {217--288},
  title = {Finding Structure with Randomness: {Probabilistic} Algorithms for Constructing Approximate Matrix Decompositions},
  volume = {53},
  year = {2011}
}

@inproceedings{papadimitriou1998latent,
  author = {Papadimitriou, Christos H. and Tamaki, Hisao and
Raghavan, Prabhakar and Vempala, Santosh},
  booktitle = pods,
  doi = {10.1145/275487.275505},
  isbn = {978-0-89791-996-8},
  month = may,
  pages = {159--168},
  title = {Latent Semantic Indexing: {A} Probabilistic Analysis},
  year = {1998}
}

@inproceedings{sarlos2006improved,
  author = {Sarl{\'o}s, Tam{\'a}s},
  booktitle = focs,
  doi = {10.1109/FOCS.2006.37},
  month = oct,
  pages = {143--152},
  title = {Improved Approximation Algorithms for Large Matrices via Random Projections},
  year = {2006}
}

@article{martinsson2011randomized,
  author = {Martinsson, Per-Gunnar and Rokhlin, Vladimir and
Tygert, Mark},
  doi = {10.1016/j.acha.2010.02.003},
  journal = acha,
  number = {1},
  pages = {47--68},
  title = {A Randomized Algorithm for the Decomposition of Matrices},
  volume = {30},
  year = {2011}
}

@article{tropp2023randomized,
  author = {Tropp, Joel A. and Webber, Robert J.},
  journal = {arXiv preprint \href{http://arxiv.org/abs/2306.12418v3}{arXiv:2306.12418v3}},
  month = sep,
  title = {Randomized Algorithms for Low-Rank Matrix Approximation: {Design}, Analysis, and Applications},
  year = {2023}
}

@book{feller1971introduction,
  author = {Feller, William},
  edition = {2nd},
  publisher = {John Wiley \& Sons},
  title = {An Introduction to Probability Theory and Its Applications},
  volume = {2},
  year = {1971}
}

@article{chen2023spectrum,
  author = {Chen, Tyler},
  doi = {10.1063/5.0166678},
  journal = {J. Chem. Phys.},
  month = sep,
  number = {11},
  pages = {114101},
  title = {A Spectrum Adaptive Kernel Polynomial Method},
  urldate = {2026-09-02},
  volume = {159},
  year = {2023}
}

@inproceedings{drineas2006sampling,
  author = {Drineas, Petros and Mahoney, Michael W. and
Muthukrishnan, S.},
  booktitle = soda,
  month = jan,
  pages = {1127--1136},
  title = {Sampling Algorithms for {$\ell^2$} Regression and Applications},
  year = {2006}
}

@article{bartan2020distributed,
  author = {Bartan, Burak and Pilanci, Mert},
  journal = {arXiv preprint \href{http://arxiv.org/abs/2002.06538v2}{arXiv:2002.06538v2}},
  month = jun,
  title = {Distributed Sketching Methods for Privacy Preserving Regression},
  year = {2020}
}

@inproceedings{woodruff2013subspace,
  author = {David P. Woodruff and Qin Zhang},
  booktitle = colt,
  month = {June 12--14},
  pages = {546--567},
  title = {Subspace Embeddings and {$\ell_p$-Regression} Using Exponential Random Variables},
  url = {https://proceedings.mlr.press/v30/Woodruff13.html},
  year = {2013}
}

@article{baumann2026vac,
  author = {Baumann, Yves and Kyng, Rasmus and Z{\"o}cklein, Gernot},
  journal = {arXiv preprint \href{http://arxiv.org/abs/2609.20241v2}{arXiv:2609.20241v2}},
  title = {{VAC:} {A} volume-sampling-based elimination rule for approximate {Cholesky} factorization},
  year = {2026}
}

@misc{hsu2026randomized,
  author = {Alex W. Hsu and Ethan N. W. Epperly},
  note = {Manuscript in preparation},
  title = {Randomized row-norm estimation and applications},
  year = {2026}
}

@misc{heidary2026rerandomized,
  author = {Heidary, Diar},
  howpublished = {Palomar, PALOMAR-2026-09-09-000003 v1},
  title = {Rerandomized Subsampled Trigonometric Transforms},
  url = {https://palomar-registry.org/entry?id=PALOMAR-2026-09-09-000003&version=1},
  year = {2026}
}

@article{mai2026nelson,
  author = {Mai, Tung and Rao, Anup},
  journal = {arXiv preprint \href{http://arxiv.org/abs/2609.22548v1}{arXiv:2609.22548v1}},
  title = {The {Nelson}--{Nguyen} conjecture via mean-to-moments concentration},
  year = {2026}
}

@article{greene2022approximating,
  author = {Greene, Samuel M. and Webber, Robert J. and Berkelbach, Timothy C. and Weare, Jonathan},
  doi = {10.1137/21M1422513},
  journal = sisc,
  number = {5},
  pages = {A3067--A3097},
  title = {Approximating matrix eigenvalues by subspace iteration with repeated random sparsification},
  volume = {44},
  year = {2022}
}

@article{ahfock2021statistical,
  author = {Ahfock, D C and Astle, W J and Richardson, S},
  doi = {10.1093/biomet/asaa062},
  journal = biometrika,
  month = jun,
  number = {2},
  pages = {283--297},
  title = {Statistical Properties of Sketching Algorithms},
  volume = {108},
  year = {2021}
}

@article{martinsson2020randomized,
  author = {Martinsson, Per-Gunnar and Tropp, Joel A.},
  doi = {10.1017/S0962492920000021},
  journal = actanum,
  language = {en},
  month = may,
  pages = {403--572},
  title = {Randomized numerical linear algebra: {Foundations} and algorithms},
  volume = {29},
  year = {2020}
}

@incollection{davidson2001local,
  author = {Davidson, Kenneth R. and Szarek, Stanislaw J.},
  booktitle = {Handbook of the Geometry of {Banach} Spaces},
  pages = {317--366},
  publisher = {North-Holland},
  title = {Local operator theory, random matrices and {Banach} spaces},
  volume = {1},
  year = {2001}
}

@phdthesis{epperly2025make,
  author = {Epperly, Ethan N.},
  doi = {10.7907/pef3-mg80},
  school = {California Institute of Technology},
  title = {Make the Most of What You Have: {Resource-efficient} Randomized Algorithms for Matrix Computations},
  year = {2025}
}

@article{drineas2011faster,
  author = {Drineas, Petros and Mahoney, Michael W. and
Muthukrishnan, S. and Sarl{\'o}s, Tam{\'a}s},
  doi = {10.1007/s00211-010-0331-6},
  journal = nm,
  month = feb,
  number = {2},
  pages = {219--249},
  title = {Faster Least Squares Approximation},
  volume = {117},
  year = {2011}
}

@article{oymak2018universality,
  author = {Oymak, Samet and Tropp, Joel A},
  doi = {10.1093/imaiai/iax011},
  journal = infinf,
  language = {en},
  month = sep,
  number = {3},
  pages = {337--446},
  title = {Universality laws for randomized dimension reduction, with applications},
  volume = {7},
  year = {2018}
}

@inproceedings{dobriban2019asymptotics,
  author = {Dobriban, Edgar and Liu, Sifan},
  booktitle = neurips,
  title = {Asymptotics for Sketching in Least Squares Regression},
  url = {https://dl.acm.org/doi/abs/10.5555/3454287.3454617},
  year = {2019}
}

@inproceedings{lacotte2020optimal,
  author = {Lacotte, Jonathan and Liu, Sifan and Dobriban, Edgar and
Pilanci, Mert},
  booktitle = neurips,
  pages = {9725--9735},
  title = {Optimal Iterative Sketching Methods with the Subsampled Randomized {Hadamard} Transform},
  url = {https://dl.acm.org/doi/10.5555/3495724.3496539},
  year = {2020}
}

@article{epperly2026sharp,
  author = {Epperly, Ethan N. and Webber, Robert J.},
  journal = {arXiv preprint \href{http://arxiv.org/abs/2605.19096v2}{arXiv:2605.19096v2}},
  title = {Sharp analysis of sketched least squares and randomized low-rank approximation},
  year = {2026}
}

@inproceedings{guzmancordero2025improving,
  author = {Guzm{\'a}n-Cordero, Andr{\'e}s and Dangel, Felix and Goldshlager, Gil and Zeinhofer, Marius},
  booktitle = neurips,
  doi = {10.52202/085713-3798},
  pages = {113870--113900},
  title = {Improving Energy Natural Gradient Descent through {Woodbury}, Momentum, and Randomization},
  year = {2025}
}

@incollection{dedamas2026randomized,
  author = {de Damas, Jean-Guillaume and Grigori, Laura and Simunec, Igor and Timsit, Edouard},
  booktitle = {Numerical Analysis and Scientific Computing: an Anthology},
  doi = {10.1007/978-3-032-33428-2_6},
  editor = {Koren, Barry and Roose, Dirk},
  pages = {151--196},
  publisher = {Springer},
  series = {Lecture Notes in Computational Science and Engineering},
  title = {Randomized Orthogonalization and {Krylov} Subspace Methods: {Principles} and Algorithms},
  year = {2026}
}

@article{dedamas2025randomized,
  author = {de Damas, Jean-Guillaume and Grigori, Laura},
  doi = {10.1137/24M1674352},
  journal = simax,
  number = {4},
  pages = {2395--2422},
  title = {Randomized implicitly restarted {Arnoldi} method for the non-symmetric eigenvalue problem},
  volume = {46},
  year = {2025}
}

@article{dong2023simpler,
  author = {Dong, Yijun and Martinsson, Per-Gunnar},
  doi = {10.1007/s10444-023-10061-z},
  journal = acom,
  number = {4},
  pages = {66},
  title = {Simpler is better: a comparative study of randomized pivoting algorithms for {CUR} and interpolative decompositions},
  volume = {49},
  year = {2023}
}

@article{melnichenko2025choleskyqr,
  author = {Melnichenko, Maksim and Balabanov, Oleg and Murray, Riley and Demmel, James and Mahoney, Michael W. and Luszczek, Piotr},
  doi = {10.1137/24M163712X},
  journal = simax,
  month = {September},
  number = {3},
  pages = {1701--1734},
  title = {{CholeskyQR} with Randomization and Pivoting for Tall Matrices {(CQRRPT)}},
  volume = {46},
  year = {2025}
}

@article{camao2025faster,
  author = {Camaño, Chris and Epperly, Ethan N. and Meyer, Raphael A. and Tropp, Joel A.},
  journal = {arXiv preprint \href{http://arxiv.org/abs/2508.21189v1}{arXiv:2508.21189v1}},
  month = aug,
  title = {Faster Linear Algebra Algorithms with Structured Random Matrices},
  year = {2025}
}

@inproceedings{chen2026gpu,
  address = {New Orleans, LA, USA},
  author = {Chen, Tyler and Niroula, Pradeep and Ray, Archan and Subrahmanya, Pragna and Pistoia, Marco and Kumar, Niraj},
  booktitle = {2026 {IEEE} International Parallel and Distributed Processing Symposium Workshops ({IPDPSW})},
  doi = {10.1109/IPDPSW71298.2026.00074},
  pages = {411--420},
  title = {{GPU}-parallelizable randomized sketch-and-precondition for linear regression using sparse sign sketches},
  year = {2026}
}

@article{tropp2017practical,
  author = {Tropp, Joel A. and Yurtsever, Alp and
Udell, Madeleine and Cevher, Volkan},
  doi = {10.1137/17M1111590},
  journal = simax,
  month = jan,
  number = {4},
  pages = {1454--1485},
  title = {Practical Sketching Algorithms for Low-Rank Matrix Approximation},
  volume = {38},
  year = {2017}
}

@inproceedings{chenakkod2026optimal,
  author = {Chenakkod, Shabarish and Dereziński, Michal and Dong, Xiaoyu},
  booktitle = soda,
  doi = {10.1137/1.9781611978971.163},
  pages = {4501--4559},
  title = {Optimal Subspace Embeddings: {Resolving} {Nelson}-{Nguyen} Conjecture up to Sub-Polylogarithmic Factors},
  year = {2026}
}

@inproceedings{nelson2013osnap,
  author = {Nelson, Jelani and Nguy\~{\^{e}}n, Huy L.},
  booktitle = focs,
  doi = {10.1109/FOCS.2013.21},
  pages = {117--126},
  title = {{OSNAP:} {Faster} Numerical Linear Algebra Algorithms via Sparser Subspace Embeddings},
  year = {2013}
}

@inproceedings{cohen2016nearly,
  author = {Cohen, Michael B.},
  booktitle = soda,
  doi = {10.1137/1.9781611974331.ch21},
  isbn = {978-1-61197-433-1},
  month = jan,
  pages = {278--287},
  title = {Nearly Tight Oblivious Subspace Embeddings by Trace Inequalities},
  year = {2016}
}

@inproceedings{chenakkod2024optimal,
  author = {Chenakkod, Shabarish and Derezi{\'n}ski, Micha{\l} and
Dong, Xiaoyu and Rudelson, Mark},
  booktitle = stoc,
  doi = {10.1145/3618260.3649762},
  isbn = {979-8-4007-0383-6},
  month = jun,
  pages = {1106--1117},
  title = {Optimal Embedding Dimension for Sparse Subspace Embeddings},
  year = {2024}
}

@inproceedings{chenakkod2025optimal,
  author = {Chenakkod, Shabarish and Derezi\'{n}ski, Micha{\l} and Dong, Xiaoyu},
  booktitle = icalp,
  doi = {10.4230/LIPIcs.ICALP.2025.55},
  isbn = {978-3-95977-372-0},
  pages = {55:1--55:20},
  title = {Optimal Oblivious Subspace Embeddings with Near-Optimal Sparsity},
  year = {2025}
}

@article{tropp2026comparisona,
  author = {Joel A. Tropp},
  doi = {10.1090/cams/72},
  journal = cams,
  pages = {844--895},
  volume = {6},
  number = {16},
  title = {Comparison Theorems for the Minimum Eigenvalue of a Random Positive-Semidefinite Matrix},
  year = {2026}
}

@incollection{nelson2014lower,
  author = {Nelson, Jelani and Nguy\~{\^{e}}n, Huy L.},
  booktitle = icalp,
  doi = {10.1007/978-3-662-43948-7_73},
  isbn = {978-3-662-43947-0 978-3-662-43948-7},
  pages = {883--894},
  publisher = {Springer},
  title = {Lower Bounds for Oblivious Subspace Embeddings},
  volume = {8572},
  year = {2014}
}

@article{huang2026well,
  author = {Huang, Han and Rudelson, Mark and
Tikhomirov, Konstantin},
  journal = {arXiv preprint \href{http://arxiv.org/abs/2607.05384v2}{arXiv:2607.05384v2}},
  month = jul,
  title = {Well-invertible column subsets of sparse matrices are rare},
  year = {2026}
}

@article{achlioptas2003database,
  author = {Achlioptas, Dimitris},
  doi = {10.1016/S0022-0000(03)00025-4},
  journal = jcss,
  number = {4},
  pages = {671--687},
  title = {Database-friendly random projections: {Johnson-Lindenstrauss} with binary coins},
  volume = {66},
  year = {2003}
}

@article{charikar2004finding,
  author = {Charikar, Moses and Chen, Kevin and
Farach-Colton, Martin},
  doi = {10.1016/S0304-3975(03)00400-6},
  journal = tcs,
  language = {en},
  month = jan,
  number = {1},
  pages = {3--15},
  series = {Automata, {Languages} and {Programming}},
  title = {Finding frequent items in data streams},
  volume = {312},
  year = {2004}
}

@inproceedings{clarkson2009numerical,
  author = {Clarkson, Kenneth L. and Woodruff, David P.},
  booktitle = stoc,
  doi = {10.1145/1536414.1536445},
  isbn = {978-1-60558-506-2},
  month = may,
  pages = {205--214},
  title = {Numerical linear algebra in the streaming model},
  year = {2009}
}

@article{clarkson2017low,
  author = {Clarkson, Kenneth L. and Woodruff, David P.},
  doi = {10.1145/3019134},
  journal = jacm,
  number = {6},
  pages = {54:1--54:45},
  title = {Low-rank approximation and regression in input sparsity time},
  volume = {63},
  year = {2017}
}

@inproceedings{meng2013low,
  author = {Meng, Xiangrui and Mahoney, Michael W.},
  booktitle = stoc,
  doi = {10.1145/2488608.2488621},
  isbn = {978-1-4503-2029-0},
  language = {en},
  month = jun,
  pages = {91--100},
  title = {Low-distortion subspace embeddings in input-sparsity time and applications to robust linear regression},
  year = {2013}
}

@inproceedings{kane2012sparser,
  author = {Kane, Daniel M. and Nelson, Jelani},
  booktitle = soda,
  pages = {1195--1206},
  title = {Sparser {Johnson}-{Lindenstrauss} transforms},
  url = {https://dl.acm.org/doi/abs/10.5555/2095116.2095210},
  year = {2012}
}

@inproceedings{andoni2015practical,
  author = {Andoni, Alexandr and Indyk, Piotr and
Laarhoven, Thijs and Razenshteyn, Ilya and
Schmidt, Ludwig},
  booktitle = neurips,
  month = dec,
  pages = {1225--1233},
  title = {Practical and optimal {LSH} for angular distance},
  url = {https://dl.acm.org/doi/10.5555/2969239.2969376},
  year = {2015}
}

@article{tropp2011improved,
  author = {Tropp, Joel A.},
  doi = {10.1142/S1793536911000787},
  journal = aada,
  month = apr,
  number = {01n02},
  pages = {115--126},
  title = {Improved Analysis of the Subsampled Randomized {Hadamard} Transform},
  volume = {3},
  year = {2011}
}

@article{sorensen1993efficient,
  author = {Sorensen, Henrik V. and Burrus, C. Sidney},
  doi = {10.1109/78.205723},
  journal = ieeetsp,
  number = {3},
  pages = {1184--1200},
  title = {Efficient computation of the {DFT} with only a subset of input or output points},
  volume = {41},
  year = {1993}
}

@article{ailon2009fast,
  author = {Ailon, Nir and Liberty, Edo},
  doi = {10.1007/s00454-008-9110-x},
  journal = dcg,
  language = {en},
  month = dec,
  number = {4},
  pages = {615--630},
  title = {Fast Dimension Reduction Using {Rademacher} Series on Dual {BCH} Codes},
  volume = {42},
  year = {2009}
}

@inproceedings{ailon2006approximate,
  author = {Ailon, Nir and Chazelle, Bernard},
  booktitle = stoc,
  doi = {10.1145/1132516.1132597},
  isbn = {978-1-59593-134-4},
  month = may,
  pages = {557--563},
  title = {Approximate Nearest Neighbors and the Fast {Johnson}-{Lindenstrauss} Transform},
  year = {2006}
}

@article{woolfe2008fast,
  author = {Woolfe, Franco and Liberty, Edo and Rokhlin, Vladimir and
Tygert, Mark},
  doi = {10.1016/j.acha.2007.12.002},
  journal = acha,
  month = nov,
  number = {3},
  pages = {335--366},
  title = {A Fast Randomized Algorithm for the Approximation of Matrices},
  volume = {25},
  year = {2008}
}

@article{tropp2019streaming,
  author = {Tropp, Joel A. and Yurtsever, Alp and
Udell, Madeleine and Cevher, Volkan},
  doi = {10.1137/18M1201068},
  journal = sisc,
  month = jan,
  number = {4},
  pages = {A2430--A2463},
  title = {Streaming Low-Rank Matrix Approximation with an Application to Scientific Simulation},
  volume = {41},
  year = {2019}
}

@article{cartis2021hashing,
  author = {Cartis, Coralia and Fiala, Jan and Shao, Zhen},
  month = may,
  journal = {arXiv preprint \href{http://arxiv.org/abs/2105.11815}{arXiv:2105.11815}},
  title = {Hashing Embeddings of Optimal Dimension, with Applications to Linear Least Squares},
  year = {2021}
}

@book{ballard2025tensor,
  author = {Grey Ballard and Tamara G. Kolda},
  doi = {10.1017/9781009471664},
  publisher = {Cambridge University Press},
  title = {Tensor Decompositions for Data Science},
  year = {2025}
}

@article{beckermann2019bounds,
  author = {Beckermann, Bernhard and Townsend, Alex},
  doi = {10.1137/19M1244433},
  journal = sirev,
  number = {2},
  pages = {319--344},
  title = {Bounds on the singular values of matrices with displacement structure},
  volume = {61},
  year = {2019}
}

@article{udell2019why,
  author = {Udell, Madeleine and Townsend, Alex},
  doi = {10.1137/18M1183480},
  journal = simods,
  month = jan,
  number = {1},
  pages = {144--160},
  title = {Why Are Big Data Matrices Approximately Low Rank?},
  volume = {1},
  year = {2019}
}

@article{webb2026low,
  author = {Webb, Marcus},
  doi = {10.1093/imanum/drag083},
  journal = imajna,
  note = {Published online 22 June 2026},
  title = {Low-rank approximation of analytic kernels},
  year = {2026}
}

@article{benner2015survey,
  author = {Benner, Peter and G{\"u}{\u{g}}ercin, Serkan and
Willcox, Karen},
  doi = {10.1137/130932715},
  journal = sirev,
  number = {4},
  pages = {483--531},
  title = {A survey of projection-based model reduction methods for parametric dynamical systems},
  volume = {57},
  year = {2015}
}

@book{jolliffe2002principal,
  author = {Jolliffe, I. T.},
  doi = {10.1007/b98835},
  edition = {second},
  publisher = {Springer-Verlag},
  title = {Principal Component Analysis},
  year = {2002}
}

@inproceedings{williams2000using,
  author = {Williams, Christopher K. I. and Seeger, Matthias},
  booktitle = neurips,
  month = jan,
  pages = {661--667},
  title = {Using the {Nyström} method to speed up kernel machines},
  url = {https://dl.acm.org/doi/10.5555/3008751.3008847},
  year = {2000}
}

@article{drineas2005nystrom,
  author = {Drineas, Petros and Mahoney, Michael W.},
  journal = jmlr,
  number = {72},
  pages = {2153--2175},
  title = {On the {Nystrom} Method for Approximating a {Gram} Matrix for Improved Kernel-Based Learning},
  url = {https://dl.acm.org/doi/10.5555/1046920.1194916},
  volume = {6},
  year = {2005}
}

@article{bucci2025numerical,
  author = {Bucci, Alberto and Nakatsukasa, Yuji and
Park, Taejun},
  journal = {arXiv preprint \href{http://arxiv.org/abs/2511.15583v2}{arXiv:2511.15583v2}},
  title = {Numerical Stability of the {Nystr\"om} Method},
  year = {2025}
}

@inproceedings{andoni2013eigenvalues,
  author = {Andoni, Alexandr and Nguy\~{\^{e}}n, Huy L.},
  booktitle = soda,
  doi = {10.1137/1.9781611973105.124},
  isbn = {9781611972511},
  pages = {1729--1737},
  title = {Eigenvalues of a matrix in the streaming model},
  year = {2013}
}

@article{meier2024fast,
  author = {Meier, Maike and Nakatsukasa, Yuji},
  doi = {10.1016/j.laa.2024.01.001},
  journal = laa,
  pages = {1--32},
  title = {Fast randomized numerical rank estimation},
  volume = {686},
  year = {2024}
}

@article{rokhlin2008fast,
  author = {Rokhlin, Vladimir and Tygert, Mark},
  doi = {10.1073/pnas.0804869105},
  journal = pnas,
  month = sep,
  number = {36},
  pages = {13212--13217},
  title = {A Fast Randomized Algorithm for Overdetermined Linear Least-Squares Regression},
  volume = {105},
  year = {2008}
}

@article{meng2014lsrn,
  author = {Meng, Xiangrui and Saunders, Michael A. and
Mahoney, Michael W.},
  doi = {10.1137/120866580},
  journal = sisc,
  language = {en},
  month = jan,
  number = {2},
  pages = {C95--C118},
  title = {{LSRN:} {A} Parallel Iterative Solver for Strongly Over- or Underdetermined Systems},
  volume = {36},
  year = {2014}
}

@article{epperly2026fasta,
  author = {Epperly, Ethan N. and Meier, Maike and
Nakatsukasa, Yuji},
  doi = {10.1002/cpa.70013},
  journal = cpam,
  number = {2},
  pages = {293--339},
  title = {Fast randomized least-squares solvers can be just as accurate and stable as classical direct solvers},
  volume = {79},
  year = {2026}
}

@article{epperly2026fastb,
  author = {Epperly, Ethan N. and Park, Taejun and
Nakatsukasa, Yuji},
  journal = {arXiv preprint \href{http://arxiv.org/abs/2602.16797v1}{arXiv:2602.16797v1}},
  title = {Fast, High-Accuracy, Randomized Nullspace Computations for Tall Matrices},
  year = {2026}
}

@article{balabanov2022randomizeda,
  author = {Balabanov, Oleg and Grigori, Laura},
  doi = {10.1137/20M138870X},
  journal = sisc,
  number = {3},
  pages = {A1450--A1474},
  title = {Randomized {Gram}--{Schmidt} Process with Application to {GMRES}},
  volume = {44},
  year = {2022}
}

@book{stewart1998matrix,
  author = {Stewart, Gilbert W.},
  doi = {10.1137/1.9781611971408},
  isbn = {978-0-89871-414-2},
  publisher = {SIAM},
  title = {Matrix algorithms. 1: {Basic} decompositions},
  year = {1998}
}

@article{yamamoto2015roundoff,
  author = {Yamamoto, Y. and Nakatsukasa, Y. and Yanagisawa, Y. and
Fukaya, T.},
  journal = etna,
  pages = {306--326},
  title = {Roundoff error analysis of the {CholeskyQR2} algorithm},
  volume = {44},
  year = {2015}
}

@article{fan2021novel,
  author = {Fan, Yuwei and Guo, Yixiao and Lin, Ting},
  journal = {arXiv preprint \href{http://arxiv.org/abs/2111.11148v2}{arXiv:2111.11148v2}},
  title = {A novel randomized {XR}-based preconditioned {CholeskyQR} algorithm},
  year = {2021}
}

@article{balabanov2022randomizedb,
  author = {Balabanov, Oleg},
  journal = {arXiv preprint \href{http://arxiv.org/abs/2210.09953v2}{arXiv:2210.09953v2}},
  title = {Randomized {Cholesky} {QR} factorizations},
  year = {2022}
}

@article{higgins2025analysis,
  author = {Andrew J. Higgins and Daniel B. Szyld and Erik G. Boman
and Ichitaro Yamazaki},
  doi = {10.1007/s00211-025-01492-5},
  journal = nm,
  pages = {1695--1737},
  title = {Analysis of Randomized {Householder}--{Cholesky} {QR} Factorization with Multisketching},
  volume = {157},
  year = {2025}
}

@article{garrison2024randomized,
  author = {Garrison, James E. and Ipsen, Ilse C. F.},
  journal = {arXiv preprint \href{http://arxiv.org/abs/2406.11751v1}{arXiv:2406.11751v1}},
  month = jun,
  title = {A randomized preconditioned {Cholesky}-{QR} algorithm},
  year = {2024}
}

@book{greenbaum1997iterative,
  author = {Greenbaum, Anne},
  doi = {10.1137/1.9781611970937},
  isbn = {978-0-89871-396-1},
  language = {en},
  publisher = {SIAM},
  series = {Frontiers in Applied Mathematics},
  title = {Iterative Methods for Solving Linear Systems},
  year = {1997}
}

@book{saad2003iterative,
  author = {Saad, Yousef},
  doi = {10.1137/1.9780898718003},
  publisher = {SIAM},
  title = {Iterative methods for sparse linear systems},
  year = {2003}
}

@book{saad2011numerical,
  author = {Saad, Yousef},
  doi = {10.1137/1.9781611970739},
  isbn = {978-1-61197-072-2},
  month = jan,
  publisher = {SIAM},
  series = {Classics in {Applied} {Mathematics}},
  title = {Numerical Methods for Large Eigenvalue Problems},
  year = {2011}
}

@book{higham2008functions,
  author = {Higham, Nicholas J.},
  doi = {10.1137/1.9780898717778},
  isbn = {978-0-89871-646-7},
  publisher = {SIAM},
  title = {Functions of Matrices: {Theory} and Computation},
  year = {2008}
}

@article{nakatsukasa2024fast,
  author = {Nakatsukasa, Yuji and Tropp, Joel A},
  doi = {10.1137/23M1565413},
  journal = simax,
  number = {2},
  pages = {1183--1214},
  title = {Fast and accurate randomized algorithms for linear systems and eigenvalue problems},
  volume = {45},
  year = {2024}
}

@article{balabanov2025randomized,
  author = {Balabanov, Oleg and Grigori, Laura},
  doi = {10.1137/22M1518396},
  journal = sisc,
  number = {1},
  pages = {A553--A585},
  title = {Randomized block {Gram}--{Schmidt} process for the solution of linear systems and eigenvalue problems},
  volume = {47},
  year = {2025}
}

@article{grigori2026restoring,
  author = {Grigori, Laura and Kressner, Daniel and Shao, Nian and
Simunec, Igor},
  journal = {arXiv preprint \href{http://arxiv.org/abs/2601.10248v1}{arXiv:2601.10248v1}},
  title = {Restoring similarity in randomized {Krylov} methods with applications to eigenvalue problems and matrix functions},
  year = {2026}
}

@article{druskin1989two,
  author = {Druskin, Vladimir L and Knizhnerman, Leonid A},
  doi = {10.1016/S0041-5553(89)80020-5},
  journal = ussrcmmp,
  number = {6},
  pages = {112--121},
  title = {Two polynomial methods of calculating functions of symmetric matrices},
  volume = {29},
  year = {1989}
}

@article{saad1992analysis,
  author = {Saad, Yousef},
  doi = {10.1137/0729014},
  journal = sinum,
  number = {1},
  pages = {209--228},
  title = {Analysis of some {Krylov} subspace approximations to the matrix exponential operator},
  volume = {29},
  year = {1992}
}

@article{cortinovis2024speeding,
  author = {Cortinovis, Alice and Kressner, Daniel and
Nakatsukasa, Yuji},
  doi = {10.1137/22M1543458},
  journal = simax,
  number = {1},
  pages = {619--633},
  title = {Speeding up {Krylov} subspace methods for computing {$f(A)b$} via randomization},
  volume = {45},
  year = {2024}
}

@article{guttel2023randomized,
  author = {G{\"u}ttel, Stefan and Schweitzer, Marcel},
  doi = {10.1137/22M1518062},
  journal = simax,
  number = {3},
  pages = {1073--1095},
  title = {Randomized sketching for {Krylov} approximations of large-scale matrix functions},
  volume = {44},
  year = {2023}
}

@article{avron2010blendenpik,
  author = {Avron, Haim and Maymounkov, Petar and Toledo, Sivan},
  doi = {10.1137/090767911},
  journal = sisc,
  language = {en},
  month = jan,
  number = {3},
  pages = {1217--1236},
  title = {{Blendenpik:} {Supercharging} {LAPACK's} Least-Squares Solver},
  volume = {32},
  year = {2010}
}

@article{paige1982algorithm,
  author = {Paige, Christopher C. and Saunders, Michael A.},
  doi = {10.1145/355993.356000},
  journal = acmtoms,
  language = {en},
  month = jun,
  number = {2},
  pages = {195--209},
  title = {Algorithm 583: {LSQR:} {Sparse} Linear Equations and Least Squares Problems},
  volume = {8},
  year = {1982}
}

@article{meier2024are,
  author = {Meier, Maike and Nakatsukasa, Yuji and Townsend, Alex and Webb, Marcus},
  doi = {10.1137/23M1551973},
  journal = simax,
  month = jun,
  number = {2},
  pages = {905--929},
  title = {Are Sketch-and-Precondition Least Squares Solvers Numerically Stable?},
  volume = {45},
  year = {2024}
}

@article{epperly2024fast,
  author = {Epperly, Ethan N.},
  doi = {10.1137/23M1616790},
  journal = simax,
  month = dec,
  number = {4},
  pages = {1782--1804},
  title = {Fast and Forward Stable Randomized Algorithms for Linear Least-Squares Problems},
  volume = {45},
  year = {2024}
}

@article{diaz2026robust,
  archiveprefix = {arXiv},
  author = {D{\'i}az, Mateo and Epperly, Ethan N. and Frangella, Zachary and Tropp, Joel A. and Webber, Robert J.},
  doi = {10.1007/s10444-026-10360-1},
  eprint = {2304.12465},
  journal = acom,
  pages = {90},
  title = {Robust, Randomized Preconditioning for Kernel Ridge Regression},
  volume = {52},
  year = {2026}
}

@article{kaminetz2026everything,
  author = {Kaminetz, Eagan and Webber, Robert J.},
  journal = {arXiv preprint \href{http://arxiv.org/abs/2603.05709v1}{arXiv:2603.05709v1}},
  month = mar,
  title = {Everything is {Vecchia:} {Unifying} low-rank and sparse inverse {Cholesky} approximations},
  year = {2026}
}

@article{frangella2023randomized,
  author = {Frangella, Zachary and Tropp, Joel A. and
Udell, Madeleine},
  doi = {10.1137/21M1466244},
  journal = simax,
  month = jun,
  number = {2},
  pages = {718--752},
  title = {Randomized {Nyström} Preconditioning},
  volume = {44},
  year = {2023}
}

@article{chen2026preconditioning,
  author = {Chen, Tyler and Huber, Caroline and Lin, Ethan and
Zaid, Hajar},
  doi = {10.1553/etna_vol65s63},
  journal = etna,
  language = {en},
  pages = {63--92},
  title = {Preconditioning without a preconditioner using randomized block {Krylov} subspace methods},
  volume = {65},
  year = {2026}
}

@article{elfving1980block,
  author = {Elfving, Tommy},
  doi = {10.1007/BF01396365},
  journal = nm,
  language = {en},
  month = mar,
  number = {1},
  pages = {1--12},
  title = {Block-iterative methods for consistent and inconsistent linear equations},
  volume = {35},
  year = {1980}
}

@article{needell2014paved,
  author = {Needell, Deanna and Tropp, Joel A.},
  doi = {10.1016/j.laa.2012.12.022},
  journal = laa,
  language = {en},
  month = jan,
  pages = {199--221},
  series = {Special {Issue} on {Sparse} {Approximate} {Solution} of {Linear} {Systems}},
  title = {Paved with good intentions: {Analysis} of a randomized block {Kaczmarz} method},
  volume = {441},
  year = {2014}
}

@inproceedings{cutajar2016preconditioning,
  author = {Cutajar, Kurt and Osborne, Michael A and
Cunningham, John P and Filippone, Maurizio},
  booktitle = icml,
  language = {en},
  pages = {2529--2538},
  title = {Preconditioning Kernel Matrices},
  url = {https://proceedings.mlr.press/v48/cutajar16.html},
  volume = {48},
  year = {2016}
}

@inproceedings{gardner2018gpytorch,
  author = {Gardner, Jacob R. and Pleiss, Geoff and Bindel, David and Weinberger, Kilian Q. and Wilson, Andrew Gordon},
  booktitle = neurips,
  title = {{GPyTorch:} {Blackbox} Matrix-Matrix {Gaussian} Process Inference with {GPU} Acceleration},
  url = {https://dl.acm.org/doi/10.5555/3327757.3327857},
  year = {2018}
}

@incollection{dereziski2025faster,
  author = {Dereziński, Michał and Musco, Christopher and
Yang, Jiaming},
  booktitle = soda,
  doi = {10.1137/1.9781611978322.62},
  month = jan,
  pages = {1972--2004},
  title = {Faster linear systems and matrix norm approximation via multi-level sketched preconditioning},
  year = {2025}
}

@techreport{girard1987un,
  author = {Girard, Didier},
  institution = {TIM3--IMAG},
  month = {may},
  number = {RR 665-M},
  title = {Un algorithme rapide pour le calcul de la trace de l'inverse d'une grande matrice},
  type = {Rapport de recherche},
  url = {https://membres-ljk.imag.fr/Didier.Girard/TR-665-M-IMAG.pdf},
  year = {1987}
}

@article{girard1989fast,
  author = {Girard, A.},
  doi = {10.1007/BF01395775},
  journal = nm,
  language = {en},
  month = jan,
  number = {1},
  pages = {1--23},
  title = {A fast {``Monte}-{Carlo} cross-validation'' procedure for large least squares problems with noisy data},
  volume = {56},
  year = {1989}
}

@article{hutchinson1989stochastic,
  author = {Hutchinson, M. F.},
  doi = {10.1080/03610918908812806},
  journal = cssc,
  language = {en},
  month = may,
  number = {3},
  pages = {1059--1076},
  title = {A stochastic estimator of the trace of the influence matrix for {Laplacian} smoothing splines},
  volume = {18},
  year = {1989}
}

@book{vershynin2026high,
  author = {Vershynin, Roman},
  doi = {10.1017/9781009490672},
  edition = {second},
  language = {en},
  publisher = {Cambridge University Press},
  title = {High-dimensional probability: {An} introduction with applications in data science},
  year = {2026}
}

@article{avron2011randomized,
  author = {Avron, Haim and Toledo, Sivan},
  doi = {10.1145/1944345.1944349},
  journal = jacm,
  language = {en},
  month = apr,
  number = {2},
  pages = {1--34},
  title = {Randomized algorithms for estimating the trace of an implicit symmetric positive semi-definite matrix},
  volume = {58},
  year = {2011}
}

@article{cortinovis2022randomized,
  author = {Cortinovis, Alice and Kressner, Daniel},
  doi = {10.1007/s10208-021-09525-9},
  journal = focm,
  language = {en},
  month = jun,
  number = {3},
  pages = {875--903},
  title = {On Randomized Trace Estimates for Indefinite Matrices with an Application to Determinants},
  volume = {22},
  year = {2022}
}

@article{lin2017randomized,
  author = {Lin, Lin},
  doi = {10.1007/s00211-016-0837-7},
  journal = nm,
  month = may,
  number = {1},
  pages = {183--213},
  title = {Randomized Estimation of Spectral Densities of Large Matrices Made Accurate},
  volume = {136},
  year = {2017}
}

@article{gambhir2017deflation,
  author = {Gambhir, Arjun Singh and Stathopoulos, Andreas and
Orginos, Kostas},
  doi = {10.1137/16M1066361},
  journal = sisc,
  month = jan,
  number = {2},
  pages = {A532--A558},
  title = {Deflation as a Method of Variance Reduction for Estimating the Trace of a Matrix Inverse},
  volume = {39},
  year = {2017}
}

@inproceedings{meyer2021hutch,
  author = {Meyer, Raphael A. and Musco, Cameron and
Musco, Christopher and Woodruff, David P.},
  booktitle = sosa,
  doi = {10.1137/1.9781611976496.16},
  isbn = {978-1-61197-649-6},
  language = {en},
  month = jan,
  pages = {142--155},
  title = {{Hutch++:} {Optimal} Stochastic Trace Estimation},
  year = {2021}
}

@article{persson2022improved,
  author = {Persson, David and Cortinovis, Alice and
Kressner, Daniel},
  doi = {10.1137/21M1447623},
  journal = simax,
  month = sep,
  number = {3},
  pages = {1162--1185},
  title = {Improved variants of the {Hutch++} algorithm for trace estimation},
  volume = {43},
  year = {2022}
}

@article{epperly2024xtrace,
  author = {Epperly, Ethan N. and Tropp, Joel A. and
Webber, Robert J.},
  doi = {10.1137/23M1548323},
  journal = simax,
  month = mar,
  number = {1},
  pages = {1--23},
  title = {{XTrace:} {Making} the Most of Every Sample in Stochastic Trace Estimation},
  volume = {45},
  year = {2024}
}

@article{bekas2007estimator,
  author = {Bekas, C. and Kokiopoulou, E. and Saad, Y.},
  doi = {10.1016/j.apnum.2007.01.003},
  journal = apnum,
  language = {en},
  month = nov,
  number = {11},
  pages = {1214--1229},
  series = {Numerical {Algorithms}, {Parallelism} and
{Applications} (2)},
  title = {An estimator for the diagonal of a matrix},
  volume = {57},
  year = {2007}
}

@article{baston2022stochastic,
  author = {Baston, Robert A. and Nakatsukasa, Yuji},
  journal = {arXiv preprint \href{http://arxiv.org/abs/2201.10684v1}{arXiv:2201.10684v1}},
  month = jan,
  title = {Stochastic diagonal estimation: probabilistic bounds and an improved algorithm},
  year = {2022}
}

@inproceedings{dharangutte2023tight,
  author = {Dharangutte, Prathamesh and Musco, Christopher},
  booktitle = sosa,
  doi = {10.1137/1.9781611977585.ch32},
  month = jan,
  pages = {353--364},
  title = {A Tight Analysis of {Hutchinson's} Diagonal Estimator},
  year = {2023}
}

@book{golub2010matrices,
  author = {Golub, Gene H. and Meurant, Gérard},
  isbn = {978-0-691-14341-5},
  language = {en},
  publisher = {Princeton University Press},
  series = {Princeton Series in Applied Mathematics},
  title = {Matrices, moments, and quadrature with applications},
  year = {2010}
}

@article{chen2024lanczos,
  author = {Chen, Tyler},
  journal = {arXiv preprint \href{http://arxiv.org/abs/2410.11090v2}{arXiv:2410.11090v2}},
  month = oct,
  title = {The {Lanczos} algorithm for matrix functions: {A} handbook for scientists},
  year = {2024}
}

@inproceedings{musco2025sharper,
  author = {Musco, Cameron and Musco, Christopher and
Rosenblatt, Lucas and Singh, Apoorv Vikram},
  booktitle = colt,
  language = {en},
  month = jul,
  pages = {4309--4358},
  title = {Sharper Bounds for {Chebyshev} Moment Matching, with Applications},
  url = {https://proceedings.mlr.press/v291/musco25a.html},
  year = {2025}
}

@article{businger1965linear,
  author = {Businger, Peter and Golub, Gene H.},
  doi = {10.1007/BF01436084},
  journal = nm,
  language = {en},
  month = jun,
  number = {3},
  pages = {269--276},
  title = {Linear least squares solutions by {Householder} transformations},
  volume = {7},
  year = {1965}
}

@article{gu1996efficient,
  author = {Gu, Ming and Eisenstat, Stanley C.},
  doi = {10.1137/0917055},
  journal = sisc,
  month = jul,
  number = {4},
  pages = {848--869},
  title = {Efficient Algorithms for Computing a Strong Rank-Revealing {QR} Factorization},
  volume = {17},
  year = {1996}
}

@article{damle2026estimating,
  author = {Damle, Anil and Glas, Silke and Townsend, Alex and
Yu, Annan},
  doi = {10.1016/j.laa.2025.11.013},
  journal = laa,
  pages = {306--342},
  title = {Estimating a matrix's singular values with interpolative decompositions},
  volume = {731},
  year = {2026}
}

@article{grigori2025randomized,
  author = {Grigori, Laura and Xue, Zhipeng},
  journal = {arXiv preprint \href{http://arxiv.org/abs/2503.18496v1}{arXiv:2503.18496v1}},
  month = mar,
  title = {Randomized strong rank-revealing {QR} for column subset selection and low-rank matrix approximation},
  year = {2025}
}

@article{deshpande2006matrix,
  author = {Deshpande, Amit and Rademacher, Luis and
Vempala, Santosh S and Wang, Grant},
  doi = {10.4086/toc.2006.v002a012},
  journal = toc,
  number = {1},
  pages = {225--247},
  title = {Matrix approximation and projective clustering via volume sampling},
  volume = {2},
  year = {2006}
}

@inproceedings{deshpande2006adaptive,
  author = {Deshpande, Amit and Vempala, Santosh},
  booktitle = approxrandom,
  doi = {10.1007/11830924_28},
  isbn = {978-3-540-38045-0},
  language = {en},
  pages = {292--303},
  series = {Lecture {Notes} in {Computer} {Science}},
  title = {Adaptive Sampling and Fast Low-Rank Matrix Approximation},
  year = {2006}
}

@article{chen2025randomly,
  author = {Chen, Yifan and Epperly, Ethan N and Tropp, Joel A and
Webber, Robert J},
  doi = {10.1002/cpa.22234},
  journal = cpam,
  number = {5},
  pages = {995--1041},
  title = {Randomly pivoted {Cholesky:} {Practical} approximation of a kernel matrix with few entry evaluations},
  volume = {78},
  year = {2025}
}

@article{dong2025robust,
  author = {Dong, Yijun and Chen, Chao and Martinsson, Per-Gunnar and
Pearce, Katherine},
  doi = {10.1137/24M1678027},
  journal = simax,
  month = sep,
  number = {3},
  pages = {1791--1815},
  title = {Robust Blockwise Random Pivoting: {Fast} and Accurate Adaptive Interpolative Decomposition},
  volume = {46},
  year = {2025}
}

@misc{meyer2026active,
  author = {Meyer, Raphael A.},
  howpublished = {\url{https://randnla.github.io/leverage-score-regression/}},
  note = {RandNLA Proof Wiki},
  title = {Active {$L_2$} Regression via Leverage Score Sampling},
  url = {https://randnla.github.io/leverage-score-regression/},
  year = {2026}
}

@article{bharadwaj2023fast,
  author = {Bharadwaj, Vivek and Malik, Osman Asif and
Murray, Riley and Grigori, Laura and Bulu{\c c}, Ayd{\i}n and
Demmel, James},
  doi = {10.52202/075280-2076},
  journal = neurips,
  pages = {47874--47901},
  title = {Fast exact leverage score sampling from {Khatri}-{Rao} products with applications to tensor decomposition},
  year = {2023}
}

@inproceedings{bharadwaj2024efficient,
  author = {Bharadwaj, Vivek and Rakhshan, Beheshteh T. and
Malik, Osman Asif and Rabusseau, Guillaume},
  booktitle = neurips,
  doi = {10.52202/079017-2345},
  pages = {73726--73744},
  title = {Efficient Leverage Score Sampling for Tensor Train Decomposition},
  year = {2024}
}

@inproceedings{musco2025provably,
  author = {Musco, Christopher and Witter, R. Teal},
  booktitle = iclr,
  month = {May},
  pages = {64532--64559},
  title = {Provably Accurate {Shapley} Value Estimation via Leverage Score Sampling},
  url = {https://proceedings.iclr.cc/paper_files/paper/2025/hash/
a224ff18cc99a71751aa2b79118604da-Abstract-Conference.html},
  year = {2025}
}

@incollection{goreinov2001maximal,
  author = {Goreinov, S. A. and Tyrtyshnikov, E. E.},
  booktitle = structmat,
  doi = {10.1090/conm/280/4620},
  isbn = {0-8218-1921-6},
  language = {English},
  pages = {47--51},
  publisher = {Providence, RI: American Mathematical Society (AMS)},
  title = {The maximal-volume concept in approximation by low-rank matrices},
  year = {2001}
}

@article{ccivril2009selecting,
  author = {{\c C}ivril, Ali and Magdon-Ismail, Malik},
  doi = {10.1016/j.tcs.2009.06.018},
  journal = tcs,
  language = {en},
  month = nov,
  number = {47-49},
  pages = {4801--4811},
  title = {On selecting a maximum volume sub-matrix of a matrix and related problems},
  volume = {410},
  year = {2009}
}

@article{derezinski2021determinantal,
  author = {Derezi{\'n}ski, Micha{\l} and Mahoney, Michael W.},
  doi = {10.1090/noti2202},
  journal = noticesams,
  month = {January},
  number = {1},
  pages = {34--45},
  title = {Determinantal Point Processes in Randomized Numerical Linear Algebra},
  volume = {68},
  year = {2021}
}

@inproceedings{kulesza2011k,
  author = {Kulesza, Alex and Taskar, Ben},
  booktitle = icml,
  language = {en},
  month = jan,
  pages = {1193--1200},
  title = {{k}-{DPPs:} {Fixed-Size} Determinantal Point Processes},
  url = {https://dl.acm.org/doi/10.5555/3104482.3104632},
  year = {2011}
}

@article{kulesza2012determinantal,
  author = {Kulesza, Alex and Taskar, Ben},
  doi = {10.1561/2200000044},
  journal = ftml,
  language = {English},
  month = dec,
  number = {2–3},
  pages = {123--286},
  title = {Determinantal Point Processes for Machine Learning},
  volume = {5},
  year = {2012}
}

@article{sood2025statistical,
  author = {Sood, Anav and Hastie, Trevor},
  doi = {10.1093/jrsssb/qkaf023},
  journal = jrssb,
  language = {English},
  number = {5},
  pages = {1382--1403},
  title = {A statistical view of column subset selection},
  volume = {87},
  year = {2025}
}

@article{martinsson2017householder,
  author = {Martinsson, Per-Gunnar and Quintana Orti, Gregorio and
Heavner, Nathan and Van De Geijn, Robert},
  doi = {10.1137/16M1081270},
  journal = sisc,
  number = {2},
  pages = {C96--C115},
  title = {{Householder} {QR} factorization with randomization for column pivoting {(HQRRP)}},
  volume = {39},
  year = {2017}
}

@article{osinsky2025close,
  author = {Osinsky, A.~I.},
  doi = {10.1016/j.laa.2025.07.016},
  journal = laa,
  pages = {359--377},
  title = {Close to optimal column approximation using a single {SVD}},
  volume = {725},
  year = {2025}
}

@article{cortinovis2026adaptive,
  author = {Cortinovis, Alice and Kressner, Daniel},
  doi = {10.1137/24M1719189},
  journal = simax,
  number = {1},
  pages = {25--47},
  title = {Adaptive randomized pivoting for column subset selection, {DEIM}, and low-rank approximation},
  volume = {47},
  year = {2026}
}

@article{fornace2026column,
  author = {Fornace, Mark and Lindsey, Michael},
  journal = {Found. Comput. Math.},
  note = {To appear},
  title = {Column and row subset selection using nuclear scores: algorithms and theory for {Nystr\"om} approximation, {CUR} decomposition, and graph {Laplacian} reduction},
  year = {2026}
}

@article{park2025accuracy,
  author = {Park, Taejun and Nakatsukasa, Yuji},
  doi = {10.1137/24M1660346},
  journal = simax,
  number = {1},
  pages = {780--810},
  title = {Accuracy and stability of {CUR} decompositions with oversampling},
  volume = {46},
  year = {2025}
}

@article{cortinovis2020low,
  author = {Cortinovis, Alice and Kressner, Daniel},
  doi = {10.1137/19M1281848},
  journal = simax,
  number = {4},
  pages = {1651--1673},
  title = {Low-rank approximation in the {Frobenius} norm by column and row subset selection},
  volume = {41},
  year = {2020}
}

@article{pritchard2025iterativecur,
  author = {Pritchard, Nathaniel and Park, Taejun and Nakatsukasa, Yuji and Martinsson, Per-Gunnar},
  journal = {arXiv preprint \href{http://arxiv.org/abs/2509.21963v2}{arXiv:2509.21963v2}},
  title = {{IterativeCUR:} {Large} Rank-Adaptive Approximation From a Small Recycled Sketch},
  year = {2025}
}

@book{martinsson2019fast,
  author = {Martinsson, Per-Gunnar},
  doi = {10.1137/1.9781611976045},
  publisher = {SIAM},
  title = {Fast direct solvers for elliptic {PDEs}},
  year = {2019}
}

@article{park2025low,
  author = {Park, Taejun and Nakatsukasa, Yuji},
  doi = {10.1137/24M1683998},
  journal = sisc,
  number = {3},
  pages = {A1858--A1887},
  title = {Low-rank approximation of parameter-dependent matrices via {CUR} decomposition},
  volume = {47},
  year = {2025}
}

@article{szyld2006many,
  author = {Szyld, Daniel B},
  doi = {10.1007/s11075-006-9046-2},
  journal = numalgs,
  number = {3-4},
  pages = {309--323},
  title = {The many proofs of an identity on the norm of oblique projections},
  volume = {42},
  year = {2006}
}

@book{unknown2006schur,
  doi = {10.1007/b105056},
  editor = {Zhang, Fuzhen},
  publisher = {Springer Science \& Business Media},
  title = {The {Schur} complement and its applications},
  year = {2005}
}

@book{golub2013matrix,
  author = {Golub, Gene H. and Van Loan, Charles F.},
  edition = {Fourth},
  publisher = {Johns Hopkins University Press},
  series = {Johns Hopkins Studies in the Mathematical Sciences},
  title = {Matrix Computations},
  year = {2013}
}

@article{townsend2013gaussian,
  author = {Alex Townsend and Lloyd N. Trefethen},
  journal = siamnews,
  month = mar,
  number = {02},
  title = {{Gaussian} Elimination as an Iterative Algorithm},
  url = {https://www.siam.org/publications/siam-news/articles/
gaussian-elimination-as-an-iterative-algorithm/},
  volume = {46},
  year = {2013}
}

@book{trefethen2019approximation,
  author = {Trefethen, Lloyd N},
  doi = {10.1137/1.9781611975949},
  publisher = {SIAM},
  title = {Approximation Theory and Approximation Practice, Extended Edition},
  year = {2019}
}

@article{townsend2013extension,
  author = {Townsend, A. and Trefethen, L. N.},
  doi = {10.1137/130908002},
  journal = sisc,
  number = {6},
  pages = {C495--C518},
  title = {An extension of {Chebfun} to two dimensions},
  volume = {35},
  year = {2013}
}

@article{jones2025subapsnap,
  author = {Jones, Eleanor and Nakatsukasa, Yuji},
  journal = {arXiv preprint \href{http://arxiv.org/abs/2510.04825v1}{arXiv:2510.04825v1}},
  title = {{SubApSnap:} {Solving} parameter-dependent linear systems with a snapshot and subsampling},
  year = {2025}
}

@article{gittens2016revisiting,
  author = {Gittens, Alex and Mahoney, Michael W},
  journal = jmlr,
  number = {1},
  pages = {3977--4041},
  title = {Revisiting the {Nystr\"om} method for improved large-scale machine learning},
  url = {https://dl.acm.org/doi/abs/10.5555/2946645.3007070},
  volume = {17},
  year = {2016}
}

@book{higham2002accuracy,
  author = {N. J. Higham},
  doi = {10.1137/1.9780898718027},
  edition = {Second},
  isbn = {0-89871-521-0},
  pages = {xxx+680},
  publisher = {SIAM},
  title = {Accuracy and Stability of Numerical Algorithms},
  year = {2002}
}

@article{epperly2025embrace,
  author = {Epperly, Ethan N and Tropp, Joel A and
Webber, Robert J},
  doi = {10.1137/24M1699048},
  journal = simax,
  number = {4},
  pages = {2527--2557},
  title = {Embrace Rejection: {Kernel} Matrix Approximation by Accelerated Randomly Pivoted {Cholesky}},
  volume = {46},
  year = {2025}
}

@inproceedings{musco2017recursive,
  author = {Musco, Cameron and Musco, Christopher},
  booktitle = neurips,
  pages = {3833--3845},
  title = {Recursive sampling for the {Nystrom} method},
  url = {https://dl.acm.org/doi/10.5555/3294996.3295140},
  volume = {30},
  year = {2017}
}

@book{rasmussen2006gaussian,
  author = {Rasmussen, Carl Edward and Williams, Christopher K.~I.},
  doi = {10.7551/mitpress/3206.001.0001},
  publisher = {MIT Press},
  title = {{Gaussian} Processes for Machine Learning},
  year = {2006}
}

@article{fowlkes2004spectral,
  author = {Fowlkes, Charless and Belongie, Serge and Chung, Fan and
Malik, Jitendra},
  doi = {10.1109/TPAMI.2004.1262185},
  journal = ieeetpami,
  number = {2},
  pages = {214--225},
  title = {Spectral grouping using the {Nystrom} method},
  volume = {26},
  year = {2004}
}

@article{zhang2010clustered,
  author = {Zhang, Kai and Kwok, James T},
  doi = {10.1109/TNN.2010.2064786},
  journal = ieeetnn,
  number = {10},
  pages = {1576--1587},
  title = {Clustered {Nystr\"om} method for large scale manifold learning and dimension reduction},
  volume = {21},
  year = {2010}
}

@article{nakatsukasa2026finding,
  author = {Nakatsukasa, Yuji and Tang, Zheng},
  journal = {arXiv preprint \href{http://arxiv.org/abs/2605.05043v1}{arXiv:2605.05043v1}},
  title = {Finding accurate eigenvalues and eigenvectors of positive semi-definite matrices given a subspace},
  year = {2026}
}

@article{wegner1981bounds,
  author = {Wegner, Franz},
  doi = {10.1007/BF01292646},
  journal = zpb,
  number = {1-2},
  pages = {9--15},
  title = {Bounds on the density of states in disordered systems},
  volume = {44},
  year = {1981}
}

@article{combes2009generalized,
  author = {Combes, Jean-Michel and Germinet, Fran\c{c}ois and
Klein, Abel},
  doi = {10.1007/s10955-009-9731-3},
  journal = jsp,
  number = {2},
  pages = {201--216},
  title = {Generalized eigenvalue-counting estimates for the Anderson model},
  volume = {135},
  year = {2009}
}

@article{minami1996local,
  author = {Minami, Nariyuki},
  journal = cmp,
  number = {3},
  pages = {709--725},
  title = {Local fluctuation of the spectrum of a multidimensional Anderson tight binding model},
  url = {http://projecteuclid.org/euclid.cmp/1104286441},
  volume = {177},
  year = {1996}
}

@article{kressner2026randomized,
  author = {Daniel Kressner and Nian Shao},
  doi = {10.1137/24M1674881},
  journal = sisc,
  title = {On a Randomized Small-Block {Lanczos} Method for Large-Scale Null Space Computations},
  year = {2026}
}

@book{parlett1998symmetric,
  author = {Parlett, Beresford N.},
  doi = {10.1137/1.9781611971163},
  isbn = {0-89871-402-8},
  note = {Corrected reprint of the 1980 original},
  pages = {xxiv+398},
  publisher = {Society for Industrial and Applied Mathematics
(SIAM), Philadelphia, PA},
  series = {Classics in Applied Mathematics},
  title = {The symmetric eigenvalue problem},
  volume = {20},
  year = {1998}
}

@inproceedings{meyer2024unreasonable,
  author = {Meyer, Raphael and Musco, Cameron and
Musco, Christopher},
  booktitle = soda,
  doi = {10.1137/1.9781611977912.32},
  isbn = {978-1-61197-791-2},
  pages = {811--845},
  title = {On the unreasonable effectiveness of single vector {Krylov} methods for low-rank approximation},
  year = {2024}
}

@techreport{stewart1990incremental,
  author = {Stewart, G. W.},
  institution = {Institute for Advanced Computer Studies, University of Maryland},
  number = {TR-90-87},
  title = {Incremental Condition Calculation and Column Selection},
  year = {1990}
}

@inproceedings{chen2026does,
  author = {Tyler Chen and Ethan N. Epperly and Raphael A. Meyer and
Christopher Musco and Akash Rao},
  booktitle = soda,
  doi = {10.1137/1.9781611978971.42},
  pages = {1026--1046},
  title = {Does block size matter in randomized block {Krylov} low-rank approximation?},
  year = {2026}
}

@techreport{shao2025structural,
  author = {Nian Shao},
  title = {A structural bound for cluster robustness of randomized small-block {Lanczos}},
  type = {{arXiv:2507.10144}},
  year = {2025}
}

@article{vandooren1981generalized,
  author = {Van Dooren, Paul M.},
  doi = {10.1109/TAC.1981.1102559},
  journal = ieeetac,
  number = {1},
  pages = {111--129},
  title = {The generalized eigenstructure problem in linear system theory},
  volume = {26},
  year = {1981}
}

@book{gantmacher1959theory,
  author = {Gantmacher, F. R.},
  publisher = {Chelsea Publishing Co.},
  title = {The theory of matrices. {Vols.} 1, 2},
  year = {1959}
}

@article{wilkinson1979kronecker,
  author = {Wilkinson, J. H.},
  doi = {10.1016/0024-3795(79)90140-X},
  journal = laa,
  pages = {285--303},
  title = {{Kronecker's} canonical form and the {$QZ$} algorithm},
  volume = {28},
  year = {1979}
}

@article{lotz2020wilkinson,
  author = {Lotz, Martin and Noferini, Vanni},
  doi = {10.1007/s10208-020-09455-y},
  journal = focm,
  number = {6},
  pages = {1439--1473},
  title = {{Wilkinson's} bus: weak condition numbers, with an application to singular polynomial eigenproblems},
  volume = {20},
  year = {2020}
}

@article{kolda2009tensor,
  author = {Tamara G. Kolda and Brett W. Bader},
  doi = {10.1137/07070111X},
  journal = sirev,
  number = {3},
  pages = {455--500},
  title = {Tensor Decompositions and Applications},
  volume = {51},
  year = {2009}
}

@article{bujanovic2021norm,
  author = {Zvonimir Bujanovi{\'c} and Daniel Kressner},
  doi = {10.1137/20M1331718},
  journal = simax,
  number = {1},
  pages = {202--223},
  title = {Norm and Trace Estimation with Random Rank-one Vectors},
  volume = {42},
  year = {2021}
}

@article{guttel2024randomized,
  author = {Stefan G{\"u}ttel and Daniel Kressner and
Bart Vandereycken},
  doi = {10.1137/22M153656X},
  journal = sisc,
  number = {5},
  pages = {A3022--A3043},
  title = {Randomized Sketching of Nonlinear Eigenvalue Problems},
  volume = {46},
  year = {2024}
}

@article{hu2026small,
  author = {Xuehan Hu and Grigoris Paouris},
  doi = {10.1016/j.jfa.2025.111309},
  journal = jfunct,
  number = {6},
  pages = {111309},
  title = {Small Ball Probabilities for Simple Random Tensors},
  volume = {290},
  year = {2026}
}

@article{vershynin2020concentration,
  author = {Roman Vershynin},
  doi = {10.3150/20-BEJ1218},
  journal = bernoulli,
  number = {4},
  pages = {3139--3162},
  title = {Concentration inequalities for random tensors},
  volume = {26},
  year = {2020}
}

@article{bamberger2022hanson,
  author = {Stefan Bamberger and Felix Krahmer and Rachel Ward},
  doi = {10.1007/s43670-022-00028-4},
  journal = stspda,
  pages = {14},
  title = {The {Hanson}--{Wright} Inequality for Random Tensors},
  volume = {20},
  year = {2022}
}

@article{meyer2026hutchinson,
  author = {Raphael A. Meyer and Haim Avron},
  doi = {10.1137/24M1720895},
  journal = simax,
  number = {1},
  pages = {353--387},
  title = {{Hutchinson's} Estimator is Bad at {Kronecker}-Trace-Estimation},
  volume = {47},
  year = {2026}
}

@inproceedings{meyer2025understanding,
  author = {Raphael A. Meyer and William Joseph Swartworth and
David P. Woodruff},
  booktitle = icml,
  pages = {43909--43933},
  title = {Understanding the {Kronecker} Matrix-Vector Complexity of Linear Algebra},
  url = {https://dl.acm.org/doi/abs/10.5555/3780338.3782102},
  volume = {267},
  year = {2025}
}

@article{oliveira2016lower,
  author = {Roberto Imbuzeiro Oliveira},
  doi = {10.1007/s00440-016-0738-9},
  journal = ptrf,
  number = {3--4},
  pages = {1175--1194},
  title = {The Lower Tail of Random Quadratic Forms with Applications to Ordinary Least Squares},
  volume = {166},
  year = {2016}
}

@article{haselby2026fast,
  author = {Cullen Haselby and Mark A. Iwen and Deanna Needell and
Elizaveta Rebrova and William Swartworth},
  doi = {10.1007/s11075-025-02101-0},
  journal = numalgs,
  number = {4},
  pages = {2567--2629},
  title = {Fast and Low-Memory Compressive Sensing Algorithms for Low {Tucker}-Rank Tensor Approximation from Streamed Measurements},
  volume = {101},
  year = {2026}
}

@article{bujanovic2026subspace,
  author = {Zvonimir Bujanovi{\'c} and Luka Grubi{\v{s}}i{\'c} and
Daniel Kressner and Hei Yin Lam},
  doi = {10.1093/imanum/draf043},
  journal = imajna,
  number = {3},
  pages = {1328--1355},
  title = {Subspace Embedding with Random {Khatri}--{Rao} Products and Its Application to Eigensolvers},
  volume = {46},
  year = {2026}
}

@incollection{martinsson2018randomized,
  author = {Martinsson, Per-Gunnar},
  booktitle = mathofdata,
  chapter = {4},
  doi = {10.1090/pcms/025},
  pages = {187--229},
  publisher = {American Mathematical Society},
  series = {IAS/Park City Mathematics Series},
  title = {Randomized methods for matrix computations},
  volume = {25},
  year = {2018}
}

@incollection{drineas2018lectures,
  archiveprefix = {arXiv},
  author = {Drineas, Petros and Mahoney, Michael W.},
  booktitle = mathofdata,
  doi = {10.1090/pcms/025},
  publisher = {AMS},
  series = {IAS/Park City Mathematics Series},
  title = {Lectures on Randomized Numerical Linear Algebra},
  volume = {25},
  year = {2018}
}

@inproceedings{ahle2020oblivious,
  author = {Ahle, Thomas D. and Kapralov, Michael and Knudsen, Jakob B. T. and Pagh, Rasmus and Velingker, Ameya and Woodruff, David P. and Zandieh, Amir},
  booktitle = soda,
  doi = {10.1137/1.9781611975994.9},
  pages = {141--160},
  title = {Oblivious Sketching of High-Degree Polynomial Kernels},
  year = {2020}
}

@inproceedings{rakhshan2020tensorized,
  author = {Beheshteh Rakhshan and Guillaume Rabusseau},
  booktitle = aistats,
  pages = {3306--3316},
  title = {Tensorized Random Projections},
  url = {https://proceedings.mlr.press/v108/rakhshan20a.html},
  volume = {108},
  year = {2020}
}

@article{sun2021tensor,
  author = {Yiming Sun and Yang Guo and Joel A. Tropp and
Madeleine Udell},
  title = {Tensor Random Projection for Low Memory Dimension Reduction},
  journal = {arXiv preprint \href{https://arxiv.org/abs/2105.00105}{arXiv:2105.00105}},
  month = apr,
  year = {2021}
}

@article{strohmer2009randomized,
  author = {Strohmer, Thomas and Vershynin, Roman},
  doi = {10.1007/s00041-008-9030-4},
  journal = fourappl,
  number = {2},
  pages = {262--278},
  title = {A randomized {Kaczmarz} algorithm with exponential convergence},
  volume = {15},
  year = {2009}
}

@article{needell2016stochastic,
  author = {Needell, Deanna and Srebro, Nathan and Ward, Rachel},
  doi = {10.1007/s10107-015-0864-7},
  journal = mathpro,
  pages = {549--573},
  title = {Stochastic gradient descent, weighted sampling, and the randomized {Kaczmarz} algorithm},
  volume = {155},
  year = {2016}
}

@article{derezinski2025randomized,
  author = {Derezi{\'n}ski, Micha{\l} and Needell, Deanna and
Rebrova, Elizaveta and Yang, Jiaming},
  doi = {10.1137/25M1727187},
  journal = simax,
  number = {4},
  pages = {2558--2588},
  title = {Randomized {Kaczmarz} methods with beyond-{Krylov} convergence},
  volume = {46},
  year = {2025}
}

@techreport{connolly2022randomized,
  author = {Connolly, Michael P and Higham, Nicholas J and
Pranesh, Srikara},
  number = {2022.10},
  title = {Randomized low rank matrix approximation: {Rounding} error analysis and a mixed precision algorithm},
  type = {MIMS EPrint},
  url = {https://eprints.maths.manchester.ac.uk/id/eprint/2863},
  year = {2022}
}

@article{nakatsukasa2020fast,
  author = {Nakatsukasa, Yuji},
  journal = {arXiv preprint \href{http://arxiv.org/abs/2009.11392v1}{arXiv:2009.11392v1}},
  title = {Fast and stable randomized low-rank matrix approximation},
  year = {2020}
}

@inproceedings{tropp2017fixed,
  author = {Tropp, Joel A. and Yurtsever, Alp and
Udell, Madeleine and Cevher, Volkan},
  booktitle = neurips,
  pages = {1225--1234},
  title = {Fixed-Rank Approximation of a Positive-Semidefinite Matrix from Streaming Data},
  url = {https://dl.acm.org/doi/10.5555/3294771.3294888},
  year = {2017}
}

@article{burke2025numerical,
  author  = {Burke, Liam and Carson, Erin and Ma, Yuxin},
  title   = {On the Numerical Stability of Sketched {GMRES}},
  journal = simax,
  volume  = {47},
  number  = {3},
  pages   = {1657--1678},
  year    = {2026},
  doi     = {10.1137/25M1744629}
}

@inproceedings{guruswami2012optimal,
  author = {Guruswami, Venkatesan and Sinop, Ali Kemal},
  booktitle = soda,
  doi = {10.1137/1.9781611973099.95},
  isbn = {978-1-61197-210-8},
  month = jan,
  pages = {1207--1214},
  title = {Optimal column-based low-rank matrix reconstruction},
  year = {2012}
}

@inproceedings{dereziski2017unbiased,
  author = {Dereziński, Michał and Warmuth, Manfred K.},
  booktitle = neurips,
  title = {Unbiased estimates for linear regression via volume sampling},
  url = {https://proceedings.neurips.cc/paper/2017/hash/54e36c5ff5f6a1802925ca009f3ebb68-Abstract.html},
  volume = {30},
  year = {2017}
}

@article{epperly2026adaptivea,
  author = {Epperly, Ethan N.},
  doi = {10.1137/25M1818503},
  journal = simax,
  number = {4},
  pages = {1785--1798},
  title = {Adaptive randomized pivoting and volume sampling},
  volume = {47},
  year = {2026}
}

@unpublished{epperly2026adaptiveb,
  author = {Epperly, Ethan N. W. and Miller, Kevin S. and Webber, Robert J.},
  note = {Manuscript in preparation},
  title = {Adaptive algorithms for data summarization},
  year = {2026}
}

@techreport{golub1976rank,
  author = {Golub, G. and Klema, V. and Stewart, G. W.},
  institution = {Dept. Comput. Sci., Stanford University},
  number = {Rep. STAN-CS-76-559},
  title = {Rank degeneracy and least squares problems},
  url = {http://infolab.stanford.edu/pub/cstr/reports/cs/tr/76/559/CS-TR-76-559.pdf},
  year = {1976}
}

@article{bischof1990incremental,
  author = {Bischof, Christian H.},
  doi = {10.1137/0611021},
  journal = simax,
  language = {en},
  month = apr,
  number = {2},
  pages = {312--322},
  title = {Incremental Condition Estimation},
  volume = {11},
  year = {1990}
}

@article{damle2026computing,
  author = {Damle, Anil},
  journal = {arXiv preprint \href{http://arxiv.org/abs/2607.13532v1}{arXiv:2607.13532v1}},
  month = jul,
  title = {Computing Strong Rank-Revealing Factorizations for Matrices with Orthonormal Rows},
  year = {2026}
}

@misc{epperly2023stochastic,
  author = {Epperly, Ethan N.},
  language = {en-US},
  month = jan,
  title = {Stochastic Trace Estimation},
  note = {Blog post, \url{https://www.ethanepperly.com/index.php/2023/01/26/stochastic-trace-estimation/}},
  year = {2023}
}

@book{wedderburn1934lectures,
  author = {Wedderburn, J.},
  doi = {10.1090/coll/017},
  isbn = {978-0-8218-4610-0 978-0-8218-3204-2},
  language = {en},
  month = dec,
  publisher = {American Mathematical Society},
  series = {Colloquium {Publications}},
  title = {Lectures on Matrices},
  volume = {17},
  year = {1934}
}

@article{chu1995rank,
  author = {Chu, Moody T. and Funderlic, Robert E. and Golub, Gene H.},
  doi = {10.1137/1037124},
  journal = sirev,
  language = {en},
  month = dec,
  number = {4},
  pages = {512--530},
  title = {A Rank–One Reduction Formula and Its Applications to Matrix Factorizations},
  volume = {37},
  year = {1995}
}

@article{li2017algorithm,
  author = {Li, Huamin and Linderman, George C. and Szlam, Arthur and Stanton, Kelly P. and Kluger, Yuval and Tygert, Mark},
  doi = {10.1145/3004053},
  journal = acmtoms,
  month = jan,
  number = {3},
  pages = {28:1--28:14},
  title = {Algorithm 971: {An} Implementation of a Randomized Algorithm for Principal Component Analysis},
  volume = {43},
  year = {2017}
}

@article{gu2015subspace,
  author = {Gu, Ming},
  doi = {10.1137/130938700},
  journal = sisc,
  month = jan,
  number = {3},
  pages = {A1139--A1173},
  title = {Subspace Iteration Randomization and Singular Value Problems},
  volume = {37},
  year = {2015}
}

@article{rokhlin2010randomized,
  author = {Rokhlin, Vladimir and Szlam, Arthur and Tygert, Mark},
  doi = {10.1137/080736417},
  journal = simax,
  language = {en},
  month = jan,
  number = {3},
  pages = {1100--1124},
  title = {A Randomized Algorithm for Principal Component Analysis},
  volume = {31},
  year = {2010}
}

@article{minster2024parallel,
  author = {Minster, Rachel and Li, Zitong and Ballard, Grey},
  doi = {10.1137/22M1540363},
  journal = sisc,
  language = {en},
  month = apr,
  number = {2},
  pages = {A1186--A1213},
  title = {Parallel Randomized {Tucker} Decomposition Algorithms},
  volume = {46},
  year = {2024}
}

@article{saibaba2025improved,
  author = {Arvind K. Saibaba and Bhisham Dev Verma and Grey Ballard},
  journal = {arXiv preprint \href{http://arxiv.org/abs/2507.23207v2}{arXiv:2507.23207v2}},
  title = {Improved Analysis of {Khatri}--{Rao} Random Projections and Applications},
  year = {2025}
}

@article{cazeaux2026linear,
  author = {Cazeaux, Paul and Dupuy, Mi-Song and
Figueroa Justiniano, Rodrigo},
  journal = {arXiv preprint \href{http://arxiv.org/abs/2603.11009v3}{arXiv:2603.11009v3}},
  title = {Linear-Scaling Tensor Train Sketching},
  year = {2026}
}

@article{camano2026linear,
  author = {Cama{\~n}o, Chris and Epperly, Ethan N. and
Meyer, Raphael A. and Tropp, Joel A.},
  journal = {arXiv preprint \href{http://arxiv.org/abs/2606.15350v1}{arXiv:2606.15350v1}},
  title = {Linear Algebra at Exponential Scale via Tensor Network Dimension Reduction},
  year = {2026}
}

@article{bujanovic2026stochastic,
  author = {Bujanovi{\'c}, Zvonimir and Kressner, Daniel and
Oli{\'c}, Hrvoje},
  journal = {arXiv preprint \href{http://arxiv.org/abs/2606.15679v1}{arXiv:2606.15679v1}},
  title = {Stochastic Trace Estimation with Tensor Train Random Vectors},
  year = {2026}
}

@article{schollwock2011density,
  author = {Schollw{\"o}ck, Ulrich},
  doi = {10.1016/j.aop.2010.09.012},
  journal = annphys,
  number = {1},
  pages = {96--192},
  title = {The density-matrix renormalization group in the age of matrix product states},
  volume = {326},
  year = {2011}
}

@article{oseledets2011tensor,
  author = {Oseledets, Ivan V.},
  doi = {10.1137/090752286},
  journal = sisc,
  number = {5},
  pages = {2295--2317},
  title = {Tensor-Train Decomposition},
  volume = {33},
  year = {2011}
}

@article{aldaas2025adaptive,
  author = {Al Daas, Hussam and Ballard, Grey and Grigori, Laura
and Martinez Aguilar, Mariana and Saibaba, Arvind K.
and Verma, Bhisham Dev},
  journal = {arXiv preprint \href{http://arxiv.org/abs/2511.03598v1}{arXiv:2511.03598v1}},
  title = {Adaptive Randomized Tensor Train Rounding using {Khatri}--{Rao} Products},
  year = {2025}
}

@incollection{huber2017randomized,
  author = {Huber, Benjamin and Schneider, Reinhold and Wolf, Sebastian},
  booktitle = compsens,
  doi = {10.1007/978-3-319-69802-1_9},
  pages = {261--290},
  publisher = {Birkh{\"a}user},
  series = {Applied and Numerical Harmonic Analysis},
  title = {A Randomized Tensor Train Singular Value Decomposition},
  year = {2017}
}

@inproceedings{ma2022cost,
  author = {Linjian Ma and Edgar Solomonik},
  booktitle = neurips,
  pages = {38980--38993},
  title = {Cost-efficient {Gaussian} Tensor Network Embeddings for Tensor-structured Inputs},
  url = {https://dl.acm.org/doi/10.5555/3600270.3603095},
  volume = {35},
  year = {2022}
}

@article{spielman2004smoothed,
  author = {Spielman, Daniel A. and Teng, Shang-Hua},
  doi = {10.1145/990308.990310},
  journal = jacm,
  number = {3},
  pages = {385--463},
  title = {Smoothed Analysis of Algorithms: {Why} the Simplex Algorithm Usually Takes Polynomial Time},
  volume = {51},
  year = {2004}
}

@article{sankar2006smoothed,
  author = {Sankar, Arvind and Spielman, Daniel A. and Teng, Shang-Hua},
  doi = {10.1137/S0895479803436202},
  journal = simax,
  number = {2},
  pages = {446--476},
  title = {Smoothed Analysis of the Condition Numbers and Growth Factors of Matrices},
  volume = {28},
  year = {2006}
}

@techreport{kireeva2024randomized,
  author = {Kireeva, Anastasia and Tropp, Joel A.},
  doi = {10.7907/7yade-5k351},
  institution = {California Institute of Technology},
  note = {Lectures at the CIME Summer School on Machine Learning,
Cetraro, Italy, July 2023},
  number = {2023-02},
  title = {Randomized Matrix Computations: {Themes} and Variations},
  type = {Caltech CMS Lecture Notes},
  year = {2024}
}

@inproceedings{shah2021smoothed,
  author = {Shah, Rikhav and Silwal, Sandeep},
  booktitle = approxrandom,
  doi = {10.4230/LIPIcs.APPROX/RANDOM.2021.40},
  pages = {40:1--40:21},
  title = {Smoothed Analysis of the Condition Number Under Low-Rank Perturbations},
  volume = {207},
  year = {2021}
}

@article{drineas2012fast,
  author = {Drineas, Petros and Magdon-Ismail, Malik and Mahoney, Michael W and Woodruff, David P},
  journal = jmlr,
  number = {1},
  pages = {3475--3506},
  title = {Fast approximation of matrix coherence and statistical leverage},
  url = {https://dl.acm.org/doi/10.5555/2503308.2503352},
  volume = {13},
  year = {2012}
}

@article{peherstorfer2020stability,
  author = {Peherstorfer, Benjamin and Drmač, Zlatko and Gugercin, Serkan},
  doi = {10.1137/19M1307391},
  journal = sisc,
  language = {en},
  month = jan,
  number = {5},
  pages = {A2837--A2864},
  title = {Stability of Discrete Empirical Interpolation and Gappy Proper Orthogonal Decomposition with Randomized and Deterministic Sampling Points},
  volume = {42},
  year = {2020}
}

@article{drineas2008relative,
  author = {Drineas, Petros and Mahoney, Michael W. and Muthukrishnan, S.},
  doi = {10.1137/07070471X},
  journal = simax,
  month = jan,
  number = {2},
  pages = {844--881},
  title = {Relative-Error {CUR} Matrix Decompositions},
  volume = {30},
  year = {2008}
}

@article{jackson1912approximation,
  author = {Jackson, Dunham},
  doi = {10.1090/S0002-9947-1912-1500930-2},
  journal = tams,
  pages = {491--515},
  title = {On approximation by trigonometric sums and polynomials},
  volume = {13},
  year = {1912}
}

@article{adcock2025optimal,
  author = {Adcock, Ben},
  doi = {10.1007/s10208-025-09738-2},
  journal = focm,
  number = {6},
  pages = {1975--2034},
  title = {Optimal sampling for least-squares approximation},
  volume = {25},
  year = {2025}
}

@article{sifuentes2015randomized,
  author = {Sifuentes, Josef and Gimbutas, Zydrunas and Greengard, Leslie},
  journal = etna,
  pages = {177--188},
  title = {Randomized Methods for Rank-Deficient Linear Systems},
  volume = {44},
  year = {2015}
}

@article{pan2010additive,
  author = {Pan, Victor Y. and Ivolgin, Dmitriy and Murphy, Brian and
Rosholt, Rhys Eric and Tang, Yuqing and Yan, Xiaodong},
  doi = {10.1016/j.laa.2009.10.020},
  journal = laa,
  number = {4},
  pages = {1070--1089},
  title = {Additive Preconditioning for Matrix Computations},
  volume = {432},
  year = {2010}
}

@article{oleary1980block,
  author = {O'Leary, Dianne P.},
  doi = {10.1016/0024-3795(80)90247-5},
  journal = laa,
  month = feb,
  pages = {293--322},
  series = {Special {Volume} {Dedicated} to {Alson} {S}. {Householder}},
  title = {The block conjugate gradient algorithm and related methods},
  volume = {29},
  year = {1980}
}

@book{dongarra1979linpack,
  author = {Dongarra, J. J. and Moler, C. B. and Bunch, J. R. and Stewart, G. W.},
  doi = {10.1137/1.9781611971811},
  isbn = {978-0-89871-172-1 978-1-61197-181-1},
  language = {en},
  month = jan,
  publisher = {SIAM},
  title = {{LINPACK} Users' Guide},
  year = {1979}
}

@incollection{higham1990analysis,
  author = {Higham, Nicholas J},
  booktitle = relnumcomp,
  doi = {10.1093/oso/9780198535645.003.0010},
  isbn = {978-0-19-853564-5},
  month = sep,
  pages = {161--185},
  publisher = {Oxford University Press},
  title = {Analysis of the {Cholesky} decomposition of a semi-definite matrix},
  year = {1990}
}

@article{fine2001efficient,
  author = {Fine, Shai and Scheinberg, Katya},
  journal = jmlr,
  month = dec,
  pages = {243--264},
  title = {Efficient {SVM} training using low-rank kernel representations},
  url = {https://dl.acm.org/doi/10.5555/944790.944812},
  volume = {2},
  year = {2001}
}

@article{nystrm1930uber,
  author = {Nyström, E. J.},
  doi = {10.1007/BF02547521},
  journal = actamath,
  language = {en},
  pages = {185--204},
  title = {Über die praktische Auflösung von Integralgleichungen mit Anwendungen auf Randwertaufgaben},
  volume = {54},
  year = {1930}
}

@article{belabbas2009spectral,
  author = {Belabbas, Mohamed-Ali and Wolfe, Patrick J.},
  chapter = {Physical Sciences},
  copyright = {© 2009 by The National Academy of Sciences of the USA},
  doi = {10.1073/pnas.0810600105},
  journal = pnas,
  language = {en},
  month = jan,
  number = {2},
  pages = {369--374},
  title = {Spectral methods in machine learning and new strategies for very large datasets},
  volume = {106},
  year = {2009}
}

@inproceedings{deshpande2010efficient,
  author = {Deshpande, A. and Rademacher, L.},
  booktitle = focs,
  doi = {10.1109/FOCS.2010.38},
  month = oct,
  pages = {329--338},
  title = {Efficient Volume Sampling for Row/Column Subset Selection},
  year = {2010}
}

@article{massei2022some,
  author = {Massei, Stefano},
  doi = {10.1007/s10543-021-00872-1},
  journal = {BIT Numer. Math.},
  month = {March},
  number = {1},
  pages = {195--220},
  title = {Some Algorithms for Maximum Volume and Cross Approximation of Symmetric Semidefinite Matrices},
  volume = {62},
  year = {2022}
}

@article{gittens2011spectral,
  author = {Gittens, Alex},
  journal = {arXiv preprint \href{http://arxiv.org/abs/1110.5305v1}{arXiv:1110.5305v1}},
  month = oct,
  title = {The spectral norm error of the naïve {Nyström} extension},
  year = {2011}
}

@inproceedings{elalaoui2015fast,
  author = {El Alaoui, Ahmed and Mahoney, Michael W.},
  booktitle = neurips,
  pages = {775--783},
  title = {Fast Randomized Kernel Ridge Regression with Statistical Guarantees},
  url = {https://dl.acm.org/doi/10.5555/2969239.2969326},
  year = {2015}
}

@inproceedings{rudi2018fast,
  author = {Rudi, Alessandro and Calandriello, Daniele and Carratino, Luigi and Rosasco, Lorenzo},
  booktitle = neurips,
  month = dec,
  pages = {5677--5687},
  title = {On fast leverage score sampling and optimal learning},
  url = {https://dl.acm.org/doi/10.5555/3327345.3327470},
  year = {2018}
}

@article{epperly2023kernel,
  archiveprefix = {arXiv},
  author = {Epperly, Ethan N. and Moreno, Elvira},
  eprint = {2306.03955},
  journal = neurips,
  pages = {65850--65868},
  title = {Kernel Quadrature with Randomly Pivoted {Cholesky}},
  url = {https://dl.acm.org/doi/10.5555/3666122.3668997},
  year = {2023}
}

@article{alben1975exact,
  author = {Alben, R. and Blume, M. and Krakauer, H. and Schwartz, L.},
  doi = {10.1103/PhysRevB.12.4090},
  journal = prb,
  pages = {4090--4094},
  title = {Exact results for a three-dimensional alloy with site diagonal disorder: {Comparison} with the coherent potential approximation},
  volume = {12},
  year = {1975}
}

@article{saibaba2025randomized,
  author = {Saibaba, Arvind K. and Mi{\k e}dlar, Agnieszka},
  doi = {10.1137/23M1593255},
  journal = simods,
  number = {1},
  pages = {136--162},
  title = {Randomized Low-Rank Approximations beyond {Gaussian} Random Matrices},
  volume = {7},
  year = {2025}
}

@inproceedings{jiang2021optimal,
  author = {Jiang, Shuli and Pham, Hai and Woodruff, David P. and Zhang, Qiuyi},
  booktitle = neurips,
  pages = {23741--23753},
  title = {Optimal Sketching for Trace Estimation},
  url = {https://dl.acm.org/doi/10.5555/3540261.3542079},
  volume = {34},
  year = {2021}
}

@article{gudmundsson1995small,
  author = {Gudmundsson, Thorkell and Kenney, Charles S. and Laub, Alan J.},
  doi = {10.1137/S0895479893243876},
  journal = simax,
  number = {3},
  pages = {776--792},
  title = {Small-Sample Statistical Estimates for Matrix Norms},
  volume = {16},
  year = {1995}
}

@article{szekely2003extremal,
  author = {Sz{\'e}kely, G{\'a}bor J. and Bakirov, Nail K.},
  doi = {10.1007/s00440-003-0262-6},
  journal = ptrf,
  number = {2},
  pages = {184--202},
  title = {Extremal Probabilities for {Gaussian} Quadratic Forms},
  volume = {126},
  year = {2003}
}

@article{roostakhorasani2015assessing,
  author = {Roosta-Khorasani, Farbod and Sz{\'e}kely, G{\'a}bor J. and Ascher, Uri M.},
  doi = {10.1137/14096311X},
  journal = siamuq,
  number = {1},
  pages = {61--90},
  title = {Assessing Stochastic Algorithms for Large Scale Nonlinear Least Squares Problems Using Extremal Probabilities of Linear Combinations of Gamma Random Variables},
  volume = {3},
  year = {2015}
}

@article{gratton2018improved,
  author = {Gratton, Serge and Titley-Peloquin, David},
  doi = {10.1137/17M1137541},
  journal = simax,
  number = {2},
  pages = {922--931},
  title = {Improved Bounds for Small-Sample Estimation},
  volume = {39},
  year = {2018}
}

@article{hallman2024extremal,
  author = {Hallman, Eric},
  journal = {arXiv preprint \href{http://arxiv.org/abs/2411.15454v1}{arXiv:2411.15454v1}},
  title = {Extremal Bounds for {Gaussian} Trace Estimation},
  year = {2024}
}

@article{epperly2026new,
  author = {Epperly, Ethan N. W.},
  journal = {arXiv preprint \href{http://arxiv.org/abs/2608.20633v1}{arXiv:2608.20633v1}},
  month = aug,
  title = {A new analysis of the randomly pivoted {Cholesky} algorithm},
  year = {2026}
}

@incollection{pearce2026randomized,
  author = {Pearce, Katherine J. and Martinsson, Per-Gunnar},
  booktitle = {Numerical Analysis and Scientific Computing: an Anthology},
  doi = {10.1007/978-3-032-33428-2_7},
  editor = {Koren, Barry and Roose, Dirk},
  pages = {197--252},
  publisher = {Springer},
  series = {Lecture Notes in Computational Science and Engineering},
  title = {Randomized Algorithms for Low-Rank Matrix and Tensor Decompositions},
  year = {2026}
}

@article{kannan2017randomized,
  author = {Kannan, Ravindran and Vempala, Santosh},
  doi = {10.1017/S0962492917000058},
  journal = actanum,
  pages = {95--135},
  title = {Randomized algorithms in numerical linear algebra},
  volume = {26},
  year = {2017}
}

@article{woodruff2014sketching,
  author = {Woodruff, David P.},
  doi = {10.1561/0400000060},
  journal = fttcs,
  language = {en},
  number = {1-2},
  pages = {1--157},
  title = {Sketching as a Tool for Numerical Linear Algebra},
  volume = {10},
  year = {2014}
}

@inproceedings{dereziski2024recent,
  author = {Dereziński, Michał and Mahoney, Michael W.},
  booktitle = kdd,
  doi = {10.1145/3637528.3671461},
  isbn = {979-8-4007-0490-1},
  language = {en},
  month = aug,
  pages = {6470--6479},
  title = {Recent and Upcoming Developments in Randomized Numerical Linear Algebra for Machine Learning},
  year = {2024}
}

@article{tropp2015introduction,
  author = {Tropp, Joel A},
  doi = {10.1561/2200000048},
  journal = ftml,
  number = {1-2},
  pages = {1--230},
  title = {An introduction to matrix concentration inequalities},
  volume = {8},
  year = {2015}
}

@inproceedings{tropp2026applied,
  author = {Tropp, Joel A.},
  booktitle = icm,
  doi = {10.1137/25M1806922},
  isbn = {978-1-61197-872-8 978-1-61197-873-5},
  language = {en},
  month = jul,
  pages = {129--149},
  title = {Applied Random Matrix Theory},
  year = {2026}
}

@techreport{vanhandel2014probability,
  author = {van Handel, Ramon},
  doi = {10.21236/ADA623999},
  institution = {Princeton University},
  language = {en},
  month = jun,
  title = {Probability in High Dimension},
  year = {2014}
}

@article{murray2023randomized,
  author = {Murray, Riley and Demmel, James and Mahoney, Michael W. and Erichson, N. Benjamin and Melnichenko, Maksim and Malik, Osman Asif and Grigori, Laura and Luszczek, Piotr and Derezi{\'n}ski, Micha{\l} and Lopes, Miles E. and Liang, Tianyu and Luo, Hengrui and Dongarra, Jack},
  journal = {arXiv preprint \href{http://arxiv.org/abs/2302.11474v2}{arXiv:2302.11474v2}},
  title = {Randomized Numerical Linear Algebra: {A} Perspective on the Field With an Eye to Software},
  year = {2023}
}

@article{anari2025optimal,
  author = {Anari, Nima and Liu, Yang P. and Vuong, Thuy-Duong},
  doi = {10.1137/22M1524321},
  journal = sicomp,
  number = {5},
  pages = {{FOCS22-93}--{FOCS22-119}},
  title = {Optimal Sublinear Sampling of Spanning Trees and Determinantal Point Processes via Average-Case Entropic Independence},
  volume = {54},
  year = {2025}
}

@inproceedings{derezinski2019exact,
  author = {Derezi{\'n}ski, Micha{\l} and Calandriello, Daniele and Valko, Michal},
  booktitle = neurips,
  title = {Exact sampling of determinantal point processes with sublinear time preprocessing},
  url = {https://dl.acm.org/doi/10.5555/3454287.3455323},
  year = {2019}
}

@inproceedings{calandriello2020sampling,
  author = {Calandriello, Daniele and Derezi{\'n}ski, Micha{\l} and Valko, Michal},
  booktitle = neurips,
  pages = {6889--6899},
  title = {Sampling from a {$k$-DPP} without looking at all items},
  url = {https://dl.acm.org/doi/10.5555/3495724.3496302},
  year = {2020}
}

@article{oseledets2010tt,
  author = {Oseledets, Ivan and Tyrtyshnikov, Eugene},
  doi = {10.1016/j.laa.2009.07.024},
  journal = laa,
  number = {1},
  pages = {70--88},
  title = {{TT}-cross approximation for multidimensional arrays},
  volume = {432},
  year = {2010}
}

@article{rudelson2015small,
  author = {Rudelson, Mark and Vershynin, Roman},
  doi = {10.1093/imrn/rnu243},
  journal = imrn,
  number = {19},
  pages = {9594--9617},
  title = {Small Ball Probabilities for Linear Images of High Dimensional Distributions},
  volume = {2015},
  year = {2015}
}

@article{rudelson2008littlewood,
  author = {Rudelson, Mark and Vershynin, Roman},
  doi = {10.1016/j.aim.2008.01.010},
  journal = advmath,
  number = {2},
  pages = {600--633},
  title = {The {Littlewood}--{Offord} Problem and Invertibility of Random Matrices},
  volume = {218},
  year = {2008}
}

@book{vandervaart1998asymptotic,
  author = {van der Vaart, Aad W.},
  doi = {10.1017/CBO9780511802256},
  publisher = {Cambridge University Press},
  series = {Cambridge Series in Statistical and Probabilistic Mathematics},
  title = {Asymptotic Statistics},
  year = {1998}
}

@article{tang2012probing,
  author = {Tang, Jok M. and Saad, Yousef},
  doi = {10.1002/nla.779},
  journal = nlaa,
  number = {3},
  pages = {485--501},
  title = {A probing method for computing the diagonal of a matrix inverse},
  volume = {19},
  year = {2012}
}

@article{stathopoulos2013hierarchical,
  author = {Stathopoulos, Andreas and Laeuchli, Jesse and Orginos, Kostas},
  doi = {10.1137/120881452},
  journal = sisc,
  number = {5},
  pages = {S299--S322},
  title = {Hierarchical probing for estimating the trace of the matrix inverse on toroidal lattices},
  volume = {35},
  year = {2013}
}

@article{laeuchli2020extending,
  author = {Laeuchli, Jesse and Stathopoulos, Andreas},
  doi = {10.1137/18M1176427},
  journal = sisc,
  number = {3},
  pages = {A1459--A1485},
  title = {Extending hierarchical probing for computing the trace of matrix inverses},
  volume = {42},
  year = {2020}
}

@article{frommer2021analysis,
  author = {Frommer, Andreas and Schimmel, Claudia and Schweitzer, Marcel},
  doi = {10.1137/20M1364461},
  journal = simax,
  number = {3},
  pages = {1290--1318},
  title = {Analysis of probing techniques for sparse approximation and trace estimation of decaying matrix functions},
  volume = {42},
  year = {2021}
}

@article{frommer2025analysis,
  author = {Frommer, Andreas and Rinelli, Michele and Schweitzer, Marcel},
  doi = {10.1090/mcom/3984},
  journal = mc,
  number = {352},
  pages = {801--823},
  title = {Analysis of stochastic probing methods for estimating the trace of functions of sparse symmetric matrices},
  volume = {94},
  year = {2025}
}

@article{saibaba2017randomized,
  author = {Saibaba, Arvind K. and Alexanderian, Alen and Ipsen, Ilse C. F.},
  doi = {10.1007/s00211-017-0880-z},
  journal = nm,
  number = {2},
  pages = {353--395},
  title = {Randomized matrix-free trace and log-determinant estimators},
  volume = {137},
  year = {2017}
}

@article{frommer2022multilevel,
  author = {Frommer, Andreas and Nasr Khalil, Mostafa and
Ramirez-Hidalgo, Gustavo},
  doi = {10.1137/21M1441894},
  journal = sisc,
  number = {4},
  pages = {A2536--A2556},
  title = {A multilevel approach to variance reduction in the stochastic estimation of the trace of a matrix},
  volume = {44},
  year = {2022}
}

@article{hallman2022multilevel,
  author = {Hallman, Eric and Troester, Devon},
  doi = {10.1016/j.laa.2021.12.010},
  journal = laa,
  pages = {125--149},
  title = {A multilevel approach to stochastic trace estimation},
  volume = {638},
  year = {2022}
}

@article{amsel2026fixed,
  author = {Amsel, Noah and Chen, Tyler and Duman Keles, Feyza and
Halikias, Diana and Musco, Cameron and Musco, Christopher},
  doi = {10.1137/25M1742710},
  journal = simax,
  number = {2},
  pages = {483--511},
  title = {Fixed-Sparsity Matrix Approximation from Matrix-Vector Products},
  volume = {47},
  year = {2026}
}

@article{chen2023krylov,
  author = {Chen, Tyler and Hallman, Eric},
  doi = {10.1137/22M1494257},
  journal = simax,
  number = {3},
  pages = {1218--1244},
  title = {{Krylov}-Aware Stochastic Trace Estimation},
  volume = {44},
  year = {2023}
}

@article{madhavan2026flextrace,
  author = {Madhavan, Madhusudan and Alexanderian, Alen and
Saibaba, Arvind K.},
  journal = {arXiv preprint \href{http://arxiv.org/abs/2603.05721v1}{arXiv:2603.05721v1}},
  title = {{FlexTrace:} {Exchangeable} Randomized Trace Estimation for Matrix Functions},
  year = {2026}
}

@article{bernardson1994monte,
  author = {Bernardson, S. and McCarty, P. and Thron, C.},
  doi = {10.1016/0010-4655(94)90004-3},
  journal = cpc,
  number = {3},
  pages = {256--264},
  title = {{Monte} {Carlo} methods for estimating linear combinations of inverse matrix entries in lattice {QCD}},
  volume = {78},
  year = {1994}
}

@article{foley2005practical,
  author = {Foley, Justin and Juge, K. Jimmy and {\'O} Cais, Alan
and Peardon, Mike and Ryan, Sin{\'e}ad M.
and Skullerud, Jon-Ivar},
  doi = {10.1016/j.cpc.2005.06.008},
  journal = cpc,
  number = {3},
  pages = {145--162},
  title = {Practical all-to-all propagators for lattice {QCD}},
  volume = {172},
  year = {2005}
}

@article{babich2012exploring,
  author = {Babich, Ronald and Brower, Richard C. and Clark, Michael A.
and Fleming, George T. and Osborn, James C.
and Rebbi, Claudio and Schaich, David},
  doi = {10.1103/PhysRevD.85.054510},
  journal = prd,
  number = {5},
  pages = {054510},
  title = {Exploring strange nucleon form factors on the lattice},
  volume = {85},
  year = {2012}
}

@article{benzi2007decay,
  author = {Benzi, Michele and Razouk, Nader},
  journal = etna,
  pages = {16--39},
  title = {Decay bounds and {$O(n)$} algorithms for approximating functions of sparse matrices},
  url = {https://etna.ricam.oeaw.ac.at/vol.28.2007-2008/pp16-39.dir/pp16-39.pdf},
  volume = {28},
  year = {2007}
}

@article{giles2015multilevel,
  author = {Giles, Michael B.},
  doi = {10.1017/S096249291500001X},
  journal = actanum,
  pages = {259--328},
  title = {Multilevel {Monte} {Carlo} Methods},
  volume = {24},
  year = {2015}
}

@article{cho2022computationally,
  author = {Cho, Taewon and Chung, Julianne and Miller, Scot M. and
Saibaba, Arvind K.},
  doi = {10.5194/gmd-15-5547-2022},
  journal = gmd,
  number = {14},
  pages = {5547--5565},
  title = {Computationally efficient methods for large-scale atmospheric inverse modeling},
  volume = {15},
  year = {2022}
}

@article{hallman2023monte,
  author = {Hallman, Eric and Ipsen, Ilse C. F. and Saibaba, Arvind K.},
  doi = {10.1137/22M1476277},
  journal = simax,
  number = {1},
  pages = {240--269},
  title = {{Monte} {Carlo} Methods for Estimating the Diagonal of a Real Symmetric Matrix},
  volume = {44},
  year = {2023}
}

@article{gebremedhin2005what,
  author = {Gebremedhin, Assefaw Hadish and Manne, Fredrik and Pothen, Alex},
  doi = {10.1137/S0036144504444711},
  journal = sirev,
  number = {4},
  pages = {629--705},
  title = {What Color Is Your {Jacobian?} {Graph} Coloring for Computing Derivatives},
  volume = {47},
  year = {2005}
}

@inproceedings{musco2026nearly,
  author = {Musco, Christopher and Ramesh, Indu},
  booktitle = esa,
  doi = {10.4230/LIPIcs.ESA.2026.85},
  pages = {85:1--85:14},
  series = {Leibniz International Proceedings in Informatics (LIPIcs)},
  title = {Nearly Instance Optimal Sparse Matrix Approximation from Matrix-Vector Products},
  volume = {388},
  year = {2026}
}

@inproceedings{sobczyk2022approximate,
  author = {Sobczyk, Aleksandros and Luisier, Mathieu},
  booktitle = neurips,
  pages = {19357--19369},
  title = {Approximate {Euclidean} Lengths and Distances Beyond {Johnson}--{Lindenstrauss}},
  url = {https://dl.acm.org/doi/10.5555/3600270.3601677},
  volume = {35},
  year = {2022}
}

@article{dereziski2018reverse,
  author = {Dereziński, Michał and Warmuth, Manfred K.},
  journal = jmlr,
  number = {23},
  pages = {853--891},
  title = {Reverse iterative volume sampling for linear regression},
  url = {https://dl.acm.org/doi/10.5555/3291125.3291148},
  volume = {19},
  year = {2018}
}

@article{hough2006determinantal,
  author = {Hough, J. Ben and Krishnapur, Manjunath and Peres, Yuval and Virág, Bálint},
  doi = {10.1214/154957806000000078},
  journal = probsurv,
  month = jan,
  pages = {206--229},
  title = {Determinantal Processes and Independence},
  volume = {3},
  year = {2006}
}

@inproceedings{barthelme2023faster,
  author = {Barthelme, Simon and Tremblay, Nicolas and Amblard, Pierre-Olivier},
  booktitle = aistats,
  pages = {5582--5592},
  title = {A Faster Sampler for Discrete Determinantal Point Processes},
  url = {https://proceedings.mlr.press/v206/barthelme23a.html},
  year = {2023}
}

@inproceedings{dereziski2019minimax,
  author = {Dereziński, Michał and Clarkson, Kenneth L. and Mahoney, Michael W. and Warmuth, Manfred K.},
  booktitle = colt,
  language = {en},
  month = jun,
  pages = {1050--1069},
  title = {Minimax experimental design: {Bridging} the gap between statistical and worst-case approaches to least squares regression},
  url = {http://proceedings.mlr.press/v99/derezinski19b.html},
  year = {2019}
}

@article{dereziski2022unbiased,
  author = {Dereziński, Michał and Warmuth, Manfred K. and Hsu, Daniel},
  journal = jmlr,
  number = {167},
  pages = {1--46},
  title = {Unbiased estimators for random design regression},
  url = {http://www.jmlr.org/papers/v23/19-571.html},
  volume = {23},
  year = {2022}
}

@inproceedings{dereziski2018leveraged,
  author = {Dereziński, Michal and Warmuth, Manfred K and Hsu, Daniel J},
  booktitle = neurips,
  pages = {2510--2519},
  title = {Leveraged Volume Sampling for Linear Regression},
  url = {https://dl.acm.org/doi/10.5555/3327144.3327176},
  volume = {31},
  year = {2018}
}

@inproceedings{higgins2025high,
  author = {Higgins, Andrew J. and Boman, Erik G. and Yamazaki, Ichitaro},
  booktitle = scworkshops,
  doi = {10.1145/3731599.3767544},
  title = {A High Performance {GPU} {CountSketch} Implementation and Its Application to Multisketching and Least Squares Problems},
  year = {2025}
}

@inproceedings{derezinski2021sparse,
  author = {Derezi{\'n}ski, Micha{\l} and Liao, Zhenyu and Dobriban, Edgar and Mahoney, Michael W.},
  booktitle = colt,
  pages = {1467--1510},
  title = {Sparse sketches with small inversion bias},
  url = {https://proceedings.mlr.press/v134/derezinski21a.html},
  year = {2021}
}

@article{tikhomirov2026level,
  author = {Tikhomirov, Konstantin},
  journal = {arXiv preprint \href{http://arxiv.org/abs/2607.23017v1}{arXiv:2607.23017v1}},
  title = {Level-set entropy and sparse randomized embeddings},
  year = {2026}
}

@article{tropp2026comparisonb,
  author = {Tropp, Joel A.},
  journal = {arXiv preprint \href{http://arxiv.org/abs/2603.04365v5}{arXiv:2603.04365v5}},
  title = {Comparison theorems for the extreme eigenvalues of a random symmetric matrix},
  year = {2026}
}

@article{simchowitz2017gap,
  author = {Simchowitz, Max and El Alaoui, Ahmed and Recht, Benjamin},
  journal = {arXiv preprint \href{http://arxiv.org/abs/1704.04548v1}{arXiv:1704.04548v1}},
  title = {On the gap between strict-saddles and true convexity: an {$\Omega(\log d)$} lower bound for eigenvector approximation},
  year = {2017}
}

@article{chu2025improved,
  author = {Chu, Ya-Chi and Cortinovis, Alice},
  doi = {10.1016/j.laa.2025.04.001},
  journal = laa,
  pages = {68--93},
  title = {Improved bounds for randomized {Schatten} norm estimation of numerically low-rank matrices},
  volume = {717},
  year = {2025}
}

@article{kong2017spectrum,
  author = {Kong, Weihao and Valiant, Gregory},
  doi = {10.1214/16-AOS1525},
  journal = annalsstats,
  number = {5},
  pages = {2218--2247},
  title = {Spectrum estimation from samples},
  volume = {45},
  year = {2017}
}

\label{lastpage}

\end{document}